\documentclass[11pt]{amsart}
\usepackage{amsmath,amssymb}
\usepackage{graphicx}

\title{Self-dual metrics and twenty-eight bitangents}
\author{Nobuhiro Honda $^{\dag}$
}  
\thanks
{$^{\dag}$Partially supported by
Research Fellowships of the 
Japan Society for the Promotion
of Science for Young Scientists.
}

\newcommand{\ol}{\overline}
\newcommand{\ra}{\rightarrow}
\newcommand{\lra}{\longrightarrow}
\newcommand{\da}{\downarrow}
\newcommand{\ua}{\uparrow}

\newcommand{\set}{\,|\,}
\newcommand{\proofend}{\hfill$\square$}
\newcommand{\linfty}{\langle l_{\infty}\rangle}

\newtheorem{prop}{Proposition}[section]
\newtheorem{lemma}[prop]{Lemma}

\newtheorem{thm}[prop]{Theorem}
\newtheorem{rmk}[prop]{Remark}
\newtheorem{cor}[prop]{Corollary}

\newtheorem{definition}[prop]{Definition}

\begin{document}

\maketitle

\begin{abstract} 

We study self-dual metrics on $3\mathbf{CP}^2$ of positive scalar curvature admitting a non-zero Killing field, but which are not conformally isometric to LeBrun's metrics.
Firstly, we determine defining equations of the twistor spaces of such self-dual metrics.
Next we prove that conversely, the complex threefolds defined by the equations always become twistor spaces of self-dual metrics on $3\mathbf{CP}^2$ of the above kind.
As a corollary, we determine a global structure of the moduli spaces of these self-dual metrics; namely we show that the moduli space  is  naturally identified with an orbifold $\mathbf{R}^3/G$, where $G$ is an involution of $\mathbf{R}^3$ having two-dimensional fixed locus.
Combined with works of  LeBrun, this settles a moduli problem of self-dual metrics on $3\mathbf{CP}^2$ of positive scalar curvature admitting a non-trivial Killing field.
In particular, it is shown that any two self-dual metrics on $3\mathbf{CP}^3$ of positive scalar curvature admitting a non-zero Killing field can be connected by deformation keeping the self-duality. 
In our proof, a key role is played by a classical result in algebraic geometry that a smooth plane quartic always possesses twenty-eight bitangents.
\end{abstract}

\tableofcontents


\section{Introduction}

A Riemannian metric on an oriented four-manifold
is called self-dual if the anti-self-dual part of the Weyl conformal
curvature of the metric identically  vanishes. 
Basic examples are
provided by the round metric on the four-sphere and the Fubini-Study
metric on the  complex projective plane. In general, one can expect that
if two four-manifolds admit self-dual metrics respectively, then their
connected sum  will also admit a self-dual metric. In fact, Y.S. Poon
\cite{P86}  constructed  explicit examples of  self-dual metrics on
$2\mathbf{CP}^2$, the connected sum of two complex projective planes.
He further showed that on $2\mathbf{CP}^2$  there are no self-dual
metrics other than his metrics, under assumption of 
the positivity of the scalar curvature.  Later, C. LeBrun \cite{LB91}
and D. Joyce \cite{J95} respectively constructed a family of self-dual metrics of
positive scalar curvature on $n\mathbf{CP}^2$ for any $n\geq 1$. These
are  called LeBrun metrics and Joyce metrics respectively, and
have nice characterizations by the (conformal) isometry group. Namely, 
A. Fujiki \cite{F00} proved that if a self-dual metric on
$n\mathbf{CP}^2$  has  an effective $U(1)\times U(1)$-isometry, then it must be a Joyce metric. LeBrun \cite{LB93} 
showed that  if a self-dual metric on $n\mathbf{CP}^2$ has a non-trivial
semi-free
$U(1)$-isometry, then the metric must be a LeBrun metric. Here, a
$U(1)$-action on a manifold $M$ is called semi-free if the isotropy
group is
$U(1)$ or identity only, at every point of $M$.
By his construction, LeBrun metrics form a connected family for each $n$.

Now drop the assumption of the semi-freeness of $U(1)$-isometry.
A work of Pedersen-Poon \cite{PP95} suggests that there should exist self-dual metrics on $n\mathbf{CP}^2$, $n\geq 3$, admitting non-semifree $U(1)$-isometry.
In fact, this is shown to be true by the author \cite{Hon02,Hon03}.
However, it seems to be quite difficult to classify {\em all}\, self-dual metrics on $n\mathbf{CP}^2$ with $U(1)$-isometry, for arbitrary $n$.
The purpose of the present paper is to settle this problem for $3\mathbf{CP}^2$.
Namely, we show the following:

\begin{thm}\label{thm-modulisp}
Let $\mathcal M$ be the moduli space of conformal classes of self-dual metrics on $3\mathbf{CP}^2$  satisfying the following conditions: (i) having a positive scalar curvature, (ii) admitting a non-zero Killing field, or equivalently, a non-trivial $U(1)$-isometry, (iii) being not conformal to LeBrun metrics. 
Then $\mathcal M$ is non-empty and is naturally identified with an orbifold $\mathbf R^3/G$, where $G$ is an involution on $\mathbf R^3$ having two-dimensional fixed locus.
\end{thm}

Thus a global structure of the moduli space is determined.
In particular, it is connected as in the LeBrun's case.
We note that the claims of the theorem involve  the uniqueness of the $U(1)$-action in (ii). The $U(1)$-action is inequivalent to LeBrun's semi-free action.

Next we mention a relationship between our self-dual metrics and LeBrun's self-dual metric.
We showed in \cite{Hon02} that some of self-dual metrics on $3\mathbf{CP}^2$ satisfying three conditions in Theorem \ref{thm-modulisp} can be obtained as a small deformation of Joyce metric with torus symmetry. 
Because Joyce metrics on $3\mathbf{CP}^2$ fall into a subfamily of LeBrun metrics, Theorem \ref{thm-modulisp} implies the following

\begin{cor}
Any two self-dual metrics on $3\mathbf{CP}^2$ satisfying (i) and (ii) in Theorem \ref{thm-modulisp} 
can be connected by deformation keeping the self-duality.
\end{cor}

Alternatively, one can also state that both of our self-dual metrics and LeBrun metrics have a Joyce metric (or LeBrun metric with torus action) as a limit.
Because the Killing field (or associated $U(1)$-action) of ours and LeBrun's are different, one must exchange the Killing filed when passing through a LeBrun metric with torus symmetry.

 Now we explain why the involution $G$ appears in Theorem \ref{thm-modulisp}.
 In the course of our proof of Theorem \ref{thm-modulisp}, we construct a family $\tilde{\mathcal M}$ of self-dual metrics on $3\mathbf{CP}^2$ parameterized by $\mathbf R^3$, which contains all metrics satisfying conditions (i)--(iii).
 But it will turned out that the family does not effectively parameterize the metrics;
nemaly we see that an involution $G$ acts on this parameter space $\mathbf R^3$ in such a way that two points exchanged by the involution represent the same conformal class.
Fixed points of $G$, which will turn out to be two-dimensional, represent  self-dual metrics having a non-trivial isometric involution, where 
non-trivial means that the involution does not belong to $U(1)$.

\vspace{2mm}

Our proof of Theorem \ref{thm-modulisp} is based on the twistor theory.
Namely, by using so called the Penrose correspondence between self-dual metrics and complex threefolds called twistor spaces, we translate the problem into that of complex algebraic geometry. 
Rather surprisingly, this translation leads us to an unexpected connection between self-dual metrics and bitangents of a plane quartic curve which is a classical but fruitful topic in algebraic geometry.
It is this connection which enable us to find {\em{arbitrary}} twistor lines in certain complex threefolds, which is equivalent to solving the self-duality equation for Riemannian metrics on a f-manifold\cite{ahs}.

In Section \ref{s-defeq} we first determine defining equations of the twistor spaces of self-dual metrics on $3\mathbf{CP}^2$ satisfying the above three conditions in Theorem \ref{thm-modulisp}.
More concretely, we prove that the twistor spaces have a structure of generically two-to-one covering branched along a quartic
\begin{equation}\label{eqn-brnch}
\left(y_2y_3+Q(y_0,y_1)\right)^2-y_0y_1(y_0+y_1)(y_0-ay_1)=0,
\end{equation}
where $Q(y_0,y_1)$ is a  quadratic form of $y_0$ and
$y_1$ with real coefficients, and $a$ is a  positive  real number
(Proposition \ref{prop-def-B}). 
Further, we show that $Q$ and $a$ must satisfy certain inequality, which we call Condition (A) (Proposition \ref{prop-necessa}).

Most of the rest of the paper is devoted to proving that, conversely, a complex threefold $Z$ having a structure of generically two-to-one covering branched along the quartic (\ref{eqn-brnch}) has a structure of a twistor space of $3\mathbf{CP}^2$. 
Since the quartic surface (\ref{eqn-brnch}) always has isolated singularities, the double covering has isolated singularities over there.
For this space, we prove the following theorem completing inversion construction of  twistor spaces:

\begin{thm}\label{thm-corresp} 
For any quartic surfaces $B$ of the form (\ref{eqn-brnch}) satisfying the condition (A), there exists a resolution of the double covering of $\mathbf{CP}^3$ branched along $B$, such that the resulting smooth threefold is a twistor space of $3\mathbf{CP}^2$.
\end{thm}

In fact, the main result of this paper (Theorem \ref{thm-inv}) implies more.
Namely, we prove that there exist  precisely two (small)  resolutions of the double cover which are actually twistor spaces, and that these two spaces are isomorphic as twistor spaces so that they determine the same self-dual structure.

It is relatively easy to derive Theorem \ref{thm-modulisp} (= Theorem \ref{thm-mod}) from  Theorem \ref{thm-inv}.
As is already mentioned, our strategy for proving this result is to find the family of twistor lines in $Z$; namely a real four-dimensional family of real smooth rational curves which foliates the whole of the threefold $Z$.
In general, it is difficult (except a few simple cases) to find arbitrary twistor lines. 
However, in the present case this is possible basically because their images onto $\mathbf{CP}^3$ (by the covering morphism)  must be, in general, conics with a very special property; namely they are conics touching the branch quartic surface (Proposition \ref{prop-image}).
We call this kind of conics in $\mathbf{CP}^3$ touching conics (following Hadan \cite{Ha00}).
Roughly speaking, every touching conics are `generated' from bitangents of the quartic, and the existence of twenty-eight bitangents guarantees the existence of twistor lines.
This is how the self-dual metrics and twenty-eight bitangents are related.

Thus a large part of  the problem of finding twistor lines in $Z$ is reduced to finding touching conics of the quartic (\ref{eqn-brnch}).
Note that a conic in $\mathbf{CP}^3$ is always contained in a unique plane.
Broadly speaking, our method of finding touching conics of the quartics consists of two parts: one is to find `general' touching conics, and the other is to find `degenerate' ones which are  limits of the `general' touching conics.
Here, `general' means that the plane on which the conic lies intersects the branch quartic (\ref{eqn-brnch}) smoothly.
It is possible to show that on this kind of planes there are 63 one-dimensional families of touching conics (Proposition \ref{prop-ftc3}; see also \cite{Ha00}).
A family which can actually be the images of twistor lines must be unique and  there are too many candidates of twistor lines, even if we take the reality into account.
This seems to be  main difficulty in proving that the threefolds are actually  twistor spaces (cf.\,\cite{KK92, Ha00}).

We are able to overcome this difficulty by considering `degenerate' touching conics, whose planes intersect the quartic (\ref{eqn-brnch}) non-smoothly.
Thanks to a very special form of the equation (\ref{eqn-brnch}), it is possible to determine which plane intersects the quartic surface non-smoothly (Lemma \ref{lemma-singsection}). 
Then it is readily seen that the sections of the quartic surface by these planes are union of two irreducible conics.
The point is that, for these planes, we can determine all families of (`degenerate') touching conics in very explicit form (Propositions \ref{prop-a}, \ref{prop-b} and \ref{prop-c}).
These in particular show that the number of families of touching conics drastically decreases for these planes.
Moreover, we can use this explicit descriptions of `degenerate' touching conics to determine which family  actually comes from twistor lines (Proposition \ref{prop-type1}). 
These are done in the latter half of  Section \ref{s-detc} and Section \ref{s-nb} by actually investigating the inverse images (in $Z$) of `degenerate' touching conics.

Thus we obtain candidates of  twistor lines lying over planes intersecting the quartic surface  non-smoothly. 
We can show that these candidates can be deformed in $Z$ to give  twistor lines whose images are `general' touching conics, and that among the above 63 families of   `general' touching conics, there is a unique family which is obtained in this way (Proposition \ref{prop-globalext}).
This is how we obtain twistor lines whose images are (touching) conics.
These form real four-dimensional family and cover an open subset of $Z$.
However, this family does not cover the whole of $Z$.
In order to cover the whole of $Z$, we need to consider another set of  twistor lines whose images becomes {\em{lines}} in $\mathbf{CP}^3$ (cf. Proposition \ref{prop-image}).
These twistor lines are investigated in Section \ref{s-line}  and will be called `twistor lines at infinity'.
A conclusion is that the inverse images of  real lines (in $\mathbf{CP}^3$) going through the unique real singular point of $B$ (coming from Condition (A)) always contain a real twistor line as its irreducible component (Proposition \ref{prop-imline}).
The parameter space of these twistor lines (at infinity) becomes a boundary for compactifying the parameter space of twistor lines whose images are touching conics.

In this way we get candidates of all twistor lines in $Z$.
It remains to show that  $Z$ is actually foliated by these candidates.
Namely, we have to show that different members of these candidates do not intersect, and that they cover the whole of $Z$.
These will be proved in Proposition \ref{prop-disj}.
The point of our proof is to consider the inverse images of arbitrary real lines, which are   smooth elliptic curves in general, and sometimes degenerate into a cycle of rational curves.
We show that these elliptic curves and their degenerations are foliated by circles which are intersection of the elliptic curve with the family of twistor lines lying over a real plane containing the original line. 
Finally we show that the parameter space of these twistor lines is actually $3\mathbf{CP}^2$, and complete the inversion construction as  Theorem \ref{thm-inv}.

Finally, we give some comments concerning the quartic surfaces (\ref{eqn-brnch}).
The equation (\ref{eqn-brnch}) itself appears in Y.\,Umezu's paper \cite[p.\,141]{Um84}, where she investigated normal quartic surfaces whose minimal resolutions are birational to elliptic ruled surfaces (which are by definition a $\mathbf{CP}^1$-bundle over an elliptic curve).
In particular, she proved that if  such a quartic surface has two simple elliptic singularities of type $\tilde{E}_7$, then one can take (\ref{eqn-brnch}) as its defining equation.
By counting the number of coefficients, these quartic surfaces form a four-dimensional family.
Our Condition (A) imposes one more condition on the coefficients and
the surface (\ref{eqn-brnch}) has one more singularity (which is an ordinary double point).
Consequently our moduli space becomes three-dimensional as in Theorem \ref{thm-modulisp}.
Also it should be noted that the equation (\ref{eqn-brnch}) is already appeared in the paper of Kreu\ss ler-Kurke \cite{KK92}, which they obtained while carefully classifying branch quartic surfaces for  {\it general}  twistor space (namely no assumptions on the Killing fields) of $3\mathbf{CP}^2$.

\vspace{3mm}
I would like to thank Shigeharu Takayama whose question to the author at the Sugadaira Symposium in the fall of 2002 was the starting point of this investigation.
Also I would like to thank Nobuyoshi Takahashi and Yoshihisa Saito for stimulating discussions and helpful  comments. 

\vspace{3mm}
\noindent{\bf{Notations and conventions.}}
Let $X$ be a complex manifold and $\sigma$ an anti-holomorphic involution on $X$. 
A subset $Y\subset X$ is said to be {\em{real}}\, if $\sigma(Y)=Y$.
The set of $\sigma$-fixed points of $Y$ is denoted by $Y^{\sigma}$.
A holomorphic line bundle $F$ on $ X$ is said to be  {\em{real}}\, if $\sigma^* F$ is isomorphic to $ \ol{F}$, where $\ol{F}$ denotes the complex conjugation of $F$.
A complete linear system $|F|$ is called real if $F$ is a real line bundle, and then $|F|^{\sigma}\subset |F|$ denotes the subset consisting of real (i.e.\,$\sigma$-invariant) members.
This is necessarily parametrized by the real projective space $\mathbf{RP}^{n-1}$, where $n$ is the complex dimension of $H^0(F)$.
Next assume that  $Z$ is  a twistor space of a self-dual four-manifold $(M,g)$. 
Let $\sigma$ be the real structure which is an anti-holomorphic involution of $Z$.
{\em A twistor line}\, is a fiber of the twistor fibration $Z\ra M$.
Next, suppose that $Z$ is a three-dimensional complex manifold with an anti-holomorphic involution $\sigma$, for which we do not assume that $Z$ is actually a twistor space. 
Then a non-singular complex submanifold $L\subset Z$ is called a {\em real line} if $L$ is biholomorphic to the complex projective line and if its normal bundle in $Z$ is isomorphic to $O(1)^{\oplus 2}$.
Needless to say, twistor line is a real line.
However, it is important to make a distinction of `twistor lines' and `real lines', because there is an example of twistor space   possessing a real line which is not a twistor lines
(as proved in Corollary \ref{cor-chatl}).

\section{Defining equations of the branch quartic surfaces and their
singularities}
\label{s-defeq} Let $g$ be a self-dual metric on $3\mathbf{CP}^2$ of positive
scalar curvature, and assume that $g$ is not conformally isometric to LeBrun
metrics. Let
$Z$  be the twistor space of $g$, and denote by $(-1/2)K_Z$  the fundamental
line bundle which is the canonical square root of the anticanonical line bundle
of $Z$.
These non-LeBrun twistor spaces of $3\mathbf{CP}^2$ are extensively studied in 
 Kreu\ss ler and Kurke  \cite{KK92} and Poon \cite{P92}, 
and it has been proved that  the fundamental
system (the complete linear system associated to the fundamental line bundle) is
free and of three-dimensional, and induces a  surjective morphism
$\Phi:Z\ra\mathbf{CP}^3$ which is generically two-to-one, and  that the branch 
divisor $B$ is a quartic surface with only isolated singularities.
Furthermore, there is the following diagram:

\begin{equation}\label{diagram001}
\begin{array}{rlll}
&Z&&\\
\mu&\da&\searrow\Phi&\\
&Z_0&\lra\mathbf{CP}^3\\
&&\Phi_0&
\end{array}
\end{equation}
where $\Phi_0:Z_0\ra \mathbf{CP}^3$
denotes the  double covering branched along
$B$, and $\mu$ is {\it{a small
resolution}}\, of the singularities of $Z_0$ over the singular points of
$B$. 

For generic non-LeBrun metric $g$, $B$ has only ordinary double points \cite{KK92, P92} and
hence is birational to a K3 surface. As a consequence, 
one can deduce that $Z$ does not
admit a non-zero holomorphic vector field.
However, the author showed
in \cite{Hon02} and \cite{Hon03}  that 
if $B$ degenerates to have non-ADE singularities, then 
$Z$ admits a non-zero holomorphic vector field,
and that such a twistor space of $3\mathbf{CP}^2$ actually exists.
Concerning the defining equation of
the branch quartic $B$ for such twistor spaces, we have the following
proposition which is the starting point of our investigation.

\begin{prop}\label{prop-def-B} 
Let $g$ be a non-LeBrun self-dual metric on 
$\mathbf{CP}^2$ of positive scalar curvature, and assume the
existence of a non-trivial Killing field. Let 
$\Phi:Z\ra\mathbf{CP}^3$ and $B\subset\mathbf{CP}^3$ be as above.
Then there exists a homogeneous coordinate
$(y_0:y_1:y_2:y_3)$ on
$\mathbf{CP}^3$ fulfilling (i)-(iii) below: 

\noindent (i) a defining equation of
$B$ is given  by
\begin{equation}\label{eqn-B}
\left(y_2y_3+Q(y_0,y_1)\right)^2-y_0y_1(y_0+y_1)(y_0-ay_1)=0,
\end{equation} where $Q(y_0,y_1)$ is a   
quadratic form of $y_0$ and
$y_1$ with real coefficients, and $a$  is a  positive  real number,

\noindent (ii) the naturally induced real structure on $\mathbf{CP}^3$ is given by
$$
\sigma(y_0:y_1:y_2:y_3)=\left(\ol{y}_0:\ol{y}_1:\ol{y}_3:\ol{y}_2\right),
$$

\noindent (iii) the naturally induced $U(1)$-action on
$\mathbf{CP}^3$ is given by
$$(y_0:y_1:y_2:y_3)\mapsto 
\left(y_0:y_1:e^{i\theta}y_2:e^{-i\theta}y_3\right),\hspace{3mm} 
e^{i\theta}\in U(1).$$
\end{prop}
 
\noindent Proof.  If the fundamental system of $Z$ 
is free, there are just four reducible members, all of which are real
\cite{P92,KK92}. 
We write $\Phi^{-1}(H_{ i})=D_{ i}+\ol{D}_{ i}$, $1\leq i\leq 4$, where
$H_i$ is a real plane in $\mathbf{CP}^3$.
 The restrictions of $\Phi$ onto $D_{i}$ and
$\ol{D}_{i}$ are obviously birational morphisms onto $H_{i}$, so that,
together with the fact that $(-1/2)K_Z\cdot L_i=2$,  it can be readily seen that
$T_{i}:=\Phi(L_{i})$ is a conic contained in $B$. This implies that
the restriction of 
$B$ onto $H_{i}$ is a  conic of multiplicity two. Namely, $T_{i}$,
$1\leq i\leq 4$, is so called {\em a trope} of $B$. 

A Killing field naturally gives rise to an isometric $U(1)$-action,
which can be canonically lifted to a holomorphic $U(1)$-action on the
twistor space.
This action  naturally goes down to 
$\mathbf{CP}^3$, and every subvarieties above are clearly preserved by
these $U(1)$-actions. In particular, $T_{i}$ is a
$U(1)$-invariant conic on a
$U(1)$-invariant plane $H_{i}$, where the $U(1)$-action
is induced by the vector field. Since the twistor fibration is 
$U(1)$-equivariant and generically one-to-one on $D_i$, and since $\Phi|_{D_i}:D_i\ra
H_i$ is also $U(1)$-equivariant and birational, the  $U(1)$-action on any
$H_{i}$ is non-trivial. Hence  $U(1)$ acts non-trivially on $T_{i}$. For
$j\neq i$, put
$l_{ij}=H_i\cap H_j$, which is clearly a real $U(1)$-invariant  line. Then 
 $C_{i}\cap
l_{ij}$ must be the two
$U(1)$-fixed points on $T_{i}$, since it is real set and since there is 
no real point on $T_i$. This implies that $l_{ij}$ is
independent of the choice of $j\neq i$. So we write
$l_{ij}=l_{\infty}$ and let $P_{\infty}$ and $\ol{P}_{\infty}$ be  the
two fixed points of the $U(1)$-action on $l_{\infty}$. Then 
$H_{i}$,
$1\leq i\leq 4$, must be  real  members of the real pencil of planes whose base locus is
$l_{\infty}$.  Since
$l_{\infty}$  is a real line, we can choose  real linear forms $y_0$ and
$y_1$ such that
$l_{\infty}=\{y_0=y_1=0\}$. Further, since any of $H_{i}$ is real, by
applying a real projective transformation (with respect to
$(y_0,y_1)$), we may assume that $\cup_{i=1}^4
H_{i}=\{y_0y_1(y_0+y_1)(ay_0-by_1)=0\}$, where $a,b\in\mathbf R$ with
$a>0$ and
$b>0$. 

 As seen above, every $T_{i}$ goes through $P_{\infty}$ and
$\ol{P}_{\infty}$.
 Let $l_{i}$, $1\leq i\leq 4$  be the tangent line of
$T_{i}$ at $P_{\infty}$. Now we claim that $l_{i}$'s  are lying on the
same plane. Let
$H$ be the plane containing $l_1$ and $l_2$.  
Then by using 
$l_1\cap B=P_{\infty}= l_2\cap B$,
we can easily deduce that  $B\cap H$  is a 
 union of lines, all of which goes through $P_{\infty}$. 
Suppose that $l_3$ is not contained in $H$. Then the line $H\cap
H_3$ is not tangent to $T_3$, so there is an intersection point  of
$T_3\cap H$ other than
$P_{\infty}$. Then the line $H\cap H_3$ is contained in $B$, because we
have already seen that $B\cap H$ is a union of lines all passing through
$P_{\infty}$. This is a contradiction since
$B\cap H_3=2T_3$. Similarly we have
$l_4\subset H$. Therefore $l_{i}\subset H$ for any $i$, as claimed.  Because
$T_{i}$'s are real, the plane
$\sigma(H)$  contains the tangent lines of
$T_{i}$'s at $\ol{P}_{\infty}$. Let $y_3$ be a linear form on
$\mathbf{CP}^3$ defining $H$, and set
$y_2:=\ol{\sigma^*y_3}$. Then  $(y_0:y_1:y_2:y_3)$ is a homogeneous
coordinate on $\mathbf{CP}^3$. By our choice, we have
$\sigma(y_0:y_1:y_2:y_3)=(\ol{y}_0:\ol{y}_1:\ol{y}_3:\ol{y}_2)$,
$P_{\infty}=(0:0:0:1)$ and $\ol{P}_{\infty}=(0:0:1:0)$.

Because the planes $\{y_i=0\}$ are $U(1)$-invariant, our $U(1)$-action can be linearized with respect to the homogeneous coordinate
$(y_0:y_1:y_2:y_3)$. 
Further, since $H_i$'s are $U(1)$-invariant, the action can be written
$(y_0:y_1:y_2:y_3)\mapsto (y_0:y_1:e^{i\alpha\theta}y_2:
e^{i\beta\theta}y_3)$ for $e^{i\theta}\in U(1)$, where $\alpha$ and $\beta$ are
relatively prime integers.
Moreover, since the conics $T_i$'s are $U(1)$-invariant, we can suppose
$\alpha=1$, $\beta=-1$. Thus the $U(1)$-action can be written in the form
(ii) of the proposition.

Let $F=F(y_0,y_1,y_2,y_3)$ be a defining equation of $B$. Since $B$ is
$U(1)$-invariant,  monomials appeared  in $F$ must be  in the ideal
$(y_2y_3, y_0^2,y_0y_1,y_1^2)^2$.
Moreover, 
$F$ contains the monomial $y_2^2y_3^2$, since otherwise the restriction onto
$\{y_1=0\}$ would be the union of two different conics, 
which contradict to the fact that $T_i$ is a trope.
We assume that its coefficient  is $1$. Then
$F$ can be written in the form $(y_2y_3+Q(y_0,y_1))^2-q(y_0,y_1)$, 
where
$Q(y_0,y_1)\in (y_0,y_1)^2$  and 
$q(y_0,y_1)\in (y_0,y_1)^4$  are uniquely determined polynomials
with real coefficients. Then it again follows from $T_i$ being a trope that 
$q(y_0,y_1)=ky_0y_1(y_0+y_1)(ay_0-by_1)$ for some constant $k\in\mathbf
R^{\times}$. 
If $k$  is negative, exchange $y_0$ and $y_1$. 
Then we get $(y_2y_3+Q(y_0,y_1))^2-ky_0y_1(y_0+y_1)(ay_0-by_1)$ where $k>0,a>0$ and $b>0$.
This can be written
$$
\left( \frac{1}{\sqrt[4]{ka}}y_2\cdot\frac{1}{\sqrt[4]{ka}}y_3+
\frac{1}{\sqrt{ka}}Q(y_0,y_1) \right)^2-
y_0y_1(y_0+y_1)\left(y_0-\frac{b}{a}y_0\right)=0.
$$
Hence if we rewrite $y_2$ and $y_3$ for $(1/\sqrt[4]{ka})y_2$ and $(1/\sqrt[4]{ka})y_3$ respectively, and $Q(y_0,y_1)$ and $a$ for $(1/\sqrt{ka})Q(y_0,y_1)$ and $a$ respectively, we obtain (\ref{eqn-B}) as an defining equation of $B$.
Finally  it is obvious that the forms of the real structure and the $U(1)$-action are invariant under the above change of coordinate.
\proofend

\vspace{3mm}\vspace{3mm}
Note that by a result of Y. Umezu \cite[p.\,141, (1)]{Um84}, every normal surface in $\mathbf{CP}^3$ which is birational to an elliptic ruled surface and which has two simple elliptic singularities of type $\tilde{E}^7$ has  (\ref{eqn-B}) as its defining equation.

\begin{figure}
\includegraphics{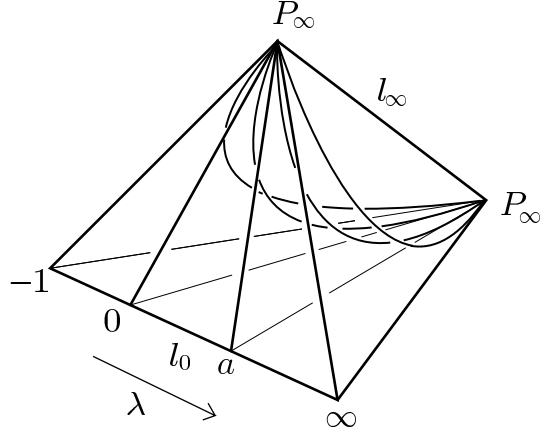}
\caption{four planes containing the tropes of $B$}
\label{fig-trope}
\end{figure}

Next we study the singular locus of $B$.

\begin{prop}\label{sing. of B} Let $B$ be a real quartic surface defined by
the equation
$$\left(y_2y_3+Q(y_0,y_1)\right)^2-y_0y_1(y_0+y_1)(y_0-ay_1)=0$$ where
$Q(y_0,y_1)$ is a real quadratic form of $y_0$ and $y_1$,
and
$a>0$. Let
$A$ be the set
$\{(y_0:y_1:0:0)\set (y_0:y_1)$ is a multiple root of the quartic
equation
$Q(y_0,y_1)^2-y_0y_1(y_0+y_1)(y_0-ay_1)=0$$\}$.  (Here we think
of this as an equation
on
$\mathbf{CP}^1=\{(y_0:y_1)\}$.) Then we have:
(i) {\em Sing}$(B)=\{P_{\infty},\ol{P}_{\infty}\}\cup A$, where we put
$P_{\infty}=(0:0:0:1)$,
(ii) $P_{\infty}$ and $\ol{P}_{\infty}$ are simple elliptic singularities of
type
$\tilde{E}_7$, and
(iii) if $Q(y_0,y_1)\neq 0$, then $(y_0:y_1:0:0)\in A$ is an ordinary double point
iff its multiplicity is two.
\end{prop}
 In particular,  every singular point of $B$ is isolated.

\vspace{3mm}
\noindent Proof. (i) First we show that (Sing$B)\cap\{y_3\neq
0\}=\{P_{\infty}\}$, by calculating the Jacobian. Let $x_i=y_i/y_3$ 
($0\leq i\leq 2$) be affine
coordinates on $y_3\neq 0$. Then the equation of $B$ becomes 
$(x_2+Q(x_0,x_1))^2-x_0x_1(x_0+x_1)(x_0-ax_1)=0$. Differentiating with
respect to 
$x_2$, we get $x_2+Q(x_0,x_1)=0$ so that we have
$x_0x_1(x_0+x_1)(x_0-ax_1)=0$. Next differentiating with respect to
$x_0$ and
$x_1$ and then substituting
$x_2+Q(x_0,x_1)=0$, we get
$x_1(x_0+x_1)(x_0-ax_1)+x_0x_1(x_0-ax_1)+x_0x_1(x_0+x_1)=0$ and
$x_0(x_0+x_1)(x_0-ax_1)+x_0x_1(x_0-ax_1)-ax_0x_1(x_0+x_1)=0$. From the
former, we obtain that  $x_0=0$ implies $x_1=0$. Then by
$x_2+Q(x_0,x_1)=0$ we have $x_2=0$. Similar argument shows that  if
$x_1, x_0+x_1$ or $x_0-ax_1$ is zero, then $x_0=x_1=x_2=0$. Conversely,
it is immediate to see that
$(x_0,x_1,x_2)=(0,0,0)$ is a double point of $B$. Thus we get
(Sing$B)\cap\{y_3\neq 0\}=\{P_{\infty}\}$. Because the given
homogeneous polynomial is symmetric with respect to $y_2$ and $y_3$, we have
(Sing$B)\cap\{y_2\neq 0\}=\{\ol{P}_{\infty}\}$.

Next we show that (Sing$B)\cap \{y_2=y_3=0\}=A$. We may suppose $y_1\neq 0$.
Putting
$v_i=y_i/y_1$ for $i=0,2,3$, the equation  of $B$ becomes
$(v_2v_3+Q(v_0,1))^2-f(v_0)=0$, where we put
$f(v_0)=v_0(v_0+1)(v_0-a)$. Substituting $v_2=v_3=0$, we get
$Q(v_0,1)^2-f(v_0)=0$. On the other hand, differentiating with respect to $v_0$ and
substituting $v_2=v_3=0$, we get $(Q(v_0,1)^2-f(v_0))'=0$, where the
prime denotes differential with respect to $v_0$. 
Thus, if $(v_0,v_2,v_3)=(\lambda_0,0,0)$ is a singular point of $B$,
$\lambda_0$ is a multiple root of $Q(v_0,1)^2-f(v_0)=0$.
 Conversely, it is easy to see that
$(v_0,v_2,v_3)=(\lambda_0,0,0)$ is a singular point for  such $\lambda_0$. Thus we get
the claim of (i).

(ii)  is obvious if one notes that we can  use
$(x_0,x_1,x_2+Q(x_0,x_1))$ instead of $(x_0,x_1,x_2)$ as a local
coordinate around $P_{\infty}$.

Finally we show (iii) by using the coordinate $(v_0,v_2,v_3)$ above. Let
$\lambda_0$ be a multiple root of $Q(v_0,1)^2-f(v_0)=0$. Then our
equation of $B$ can be written
$(v_2v_3+2Q(v_0,1))v_2v_3+g(v_0)(v_0-\lambda_0)^2=0$, where $g(v_0)$
is a  polynomial of degree two.
Clearly $\lambda_0$ is a double root iff $g(\lambda_0)\neq 0$.
Suppose $g(\lambda_0)\neq 0$ and define
\begin{equation}\label{coord}
\begin{array}{l}
w_1=\sqrt{g(v_0)}\cdot(v_0-\lambda_0),\\
w_2=\sqrt{2Q(v_0,1)+v_2v_3}\cdot v_2,\\
w_3= \sqrt{2Q(v_0,1)+v_2v_3}\cdot v_3.
\end{array}
\end{equation}
Because $g(\lambda_0)\neq 0$ and  $Q(\lambda_0,1)\neq0$,
$(w_1,w_2,w_3)$ is a local coordinate around $(\lambda_0,0,0)$.
 Then our equation of $B$ becomes
$w_2w_3+w_1^2=0$. Thus the singularity is an ordinary double point.
Conversely, if $g(\lambda_0)=0$, it is immediate to see that our equation of $B$
can be written of the form $w_2w_3+w_1^3=0$ or $w_2w_3+w_1^4=0$ depending on
whether the multiplicity of $\lambda_0$ is three or four.
This implies that $(\lambda_0,0,0)$ is not an ordinary double point.
\proofend

\begin{prop}\label{prop-pos} Let $B$ be as in Proposition \ref{sing. of
B}. Put $f(\lambda)=\lambda(\lambda+1)(\lambda-a)$. Let
$Z_0\ra\mathbf{CP}^3$ be the double covering branched along $B$. Then if
$Z_0$ admits a small resolution  $Z\ra Z_0$ such that $Z$ is a twistor
space of
$3\mathbf{CP}^2$, then 
$Q(\lambda,1)^2-f(\lambda)\geq 0$ for any $\lambda\in\mathbf R$ and the
equality holds for a unique $\lambda_0\in\mathbf R$. Further, in this case,
the multiplicity of $\lambda_0$ is two.
\end{prop}

\noindent Note that it follows from this proposition that  
$f(\lambda_0)>0$ holds, because we have $Q(\lambda_0,1)^2=f(\lambda_0)$ and
$Q(\lambda_0,1)$ is a real number which is non-zero because otherwise
the restriction $B|_{H_{\lambda_0}}$ would be  $y_0-\lambda_0 y_1=(y_2y_3)^2=0$ that
yields another reducible fundamental divisor.

\vspace{3mm}
\noindent Proof. By  results of Kreu\ss ler \cite{Kr89} and  Kreu\ss
ler-Kurke
\cite{KK92},  we have $\sum(\mu(x)+c(x))=26$ for $Z$ to be a twistor
space of
$3\mathbf{CP}^2$ for a topological reason, where
$\mu(x)$ is the Milnor number of the  singularity $x$ of $B$ and $c(x)$
is the number of irreducible components of a small resolution $Z\ra
Z_0$.  Because elliptic singularity of type
$\tilde{E}_7$ has $\mu=9$ and $c=3$, we get $\sum(\mu(x)+c(x))=2$ for 
other remaining singularities. This implies that there is only one
singularity remaining, and that it must be an ordinary double  point, 
which will be denoted by $P_0$.
Therefore, by
Proposition
\ref{sing. of B} (i), we have
$A=\{P_0\}$. Namely, $Q^2-f=0$ has a unique  multiple root $\lambda_0$.
 The multiplicity
is two by Proposition \ref{sing. of B} (iii).  It is obvious from the uniqueness that
this ordinary double point is  real. Namely, $\lambda_0$ is real.

Next we show that other  solutions of $Q(\lambda,1)^2-f(\lambda)=0$
are not real.  
Assume $\lambda\in\mathbf R$, $\lambda\neq\lambda_0$ is a
solution. Then by restricting $B$ to the plane $y_0=\lambda y_1$, we get
$(y_2y_3+Q(\lambda,1)y_1^2)^2-f(\lambda)y_1^4=
(y_2y_3)^2+2Q(\lambda,1)y_2y_3y_1^2+(Q(\lambda,1)^2-f(\lambda))y_1^4=
y_2y_3(y_2y_3+2Q(\lambda,1)y_1^2)$ $(=0)$. 
Therefore, the point
$(\lambda:1:0:0)$ is a real point of $B$. 
Since the multiplicity of the solution $\lambda$ is one, 
 Proposition \ref{sing. of B} shows that this is a smooth point of
$B$. 
This implies that $Z$ has a real point, contradicting to the
absence of real points on any twistor spaces. Hence the equation
$Q(\lambda,1)^2-f(\lambda)=0$ has no real solution other than
$\lambda_0$. Because $\lambda_0$ is a solution whose multiplicity is
 two, this implies that the polynomial
$Q(\lambda,1)^2-f(\lambda)$ has constant sign on
$\mathbf R\backslash \{\lambda_0\}$.  This sign must be clearly positive.
\proofend

\vspace{3mm}
To investigate the real locus of $B$,  we need the
following elementary 
\begin{lemma}\label{lemma-element} Let $C_{\alpha}=\{y_2y_3=\alpha
y_1^2\}$,
$\alpha\in\mathbf R$ be a real conic in
$\mathbf{CP}^2$, where the real structure is given by
$(y_1:y_2:y_3)\mapsto (\ol{y}_1:\ol{y}_3:\ol{y}_2).$ Then $C_{\alpha}$
has no real point iff
$\alpha<0$.
\end{lemma}

\noindent Proof.
It is immediate to see that the real locus of $C_{\alpha}$ is 
$$\left\{(1:v:\ol{v})\in\mathbf{CP}^2\set |v|=\sqrt{\alpha}\right\}.$$
This is empty iff $\alpha<0$.
\endproof

\begin{prop}\label{prop-realpoint} Let $B$ be as in Proposition
\ref{sing. of B} and suppose that the inequality
$Q(\lambda,1)^2-f(\lambda)\geq 0$ holds on $\mathbf R$ with the
equality holding iff
$\lambda=\lambda_0$ as in Proposition
\ref{prop-pos}. Put $P_0:=(\lambda_0:1:0:0)$, which is clearly a real
point of $B$. Then we have: (i) there is no real point on $B$ other than 
$P_0$ iff the following condition is satisfied:
if $f(\lambda)\geq 0$ and $\lambda\neq\lambda_0$, then
$Q(\lambda,1)>\sqrt{f(\lambda)}$
$\cdots (*)$,
  (ii) if $ (*)$ is satisfied, then there is no real point on any
small resolutions of
$Z_0$.

\end{prop}
\noindent

\vspace{3mm}
\noindent Proof. It is immediate to see that 
 any real point of $B$ is contained in some
real plane
$H_{\lambda}:=\{y_0=\lambda y_1\}$, $\lambda\in\mathbf R\cup\{\infty\}$.
An equation of the restriction
$B_{\lambda}:=B\cap H_{\lambda}$ is given by (as in the proof of
Proposition
\ref{prop-pos}) 
$\left(y_2y_3+Q(\lambda,1)y_1^2\right)^2-f(\lambda)y_1^4=0$. This can be rewritten as
$$B_{\lambda}:
\left\{y_2y_3+\left(Q(\lambda,1)-\sqrt{f(\lambda)}\right)y_1^2\right\}
\left\{y_2y_3+\left(Q(\lambda,1)+\sqrt{f(\lambda)}\right)y_1^2\right\}=0.
$$ Namely, $B_{\lambda}$ is a union of two conics.
If $\lambda\neq\lambda_0$, we have $Q(\lambda,1)^2-f(\lambda)> 0$ by 
our assumption, so both of 
$Q(\lambda,1)-\sqrt{f(\lambda)}$
and $Q(\lambda,1)+\sqrt{f(\lambda)}$ are non-zero.

  Recall that our real structure is given by
$\sigma(y_1:y_2:y_3)=(\ol{y}_1:\ol{y}_3:\ol{y}_2)$ on $H_{\lambda}$
(Proposition \ref{prop-def-B}, (ii)).
Thus each component of $B_{\lambda}$ is real iff the coefficients are
real; namely
$f(\lambda)\geq 0$. Further, the intersection of these two conics are
$\{P_{\infty},\ol{P}_{\infty}\}$. Therefore, there is no real point on
$B_{\lambda}$ if $f(\lambda)<0$. So suppose $f(\lambda)\geq 0$. In this
case, each of the two conics are real, and by Lemma \ref{lemma-element},
both components have no real point iff $Q(\lambda,1)>\sqrt{f(\lambda)}$.
On the other hand, we have $B_{\infty}=\{(y_2y_3+Q(y_0,0))^2=0\}$.
Hence again by Lemma \ref{lemma-element}, we have $Q(1,0)>0$ if
$B_{\infty}$ has no real point. But this   follows from the first
condition.  If $\lambda=\lambda_0$, we have 
$Q(\lambda_0,1)=\sqrt{f(\lambda_0)}$ and hence one of the
components of
$B_{\lambda_0}$ degenerates into a union of two lines whose intersection is $P_0$.
And the other component has no real point since
$Q(\lambda_0,1)+\sqrt{f(\lambda_0)}=2\sqrt{f(\lambda_0)}>0$. Thus we get (i).
 
Next we show that $Z_0$ has no real point other than $p_0:=\Phi_0^{-1}(P_0)$,
under the condition
$(\ast)$. Suppose $y_1\neq 0$ and use the coordinate $(v_0,v_2,v_3)$
defined in the proof of  Proposition \ref{sing. of B}. Then  $Z_0$ is
given by the equation
$z^2+(v_2v_3+Q(v_0,1))^2-f(v_0)=0$, where $z$ is a fiber coordinate on
the line bundle $ O(2)$ over $\mathbf{CP}^3$.  Thus to prove that 
the point $p_0$ is the only real locus of
$Z_0\cap\Phi_0^{-1}(\{y_1\neq 0\})$,  it suffices to show
that 
\begin{equation}\label{ineq} (v_2v_3+Q(v_0,1))^2-f(v_0)>0
\end{equation} for any real $(v_0, v_2,v_3)\neq (\lambda_0,0,0)$. Recall
that
$(v_0, v_2,v_3)$ is real iff
$v_0\in\mathbf R$ and
$v_2=\ol{v}_3$. Hence (\ref{ineq}) is obvious for real $(v_0, v_2,v_3)$ 
with
$f(v_0)<0$. Assume $f(v_0)\geq 0$. Using the reality condition, we have
$(v_2v_3+Q(v_0,1))^2-f(v_0)=|v_2|^4+2Q(v_0,1)|v_2|^2+
(Q(v_0,1)^2-f(v_0))$.
By our assumption we have $Q(v_0,1)^2-f(v_0)> 0$ for any 
$v_0\in\mathbf R$ with $v_0\neq \lambda_0$. Further, by the condition
$(\ast)$  we have $Q(v_0,1)>\sqrt{f(v_0)}$ for $v_0\neq \lambda_0$
with $f(v_0)\geq 0$. Therefore we have (\ref{ineq}) also for  real
$(v_0, v_2,v_3)\neq (\lambda_0,0,0)$ with $f(v_0)\geq 0$. Thus we get that
on $\Phi_0^{-1}(\{y_1\neq 0\})$, there are no real locus other than $p_0$. 
In the same way we can see that, over $y_0\neq 0$, there are no real points
other than $p_0$.
So it remains to show that there is no real point over the line
$l_{\infty}=\{y_0=y_1=0\}$.
To check this, we introduce a new homogeneous coordinate
$(y_0,y_1,y_2-y_3,y_2+y_3)$ on $\mathbf{CP}^3$.
Then the two subsets $y_2-y_3\neq 0$ and $y_2+y_3\neq 0$ are real,
and it can be easily seen that over the line $l_{\infty}$,
 the equation of $Z_0$ is
of the form $z^2+q^2=0$ on these two  open subset, where $q$ is  non-zero
real valued on the real set $l_{\infty}^{\sigma}$.
 Hence $Z_0$ does not have real
point over
$l_{\infty}$. Thus $p_0$ is the unique real point on the whole of $Z_0$.

Finally we show that $\Gamma_0=\Phi^{-1}(P_0)$ has no real point. 
To show this, we use a coordinate
$(w_0,w_2,w_3)$ (around $P_0$) defined in (\ref{coord}). Then $B$
is  given by  $w_1^2+w_2w_3=0$. Further, it is easy to
see that  the real structure is also given by
$\sigma(w_1,w_2,w_3)=(\ol{w}_1,\ol{w}_3,\ol{w}_2)$.  
Now because $\Gamma_0$ is the exceptional curve of a small resolution of an ordinary double point, $\Gamma_0$ can be canonically identified with the set of lines contained in the cone $w_1^2+w_2w_3=0$. 
If $\{(w_1:w_2:w_3)=(a_1:a_2:a_3)\}$ is a real line, we can suppose
$a_1\in\mathbf R$, $a_3=\ol{a}_2$.
It follows that it cannot be 
contained in the cone. This implies that $\Gamma_0$ has no real point.
On the other hand, resolutions of the singularities over $P_{\infty}$ and 
$\ol{P}_{\infty}$ do not yield real points.
Thus we can conclude that $Z$ has no real point.
 \proofend

\vspace{3mm}
Here we summarize necessary conditions for our threefolds to be (birational to) a 
twistor space:

\begin{prop}\label{prop-necessa}
Let $B$ be a quartic surface in $\mathbf{CP}^3$ defined by
$$
\left(y_2y_3+Q(y_0,y_1)\right)^2-y_0y_1(y_0+y_1)(y_0-ay_1)=0,$$
where $Q(y_0,y_1)$ is a quadratic form of $y_0$ and $y_1$ with real coefficients, and  $a>0$.
Let $Z_0\ra \mathbf {CP}^3$ be the double covering branched along $B$.
If there exists a resolution $Z\ra Z_0$ such that $Z$ is a twistor space,  then  $Q$ and $a$  satisfy the following condition:

\noindent
{\textbf{Condition (A)}} If $\lambda$ satisfies $\lambda(\lambda+1)(\lambda-a)\geq 0$, then 
 \begin{equation}\label{eqn-cond}
 Q(\lambda,1)\geq\sqrt{\lambda(\lambda+1)(\lambda-a)}
 \end{equation}
holds. 
Moreover, there exists a unique $\lambda_0$ such that  the equality of (\ref{eqn-cond}) holds for $\lambda=\lambda_0$.

\end{prop}

Proposition \ref{prop-def-B}
(combined with some calculations given in \S
\ref{ss-special})  has the following consequence:

\begin{prop}\label{prop-class} Let $g$ be a self-dual metric on
$3\mathbf{CP}^2$ of positive scalar curvature with a non-trivial Killing
field, and assume that 
$g$ is not conformally isometric to LeBrun metric. Then the naturally
induced
$U(1)$-action   on
$3\mathbf{CP}^2$ is uniquely determined up to diffeomorphisms.
\end{prop}

\noindent Proof. Let $Z$ be the twistor space of $g$. Then $Z$ is as in
Proposition
\ref{prop-def-B}. Let $H_i$ and $T_i$ ($1\leq i\leq 4$) be  as in the
proof of Proposition
\ref{prop-def-B}. Namely, $H_i$ is a real $U(1)$-invariant plane
such that 
$B|_{H_i}$ is a trope  whose reduction is denoted by
$T_i$  as before.  Then $\Phi^{-1}_0(H_i)$ consists of two irreducible components,
both of which are biholomorphic to $H_i$ ($\simeq\mathbf{CP}^2$). Since
$\mu:Z\ra Z_0$ is small,
$\Phi^{-1}(H_i)$ also consists of two irreducible components, which are
denoted by
$D_i$ and $\ol{D}_i$. These are $U(1)$-invariant.
Then as will be proved in Lemma \ref{lemma-promise-1}, there is 
a smooth rational surface $D$ with $U(1)$-action such that
$D$ is $U(1)$-equivariantly biholomorphic to $D_i$ for some $1\leq i\leq 4$.
%
%
Since
$D_i+\ol{D}_i$ is a fundamental divisor, and since we  have
$-(1/2)K_Z\cdot L=2$ ($L$ is a twistor line), we have
$D_i\cdot L=\ol{D}_i\cdot L=1$. Hence by a result of Poon
\cite{P92},  $L_i:=D_i\cap \ol{D}_i$ is a twistor line which is
obviously $U(1)$-invariant, and that $L_i$ is contracted to a
point by the twistor fibration $Z\ra 3\mathbf{CP}^2$ which is
$U(1)$-equivariant.
Hence the $U(1)$-action on $3\mathbf{CP}^2$ can be read  from that on
any one of 
$D_i$'s. Therefore the conclusion of the proposition follows.
\proofend

\vspace{3mm} The proposition implies that, up to
diffeomorphisms,  there are only two effective $U(1)$-actions on
$3\mathbf{CP}^2$ which can be the identity component of the isometry
group of a self-dual metric whose scalar curvature is positive. One is
the semi-free $U(1)$-action, which is the identity component of generic
LeBrun metric, and the other is the action  obtained in
Proposition \ref{prop-class}. Of course, there are many other
differentiable $U(1)$-actions on $3\mathbf{CP}^2$ in general: for example, we can get an infinite
number of mutually inequivariant $U(1)$-actions  by first taking an
effective
$U(1)\times U(1)$-action on $3\mathbf{CP}^2$ and then choosing 
many different
$U(1)$-subgroup of $U(1)\times U(1)$.

In the rest of this section we prove propositions which will be used in the subsequent sections.

\begin{prop}\label{prop-pic}
Let $B$ be the quartic surface defined by the equation (\ref{eqn-B}) and suppose 
that $Q$ and $f$ satisfy the assumption in Proposition \ref{prop-realpoint}.
Let $\Phi_0:Z_0\ra\mathbf{CP}^3$ be the double covering branched along $B$, and
$\mu:Z\ra Z_0$ any small resolution preserving the real structure.
 Put $\Phi=\mu\cdot\Phi_0$. 
 Then we have (i) $K_Z\simeq \Phi^*O(-2)$, (ii) the line bundle $(-1/2)K_Z$ is 
uniquely determined, (iii) $\Phi$ is the induced morphism associated to the complete linear system $|(-1/2)K_Z|$.
\end{prop}

\noindent Proof.
Let $K_{Z_0}$ denote the canonical sheaf of $Z_0$. 
Then we have $K_{Z_0}\simeq
\Phi_0^*(K_{\mathbf{CP}^3}+(1/2)O(B))\simeq\Phi_0^*O(-2)$. Moreover,
since $\mu$ is small, we have $K_Z\simeq \mu^* K_{Z_0}$. Hence we get
$K_Z\simeq \Phi^*O(-2)$. 
For (ii) it suffices to show that $H^1(O_Z)=0$.
Since the singularities of $Z_0$ are normal, and since the exceptional curves
of $\mu$ are rational, 
we get by Leray spectral sequence $H^1(O_Z)\simeq H^1(O_{Z_0})$.
Then applying the spectral sequence to $\Phi_0$, and using $\Phi_{0*}O_{Z_0}
\simeq O\oplus O(-2)$ and $R^q\Phi_{0*}O_{Z_0}=0$ for $q\geq 1$, we get
$H^1(O_{Z_{0}})=0$ and we get (ii). 
(iii) is easily obtained from (i).
\proofend
 
 \vspace{3mm}
  The following result will be also needed in the next section:
 
 \begin{prop}\label{prop-noline2}
 Let $B$ be a quartic surface defined by (\ref{eqn-B}). 
 Then $B$ does not contain real lines.
 \end{prop}
 
 \noindent Proof.
 Suppose that $l$ is a real line lying on $B$.
 Because  $B$ is $\mathbf C^*$-invariant by the action $\rho_t:(y_0,y_1,y_2,y_3)\mapsto(y_0,y_1,ty_2,t^{-1}y_3)$ for $t\in\mathbf C^*$ (which is the complexification of the $U(1)$-action  in Proposition \ref{prop-def-B} (iii)),  $\rho_t(l)$ is a line contained in $B$.
 Hence if $l$ is not $\mathbf C^*$-invariant, then it yields a one-parameter family of rational curves in $B$. 
 This implies that $B$ is a rational surface.
 However, $B$ is birational to an elliptic ruled surface \cite{Um84, Hon03}.
 Therefore $l$ must be $\mathbf C^*$-invariant.
 It is readily verified from the above explicit form of the $\mathbf C^*$-action that $\mathbf C^*$-invariant real lines on $\mathbf{CP}^3$ are $l_0=\{y_2=y_3=0\}$ and $l_{\infty}=\{y_0=y_1=0\}$ only.
 Both of these lines are clearly not contained in $B$.
 Hence there is no real line contained in $B$.
 \proofend
 
\begin{figure}[htbp]
\includegraphics{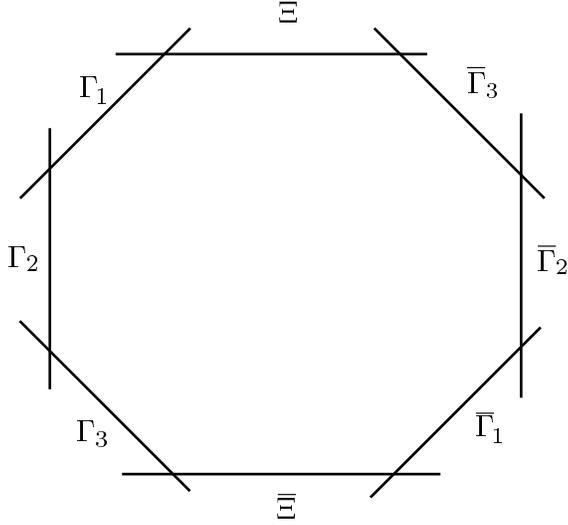}
\caption{$\Phi^{-1}(l_{\infty})$ becomes a cycle of 8 smooth rational curves}
\label{fig-cycle1}
\end{figure}

 \vspace{3mm}
 The following `octagon' will play a significant role throughout our investigation:
 
  \begin{prop}\label{prop-octagon}
 Let $B$, $Z_0$, $\Phi_0$, $\mu$ and $\Phi$ be as in Proposition \ref{prop-pic}. 
 Let  $l_{\infty}$ be a real line defined by $y_0=y_1=0$.  
 Then $\Phi^{-1}(l_{\infty})$ is a cycle of eight smooth rational curves 
 as in Figure \ref{fig-cycle1}, where  $\Gamma=\Gamma_1\cup\Gamma_2\cup\Gamma_3$ and  $\ol{\Gamma}=\ol{\Gamma}_1\cup\ol{\Gamma}_2\cup\ol{\Gamma}_3$ are the exceptional curves of the conjugate pair of  singularities of $Z_0$, and $\Xi$ and $\ol{\Xi}$ are conjugate pair of curves which are mapped biholomorphically onto $l_{\infty}$. 
 \end{prop}
 
 \noindent Proof.
 Let $P_{\infty}$ and $\ol{P}_{\infty}$ be the conjugate pair of elliptic singularities of $B$  as before and put $p_{\infty}=\Phi_0^{-1}(P_{\infty})$ and $\ol{p}_{\infty}=\Phi_0^{-1}(\ol{P}_{\infty})$ be the corresponding points of $Z_0$.
Then since $p_{\infty}$ and $\ol{p}_{\infty}$ are compound $A_3$-singularities of $Z_0$,  the exceptional curve of a small resolution is a chain of three smooth rational curves.
(See Section \ref{ss-special}, where small resolutions are explicitly given.
Another way to see this is that, cutting $B$ by a generic plane containing $l_{\infty}$ and taking the inverse image by $\Phi_0$, one obtains a surface which has $p_{\infty}$ and $\ol{p}_{\infty}$ as  $A_3$-singularities.
Then any small resolution of $p_{\infty}$ and $\ol{p}_{\infty}$ necessarily gives the minimal resolution of the surface, by a well-known property of rational double points; cf. Proof of Lemma \ref{lemma-degbitan3}.)
It remains to see how about the remaining components.
It is readily seen that  $\Phi_0^{-1}(l_{\infty})$ consists of two smooth rational curves both of which are mapped biholomorphically onto $l_{\infty}$.
Moreover, as seen in the proof of Proposition \ref{prop-realpoint}, there is no real point on $\Phi_0^{-1}(l_{\infty})$. 
Therefore the two components of $\Phi_0^{-1}(l_{\infty})$ are conjugate of each other.
Strictly speaking, we have to verify that the strict transforms in $Z$ of these two components intersect $\Gamma_1$ and $\Gamma_3$ (and $\ol{\Gamma}_1$ and $\ol{\Gamma}_3$) as in Figure \ref{fig-cycle1}, and do not intersect $\Gamma_2$ and $\ol{\Gamma}_2$.
But this can be checked by explicitly giving small resolutions as in Section \ref{ss-special} and chasing the inverse images.
(Indeed this is essentially done in the proof of Lemma \ref{lemma-intersection3}.)
Thus we get the claim of the proposition.
\proofend

 \section{How to find twistor lines}\label{s-how}
 In the previous section we obtained necessary conditions (Proposition \ref{prop-necessa}) of a quartic surface  $B$ to be the branch divisor of a twistor space of non-LeBrun self-dual metric on $3\mathbf{CP}^2$ with a non-trivial Killing field. 
 The goal of this paper is to show that these are also  sufficient conditions. 
 We solve this problem by finding the family of twistor lines. 
 In this section we will describe how to find twistor lines. 
 The basic role is played by a classical result in algebraic geometry   that  a smooth plane quartic always has twenty-eight bitangents. 
 
 Let $B$ be a quartic surface defined by the equation
\begin{equation}\label{eqn-B2}
\left(y_2y_3+Q(y_0,y_1)\right)^2-y_0y_1(y_0+y_1)(y_0-ay_1)=0
\end{equation}
 and suppose that Conditions (A) for $Q$ and $a$ of Proposition \ref{prop-necessa} is satisfied. 
 Let $P_0=(\lambda_0,1,0,0)$, $P_{\infty}=(0,0,0,1)$ and $\ol{P}_{\infty}=(0,0,1,0)$ be the singular locus of $B$, where $P_0$ is a real ordinary double point, and $P_{\infty}$ and $\ol{P}_{\infty}$ are simple elliptic singularities of type $\tilde{E}_7$ (Proposition \ref{sing. of B}).
 Let $\Phi_0:Z_0\ra \mathbf{CP}^3$ be the double covering branched along $B$.
 Let $\mu:Z\ra Z_0$ be any one of small resolution of the three singular points of $B$, preserving the real structure of $Z_0$. 
 (There are many choices of such $\mu$. 
 Later (in Sections \ref{s-nb} and \ref{s-line}) we will determine which resolution must be taken.
 Here we only suppose that $\mu$ preserves the real structure.)
 We do not assume that $Z$ is a twistor space.
 Further we put $\Phi=\Phi_0\cdot\mu$, and $\Gamma_0=\Phi^{-1}(P_0),\Gamma=\Phi^{-1}(P_{\infty})$ and $\ol{\Gamma}=\Phi^{-1}(\ol{P}_{\infty})$ which are  the exceptional curves of the small resolution of $Z_0$.
 
 Set $H_1=\{y_0+y_1=0\}, H_2=\{y_1=0\}, H_3=\{y_0-ay_1=0\}$ and $H_4=\{y_0=0\}$. These are planes such that $H_i\cap B$ is a double conic (= trope), which is denoted by $T_i\, (\subset H_i)$. Every $H_i$ and  $T_i$ are real, where the real structure on $\mathbf{CP}^3$ is defined by 
 \begin{equation}\label{eqn-rs2}
 \sigma(y_0:y_1:y_2:y_3)=(\ol{y}_0:\ol{y}_1:\ol{y}_3:\ol{y}_2)
\end{equation}
 (cf. Proposition \ref{prop-def-B}).
 Further, we write $\Phi^{-1}(H_i)=D_i+\ol{D}_i$ for $1\leq i\leq 4$.
 
 \begin{definition}\label{def-tc}{\em{ An irreducible conic $C$ in $\mathbf{CP}^3$ is
called a {\it{touching conic}} of $B$ if $C\subset B$ or otherwise if the
intersection number with
$B$ is even at any intersection points.}}
\end{definition}

The following proposition motivates the definition:

\begin{prop}\label{prop-image} Let $L$ be a
real line in $Z$.
 Then $\Phi(L)$ is a line in $\mathbf{CP}^3$ iff
$L\cap \Gamma_0\neq\phi$. Otherwise $\Phi(L)$ is a real touching conic of $B$.
\end{prop}

Note that any irreducible conic in $\mathbf{CP}^3$ is contained in a unique plane, and that the plane is real if the conic is real. 
In particular, it follows from Proposition \ref{prop-image}  that  any real line in $Z$ always lies on the inverse image of some real plane.

\vspace{3mm}
\noindent Proof of Proposition \ref{prop-image}.
 By adjunction formula, we have $-2=K_Z\cdot L+\deg N_{L/Z}=
K_Z\cdot L+2$. 
Hence by Proposition \ref{prop-pic} we have $(-1/2)K_Z\cdot L=2$. 
Therefore $\Phi(L)$ is a curve whose  degree is at most two,
and $\Phi(L)$ is a line iff $\Phi|_L:L\ra\Phi(L)$ is two-to-one. 

Assume that $\Phi(L)$ is a line, which is necessarily real. 
By Proposition \ref{prop-noline2}, $\Phi(L)$ is not contained in $B$.
To prove $P_0\in\Phi(L)$, it suffices to show that $\Phi(L)\cap B$ consists of three or  one point,
since by Proposition \ref{prop-realpoint}, $P_0$ is the unique real point of $B$.
Since $B$ is a quartic, $\Phi(L)\cap B$ consists of at most four points.
If $\Phi(L)\cap B$ consists of four points, the intersection points are smooth points of $B$ and $\Phi^{-1}(\Phi(L))$ must be a smooth elliptic curve, contradicting our assumption.
If $\Phi(L)\cap B$ consists of two points, the intersection points must be singular points of $B$, since otherwise $\Phi^{-1}(\Phi(L))$ would split into two rational curves, which contradicts that $\Phi|_L:L\ra\Phi(L)$ is two to one.
Therefore the two intersection points must be singular points of $B$ and it follows $\Phi(L)=l_{\infty}$.
However, by Proposition \ref{prop-octagon}, $\Phi^{-1}(l_{\infty})$ does not contain real components. 
(Or more strongly any component of $\Phi^{-1}(\Phi(L))$ does not mapped two-to-one  onto its image.)
Hence $\Phi(L)\cap B$ cannot consist of two points.
Therefore it must consists of three or one point
(although the latter cannot happen as is readily seen), and 
 we have $P_0\in\Phi(L)$.
It follows that $L\cap \Gamma_0\neq\phi$.
Conversely assume that $L$ is a real line intersecting $\Gamma_0$.
Then since there are no real points on $\Gamma_0$
(Proposition \ref{prop-realpoint} (ii)), 
the intersection is not one point. Because
$\Phi(\Gamma_0)=P_0$, this implies that $\Phi$ is not one-to-one on $L$.
Hence $\Phi(L)$ must be a line.

Finally suppose that $\Phi(L)$ is a conic. If $(\Phi(L), B)_P=1$ for
some $P\in \Phi(L)\cap B$, then $P$ is a smooth point of $B$
and  the intersection is transversal. Therefore
$\Phi^{-1}(\Phi(L))$ is locally irreducible near $\Phi^{-1}(P)$.
This contradicts to the fact that $\Phi|_L$ is bijective.
Therefore we have $(\Phi(L), B)_P\geq 2$ for any $P\in \Phi(L)\cap B$.
Moreover, $(\Phi(L), B)_P$ must be even, since otherwise $\Phi^{-1}(\Phi(L))$ would be locally irreducible. This implies that
$\Phi(L)$ is a touching conic of $B$, which is necessarily real.
\hfill $\square$

\vspace{3mm}
In the sequel $(\mathbf{CP}^3)^{\vee}$ denotes the dual complex projective space, which can be viewed as the set of planes in $\mathbf{CP}^3$.
Let  $(\mathbf{RP}^3)^{\vee}\subset(\mathbf{CP}^3)^{\vee}$ be the dual projective space of $\mathbf{RP}^3$,  which is the set of  real planes in $\mathbf{CP}^3$, where we are assuming that the real structure on $\mathbf{CP}^3$ is given by (\ref{eqn-rs2}).
 According to Proposition \ref{prop-pic}, for $H\in(\mathbf{RP}^3)^{\vee}$, $S_H:=\Phi^{-1}(H)$ is a real member of $|(-1/2)K_Z|$. 
 If $Z$ is actually a twistor space, the following basic result of Pedersen-Poon  \cite{PP94} concerning the structure of real members of this complete system  is known:

\begin{prop}\label{prop-PP}
Any real irreducible member $S$ of $|(-1/2)K_Z|$ of a compact twistor space is a smooth surface. 
Moreover, such an $S$ always contains a real pencil whose real members (which are automatically parametrized by a circle) are just the set of twistor lines of $Z$ contained in $S$. 
Further, each member has the trivial normal bundle in $S$.
\end{prop}
 
This provides  a necessary condition  for a threefold to be a twistor space. 
However,  for a general twistor space, it is not true that  generic twistor lines are contained in some $S\in|(-1/2)K|^{\sigma}$.
An important point in the present case is that this is true, due to the remark after Proposition \ref{prop-image}.
We will see in Proposition \ref{prop-str_of_S} that $S_H=\Phi^{-1}(H)$ is smooth for any $H\in (\mathbf{RP}^3)^{\vee}$, unless  $\Phi^{-1}(H)$ is reducible (i.e. $H\neq H_i,\, 1\leq i\leq 4$).

By using the $\mathbf C^*$-action, we can determine real planes $H$ for which $B\cap H$ is singular:

\begin{lemma}\label{lemma-singsection}
For a real plane $H\in (\mathbf{RP}^3)^{\vee}$, $B\cap H$ is a singular quartic precisely when $H$ goes through a singular point of $B$. 
Moreover, any singular point  of $B\cap H$ must be a singular point of $B$.
\end{lemma}
 
 \noindent Proof. 
 Recall that $B$ is invariant by the $\mathbf C^*$-action $(y_0,y_1,y_2,y_3)\mapsto (y_0,y_1,ty_2,t^{-1}y_3)$, $t\in\mathbf{C}^*$, and that $B$ has only isolated singular points, which are $P_0, P_{\infty}$ and $\ol{P}_{\infty}$.
It follows that any singular point of  $H\cap B$ is a $\mathbf C^*$-fixed point, since otherwise $B$ would have singularities along the $\mathbf C^*$-orbit. 
It is immediate to see that the $\mathbf C^*$-fixed set on $\mathbf{CP}^3$ is  a line $l_0=\{y_2=y_3=0\}$ and two points $P_{\infty}$ and $\ol{P}_{\infty}$. 
But it is readily seen  that  $B\cap l_0$ consists of three points, one of which is $P_0$ and the other two are a pair of conjugate points (cf. Proof of Proposition \ref{prop-imline}). 
It follows from the reality of $H$ that if $H$ contains one of the conjugate pair of points, then  $H$ contains $l_0$ and in particular contains $P_0$. 
Thus we have seen that  $H\cap B$ is singular precisely when $H$ goes through $P_0, P_{\infty}$ or $\ol{P}_{\infty}$ and we obtain the former claim of the lemma.
Moreover, the conjugate pair of points of $B\cap l_0$ cannot be singular points of $B\cap H$ ($H\supset l_0$), since $B$ is a quartic and $P_0$ is a double point.
Thus any singular point of $H\cap B$ must be  $P_0,P_{\infty}$ or $\ol{P}_{\infty}$ and we get the latter claim of the lemma.\proofend

\vspace{3mm}
Next we introduce some notations that will be frequently used in the sequel. 
We set $U:=\{H\in (\mathbf{CP}^3)^{\vee}\set H\cap B$ is a non-singular quartic$\}$, which is a Zariski-open subset of $(\mathbf{CP}^3)^{\vee}$. Set $U^{\sigma}=U\cap(\mathbf{RP}^3)^{\vee}$.
By Proposition \ref{lemma-singsection}, $U^{\sigma}$ is the set of planes not going through $P_0,P_{\infty}$ and $\ol{P}_{\infty}$. 
From this, it is immediate to see that 
$(\mathbf{RP}^3)^{\vee}\backslash U^{\sigma}$ consists of two components: 
one is the set of real planes going through $P_0$ and another is the set of real planes containing $l_{\infty}=\{y_0=y_1=0\}$. 
We denote by $\mathbf{RP}^2_{\infty}$ for the former component and 
$\langle l_{\infty}\rangle^{\sigma}$ for the latter component respectively. 
Note that  $\mathbf{RP}^2_{\infty}\cap\linfty^{\sigma}$ is a single plane, which we will denote by $H_{\lambda_0}$ (since it is defined by $y_0=\lambda_0y_1$, where $\lambda_0$ appears in Condition (A) of Proposition \ref{prop-necessa}) .
Members of $\linfty^{\sigma}$ can be characterized by the $\mathbf C^*$-invariance, where the $\mathbf C^*$-action is as in the proof of Lemma \ref{lemma-singsection}.  
$\mathbf{RP}^2_{\infty}$ can be considered as a plane at infinity on $(\mathbf{RP}^3)^{\vee}$, and by removing it, we get an Euclidean space $\mathbf{R}^3$.
Then $\langle l_{\infty}\rangle^{\sigma}$ is a line in this $\mathbf R^3$, and $U^{\sigma}$ is isomorphic to $\mathbf R^3\backslash\mathbf R\simeq \mathbf C^*\times\mathbf R$. 
In particular we have $\pi_1(U^{\sigma})\simeq\mathbf Z$.
This fact will play an essential role in our global construction of twistor lines given in  Section \ref{s-disj}.

The following proposition implies that our threefold $Z$ actually satisfies a necessary condition to be a twistor space imposed by Proposition \ref{prop-PP}. 

\begin{prop}\label{prop-str_of_S}
If $H$ is a real plane different from $H_i$ ($1\leq i\leq 4$), then $S_H=\Phi^{-1}(H)$ is a smooth rational surface satisfying $c_1^2(S_H)=2$.
\end{prop}

\noindent Proof.
Recall that $\Phi_0:Z_0\ra\mathbf{CP}^3$ is the double covering branched along $B$ and $\Phi:Z\ra\mathbf{CP}^3$ is the composition of $\Phi_0$ with a small resolution of $Z_0$. 
Thus if $H$ is a real plane not going the singular points of $B$, $S_H\ra H$ is simply a double covering   whose branch is $B\cap H$. 
By Lemma \ref{lemma-singsection}, this kind of $H$ is exactly real planes such that $B\cap H$ is a smooth quartic; namely $H\in U^{\sigma}$. 
Hence $S_H\ra H$ is a double covering branched along a smooth quartic. 
It follows that $S_H$ is a smooth rational surface with $c_1^2=2$
(which is a del-Pezzo surface). 
Thus the claim is proved for $H\in U^{\sigma}$. 
Next suppose that $H\in \mathbf{RP}_{\infty}^2$ and $H\neq H_{\lambda_0}$.
Then $H\cap B$ is a quartic curve whose singular points is $P_0$ only by Lemma \ref{lemma-singsection}. 
 $\Phi_0^{-1}(H)$ has a unique ordinary point over $P_0$. 
But this is always resolved through any small resolution $Z\ra Z_0$.
Hence $S_H$ is again smooth. 
Moreover, $S_H$ must be a rational surface with $c_1^2=2$, since by moving the plane $H$ to be a member of $U^{\sigma}$,  $S_H$ can be considered as a deformation of $S_{H'}$ for $H'\in U^{\sigma}$ and since $S_{H'}$ is a rational surface with $c_1^2=2$ as is already mentioned. 
(However, this time $S_H$ is not a del-Pezzo surface any more.) 
Thus the claim of the proposition is also proved for $H\in \mathbf{RP}^2_{\infty}$.
 Finally suppose that $H\in\linfty^{\sigma}$. 
 Since $H$ is $\mathbf C^*$-invariant, $H\cap B$ must be also $\mathbf C^*$-invariant.
  It is readily seen that $H\cap B$ is a union of two $\mathbf C^*$-invariant conics going through $P_{\infty}$ and $\ol{P}_{\infty}$, and the two conics are mutually different iff $H\neq H_i$ ($1\leq i\leq 4$). 
  It follows that $\Phi_0^{-1}(H)$ has $A_3$-singularities over $P_{\infty}$ and $\ol{P}_{\infty}$. 
  These singularities are also resolved by any small resolution of the conjugate pair of singular points of  $Z_0$, and therefore $S_H$ is smooth. 
  (These are proved Lemma \ref{lemma-promise-2} by direct calculation using local coordinates.) 
  The structure of $S_H$ is again as in the claim of the proposition by the same reason for the previous case, and we have obtained all of the claims.\proofend

\vspace{3mm}
Our next task is to find twistor lines lying on $S_H$.  
Again by Proposition \ref{prop-PP}  twistor lines on $S_H$ must be the real part of a real pencil on $S_H$ whose corresponding line bundle $\mathcal L_H\ra S_H$  satisfies $(\mathcal L_H)^2_S=0$. 
There are many real line bundles $\mathcal L$ on $S$ satisfying $(\mathcal L)^2_S=0$. 
However, we know by Proposition \ref{prop-image} that the image of general lines is a touching conic. 
Thus we need to know about touching conics lying on  real planes.

First we consider the case that the quartic curve is smooth, and forget the reality for a while (until the proof of Proposition \ref{prop-ftc3}). 
We recall the following classical result in algebraic geometry, which follows from the Pl\"ucker formula.

\begin{prop}\label{prop-28bitan}
Any smooth plane quartic has 28 bitangents.
\end{prop}

Next we see how the bitangents generate touching conics. 
For this purpose, let $B_1$ be a smooth plane quartic and $\phi:S\ra\mathbf{CP}^2$ the double covering branched along $B_1$. 
It is easily seen that $-K_S\simeq\phi^{*}O(1)$ and that $S$ is a del-Pezzo surface of degree two.
The following result is also well-known and easy to prove:

\begin{prop}\label{prop-exceptionalcurve}
(i) The inverse image of any bitangent of a smooth plane quartic consists of two $(-1)$-curves $E_1$ and $E_2$ on $S$, both of which are mapped biholomorphically onto the original bitangent, (ii) $E_1\cap E_2$ consists of two points and the intersections are transversal, (iii) the image of any $(-1)$-curve on $S$ is a bitangent.
\end{prop}

It follows from Propositions \ref{prop-28bitan} and \ref{prop-exceptionalcurve}  that the number of  $(-1)$-curves on $S$ is $28\times 2=56$.

\begin{lemma}\label{lemma-ftc}
Let $B_1$ and $\phi:S\ra\mathbf{CP}^2$ be as above and $C$ a touching conic of $B_1$.
 Then  (i) $\phi^{-1}(C)$ is a reducible curve whose irreducible components are two smooth rational curves which are mapped biholomorphically onto $C$ by $\phi$, 
 (ii)  the two irreducible components of $\phi^{-1}(C)$ have  trivial normal bundles in $S$, and they define two pencils on $S$,
 (iii) the images of general members of the two pencils in (ii) are touching conics of $B_1$, and they define a one-dimensional family of touching conics.
\end{lemma}

\noindent Proof.  Let $P\in B_1\cap C$. 
Then by definition of touching conic (Definition \ref{def-tc}),
we can suppose that in a neighborhood of $P$, $B_1$ and $C$ are defined by  $x^{2m}-y=0$ and $y=0$ respectively for some $m\geq 1$. 
Then on the inverse image $S$ is locally defined by $z^2=x^{2m}-y$ and hence $\phi^{-1}(C)$ is locally defined by $z^2=x^{2m}$. 
It follows that $\Phi^{-1}(C)$ is locally reducible at $\phi^{-1}(P)$. 
Let $(\phi^{-1}(C))'\ra\phi^{-1}(C)$ be the normalization. 
Then the naturally induced map $(\phi^{-1}(C))'\ra C$ is an unramified double covering of $C$. 
Therefore $(\phi^{-1}(C))'$ is a disjoint union of two rational curves, since $C$ is simply connected. 
This implies that $\phi^{-1}(C)$ is reducible. 
We write $\phi^{-1}(C)=L_1+L_2$. 
It is obvious $L_1$ and $L_2$ are smooth rational curves  which are mapped biholomorphically onto $C$. 
Thus we have (i). 
Since $L_1+L_2\in|\phi^*O(2)|$ and since $(\phi^*O(1))^2=c_1^2(S)=2$, we have $(L_1+L_2)^2=8$ on $S$. 
Moreover, by the above local description of $\phi^{-1}(C)$ near $P\in B_1\cap C$, and  recalling that $(C, B_1)=2\cdot 4=8$ on $\mathbf{CP}^2$, it follows that $(L_1,L_2)=4$ on $S$.
 Therefore we have $L_1^2+L_2^2=0$ on $S$. 
 Furthermore, because $S$ is a del-Pezzo surface, there is no smooth rational  curve whose self-intersection number on $S$ is less than $-2$. 
 Further, since the image of $(-1)$-curves are bitangents (Proposition \ref{prop-exceptionalcurve}), $L_1$ and $L_2$ are not $(-1)$-curves. 
 Thus we have deduced $L_1^2=L_2^2=0$.
 It is obvious that $|L_1|$ and $|L_2|$ are pencils on $S$, and $|L_1|\neq|L_2|$ since $L_1L_2=4$.
 Thus we get (ii). 
 Finally $L_1$ and $L_2$ respectively define pencils whose self-intersections are zero, and $L_1$ and $L_2$ are a smooth member of the pencil respectively. 
 By taking the image of the members of the pencils, we get a one-dimensional  family of conics parametrized by $\mathbf{CP}^1$. 
 It remains to see that these are touching conics. 
Let $L\subset S$ be a generic member of   the pencil $|L_1|$. 
Then since the property that $\phi|_L:L\ra\phi(L)$ is isomorphic is an open condition, $\phi|_L$ is isomorphic onto the image. 
Therefore as in the final step in the proof of Proposition \ref{prop-image}, $\phi(L)$ must be a touching conic. 
\proofend

\vspace{3mm}
Next we investigate the number of (one-dimensional) families of touching conics appeared in Lemma \ref{lemma-ftc}. 
So let $C$ be a touching conic and write $\phi^{-1}(C)=L_1+L_2$ as in the proof. 

\begin{lemma}\label{lemma-ftc2}
(i) Both of the pencils $|L_1|$ and $|L_2|$ contain exactly 6 reducible members, any of which are the sum of two $(-1)$-curves intersecting transversally at a point. 
(ii) Each one-dimensional family of touching conics contains exactly 6 reducible members, which are the sum of two bitangents.
\end{lemma}

\noindent Proof. 
(i) Because $S$ is a del-Pezzo surface, there is no smooth rational curve whose self-intersection is less than $-1$. 
Therefore any reducible member of $|L_1|$ is of the form in the lemma. 
Since we have $c_1^2(S)=2$, the number of such reducible members must be 6. 
(ii) is immediate from (i) and Proposition \ref{prop-exceptionalcurve}.\proofend

\begin{prop}\label{prop-ftc3}
A smooth plane quartic has precisely 63 (one-dimensional) families of touching conics.
\end{prop}

\begin{figure}[htbp]
\includegraphics{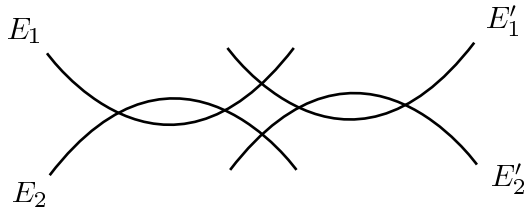}
\caption{the inverse image of a pair of bitangents}
\label{fig-bitan1}
\end{figure}

\noindent Proof. 
Let $l$ and $l'$ be any two different bitangents. 
Then by Proposition \ref{prop-exceptionalcurve} we can write $\phi^{-1}(l)=E_1+E_2$ and $\phi^{-1}(l')=E_1'+E_2'$, where $E_1, E_2,E_1',E_2'$ are $(-1)$ curves, and $E_1$ and $E_2$ (and $E_1'$ and $E_2'$ also) intersect transversally at two points.
 Because $l\cap l'$ is not on the branch quartic,  we can suppose that $E_1E_1'=E_2E'_2=1$ and $E_1E_2'=E_2E_1'=0$ 
 (Figure \ref{fig-bitan1}).
Then we have $(E_1+E_1')^2=(E_2+E'_2)^2=0$.
From this it can be easily shown that  $|E_1+E_1'|$ and $|E_2+E'_2|$ are pencils on $S$ whose general members are smooth rational curves.
  Taking the image by $\phi$, we get a one-dimensional family of touching conics, which contains $l+l'$ as one of the 6 reducible members. 
  In this way a choice of two bitangents defines a (one-dimensional) family of touching conics. 
  By Proposition \ref{prop-28bitan} there are $28!/2!26!=378$ such choices. 
  Moreover we know by Lemma \ref{lemma-ftc2} that each family contains just 6 reducible members. 
  Therefore, the number of the families must be $378/6=63$.\proofend

\vspace{3mm}
So far we have considered bitangents and touching conics of a fixed smooth plane quartic.
From now on we go back to the situation where $B$ is a real quartic surface defined by (\ref{eqn-B2}) and 
consider the space of touching conics of $B$. 
Let $\mathcal P\ra(\mathbf{CP}^3)^{\vee}$ be the $\mathbf{CP}^5$-bundle whose fiber over $H\in(\mathbf{CP}^3)^{\vee}$  is the space of conics on $H$.
 Let $\mathcal P_U\ra U$ be the restriction onto $U$ and $\mathcal C\subset \mathcal P_U$ the closed subset formed by touching conics. 
Since $B$ has a real structure, $\mathcal C$ also has a real structure, and real members of $\mathcal C$ are candidates of the image of twistor lines. 
Then the following proposition is obvious from Proposition \ref{prop-ftc3}:

\begin{prop}\label{prop-ftc4}
The natural projection 
 $\mathcal C\ra U$ is a fiber bundle whose fibers are  disjoint union of 63 copies of $\mathbf{CP}^1$.
\end{prop}

Now it is ready to explain how we will obtain all the twistor lines in $Z$. 
(A complete proof will be given in Sections \ref{s-disj} and \ref{s-thms}.)
By Lemma \ref{lemma-ftc}, $\mathcal C$ admits a natural double covering $\tilde{\mathcal C}\ra\mathcal C$ whose fiber over $C\in\mathcal C$ represents the two irreducible components of $\Phi^{-1}(C)$.
The composition $\tilde{\mathcal C}\ra\mathcal C\ra U$ is a fiber bundle whose fibers are 126 copies of $\mathbf{CP}^1$.
Because $Z_0$ has a real structure, the real structure on $\mathcal C$ naturally lifts on $\tilde{\mathcal C}$.
Let $\tilde{\mathcal C}_{U^{\sigma}}\ra U^{\sigma}$ be the restriction of $\tilde{\mathcal C}\ra U$ onto $U^{\sigma}$.
Because $U^{\sigma}\simeq\mathbf C^*\times\mathbf R$, a monodromy problem naturally arises  for this bundle, and generic twistor lines in $Z$ must be  real members of $\tilde{\mathcal C}_{U^{\sigma}}$ which are monodromy invariant.
Hence studying monodromy of $\tilde{\mathcal C}_{U^{\sigma}}\ra U^{\sigma}$ is one way to detect generic twistor lines in $Z$.
However, we do not take this direction and instead investigate what happens as a plane $H\in U^{\sigma}$ moves to be in $\linfty^{\sigma}$.
For any $H\in \linfty^{\sigma}$, we will determine real lines lying on $S_H=\Phi^{-1}(H)$ in explicit form. 
This is accomplished at length in Section \ref{s-detc}--\ref{s-nb}.
Then we will show that there exists a  neighborhood $W$ of $\linfty^{\sigma}$ in $(\mathbf{RP}^3)^{\vee}$ such that for any $H'\in W$, $S^1$-family of real lines on $S_{H'}$ is obtained as a unique extension of the $S^1$-family of real lines on $S_H$, $H\in\linfty^{\sigma}$.
Real lines obtained in this ways become automatically invariant by the monodromy action, since a loop around $\linfty^{\sigma}$ in $W$ generates $\pi_1(U^{\sigma})\simeq\mathbf Z$. 
Moreover, these real lines  extend in a unique way to give a $S^1$-family of real smooth rational curves on $S_H$ for any $H\in U^{\sigma}$.
On the other hand, we will determine in Section \ref{s-line}, $S^1$-family of real lines lying on $S_H$ for $H\in \mathbf{RP}^2_{\infty}$.
Further, we will show that these real lines on $S_H$, $H\in\mathbf{RP}^2_{\infty}$ is a deformation of real rational curves obtained for $H\in U^{\sigma}$.
In this way we will get a connected four-dimensional family of real rational curves in $Z$.
Once  candidates of twistor lines are determined, we next have to show that these curves foliate $Z$.
In Section \ref{s-disj}, we show that this is the case; 
namely we prove that, for any point of $Z$, there uniquely exists a member of the family going through the point.
We have to be careful at this point since as will be proved in Section \ref{ss-generic},  
there exists a connected family of real  lines in $Z$ such that different members of the family really intersect (Proposition \ref{prop-break}). 
In our proof of the disjointness of lines, we will find the following interesting  geometric obstruction for $B$ to define a twistor space. 
In order to explain that, recall that   $H_i$ ($1\leq i\leq 4$) are the real planes on which the tropes $T_i$ lie. 
Then any line on $H_i$ is automatically a bitangent.

\begin{definition}\label{def-tb}
\em{We call lines on $H_i$ {\it{trivial bitangents}}. 
Bitangent of $B$ which is not trivial is called a {\it{non-trivial bitangent}}.}
\end{definition}
Evidently, the space of trivial bitangents consists of four irreducible components,  which are $H_i^{\vee}$ (the dual projective plane). 
By configuration of $\{H_i\}$ we have $(H_i)^{\vee}\cap (H_j)^{\vee}=\{l_{\infty}\}$ for any $i\neq j$. 
On  any plane other than $H_i$ there are just four trivial bitangents which are the intersection of the plane with $H_i$. 
If the plane is real, these trivial bitangents are obviously real. 
The following proposition implies that non-trivial real bitangent becomes an obstruction to define a twistor space:

\begin{prop}\label{prop-bitanobs}
If $B$ has a non-trivial real bitangent, then there is no resolution $Z\ra Z_0$, such that $Z$ is a twistor space.
\end{prop}

\noindent Proof. 
Let $l$ be a non-trivial real bitangent. 
Since $l_{\infty}$ is a trivial bitangent, $l\neq l_{\infty}$. 
Further, a real line going through $P_0$ cannot be a bitangent since by Proposition \ref{prop-realpoint}, $P_0$ is the unique real point of $B$ and hence the other touching point cannot be real. 
Thus $l$ does not go through the singular points $P_0$, $P_{\infty}$ and $\ol{P}_{\infty}$. 
Let $P$ and $\ol{P}$ be the touching points of $l$ with $B$, and $p$ and $\ol{p}$ the corresponding points of $Z$ respectively. 
Suppose that $Z$ is a twistor space and take the twistor line $L\subset Z$  joining $p$ and $\ol{p}$. 
Then $\Phi(L)$ must be a conic, since $\Phi^{-1}(l)$ cannot contain a twistor line. 
Hence by Proposition \ref{prop-image}, $\Phi(L)$ must be a real conic and there exists a unique real plane $H$ containing $\Phi(L)$.
 This $H$ does not go through $P_0$ since the image of a twistor line lying on the inverse image of such a real plane must be a line by Proposition \ref{prop-image}.
  Assume that $H$ goes through $P_{\infty}$. 
  Then $H$ contains $l_{\infty}$ by reality. 
  Namely, $H\in\linfty^{\sigma}$. 
  Hence $H\cap B$ is a union of two irreducible conics which are the closure of $\mathbf C^*$-orbits. 
  Moreover, since $l$ is not a trivial bitangent, $H\neq H_i$ for any $1\leq i\leq 4$.
   Hence the two conics are mutually different. 
   However, as will be seen in Proposition \ref{prop-3bitan}, a quartic of this kind
   (a union of these two conics) does not have real bitangent other than $l_{\infty}$. 
   This contradicts our choice of $l$. 
   Hence we can suppose that $H$ does not go through $P_{\infty}$ (and $\ol{P}_{\infty}$). 
   Then by Lemma \ref{lemma-singsection}, $H\cap B$ is smooth and $S_H$  (the double cover of $H$ branched along $H\cap B$)  can be considered as a smooth surface in $Z$. 
   By our construction $S_H$ contains both of $L$ and $\Phi^{-1}(l)$, and the latter is an anticanonical curve of $S_H$ which is a sum of two $(-1)$-curves intersecting transversally at $p$ and $\ol{p}$. 
   Thus the intersection number of $L$ and $\Phi^{-1}(L)$ on $S_H$ must be at least 4
   (since $L$ goes through $P$ and $\ol{P}$).
    On the other hand since $L^2=0$ on $S_H$ (see Proposition \ref{prop-PP}), we have $-2=K_SL+L^2=K_SL$ and hence $-K_SL=2$. 
    This is a contradiction, and hence there cannot exist non-trivial real bitangents, if $Z$ is a twistor space.\proofend

\vspace{3mm}
But fortunately, we have the following 
\begin{prop}\label{prop-norealbitan}
There is no non-trivial real bitangent of $B$.
\end{prop}

In the next section we prove this proposition   by studying  degeneration of quartic curves and their bitangents.

\section{Behavior of bitangents as a plane quartic degenerates}\label{s-limit}
In this section we investigate the behavior of  bitangents of plane quartics when the planes move to contain the line $l_{\infty}$. 
As a result  we prove Proposition \ref{prop-norealbitan} saying that our branch quartic surface $B$ does not have non-trivial bitangents (see Definition \ref{def-tb}). 
As seen in Proposition \ref{prop-bitanobs}, non-trivial real bitangents are obstruction for the threefold $Z$ to be a twistor space and non-existence is an important step in our proof of the main theorem.

Let us first consider the simplest situation when a plane quartic degenerates into a quartic which has a node as its only singularity. 
Although we do not need this case itself, clarifying the behavior of bitangents in this simplest case will help to understand the case we really want.

As in Proposition \ref{prop-28bitan}, a smooth plane quartic has 28 bitangents.
Pl\"ucker formula also implies that a quartic with a unique node has 16 bitangents. 
However, the formula only counts the number of bitangents which do not go through the node; namely it does not count lines going through the node and touching the quartic at a different point.
We can count the number of such bitangents:

\begin{prop}\label{prop-6bitan}
A plane quartic $B_0$ having a node as its only singularity possesses precisely 6 bitangents going through the node.
\end{prop}

\begin{figure}[htbp]
\includegraphics{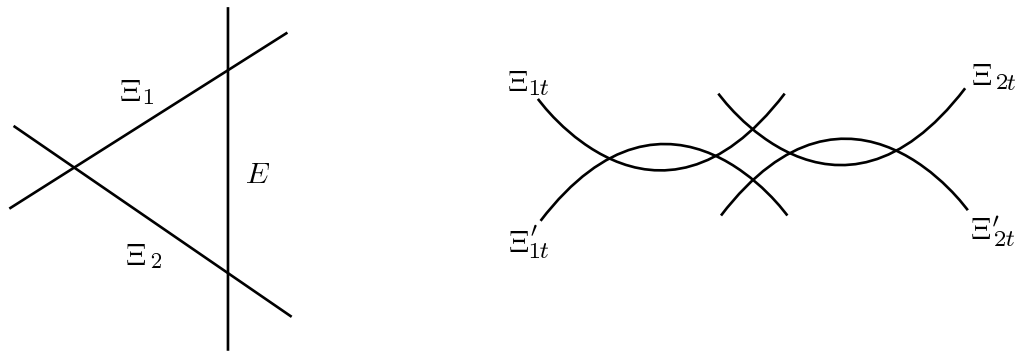}
\caption{the inverse image of a bitangent through the node (left) and its small deformation (right)}
\label{fig-bitan2}
\end{figure}

\noindent Proof. 
Let $P_0$ be the unique node of $B_0$. 
Let $\phi_0:S_0\ra \mathbf{CP}^2$ be the double covering whose branch is $B_0$. $S_0$ has an ordinary double point over $P_0$. 
Let $S_0'\ra S_0$ be the minimal resolution and $E$ the exceptional curve which is a smooth rational curve on $S'_0$ whose self-intersection is $-2$.
We denote by $\psi:S'_0\ra\mathbf{CP}^2$ the composition of $\phi_0$ with the resolution. 
$S'_0$ is a smooth rational surface with $c_1^2=2$, and we have $-K\simeq\psi^*O(1)$. 
Hence we have $(\psi^*O(1)-E)^2=c_1^2+E^2=0$ on $S'_0$. 
If $l$ is a  line going through $P_0$ and intersecting $B_0$ at other two different points, $\psi^{-1}(l)-E$ is a smooth member of the system $|\psi^*O(1)-E|$ and it follows that $|\psi^*O(1)-E|$ is a pencil whose general members are smooth rational curves. Moreover, even if $l$ coincides with one of the two tangent lines of $B_0$ at $P_0$, $\psi^{-1}(l)-E$ is still a smooth rational curve (in this case $\psi^{-1}(l)-E$ touches $E$), and these two tangent lines also give  smooth members of $|\psi^*O(1)-E|$. Remaining is the case that $l$ being a line through $P_0$ but touches $B_0$ at a different point; namely  $l$ is a bitangent through $P_0$.
 If $l$ is this kind of bitangent, $\phi_0^{-1}(l)-E$ consists of two smooth rational curves intersecting at two points, one of which is the ordinary double point of $S_0$ and the other is over the touching point. 
 Thus $\psi^{-1}(l)$ is a triangle of smooth rational curves such that every three intersection is transversal (Figure \ref{fig-bitan2}).
  Write $\Xi_1$ and $\Xi_2$ for the two components other than $E$. 
$\Xi_1+\Xi_2$ is a reducible member of $\psi^{-1}(l)-E$ and we have $\Xi_1^2=\Xi_2^2=-1$ on $S'_0$. 
Thus $|\psi^*O(1)-E|$ is a pencil whose reducible members consist of two irreducible components.
Because we know $c_1^2(S'_0)=2$, it follows that the number of reducible members must be six. 
This implies that $B_0$ has exactly 6 bitangents going through the node.\proofend

\vspace{3mm}
Thus $B_0$ has $16+6=22$ bitangents in all.
We next investigate the relationship of these bitangents with 28 bitangents of a smooth quartic.
Let $\Delta$ be a unit disk in $\mathbf C$ and  $\mathcal B_1\subset \mathbf{CP}^2\times\Delta$ a complex subspace such that $B_t:=\mathcal B_1\cap (\mathbf{CP}^2\times\{t\})$ is a smooth quartic for any $t\neq 0$ and such that $B_0=\mathcal B_1\cap (\mathbf{CP}^2\times\{0\})$ is a quartic which has a node as its only singularity. 
Then $\mathcal B_1$ has a unique singular point which is the node of $B_0$. 
By utilizing a projection $\mathcal B_1\ra \Delta$ we consider $\mathcal B_1$ as a one-dimensional family of plane quartics.
Then 16 bitangents of $B_0$ not going through $P_0$ clearly extend in a unique way to give 16 bitangents of $B_t$ for $t\neq 0$. 
The following lemma describes how the remaining $12$ $(=28-16)$ bitangents of $B_t$ are obtained from the 6 bitangents of $B_0$ going through the node.
  
 \begin{lemma}\label{lemma-degbitan1}
 Among 28 bitangents of a smooth plane quartic $B_t$ ($t\neq 0$), there are 12 bitangents with the following property:
 the 12 bitangents form 6 pairs and the two members of the same pair will converge, as $t\ra 0$, to the same bitangents of $B_0$ going through the node. 
 \end{lemma}
  Of course, the 12 bitangents  in the proposition depend on a choice of degeneration $\mathcal B_1\ra\Delta$.

 \vspace{3mm}
 \noindent Proof of Lemma \ref{lemma-degbitan1}.
 Let $l$ be a bitangent of $B_0$ going through the node.
Let $\mathcal S_1\ra\mathbf{CP}^2\times\Delta$ be the double covering branched along $\mathcal B_1$ and $\mathcal S_1'$ a small resolution of the ordinary double point of $\mathcal S_1$. 
By composition $\mathcal S'_1\ra\mathcal S_1\ra\mathbf{CP}^2\times\Delta\ra\Delta$, we get a projection $\mathcal S'_1\ra\Delta$ whose fibers are smooth rational surfaces with $c_1^2=2$. 
This can be viewed as a degeneration of del-Pezzo surface of degree two. 
Let $\psi:S'_0\ra \mathbf{CP}^2$ denote the restriction of $\mathcal S'_1\ra\mathbf{CP}_2\times\Delta$ onto the central fiber (=$\mathbf{CP}^2\times\{0\}$).
Then we can write $\psi^{-1}(l)=E+\Xi_1+\Xi_2$ which is a triangle of smooth rational curve, where $E$ is the exceptional $(-2)$-curve, and $\Xi_1$ and $\Xi_2$ are $(-1)$-curves as in the proof of Proposition \ref{prop-6bitan}.
Then since $(-1)$-curves survive under small deformation of  surface, $\Xi_1$ and $\Xi_2$ extend to be $(-1)$-curves on $S_t$ (= the fiber of $\mathcal S'_1\ra\Delta$ over $t$) for sufficiently small $t\in\Delta$. 
Moreover, such extensions are obviously unique. 
Let $\Xi_{1t}$ and $\Xi_{2t}$ be the $(-1)$-curves on $S_t$ obtained in this way. 
The stability of $(-1)$-curves only guarantees the existence of such $\Xi_{1t}$ and $\Xi_{2t}$ for small $t$, but these extend to any $t\in\Delta$ since these are irreducible components of the inverse image of bitangents.
Since $\Xi_1\Xi_2=1$ on $S'_0$, we have $\Xi_{1t}\Xi_{2t}=1$ on $S_t$. 
 Taking the image onto $\mathbf{CP}_2$, we get a pair of bitangents of $B_t$, $t\neq 0$. 
These bitangents are mutually different since if they coincide, the inverse image must be a pair of $(-1)$-curves intersecting transversally at two points  over the touching points (Proposition \ref{prop-exceptionalcurve}), which contradicts $\Xi_{1t}\Xi_{2t}=1$.  
Thus we have seen  that a bitangent through the node generates two bitangents of a smooth quartic.
Moreover it is obvious that different bitangents cannot generate the same bitangent.
Hence we get $2\times 6=12$ bitangents of a smooth quartic.
\proofend

\vspace{3mm}
As in the above proof of Lemma \ref{lemma-degbitan1}, we have two  $(-1)$-curves $\Xi_{1t}$ and $\Xi_{2t}$ of $S_t$ which converge to $\Xi_1$ and $\Xi_2$ as $t\ra 0$ respectively.
One may wonder how about the limits of the other components $\phi_t^{-1}(\phi_t(\Xi_{1t}))-\Xi_{1t}$ and $\phi_t^{-1}(\phi_t(\Xi_{2t}))-\Xi_{2t}$,  where $\phi_t:S_t\ra\mathbf{CP}_2$ is the natural projection. 
At this point, we have the following

\begin{lemma}\label{lemma-degbitan2}
When $t$ goes to $0$ in $\Delta$, then $\phi_t^{-1}(\phi_t(\Xi_{1t}))-\Xi_{1t}$ and $\phi_t^{-1}(\phi_t(\Xi_{2t}))-\Xi_{2t}$ respectively converge to  reducible curves $E+\Xi_2$ and $E+\Xi_1$ in $S'_0$. 
\end{lemma}

\noindent Proof. 
This is also proved by standard argument in deformation theory.  
In this proof we write $S=S'_0$ for simplicity.
Let $F$ be the line bundle on $S$ associated to the divisor $E+\Xi_2$. 
We readily have $F^2=-1$ on $S$. 
The exact sequence $0\ra O_S\ra O_S(\Xi_2)\ra O_{\Xi_2}(\Xi_2)\ra 0$ and $\Xi_2^2=-1$ and the rationality of $S$ imply $H^i(O_S(\Xi_2))=0$ for $i=1,2$. 
Then by the exact sequence $0\ra O_S(\Xi_2)\ra F\ra O_E(\Xi_2+E)\ra 0$ and $E\Xi_2=1$ and $E^2=-2$, we get $H^0(S, F)$ is one-dimensional and  $H^i(S,F)=0$ for $i=1,2$.  
Since $S$ is a rational surface, the line bundle $F\ra S$ extends in a unique way to give a line bundle $F_t\ra S_t$ for  $t\neq 0$. 
Then we have $H^i(S_t,F_t)=0$ for $i=1,2$ by the upper-semicontinuity and  $\dim H^0(S_t, F_t)=1$ by the invariance of the Euler characteristic $\chi(F_t)$ under deformation.
In particular $|F_t|$ still consists of a unique member. 
Obviously $F_t^2=-1$ on $S_t$. 
Then because $S_t$ is a del-Pezzo surface for $t\neq 0$, there does not exist a smooth rational curve on $S_t$ whose self-intersection is less than $-1$. 
This implies that the unique member of $|F_t|$ ($t\neq 0$) must be irreducible and it must be a $(-1)$-curve on $S_t$. 
Thus we get a $(-1)$-curve $\Xi_{1t}'$ on $S_t$ as a deformation of $E+\Xi_2$. 
To complete a proof of the lemma it suffices to show that $\Xi'_{1t}=\phi_t^{-1}(\phi_t(\Xi_{1t}))-\Xi_{1t}$. 
It is immediate to verify that $\Xi_1(E+\Xi_2)=2$ on $S$. 
Hence we have $\Xi_{1t}\Xi'_{1t}=2$ on $S_t$.
Namely $\{\Xi_{1t},\Xi'_{1t}\}$ is a pair of $(-1)$-curves satisfying $\Xi_{1t}\Xi'_{1t}=2$.
On the other hand, in view of Proposition \ref{prop-exceptionalcurve},
the intersection number  of different $(-1)$-curves on $S_t$ is $0,1$ or 2, and it becomes two exactly when their images are the same bitangent (see Figure \ref{fig-bitan2}).
Thus we have $\phi_t(\Xi_{1t})=\phi_t(\Xi_{1t}')$.  
Therefore we obtain $\Xi'_{1t}=\phi_t^{-1}(\phi_t(\Xi_{1t}))-\Xi_{1t}$, as claimed. 
The claim for $E+\Xi_1$ is completely parallel. \proofend

\vspace{3mm}
The following lemma is obvious from the proofs of Lemmas \ref{lemma-degbitan1} and \ref{lemma-degbitan2}.

\begin{lemma}\label{lemma-limit1}
Let $S'_0$ be as in the proof of Proposition \ref{prop-6bitan} (or  Lemmas \ref{lemma-degbitan1} and \ref{lemma-degbitan2}) 
and $\mathcal S'_1\ra\Delta$ a degeneration of del-Pezzo surfaces  introduced in the proof of Lemma \ref{lemma-degbitan1}, having $S'_0$ as the central fiber. 
Then there exists (in general reducible) 56 curves on $S'_0$ satisfying the following properties:
(i) the self-intersection numbers of the 56 curves on $S'_0$ are $(-1)$,
(ii) the  56 curves can be naturally extended to $(-1)$-curves on $S_t$ for any $t\in\Delta$, $t\neq 0$,
(iii) the $(-1)$-curves obtained in (ii) are the set of $(-1)$-curves on $S_t$.
\end{lemma}

\vspace{3mm}
So far we have considered the simplest situation that the quartic degenerates to have a unique node. 
From now on we consider the situation it actually happens for our quartic surface. 
Recall that if $H\in \linfty^{\sigma}$ (and $H\neq H_i, H_{\lambda_0}$), then $H\cap B$ is a union of two $\mathbf C^*$-invariant conics
and that the two conics touches each other at two $\mathbf C^*$-fixed points. 
So let $B_0$ be  such a $\mathbf C^*$-invariant quartic. 
$B_0$ has just two singular points, both of which are $A_3$-singularities of a curve.

\begin{prop}\label{prop-3bitan}
$B_0$ has three bitangents. 
One is the line connecting the two singular points of $B_0$ and the other two are the two common tangents  at the singular points of $B_0$.
\end{prop}

\noindent
Proof. 
It suffices to show that there are not bitangents other than the three bitangents in the proposition. 
So assume that  $l$ is such a bitangent of $B_0$.
Then it is immediate to see that $l$ does not pass through the two singular points of $B_0$. 
This implies that $l$ is not $\mathbf C^*$-invariant. 
Moreover, being bitangent is an invariant condition under $\mathbf C^*$-action since the branch curve $B_0$ is assumed to be $\mathbf C^*$-invariant.
Hence we get a one-dimensional family of bitangents of $B_0$.
However, bitangents always must be isolated, since the inverse image contains $(-1)$-curve as an irreducible component, which cannot be moved in the double covering. 
Hence there are no bitangents other than the three ones.\proofend

\vspace{3mm}
In the following we investigate the behavior of 28 bitangents when a smooth plane quartic degenerates into the above $B_0$. 
We consider a family of plane quartics $\mathcal B_2\ra\Delta$ such that $B_t$ (= the fiber over $t$) is smooth for $t\neq 0$ and the central fiber is $B_0$. 
 The following lemma corresponds to Lemma \ref{lemma-degbitan1} in the simplest case.
In the subsequent statement $l_1$ means the line through the two ($A_3$-) singularities of $B_0$, and $l_2$ and $l_3$ mean the remaining two bitangents (which are common tangents  of the two irreducible components of $B_0$). 

\begin{lemma}\label{lemma-degbitan3}
The 28 bitangents of $B_t$, $t\neq 0$, can be grouped into three subsets $\mathcal G_1$, $\mathcal G_2$ and $\mathcal G_3$ which behave as follows when $t$ goes to $0$ in $\Delta$: every bitangent in $\mathcal G_i$ converges to $l_i$ for $i=1,2$ and $3$. Furthermore, $\mathcal G_1$ consists of $16$ bitangents, whereas both of $\mathcal G_2$ and $\mathcal G_3$ consist of 6 bitangents.
\end{lemma}

\begin{figure}[htbp]
\includegraphics{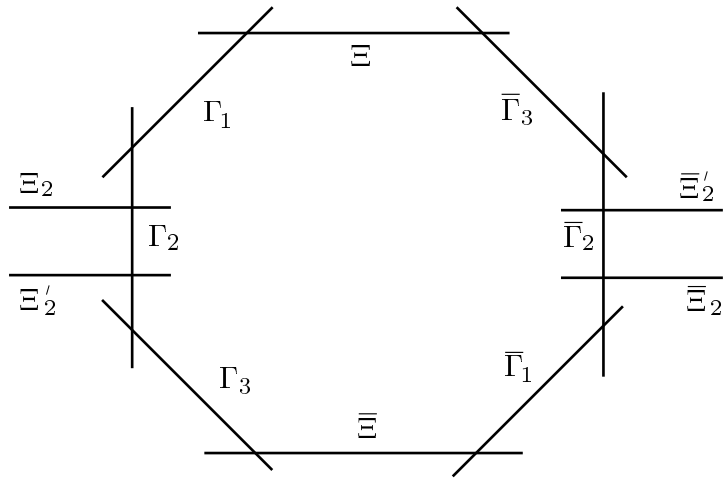}
\caption{the inverse image of $l_1+l_2+l_3$}
\label{fig-cycle2}
\end{figure}

\noindent Proof.
Take a double cover $\mathcal S_2\ra\mathbf{CP}^2\times\Delta$ branching along $\mathcal B_2$.
Then $\mathcal S_2$ has two compound $A_3$-singularities, which are of course over the singular points of $B_0$.
By using the natural projection $\mathcal S_2\ra\Delta$, $\mathcal S_2$ can be viewed as a degeneration of del-Pezzo surfaces (of degree two) into a surface $S_0$ which has two $A_3$-singularities. 
($S_0$ is the double cover branched along $B_0$.)
Hence by a well-known property of rational double points of surface, $\mathcal S_2\ra\Delta$ admits a simultaneous resolution $\mathcal S'_2\ra\mathcal S_2$, which can be also regarded a small resolution of the compound singularities of $\mathcal S_2$. 
(There are many small resolutions of $\mathcal S_2$ as will be explicitly given in Section \ref{ss-special}.)
It is also well know that over the central fiber the simultaneous resolutions give the minimal resolution of the $A_3$-singularities of $S_0$. 
We again denote $\psi:S'_0\ra\mathbf{CP}^2$ for the composition of the minimal resolution with the covering map $S_0\ra\mathbf{CP}^2$ as in the proof of Lemma \ref{lemma-degbitan1}.
The exceptional curve of the minimal resolution  of a $A_3$-singularity is a chain of three $(-2)$-curves which we write $\Gamma=\Gamma_1+\Gamma_2+\Gamma_3$, where we may suppose $\Gamma_1\Gamma_2=\Gamma_2\Gamma_3=1$ and $\Gamma_3\Gamma_1=0$ on $S'_0$.
We denote the exceptional curve over the other $A_3$-singularity by $\ol{\Gamma}=\ol{\Gamma}_1+\ol{\Gamma}_2+\ol{\Gamma}_3$ (since we keep in mind the case that $B_0$ and $S_0$ have real structures). 
These are also $(-2)$-curves on $S'_0$ satisfying $\ol{\Gamma}_1\ol{\Gamma}_2=\ol{\Gamma}_2\ol{\Gamma}_3=1$ and $\ol{\Gamma}_3\ol{\Gamma}_1=0$. 
$l_1$ was the line connecting the two $A_3$-singularities of $B_0$. 
Then $\psi^{-1}(l_1)$ is a cycle of 8 smooth rational curves, six of which are $\Gamma_1,\cdots,\ol{\Gamma}_3$ and the other two are   smooth rational curves which are mapped biholomorphically onto $l_1$. 
We still denoted these curves by $\Xi$ and $\ol{\Xi}$. 
(See Figure \ref{fig-cycle2}.) 
We claim that $\Xi^2=\ol{\Xi}^2=-1$ on $S_0'$. 
Since we still have $-K=\psi^*O(1)$ for the anticanonical class of $S'_0$, and since we have $c_1^2(S'_0)=2$, we have $(\psi^{-1}(l_1))^2=2$ on $S'_0$. 
Because the incidence relations of the irreducible components of $\psi^{-1}(l_1)$ are as in Figure \ref{fig-cycle2}, we easily get $\Xi^2+\ol{\Xi}^2=-2$. 
Suppose that one of $\Xi$ and $\ol{\Xi}$, say $\Xi$ satisfies $\Xi^2\geq 0$ on $S'_0$. Then by deformation theory, $\Xi$ has a non-trivial deformation in $S'_0$. 
By taking the image by $\psi$, we get a non-trivial family of lines. 
This cannot happen because the curves obtained as deformation of $\Xi$ in $S'_0$ always intersect $\Gamma$ and $\ol{\Gamma}$ (because $\Gamma\Xi=\ol{\Gamma}\Xi=1$) and therefore the image (by $\psi$) must go through the two singular points of $B_0$.
Thus we get the claim $\Xi^2=\ol{\Xi}^2=-1$ on $S'_0$.
Then as in the proof of Lemma \ref{lemma-degbitan2}, a connected reduced curve $\Omega\subset\psi^{-1}(l_1)$ containing $\Xi$ but not containing $\ol{\Xi}$ as its irreducible component will uniquely extended to be a $(-1)$-curve on  $S_t$ for  small $t\in \Delta$.
It is easy to list up all such curves $\Omega\subset\psi^{-1}(l_1)$ in explicit forms and conclude that the number of such $\Omega$ is $4\times 4=16$. 
Moreover, the intersection numbers of any two different curves among these 16 curves can be verified to be 0 or 1. 
It then follows, as in the final part of the  proof of Lemma \ref{lemma-degbitan2}, that any of the 16 $(-1)$-curves on $S_t$ ($t\neq 0$) are mapped (by $\phi_t:S_t\ra \mathbf{CP}^2$) to mutually different bitangents. 
By construction, all these 16 bitangents clearly converge to $l_1$ when $t\ra 0$. 
These bitangents of $B_t$ give the members of $\mathcal G_1$ in the lemma.

Next we look for 6 bitangents in $\mathcal G_2$ in a similar way.
$\psi^{-1}(l_2)$ consists of 5 smooth rational curves, three of which are $\Gamma_1,\Gamma_2$ and $\Gamma_3$, and the remaining two are mapped biholomorphically onto $l_2$.
Let $\Xi_2$ and $\Xi_2'$ be the latter two components, which are $(-1)$-curves on $S'_0$. 
In this case, one should note that $\psi^{-1}(l_2)$ is non-reduced and contains $\Gamma_2$ as a component of multiplicity two.
Moreover, we have  $\Xi_2\Gamma_1=\Xi_2\Gamma_3=\Xi_2'\Gamma_1=\Xi_2'\Gamma_3=0$ and $\Xi_2\Gamma_2=\Xi'_2\Gamma_2=1$ on $S'_0$ (see Figure \ref{fig-cycle2}).
 (These can be easily seen by calculations using local coordinates. 
 We do not explain it here since in Lemma \ref{lemma-intersection3}, we calculate the intersection numbers of $\Gamma_i$'s and the inverse image of a conic which has $l_2$ and $l_3$ as a tangent line at the two singular points of $B_0$. )
  Then we take 6 curves
\begin{equation}\label{eqn-6}
\Xi_2,\hspace{2mm}\Xi_2+\Gamma_2,\hspace{2mm}\Xi_2+\Gamma_2+\Gamma_1,\hspace{2mm}\Xi_2+\Gamma_2+\Gamma_3,\hspace{2mm}
\Xi_2+\Gamma_2+\Gamma_1+\Gamma_3,\hspace{2mm}\Xi_2+2\Gamma_2+\Gamma_1+\Gamma_3\end{equation}
as candidates of curves which extend to be $(-1)$-curves on  $S_t$,
$t\neq 0$. 
It is easily verified that the self-intersections of these 6 curves are $-1$. 
Then again an argument in the proof of Lemma \ref{lemma-degbitan2} with a slight modification shows that these curves (\ref{eqn-6}) also can be deformed to be $(-1)$-curves on $S_t$. 
Moreover the intersection numbers of the different two curves among (\ref{eqn-6}) are either 0 or 1 and not 2.
Hence the images of the $(-1)$-curves on $S_t$ must be mutually different. 
Thus we get 6 bitangents whose limit as $t\ra 0$ are $l_2$. 
These give the members of $\mathcal G_2$ in  the lemma. 

Finally, the remaining 6 bitangents in $\mathcal G_3$ are similarly obtained by considering $l_3$, instead of $l_2$ in the above argument.
\proofend

\vspace{3mm}
Analogous to Lemma \ref{lemma-limit1}, we have the following

\begin{lemma}\label{lemma-limit2}
Let  $\mathcal S'_2\ra\Delta$ be the degeneration of del-Pezzo surfaces of degree two introduced in the proof of Lemma \ref{lemma-degbitan3} and $S'_0$ the fiber over the origin.
Then there exist (generally reducible) 56 curves on $S'_0$ satisfying the following properties:
(i) the self-intersection numbers of the 56 curves on $S'_0$ are $(-1)$,
(ii) the  56 curves can be naturally extended to $(-1)$-curves on $S_t$ for $t\in\Delta$, $t\neq 0$,
(iii) the $(-1)$-curves obtained in (ii) are the set of $(-1)$-curves on $S_t$.
\end{lemma}

\noindent Proof.
We continue to use the notations in  the proof of Lemma \ref{lemma-degbitan3}.
In the proof, we have obtained sixteen connected curves contained in $\psi^{-1}(l_1)$ containing $\Xi$ but not containing $\ol{\Xi}$ as their irreducible components.
These give sixteen $(-1)$-curves on $S_t$, $t\neq 0$. 
Interchanging the role of $\Xi$ and $\ol{\Xi}$, we obtain other sixteen $(-1)$-curves on $S_t$.
Thus we get $16\times 2=32$ $(-1)$-curves so far.
On the other hand, the six curves (\ref{eqn-6}) yield six $(-1)$-curves on $S_t$.
Exchanging $\Xi_2$ and $\Xi_3$, we get other six $(-1)$-curves on $S_t$.
Thus we get $6\times 2=12$ $(-1)$-curves on $S_t$  from the curves contained in $\psi^{-1}(l_2)$. 
Then by exchanging the role of $l_2$ and $l_3$, we get  twelve $(-1)$-curves.
Thus we get $32+2\times12=56$ $(-1)$-curves on $S_t$.
These must be all of the $(-1)$-curves on $S_t$ by Proposition \ref{prop-exceptionalcurve}.
\proofend

\vspace{3mm}
Next we take real structures into consideration and prove Proposition \ref{prop-norealbitan} (non-existence of non-trivial real bitangent).
Let $\mathcal S_2\ra \Delta$ be the degeneration of del-Pezzo surface constructed from $\mathcal B_2\ra\Delta$ as above and suppose that  $\mathcal B_2$ is invariant under the real structure on $ \mathbf{CP}^2\times\Delta$ which is the product of standard real structures on $\mathbf{CP}^2$ and $\Delta\subset\mathbf C$
(i.e. complex conjugations). 
Then $\mathcal S_2$ admits a natural  real structure and a fiber $S_t$ is real if $t\in\Delta^{\sigma}$. 
Suppose moreover that on the real fibers of $\mathbf{CP}^2\times\Delta\ra\Delta$ the real structure keeps the line $l_1$ invariant and exchanges $l_2$ and $l_3$. 
Then the two $A_3$-singularities of $B_0$ are necessarily a conjugate pair and $l_2\cap l_3$ is a real  point.
Moreover there are simultaneous resolutions $\mathcal S'_2\ra\Delta$ of $\mathcal S_2\ra\Delta$ such that the real structure on $\mathcal S_2$ lifts onto $\mathcal S'_2$,
since once one of the singularities is resolved, the other singularity is automatically resolved by reality.
We still denote by $\psi:S'_0\ra\mathbf{CP}^2$ for the composition of the minimal resolution $S'_0\ra S_0$ and the covering map $S_0\ra\mathbf{CP}^2$ whose branch is $B_0$.

\begin{lemma}\label{lemma-4realbitan}
Let $\mathcal B_2\ra\Delta$ and $\mathcal S'_2\ra\mathcal S_2\ra\Delta$ be as above and suppose that the real structure acts on $\psi^{-1}(l_1)$ in such a way that it interchanges  irreducible components on the opposite side. 
Then for $t\in\Delta^{\sigma}$ and $t\neq 0$ the number of real bitangents of $B_t$ is four.
\end{lemma}

\noindent Proof.
By Lemma \ref{lemma-limit2}, 28 bitangents of $B_t$ ($t\neq 0$) are obtained from (generally reducible) 56 curves on $S'_0$ by  smoothing and then taking the image by the double covering map. 
Two members $E_1$ and $E_2$ among the 56 curves on $S'_0$ generate $(-1)$-curves that are mapped the same bitangent of $B_t$ precisely when $E_1E_2=2$ on $S'_0$.
Thus it suffices to show that there are just four such pairs $\{E_1,E_2\}$ satisfying $E_2=\ol{E}_1$.
This is possible because the 56 curves have been explicitly determined.
Actually, if $E_1$ is contained $\psi^{-1}(l_2)$ or $\psi^{-1}(l_3)$, then $\ol{E}_1$ is contained in $\psi^{-1}(l_3)$ or $\psi^{-1}(l_2)$ respectively by our assumption on the real structure.
Thus we always have $E_1\ol{E}_1=0$.
So we can suppose that $E_1$ is contained in $\psi^{-1}(l_1)$.
In this case, $E_2=\psi^{-1}(l_1)-E_1$ is the unique curve (among the 56 curves on $S'_0$) which satisfies $E_1E_2=2$. 
It is easily seen from the action of the real structure on the irreducible components on $\psi^{-1}(l_1)$ that there are just four $E_1$'s satisfying $\ol{E}_1=\psi^{-1}(l_1)-E_1$.
Thus we have obtained the claim of the proposition.
\proofend

\vspace{3mm}\noindent
{\textbf{Proof of Proposition \ref{prop-norealbitan}}}.
It is readily seen by Lemma \ref{lemma-singsection} that real lines on $\mathbf{CP}^3$ which are not contained in some $H\in U^{\sigma}$ are $l_{\infty}$ and real lines going through $P_0$.
$l_{\infty}$ is a trivial bitangent and real lines through $P_0$ cannot be a bitangent as we will see in the proof of  Proposition \ref{prop-imline} (easy to show).
Thus to prove the non-existence of  non-trivial real bitangent of $B$, it suffices to show that any $H\in U^{\sigma}$ has just four real bitangents (which are of course trivial bitangents).

Fix any real plane $H_0\in\linfty^{\sigma}$ which is different from $H_i$ ($1\leq i\leq 4$) and $H_{\lambda_0}$. 
Then $H\cap B$ is a union of two  $\mathbf C^*$-invariant conics.
Next take a real line $\gamma$ in $(\mathbf{CP}^3)^{\vee}$ going through $H_0\in(\mathbf{CP}^3)^{\vee}$.
Then $\gamma$ defines a real pencil of planes containing $H_0$ as a (real) member, and by considering their intersections with $B$, we get a one-dimensional family of plane quartics.
By taking $\gamma$ sufficiently general, we can suppose that  general members of this family is smooth.
%
Let $\mathcal B\ra\gamma$ be a family of plane quartics thus obtained.
By construction, this family  enjoys the same properties as $\mathcal B_2\ra\Delta$ in Lemma \ref{lemma-4realbitan}.
Therefore fibers of $\mathcal B\ra\gamma$ has just four real bitangents, at least in a neighborhood of $H_0\in\gamma$.
But since we know that $U^{\sigma}$ is connected, we get that $H\cap B$ has just four real bitangents for any $H\in U^{\sigma}$. 
On the other hand, $H\cap H_i$ ($1\leq i\leq 4$) are actually bitangents of $B\cap H$.
This implies that every bitangent of $B\cap H$ is trivial, and   we get the claim of Proposition \ref{prop-norealbitan}.\proofend



\section{Defining equations of real touching conics}\label{s-detc}
In Section \ref{s-how} we have shown that  real lines in our threefold $Z$ are mapped biholomorphically (by $\Phi:Z\ra\mathbf{CP}^3$)  onto real touching conics  in general (Proposition \ref{prop-image}).
Then we studied touching conics lying on general real planes (i.e.\,planes in $U^{\sigma}$) and showed that they form 63 one-dimensional families (Proposition \ref{prop-ftc3}).
We also explained that in order to select out the right family which will actually be the image of twistor lines, it is important to investigate real lines contained in $S_H=\Phi^{-1}(H)$ where $H$ is a real $\mathbf C^*$-invariant planes (i.e.\,$H\in\linfty^{\sigma}$), or equivalently, real touching conics lying on these planes.
In Sections \ref{s-detc}--\ref{s-nb} we investigate these special real lines and touching conics.
These lines and conics can be viewed as a limit of general twistor lines and their images.
Since the plane sections of $B$ by the $\mathbf C^*$-invariant planes are singular (splitting into two $\mathbf C^*$-invariant conics in general), we cannot apply the results in Section \ref{s-how}.
In this section we determine real touching conics lying on these real $\mathbf C^*$-invariant planes.
Thanks to the simple form of the quartics, it is possible to classify them  and write down their defining equation.

Let $B, P_0,P_{\infty},\ol{P}_{\infty}, \Phi_0:Z_0\ra\mathbf{CP}^3, \mu:Z\ra Z_0$, $\Phi:Z\ra\mathbf{CP}^3,\Gamma_0,\Gamma,\ol{\Gamma}$ and $\sigma$ have the meaning as in the beginning of Section \ref{s-how}.
We fix a $\mathbf C^*$-action defined by $$(y_0,y_1,y_2,y_3)\mapsto (y_0,y_1,ty_2,t^{-1}y_3),\hspace{3mm} t\in\mathbf C^*$$ 
(cf. Proposition \ref{prop-def-B} (iii)).
Then $\mathbf C^*$-invariant real planes must be of the form $H_{\lambda}=\{y_0=\lambda y_1\}$, $\lambda\in \mathbf R$ or $H_{\infty}=\{y_1=0\}$.
($H_{\lambda_0}$ is the unique plane in $\linfty^{\sigma}$ going through $P_0$.)
Moreover, as in Section \ref{s-defeq}, we put 
$$
f(\lambda)=\lambda(\lambda+1)(\lambda-a).
$$
(Note that $a>0$.)
The sign of $f(\lambda)$ will be important in the sequel.
We are assuming that $Q$ and $f$ satisfy Condition (A) in Proposition \ref{prop-necessa}.

\begin{figure}[htbp]
\includegraphics{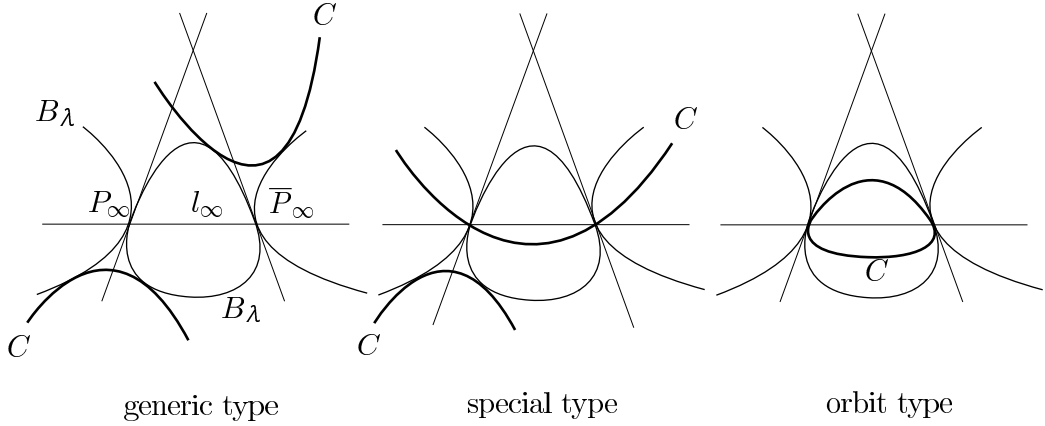}
\caption{real touching conics ($C$ is irreducible)}
\label{fig-conics}
\end{figure}

As seen in the proof of Proposition \ref{prop-realpoint}, $B_{\lambda}=B\cap
H_{\lambda}$ is a union of two 
$\mathbf C^*$-invariant conics, and their
intersection is
$P_{\infty}$ and $\ol{P}_{\infty}$.
Note that if $f(\lambda)=0$ (namely if $\lambda=-1,0$ or $a$) or $\lambda=\infty$, $B_{\lambda}$ becomes a trope.
Let $C$ be a real touching conic contained in $H_{\lambda}$, $\lambda\neq -1,0,a,\infty$.
Since  the two components of $B_{\lambda}$ have the same tangent lines at
$P_{\infty}$ and $\ol{P}_{\infty}$,  the local intersection number
$(C, B_{\lambda})_{P_{\infty}}$ is
zero, two, or four, and by reality, the same holds for  $(C,
B_{\lambda})_{\ol{P}_{\infty}}$. 
Correspondingly, real touching conics lying on some $H_{\lambda}$ can be
classified into the following three types:

\begin{definition}{\em{ Let $C$  be a real touching conic in  $H_{\lambda}$, $\lambda\neq -1,0,a,\infty$.
Then 
(i)  $C$  is called {\em generic type}\,
if $(C, B_{\lambda})_{P_{\infty}}=(C,B_{\lambda})_{\ol{P}_{\infty}}=0$.
(ii)  $C$  is called {\em special type}
if $(C, B_{\lambda})_{P_{\infty}}=(C,B_{\lambda})_{\ol{P}_{\infty}}=2$.
(iii)  $C$  is called {\em orbit type}
if $(C, B_{\lambda})_{P_{\infty}}=(C,B_{\lambda})_{\ol{P}_{\infty}}=4$.
(See Figure \ref{fig-conics}.)}}
\end{definition}

If $C\subset H_{\lambda}$ is a real touching conic of generic type, then $P_{\infty},
\ol{P}_{\infty}\not\in C$ and 
 $C\cap B_{\lambda}$ consists of just four points, all satisfying
$(C,B_{\lambda})_P=2$. 
If $C$ is a real touching conic of special type, then $C$ goes through  $P_{\infty}$
and $
\ol{P}_{\infty}$ but the tangent lines  at
$P_{\infty}$ and $\ol{P}_{\infty}$ are different from the common bitangents  of $B_{\lambda}$.
Further,  there  are other intersection $P$ and $\ol{P}$ satisfying
$(C,B_{\lambda})_P=(C,B_{\lambda})_{\ol{P}}=2$. 
If $C$ is a real touching conic of orbit type, there are no other intersection
points. 
 In this case, $C$ is the closure of $\mathbf C^*$-action,
where $\mathbf C^*$-action is the complexification of $U(1)$-action.

First we classify real touching conics of generic type:

\begin{prop}\label{prop-a} (i) If $\lambda\in\mathbf R$ satisfies $f(\lambda)>0$ and 
$Q(\lambda,1)^2>f(\lambda)$ (i.e.\,$\lambda\neq\lambda_0)$, there exists a family of real touching conics
of generic type  on
$H_{\lambda}$, parametrized by a circle. 
Their  defining equations are
explicitly given by
\begin{equation}\label{eqn-a-0}
2(Q^2-f)y_1^2+\sqrt{f}e^{i\theta}y_2^2+2Qy_2y_3+\sqrt{f}e^{-i\theta}y_3^2=0,
\end{equation} where we put  $Q=Q(\lambda,1)$ and $f=f(\lambda)$, and
$\theta\in\mathbf R$. 
Further, every real touching conic of generic type in
$H_{\lambda}$ is a member of this family. 
(ii) If $f(\lambda)<0$ or
$Q(\lambda,1)^2=f(\lambda)$, there is no real touching conic of generic
type on
$H_{\lambda}$. 
(iii) For the case (i), the conic (\ref{eqn-a-0}) has no real point for any
$\theta\in\mathbf R$.
\end{prop}

We note that  $U(1)$ acts transitively (but non-effectively) on the space of these
touching  conics. 
To prove the proposition, we need the following 

\begin{lemma}\label{lemma-2dr} An equation
$x^4+a_1x^3+a_2x^2+a_3x+a_4=0$ has two double roots iff (i) if $a_1\neq
0$, then $4a_1a_2=a_1^3+8a_3$ and
$a_1^2a_4=a_3^2$ hold, (ii) if $a_1=0$, then  $a_3=0$ and $4a_4=a_2^2$
hold. In these cases, the double roots are given by 
$$-\frac{a_1}{4}\pm\sqrt{\frac{a_1^2}{16}-\frac{a_3}{a_1}}$$
for the case (i) and 
$$
\pm\sqrt{\frac{-a_2}{2}}$$
for the case (ii).
\end{lemma}

We omit a proof since it is elementary.

\vspace{3mm}
\noindent Proof of Proposition \ref{prop-a}. 
Let 
$ay_1^2+by_1y_2+cy_2^2+dy_1y_3+ey_2y_3+hy_3^2=0$ be 
a defining equation of $C$ in $H_{\lambda}$. 
Then since
$B_{\lambda}\cap\{y_3=0\}=\{\ol{P}_{\infty}\}$ and since $C$ is assumed to be
generic type, all of the touching points are  on  $\{y_3\neq0\}$.
Putting $x_1=y_1/y_3$ and $x_2=y_2/y_3$ on $\{y_3\neq 0\}$ as before, 
$C$ is defined by
\begin{equation}\label{eqn-a}
 ax_1^2+bx_1x_2+cx_2^2+dx_1+ex_2+h=0
\end{equation} 
and $B_{\lambda}$ is defined by (as in the proof of Proposition
\ref{prop-realpoint})
$$\left(x_2+\left(Q-\sqrt{f}\right)x_1^2\right)
\left(x_2+\left(Q+\sqrt{f}\right)x_1^2\right)=0.$$
 Let $g$ denote $g_+:=Q+\sqrt{f}$ or
$g_-:=Q-\sqrt{f}$. 
Substituting $x_2=-gx_1^2$ into (\ref{eqn-a}), we get
\begin{equation}\label{eqn-a-1} g^2cx_1^4-gbx_1^3+(a-ge)x_1^2+dx_1+h=0.
\end{equation} If $c=0$, (\ref{eqn-a-1}) cannot have two double roots,
so we have $c\neq 0$. 
Suppose $b\neq 0$. 
Then by Lemma \ref{lemma-2dr}
(i), (\ref{eqn-a-1}) has two double roots iff 
$$-\frac{4b}{cg}\cdot\frac{a-eg}{cg^2}=-\frac{b^3}{c^3g^3}+\frac{8d}{cg^2}
\hspace{3mm}{\rm{and}}\hspace{2mm}
\left(-\frac{b}{cg}\right)^2\cdot\frac{h}{cg^2}=\frac{d^2}{c^2g^4}
$$ hold. 
These can be respectively written
\begin{equation}\label{eqn-a-2} 4bc(a-eg)=b^3-8gc^2d,\,\,\,\,\,b^2h=cd^2.
\end{equation} Namely, a conic (\ref{eqn-a}) with $b\neq 0$ is a
touching conic of generic type iff (\ref{eqn-a-2}) is satisfied for both of
$g=g_+$ and $g=g_-$. 
In this case, simple calculations show that 
$4ac=b^2$ and $4ah=d^2$ and $2ae=bd$. 
From these we get
$a=b^2/4c, e=2cd/b$ and $h=cd^2/b^2$. 
Substituting these into (\ref{eqn-a}),
we get $(bx_1+2cx_2+2cd/b)^2=0$. 
This implies that $C$ is a double line.
 Thus
contradicting our assumption and we get $b=0$. 
Then by Lemma
\ref{lemma-2dr} (ii), we have $d=0$ and 
\begin{equation}\label{eqn-a-3} 4g^2ch=(a-ge)^2.
\end{equation} If we regard (\ref{eqn-a-3}) as a homogeneous equation
about
$(a:c:e:h)\in\mathbf{CP}^3$, (\ref{eqn-a-3}) is a quadratic cone whose
vertex is
$(a:c:e:h)=(g:0:1:0)$.
We need to get the intersection of these two quadrics. 
Restricting
(\ref{eqn-a-3}) onto the plane
$a=\kappa e$, we get
\begin{equation}\label{eqn-a-4} 4g^2ch=(\kappa-g)^2e^2,\,\, g=g_{\pm}.
\end{equation} It is readily seen that these two conics (for the case
$g=g_+$ and
$g=g_-$) coincide iff
$\kappa=0$ or $\kappa=(Q^2-f)/Q$. 
If $\kappa=0$, we have $a=0$, so 
(\ref{eqn-a}) becomes
$cx_2^2+ex_2+h=0$, where the coefficients are subjected to $4ch=e^2$. 
This
implies that
$C$ is a union of two lines, contradicting  our assumption. 
Hence we have
$\kappa=(Q^2-f)/Q$. 
 Then (\ref{eqn-a-4}) becomes 
\begin{equation}\label{eqn-a-5} 4Q^2ch=fe^2.
\end{equation} If $e=0$, it follows that $h=0$, and hence
(\ref{eqn-a}) will again be a union of lines. 
Therefore we have $e\neq 0$. 
It is
easily seen that the real structure on the space of coefficients is
given by
$(a:c:e:h)\mapsto (\ol{a}:\ol{h}:\ol{e}:\ol{c})$.  
 Hence if $f=f(\lambda)<0$,
 (\ref{eqn-a-5})
cannot hold for real $(a:c:e:h)$. 
Namely, on $H_{\lambda}$, there does not
exist a real touching conic of generic type if
$f(\lambda)<0$. 
This implies (ii) of the proposition for the case
$f(\lambda)<0$. 
If $f=f(\lambda)>0$, putting $e=1$ and $h=\ol{c}$ in
(\ref{eqn-a-5}) we get
$4Q^2|c|^2=f$. 
Hence we can write $$c=\frac{\sqrt{f}}{2Q}e^{i\theta}$$
for some
$\theta\in\mathbf R$. 
Further we have $a=(Q^2-f)/Q$. 
Substituting these
into (\ref{eqn-a}), we get (\ref{eqn-a-0}).  
Then it is immediate to see
that the determinant of the matrix defining (\ref{eqn-a-0}) is
$-2(Q^2-f)^2$. 
Therefore the conic (\ref{eqn-a-0}) is irreducible iff
$Q^2-f\neq 0$.
 Thus we get 
 (i), and also (ii) for the case $Q(\lambda)^2=f(\lambda)$.

Finally we show (iii). 
Recall that the real structure on $H_{\lambda}$
is given by
$(y_1:y_2:y_3)\mapsto (\ol{y}_1:\ol{y}_3:\ol{y}_2)$. 
Hence if
$(y_1:y_2:y_3)\in H_{\lambda}$ is a real point, we can suppose $y_1\in\mathbf
R$ and $y_3=\ol{y}_2$. 
Substituting these into (\ref{eqn-a-0}), we get  
\begin{equation}\label{eqn-a-6}
(Q^2-f)y_1^2+Q|y_2|^2+\sqrt{f}\cdot{\rm{Re}}(e^{i\theta}y_2^2)=0,
\end{equation} where ${\rm{Re}}(z)$ denotes the real part of $z$.
From this it follows that $y_2=0$ implies $y_3=y_1=0$, so we have
$y_2\neq 0$. 
Then recalling that $Q>\sqrt{f}$ (Proposition \ref{prop-realpoint}
(i)) we have
$Q|y_2|^2>\sqrt{f}|y_2|^2$. 
Also we have
${\rm{Re}}(e^{i\theta}y_2^2)\leq |y_2|^2$. 
Therefore we have
\begin{eqnarray*}
(Q^2-f)y_1^2+Q|y_2|^2+\sqrt{f}\cdot{\rm{Re}}(e^{i\theta}y_2^2)& >&
(Q^2-f)y_1^2+\sqrt{f}|y_2|^2-\sqrt{f}|y_2|^2\\ &=& (Q^2-f)y_1^2\,\,\geq 0
\end{eqnarray*} This implies that (\ref{eqn-a-6}) does not hold for any
real $(y_1:y_2:y_3)$ and
$\theta\in\mathbf R$. 
Thus there is no real point on the conic (\ref{eqn-a-0})
provided $Q^2>f>0$, and we get (iii).
\hfill $\square$

\vspace{3mm}
Next we classify real touching conics of special type:

\begin{prop}\label{prop-b} (i) If $\lambda\in\mathbf R$ satisfies $f(\lambda)>0$, there is no real
touching conic of special type on
$H_{\lambda}$.
 (ii) If $f(\lambda)<0$, there exists a family of 
 real touching conics of special type, parametrized by a circle. 
 Their
defining equations are given by
\begin{equation}\label{eqn-b-0}
\sqrt{Q^2-f}\cdot y_1^2+\sqrt{\frac{\sqrt{Q^2-f}-Q}{2}}\cdot
e^{i\theta}y_1y_2+
\sqrt{\frac{\sqrt{Q^2-f}-Q}{2}}\cdot e^{-i\theta}y_1y_3+y_2y_3=0,
\end{equation} where we put  $Q=Q(\lambda,1)$ and $f=f(\lambda)$,
and $\theta\in \mathbf R$ as before.
Further, every
real touching conic of special type  in $H_{\lambda}$ is a member of 
this family. 
 (iii) The conic (\ref{eqn-b-0}) has no real point for any
$\theta\in\mathbf R$.
\end{prop}

Note again that $U(1)$ acts transitively on the parameter space of
these touching conics.
Also note that if $f<0$ we have $Q^2-f>0$ and $\sqrt{Q^2-f}-Q>0$. 
Hence
every square root in the equation make a unique sense (i.e. we always take the positive
root).

\vspace{3mm}
\noindent Proof. 
Let $C$ be a real touching conic of special type on
$H_{\lambda}$. 
Then since $C$ goes through $P_{\infty}$ and $\ol{P}_{\infty}$, 
the other two touching points
belong to mutually different irreducible component of $B_{\lambda}$.
On the other hand, as shown in the proof of Proposition \ref{prop-realpoint},
each irreducible components of $B_{\lambda}$ is real iff $f(\lambda)\geq 0$.
Therefore, on $H_{\lambda}$, there does not exist real touching conic of special type 
if $f(\lambda)>0$. 
Thus we get (i). 
So in the sequel we suppose $f(\lambda)<0$.

It is again readily seen that touching points are not on the line
$\{y_3=0\}$. 
 So we still use $(x_1,x_2)$ as a non-homogeneous coordinate
 on
$\mathbf C^2=\{y_3\neq 0\}\subset H_{\lambda}$. 
Then because
$C$ contains $P_{\infty}$ and
$\ol{P}_{\infty}$, an equation of a touching conic $C$ of special
type is of the form
\begin{equation}\label{eqn-b} ax_1^2+bx_1x_2+dx_1+ex_2=0.
\end{equation} Substituting $x_2=-gx_1^2$ into (\ref{eqn-b}), we get
\begin{equation}\label{eqn-b-1}
x_1\cdot\left(gbx_1^2+(ge-a)x_1-d\right)=0.
\end{equation} If $d=0$, $x_1=0$ is a double root of (\ref{eqn-b-1}).
Then the tangent line of $C$ at $P_{\infty}=(0,0)$ becomes $x_2=0$,  as
is obvious from (\ref{eqn-b}). 
This implies that $C$ is a touching conic of
orbit type, contradicting our assumption.
Therefore, we have $d\neq 0$.
 (\ref{eqn-b-1}) has a
double root other than $x_1=0$ iff
\begin{equation}\label{eqn-b-2} (ge-a)^2+4gbd=0.\end{equation}
 Namely
(\ref{eqn-b}) is a touching conic of special type iff (\ref{eqn-b-2}) is
satisfied for both of $g=g_+$ and $g=g_-$.
 (Note that $g_+\neq g_-$.)
  If we
regard (\ref{eqn-b-2}) as a  homogeneous equation of
$(a:b:d:e)\in\mathbf{CP}^3$, (\ref{eqn-b-2}) is a quadratic cone whose
vertex is 
$(a:c:d:e)=(g:0:0:1)$. 
That is, the parameter space of  touching conics
of special type is the intersection of two quadratic cones in
$\mathbf{CP}^3$.
 Restricting (\ref{eqn-b-2}) onto the plane $a=\kappa
e$, we get
$$ (g-\kappa)^2e^2+4gbd=0,\,\, g=g_{\pm}.$$
 It is readily seen that these
two conics coincide iff $\kappa=\pm\sqrt{g_+g_-} =\pm\sqrt{Q^2-f}$.
 Therefore $C$ is a
touching conic of special type iff
\begin{equation}\label{eqn-b-3} a=\sqrt{Q^2-f}\cdot e,\,\,\,
\left(Q-\sqrt{Q^2-f}\right)e^2+2bd=0
\end{equation} or
\begin{equation}\label{eqn-b-4} a=-\sqrt{Q^2-f}\cdot e,\,\,\,
\left(Q+\sqrt{Q^2-f}\right)e^2+2bd=0
\end{equation} hold. 
 It is easily seen that the real structure on the space of coefficients
is given by
$(a:b:d:e)\mapsto (\ol{a}:\ol{d}:\ol{b}:\ol{e})$.  

Since we have assumed  $f<0$, we have $Q+\sqrt{Q^2-f}>0$ 
and there is no real conic satisfying (\ref{eqn-b-4}). 
On the other
hand, we have
$Q-\sqrt{Q^2-f}<0$.
 Hence by (\ref{eqn-b-3}) we have
$$ 2|b|^2=\left(\sqrt{Q^2-f}-Q\right)e^2,
$$ where $b\in \mathbf C$ and $e\in\mathbf R$. 
If $e=0$, then $b=a=0$,
contradicting  the assumption that $C$ is a conic. 
Hence $e\neq 0$, and
we may put $e=1$. 
Then we can write
$$b=\sqrt{\frac{\sqrt{Q^2-f}-Q}{2}}\cdot e^{i\theta}$$ for some
$\theta\in\mathbf R$. 
Also we have $d=\ol{b}$. 
Thus we obtain
(\ref{eqn-b-0}) of the proposition.

Finally we show (iii).
 If $(y_1:y_2:y_3)$ is a real point of
$H_{\lambda}$, we can suppose $y_1\in\mathbf R$ and $y_3=\ol{y}_2$.
 Substituting
these into (\ref{eqn-b-0}), we get
\begin{equation}\label{eqn-b-5}
\sqrt{Q^2-f}\cdot y_1^2+\sqrt{2\left(\sqrt{Q^2-f}-Q\right)}\cdot
y_1\cdot{\rm{Re}}(e^{i\theta}y_2)+|y_2|^2=0.
\end{equation} If $y_1=0$, it follows $y_2=y_3=0$. 
Hence $y_1\neq 0$
 and we can suppose $y_1>0$. 
 Then
 we have $y_1{\rm{Re}}(e^{i\theta}y_2)\geq -y_1|y_2|$. 
 Hence we
have
\begin{eqnarray*}
\lefteqn{
\sqrt{Q^2-f}\cdot y_1^2+\sqrt{2\left(\sqrt{Q^2-f}-Q\right)}\cdot
y_1{\rm{Re}}(e^{i\theta}y_2)+|y_2|^2}\\ &\geq&\sqrt{Q^2-f}\cdot
y_1^2-\sqrt{2\left(\sqrt{Q^2-f}-Q\right)}\cdot y_1|y_2|+|y_2|^2\\
&=&\left(|y_2|-\sqrt{\frac{\sqrt{Q^2-f}-Q}{2}}\cdot y_1\right)^2+
\frac{\sqrt{Q^2-f}+Q}{2}\cdot y_1^2
\end{eqnarray*} Because $y_1\neq 0$ and $f<0$, we have $(\sqrt{Q^2-f}+Q)y_1^2>0$.
Therefore, the left hand side of (\ref{eqn-b-5}) is strictly positive. 
Thus 
(\ref{eqn-b-5}) does not hold for any real $(y_1:y_2:y_3)\in H_{\lambda}$   and any
$\theta\in\mathbf R$. 
Therefore the conic (\ref{eqn-b-0}) has no real point for any
$\theta\in\mathbf R$.

It is immediate to see that the determinant of the matrix defining
(\ref{eqn-b-0}) is $-(Q+\sqrt{Q^2-f})/8$ and this is negative if $f<0$.
Hence the conic (\ref{eqn-b-0}) is irreducible
\hfill$\square$

\vspace{3mm} The case of orbit type is straightforward and need no assumption
on the sign of $f(\lambda)$:
\begin{prop}\label{prop-c}  There exists a family of real touching
conics of orbit type, parametrized by non-zero real numbers. 
Their 
defining equations are 
\begin{equation}\label{eqn-c} y_2y_3=\alpha y_1^2,\,\,\alpha\in \mathbf
R^{\times}.
\end{equation} Further, every real touching conic of orbit type in
$H_{\lambda}$ is contained in this family. 
\end{prop}

Note that by Lemma \ref{lemma-element}, the conic (\ref{eqn-c})
has no real point iff $\alpha<0$.

Combining Propositions \ref{prop-a}, \ref{prop-b} and \ref{prop-c}, we get the
following

\begin{prop}\label{prop-type1} Let
$\{S_{\lambda}:=\Phi^{-1}(H_{\lambda})\set H_{\lambda}\in\langle
l_{\infty}\rangle^{\sigma}\}$ be the real members of the pencil of
$U(1)$-invariant  divisors on $Z$. 
Then  (i) if
$f(\lambda)>0$,  the images of real lines in $S_{\lambda}$ are real
touching conics of generic type or orbit type, (ii) if $f(\lambda)<0$, 
the images of real lines in
$S_{\lambda}$ are real touching conics of special type or orbit type.
\end{prop}

In Section \ref{s-nb} we will determine which one of the two candidates in the proposition must be the images of twistor lines (Proposition \ref{prop-type2}).

\section{The inverse images of real touching conics}\label{s-inv} According
to the previous section, real touching conics in a $\mathbf C^*$-invariant real plane form
 families parametrized by a circle for generic and special types, or
$\mathbf R^{\times} (\ni\alpha$) for orbit type. 
In this section we study the
inverse images  in $Z$ of these touching conics, which are candidates of  twistor lines. 
We continue to use the same notations and assumptions.
Recall that if $\lambda\neq -1,0,a,\infty$, $S_{\lambda}=\Phi^{-1}(H_{\lambda})$ is a smooth rational surface by Proposition \ref{prop-str_of_S}. 
Clearly $S_{\lambda}$ is $\mathbf C^*$-invariant real member of $|(-1/2)K_Z|$
(Proposition \ref{prop-pic}).

First we investigate the inverse images of touching conics of generic type.

\begin{prop}\label{prop-inv-a} Suppose that $\lambda\in\mathbf R$ satisfies $f(\lambda)>0$ and $\lambda\neq\lambda_0$
(namely $Q^2-f>0$), and let
$C_{\theta}\subset H_{\lambda}$ be a real touching conic of generic type
defined by the equation (\ref{eqn-a-0}). 
 Then the following (i)--(iv)
hold:  (i) 
$\Phi^{-1}(C_{\theta})$ has just two irreducible components, both of
which are smooth rational curves that are mapped biholomorphically onto
$C_{\theta}$, (ii) each irreducible component of $\Phi^{-1}(C_{\theta})$
has a trivial normal bundle in $S_{\lambda}$,   (iii)  these two
irreducible components of
$\Phi^{-1}(C_{\theta})$ belong to mutually different pencils on
$S_{\lambda}$, (iv) each  irreducible component of 
$\Phi^{-1}(C_{\theta})$ is  real.
\end{prop}

\noindent Proof. 
Since $C_{\theta}$  and the branch quartic 
$B_{\lambda}$  have the same
tangent line at any intersection points, it is obvious that 
$\Phi_0^{-1}(C_{\theta})$  splits into two irreducible  components
$L_1$ and $L_2$ which are
mapped biholomorphically onto $C_{\theta}$. 
 Thus we get (i).
 For (ii) first note that
$\Phi^{-1}(C_{\theta})=L_1+L_2$ belongs to $|-2K|$ of
$S_{\lambda}$ since we have $\Phi^*O_{H_{\lambda}}(1)\simeq -K$.
 Hence we have
$(-2K)^2=(L_1+L_2)^2=L_1^2+L_2^2+2L_1L_2$ on $S_{\lambda}$.
 On the other hand, we
have
$4c_1^2=8$  by Proposition
\ref{prop-str_of_S}. 
Hence we get $L_1^2+L_2^2+2L_1L_2=8$. 
Further, since $L_1$
and $L_2$ intersect transversally at four points (over the
touching points of
$C_{\theta}$  with $B_{\lambda}$), we have $L_1L_2=4$. 
Therefore
we get
$L_1^2+L_2^2=0$. 
Moreover, by (\ref{eqn-a-0}), $C_{\theta}$ actually moves in
a holomorphic family of curves on
$H_{\lambda}$. 
Hence we have $L_1^2\geq 0$ and $L_2^2\geq 0$. 
Therefore we get
$L_1^2=L_2^2=0$. 
Namely we have (ii).
 (iii) immediately follows from
(ii), since we have $L_1L_2=4$ on $S_{\lambda}$. 
(iv) is harder than  one may think
at first glance,  since there is no real point on
$C_{\theta}$. 
First we note that it suffices to prove the claim for
$C_0$ (= the curve obtained by setting $\theta=0$ for $C_{\theta}$), since
$U(1)$ acts transitively on the parameter space of real  touching conics
of generic type (see Proposition \ref{prop-a}).
 The idea of our proof  of the reality is as follows:   the map
$\Phi^{-1}(C_0)\ra C_0$ is  finite, two sheeted covering whose branch
consists of four points. 
 We choose a real simple closed curve $\mathcal
C$ in $C_0$ containing all of these branch points,  in such a way that
over $\mathcal C$  we can
distinguish two sheets, so that 
 we can explicitly see the reality of each irreducible components.
  To
this end, we still use $(y_1:y_2:y_3)$ as a homogeneous coordinate on $H_{\lambda}$
and set $V:=\{y_1\neq 0\}$ ($\simeq \mathbf{C}^2$), which is clearly a
real subset of
$H_{\lambda}$, and use $(v_2,v_3)=(y_2/y_1,y_3/y_1)$ as an affine coordinate on
$V$.
 Then $Z_0|_{V}=\Phi^{-1}_0(V)$ is
defined by the equation
\begin{equation}\label{eqn-Z-V} z^2+\left(v_2v_3+Q\right)^2-f=0,
\end{equation} where  $z$ is a fiber coordinate of $ O(2)$, and
$Q=Q(\lambda,1)$, $f=f(\lambda)$ as in the previous section. 
 The real
structure is given by 
$(v_2,v_3, z)\mapsto (\ol{v}_3,\ol{v}_2,\ol{z})$. 
  Then on $V$, our equation (\ref{eqn-a-0}) of $C_0$ becomes 
$2(Q^2-f)+\sqrt{f}v_2^2+2Qv_2v_3+\sqrt{f}v_3^2=0$.
  Now we introduce a
new coordinate $(u,v):=(v_2+v_3,v_2-v_3)$ which is valid on $V$.
 Then our real
structure is given by $(u,v)\mapsto (\ol{u},-\ol{v})$, and the defining
equation of $C_0$ becomes
$4(Q^2-f)+\sqrt{f}(u^2+v^2)+Q(u^2-v^2)=0$. From this we immediately have 
\begin{equation}\label{eqn-uv-1}
C_0:\hspace{3mm}v^2=4\left(Q+\sqrt{f}\right)+\frac{Q+\sqrt{f}}{Q-\sqrt{f}}u^2.
\end{equation} 
We put $V':=\{(u,v)\in  V\set u\in i\mathbf{R},
v\in\mathbf R\}$ which  is clearly a real subset of $V$.
 Then $\mathcal
C:=V'\cap C_0$ is a real simple closed curve (an ellipse) in $V'\simeq\mathbf R^2$.
 Indeed,
putting $u=iw$ ($w\in\mathbf R$), we get from (\ref{eqn-uv-1})  
$$C_0:\hspace{3mm}v^2+\frac{Q+\sqrt{f}}{Q-\sqrt{f}}w^2=4\left(Q+\sqrt{f}\right).$$
(Note that  $Q-\sqrt{f}>0$ by our assumption $f>0$ and Proposition
\ref{prop-necessa}.)
  Substituting
$v_2v_3=(u^2-v^2)/4$ into (\ref{eqn-Z-V}), and then using (\ref{eqn-uv-1}), 
we get
\begin{equation}\label{eqn-inv-C}
\Phi^{-1}(C_0):\hspace{3mm}4\left(Q-\sqrt{f}\right)^2z^2
+fu^2\left(u^2+4\left(Q-\sqrt{f}\right)\right)=0,
\end{equation} or,  using $w$ above, 
\begin{equation}\label{eqn-inv-C2}
\Phi^{-1}(C_0):\hspace{3mm}4\left(Q-\sqrt{f}\right)^2z^2
=fw^2\left(4\left(Q-\sqrt{f}\right)-w^2\right).
\end{equation}
 Here note that in (\ref{eqn-uv-1}) $u$ can be used as a coordinate on
$C_0$, only outside the two branch points of $C_{\lambda}\cap V\ra\mathbf C$
defined by
$(u,v)\mapsto u$.
 In a neighborhood of these branch points,  we
have to use $v$ instead of $u$ as a local coordinate on $C_0$.
 Then we can
see that the inverse image of a neighborhood of
$u=\pm 2i(Q-\sqrt{f})^{\frac{1}{2}}$ \,(i.e. the branch points) also
splits into two irreducible components, which is of course as expected.
From (\ref{eqn-inv-C}), we easily deduce that the branch points of
$\Phi^{-1}(C_0)\ra C_0$ are $(u,v)=(0,\pm2(Q+\sqrt{f})^{\frac{1}{2}})$
and
$(u,v)=(\pm 2i(Q-\sqrt{f})^{\frac{1}{2}},0)$.
 All of these four points
clearly lie on
$\mathcal C$, and
$\mathcal C$ is divided into four segments. 
It immediately follows from
(\ref{eqn-inv-C2}) that 
 $z$ always takes  real value over $\mathcal C$.
  Moreover, it is clear
that  the sign of $z$ is constant on each of the four segments in
$\mathcal C$, and that the sign changes when passing though the branch
points. 
On the other hand, since the real structure is given by
$(w,v)\mapsto (-w, -v)$ on $V'$, the real structure on
$\mathcal C$ sends each segment to another segment which is not adjacent
to the original one. 
From these, and because the real structure on
$\Phi^{-1}(\mathcal C)$ is given by $(w,z)\mapsto (-w,\ol{z})=(-w,z)$,
it follows that each of the two irreducible component of
$\Phi^{-1}(\mathcal C)$ is real. 
Hence the same is true for $\Phi^{-1}(C_0)$.
Thus we have proved (iv) of the proposition.
\hfill
 \proofend

\vspace{3mm}
We have similar statements for touching conics of special type:

\begin{prop}\label{prop-inv-b} Assume $f(\lambda)<0$ and let
$C_{\theta}\subset H_{\lambda}$ be a real touching conic of special type
given by the equation (\ref{eqn-b-0}). 
Then  (i)--(iv) of Proposition
\ref{prop-inv-a} hold if we replace $\Phi^{-1}(C_{\theta})$ by
$\Phi^{-1}(C_{\theta})-\Gamma-\ol{\Gamma}$, where  we set
$\Gamma:=\Phi^{-1}(P_{\infty})$ and
$\ol{\Gamma}:=\Phi^{-1}(\ol{P}_{\infty})$.
\end{prop}

\noindent Proof.
 (i) can be proved in the same way as in Proposition
\ref{prop-inv-a}.
 (But in this case, any small resolution $Z\ra Z_0$
gives the resolution of 
$\Phi^{-1}_0(C_0)$ at the points over $P_{\infty}$ and $\ol{P}_{\infty}$, as will be
mentioned below.)
 For (ii) first note that we have
$\Phi^{-1}(C_{\theta})=L_1+L_2+\Gamma+\ol{\Gamma}\in |-2K|$ this time,
where
$L_1$ and
$L_2$ are irreducible components of
$\Phi^{-1}(C_{\theta})-\Gamma-\ol{\Gamma}$. 
$\Gamma$ and $\ol{\Gamma}$ are chains of  three $(-2)$-curves on
$S_{\lambda}=\Phi^{-1} (H_{\lambda})$, since they are exceptional curves
of the minimal resolution of $A_3$-singularities of surface.
We write $\Gamma=\Gamma_1+\Gamma_2+\Gamma_3$, where $\Gamma_i$'s are 
smooth rational curves satisfying $\Gamma_1\Gamma_2=\Gamma_2\Gamma_3=1$ and
$\Gamma_1\Gamma_3=0$ on $S_{\lambda}$.
We then have 
$\Gamma^2=\ol{\Gamma}^2=-2$. 
Furthermore, as we will state and prove in Lemma 
\ref{lemma-interse},
 we have (or more precisely, can suppose)
 $L_1\Gamma_1=1$, $L_1\Gamma_i=0$ for $i=2,3$, and
 $L_3\Gamma_3=1$, $L_3\Gamma_i=0$ for $i=1,2$.
 (In that lemma we write $L^+_{\theta}$ and $L^-_{\theta}$ instead of 
 $L_1$ and $L_2$ respectively.)
 It follows that 
  $L_1\Gamma=L_1\ol{\Gamma}=L_2\Gamma=L_2\ol{\Gamma}=1$. 
  Therefore
again by Proposition \ref{prop-str_of_S}, we get
$8=(-2K)^2= (L_1+L_2+\Gamma+\ol{\Gamma})^2=L_1^2+L_2^2+2L_1L_2+4$.
  But
because
$L_1$ and $L_2$ intersect transversally at two points 
(over the touching points of $C_{\theta}$ and $B_{\lambda}$), and
because $L_1$ and $L_2$ do not intersect on $\Gamma\cup\ol{\Gamma}$,
we have
$L_1L_2=2$.
 Therefore we have $L_1^2+L_2^2=0$.
  Hence by the same reason in 
the proof of the previous proposition, we again have
$L_1^2=L_2^2=0$.
 This implies (ii).
  (iii) follows from (ii),
because we have $L_1L_2=2$ as is already seen.
 (iv) can be proved by the
same idea as in the previous proposition: first we may assume
$\theta=0$.
 Then by (\ref{eqn-b-0}) the equation of
$C_0$ on $V=\{y_1\neq 0\}=\{(v_2,v_3)\}$ is given by
$$
\sqrt{Q^2-f} +\sqrt{\frac{\sqrt{Q^2-f}-Q}{2}}\cdot
(v_2+v_3)+v_2v_3=0.
$$   If we
use another coordinate $(u,v)$ defined in the proof of the previous
proposition, this can be written as
\begin{equation}\label{eqn-circl1}
 C_0:\,\,v^2=u^2+2\sqrt{2}\sqrt{\sqrt{Q^2-f}-Q}\cdot
u+4\sqrt{Q^2-f}.
\end{equation}
 Next we introduce a new variable $w$ by setting 
$u=-\sqrt{2}(\sqrt{Q^2-f}-Q)^{\frac{1}{2}}+iw$.
 Then the equation becomes
\begin{equation}\label{eqn-bfdj}
C_0:\hspace{3mm}v^2+w^2=2\left(\sqrt{Q^2-f}+Q\right).
\end{equation}
Put $\mathcal C:=C_0\cap\mathbf R^2$, where $\mathbf
R^2=\{(w,v)\set w\in\mathbf R,v\in\mathbf R\}$.
 Then since 
$\sqrt{Q^2-f}+Q>0$, 
$\mathcal C$ is a real circle in $\{(w,v)\in\mathbf R^2\}$. 

By (\ref{eqn-Z-V}) we have 
\begin{equation}\label{eqn-abrcbr}
\Phi_0^{-1}(V):\hspace{2mm}
z^2=f-\left(\frac{u^2-v^2}{4}+Q\right)^2.
\end{equation}
On the other hand,  by (\ref{eqn-circl1}), we have 
$$C_0: \hspace{2mm}u^2-v^2=-2\sqrt{2}\sqrt{\sqrt{Q^2-f}-Q}\cdot
u-4\sqrt{Q^2-f}
$$ on $C_0$.
 Substituting this into (\ref{eqn-abrcbr}), we get
$$
\Phi_0^{-1}(C_0):\,\, z^2=f+\frac{\sqrt{Q^2-f}-Q}{2}w^2.
$$
Then by using (\ref{eqn-bfdj}), we get
\begin{equation}\label{eqn-circdl1}
\Phi_0^{-1}(C_0):\,\, z^2=-\frac{\sqrt{Q^2-f}-Q}{2}v^2.
\end{equation}
Therefore, $z$ is pure imaginary over $\mathcal C$,
so that we can distinguish  two sheets by looking the sign of $z/i$.
By (\ref{eqn-circdl1}) and (\ref{eqn-bfdj}), the branch points of
$\Phi^{-1}(C_0)\ra C_0$ is the two points
$(w,v)=(\pm (2(\sqrt{Q^2-f}+Q))^{\frac{1}{2}},0)$
 which lie on $\mathcal C$.
  The real structure is given by
$(w,v)\mapsto (-\ol{w},-\ol{v})$ and this is equal to $(-w,-v)$ on
$\mathcal C$.
 Thus the real structure on $\mathcal C$ interchanges the
two segments separated by the two branch points.     
Moreover, the real structure on the fiber coordinate is given by $z\mapsto \ol{z}$. 
Therefore, it changes the sign of $z/i$ over $\mathcal C$.
This implies that each component of $\Phi^{-1}(\mathcal C)$
is real. Therefore that of $\Phi^{-1}(C_0)$ is also real.
\proofend

\vspace{3mm} The situation slightly changes for touching conics of orbit
type:

\begin{prop}\label{prop-inv-c} Suppose $f(\lambda)\neq 0$ and let
$C_{\alpha}\subset H_{\lambda}$ be a real touching conic of orbit type
given by the equation (\ref{eqn-c}).
  Then we have: (i) $C_{\alpha}$ is
contained in $B$ iff $\alpha=-Q\pm \sqrt{f}$. 
In the following (ii)--(vi) suppose that $\alpha\neq -Q\pm \sqrt{f}$.
(ii)
$\Phi^{-1}(C_{\alpha})-\Gamma-\ol{\Gamma}$ has just two irreducible
components, both of which are smooth rational curves that are mapped
biholomorphically onto
$C_{\alpha}$.
 (Here $\Gamma$ and $\ol{\Gamma}$ are as in Proposition
\ref{prop-inv-b}.)
 (iii) Each irreducible component of
$\Phi^{-1}(C_{\alpha})-\Gamma-\ol{\Gamma}$ has a trivial normal bundle in
$S_{\lambda}=\Phi^{-1}(H_{\lambda})$.
  (iv) The two irreducible components
of
$\Phi^{-1}(C_{\alpha})-\Gamma-\ol{\Gamma}$ belong to one and  the same
pencil on
$S_{\lambda}$.
  (v) Each irreducible component  of
$\Phi^{-1}(C_{\alpha})-\Gamma-\ol{\Gamma}$ is real iff $f>0$ and 
$-Q-\sqrt{f}<\alpha<-Q+\sqrt{f}$ are satisfied.
(vi) There is no real point on $C_{\alpha}$ if $f$ and $\alpha$ satisfies
the inequalities of (v).
(vii) If $f>0$ and $\alpha=-Q\pm\sqrt{f}$, then $C_{\alpha}$ has no real point.
\end{prop}

Note that $f>0$ implies $Q\geq\sqrt{f}$ by Proposition \ref{prop-realpoint}
(ii). 

\vspace{3mm}
\noindent Proof.
 (i) By substituting $y_2y_3=\alpha y_1^2$ into the defining
equation of $B_{\lambda}$, we obtain
$\left((\alpha+Q)^2-f\right)y_1^4=0$.
 Thus if
$C_{\alpha}$ is contained in $ B$ iff 
$(\alpha+Q)^2=f$, which implies $\alpha=-Q\pm \sqrt{f}$, as desired.
 (ii)
can be seen in the same way as in (i) of Proposition \ref{prop-inv-a}.
(This time, any small  resolution $Z\ra Z_0$ gives the
normalization of 
$\Phi^{-1}_0(C_0)$.)
  Next we prove (iii).
   Let
$\Gamma=\Gamma_1+\Gamma_2+\Gamma_3$, $\ol{\Gamma}=
\ol{\Gamma}_1+\ol{\Gamma}_2+\ol{\Gamma}_3$, $L_1$ and
$L_2$ have the same meaning as in the proof of the last proposition. 
It is obvious that  $\Gamma$ and $\ol{\Gamma}$ are 
disjoint.
 In Lemma \ref{lemma-intersection3} we will show that 
$L_1$ is disjoint from $\Gamma_1$ and $\Gamma_3$, 
and the same for $L_2$, and both of $L_1$ and $L_2$ intersect transversally
at a unique point on $\Gamma_2$.
(In the lemma we write $L^+_{\alpha}$ and $L^-_{\alpha}$ instead of
$L_1$ and $L_2$.)
 By reality,
the same is true for the intersection of $L_1$ and $L_2$ with
$\ol{\Gamma}$.
 Moreover, it will also be shown in Lemma \ref{lemma-intersection3} that 
if $L_1$ and $L_2$ are different, these two curves are disjoint.
Thus 
we have
$L_1L_2=\Gamma\ol{\Gamma}=L_1\Gamma_1=L_1\ol{\Gamma}_1=L_2\Gamma_1=L_2\ol{\Gamma}_1=
L_1\Gamma_3=L_1\ol{\Gamma}_3=L_2\Gamma_3=L_2\ol{\Gamma}_3=0$
 on $S_{\lambda}$, while
$L_1\Gamma_2=L_1\ol{\Gamma}_2=L_2\Gamma_2=L_2\ol{\Gamma}_2=1$
(also on $S_{\lambda}$).
 In
particular, we have
$L_1\Gamma=L_1\ol{\Gamma}=L_2\Gamma=L_2\ol{\Gamma}=1$.
   On the other
hand, we still have $8=(-2K)^2= (L_1+L_2+\Gamma+\ol{\Gamma})^2$
and $L_1^2\geq 0$ and $L_2^2\geq 0$.
Combining these, we get
$L_1^2=L_2^2=0$ on $S_{\lambda}$.
 Thus we have  (iii).
  (iv) easily
follows if we consider  the linear systems $|L_1|$ and $|L_2|$, and if we
note that $L_1 L_2=0$.
 Next we show (v).
  Substituting $v_2v_3=\alpha$ into
(\ref{eqn-Z-V}), we get
$ z^2+(\alpha+Q)^2-f=0.
$ From this, the equations of irreducible components of 
$\Phi_0^{-1}(C_{\alpha})$ can be  calculated to be
\begin{equation}\label{eqn-inv-gen} z=\pm \sqrt{f-(\alpha+Q)^2}
\end{equation} 
Recalling
that the real structure is given by $z\mapsto \ol{z}$,
these curves
are real iff $f-(\alpha+Q)^2\geq 0$.
  In particular, $f\geq 0$ follows.
Then we have
$-Q-\sqrt{f}<\alpha< -Q+\sqrt{f}$ under our assumption, and we get (v)
and (vii). 
Finally, (vi) immediately follows from Lemma \ref{lemma-element}.
\proofend

\vspace{3mm}
Proposition \ref{prop-inv-c} implies that  not all  touching conics of
orbit type can be the image of a twistor line:   $f(\lambda)>0$
is needed, and further, 
$-Q-\sqrt{f}\leq\alpha\leq-Q+\sqrt{f}$ must be satisfied.
But once we know that one of the two irreducible components is actually a twistor
line, it follows that  the other component is also a twistor line,  since
by (iv) these two components can be connected by deformation in
$S_{\lambda}$ (hence also in
$Z$) preserving the real structures.
 This is not true for touching conics of
generic type and special type, because the two irreducible components
of $\Phi^{-1}(C_{\theta})$ or $\Phi^{-1}(C_{\theta})-\Gamma-\ol{\Gamma}$
 intersect.

\section{Twistor lines lying on $\mathbf C^*$-invariant fundamental divisors}\label{s-nb}
 In this section we calculate the normal bundle of $L$ in $Z$, where $L$ is
a real irreducible component of the inverse images of the real
touching conics which are determined in Section \ref{s-detc}.
We continue to use the notations in the previous section.
Because touching conics in $H_{\lambda}\in\linfty^{\sigma}$ generally go through the singular points ($P_{\infty}$ and $\ol{P}_{\infty}$) of $B$, the normal bundles generally depend on the choice of a small resolution of $Z_0$.
So in this section (especially in Section \ref{ss-special}--\ref{ss-conn}), we need to discuss the choice of small resolution of the conjugate pair of singularities of $Z_0$.
(But we do not discuss resolution of the unique real ordinary double point of $Z_0$.
This is discussed in the next section.)

Let us briefly describe the content of this section.
In Section \ref{ss-prel} we give a simple lemma which will be used to determine the normal bundles. 
Some notations are also introduced.
In Sections \ref{ss-generic}--\ref{ss-orbit}  we determine the normal bundles, according to the three types of  real touching conics.
 These subsections are organized as follows:
 first we explicitly calculate the intersection of irreducible components (of  the inverse images of touching conics)   with some curves. 
Consequently we get a function of
$\lambda$ (which will be written $h_i$). 
Second we show that the normal bundle in problem degenerates into 
$O\oplus O(2)$ precisely when
$\lambda$ is a critical point of this function.
Then we determine the critical points.
In Section \ref{ss-con} we use the results in Sections \ref{ss-generic}--\ref{ss-orbit} to determine which type of touching conics actually come from twistor lines (Proposition \ref{prop-type2}).
We also prove that among various ways there are only two small resolutions of the conjugate pair of singularities of $Z_0$  which can yield twistor space (Proposition \ref{prop-sr}).
Then finally in Section \ref{ss-conn} we show that our candidates of twistor lines obtained in Section \ref{ss-con} actually form a connected family (Proposition \ref{prop-limit02}); namely they can be connected by deformations in the threefold, by adding the inverse image of four tropes of $B$.
This result is used in Section \ref{s-disj} to construct  {\em{arbitrary}} twistor lines in $Z$.

%

\subsection{Preliminary lemma and notations}\label{ss-prel}
 In order to determine which one
of
$O(1)^{\oplus 2}$ and $O\oplus O(2)$ actually occurs,
 we use the following elementary criterion:

\begin{lemma}\label{lemma-nb} Let $E\ra\mathbf{CP}^1$ be a holomorphic
line bundle of rank two, and assume that 
$E$ is isomorphic to $ O(1)^{\oplus 2}$ or $O\oplus O(2)$.
 Let $s$ and
$t$ be global sections of $E$ which are  linearly independent as global
sections.
 Then $E\simeq O\oplus O(2)$ iff there are constants
$a,b\in\mathbf C$ such that 
$as+bt$ has two zeros. \end{lemma}

\noindent Proof.
 It is immediate to see that non-zero sections of $O(1)^{\oplus 2}$ 
have at most one zero.
 So sufficiency follows.
  Conversely let
$s=(s_1,s_2)$ and
$t=(t_1,t_2)$ be any linearly independent sections of $O\oplus O(2)$,
where
$s_1,t_1\in\Gamma(O)=\mathbf C$ and
$s_2,t_2\in\Gamma (O(2))$.
 Then take $a,b\in\mathbf C$ such that
$as_1+bt_1=0$.
 Then  $as+bt$ can be regarded as a non-zero section of $O(2)$ so
that it has two zeros.
\proofend

\vspace{3mm} Next we introduce some notations which will be used 
throughout this section.
As in the previous sections, $\lambda\in\mathbf R$ denotes a parameter on 
the space of real $U(1)$-invariant planes.
In other words, $\lambda$ is a parameter on the real locus $l_0^{\sigma}$
of the real line $l_0:=\{y_2=y_3=0\}$.
 The
function $f(\lambda)=\lambda(\lambda+1)(\lambda-a)$, ($a>0$) defines
four open intervals in the circle $l_0^{\sigma}$:
$$I_1=(-\infty,-1),\, I_2=(-1,0),\, I_3=(0,a)\hspace{2mm}
{\rm{and}}\hspace{2mm} I_4=(a,+\infty).
 $$ Namely,
$I_1\cup I_3=\{\lambda\in\mathbf R\set f(\lambda)<0\}$ and $I_2\cup
I_4=\{\lambda\in\mathbf R\set f(\lambda)>0\}$.
 By Proposition
\ref{prop-pos}, the equation $Q(\lambda,1)^2-f(\lambda)=0$ has a
unique real solution
$\lambda=\lambda_0$ which is necessarily a   double root.
  Since
we have
$f(\lambda_0)=Q(\lambda_0,1)^2>0$, $\lambda_0\in I_2\cup I_4$. 
By applying a projective transformation $(y_0,y_1)\mapsto (ay_1,-y_0)$ (interchanging $I_2$ and $I_4$) if necessary, we may suppose that $\lambda_0\in I_4$.
  Then we set
$I_{4}^-=(a,\lambda_0)$ and $I_{4}^+=(\lambda_0, +\infty)$. 

Next suppose $\lambda\in I_2\cup I_4^-\cup I_4^+$, and let
$$\mathcal C_{\lambda}^{\,\rm{gen}}=\{C_{\theta}\subset H_{\lambda}
\set C_{\theta} {\mbox{ is defined by (\ref{eqn-a-0})}}\}$$ be
the set of real touching conics of generic type on $H_{\lambda}$.
 Note that if
$\lambda=\lambda_0$ or if $\lambda\in I_1\cup I_3$ there is no real
touching conic of generic type on
$H_{\lambda}$ by Proposition \ref{prop-a} (ii).
Similarly, for $\lambda\in I_1\cup I_3$, 
let $$\mathcal C_{\lambda}^{\,\rm{sp}}=\{C_{\theta}\subset
H_{\lambda}\set C_{\theta} {\mbox{ is defined by (\ref{eqn-b-0})}}\}\hspace{2mm} 
$$ be
the set of real touching conics of special type  on
$H_{\lambda}$.
Note  that if $\lambda\in I_2\cup I_4$ there is no real touching
conic of generic type by Proposition \ref{prop-b} (i).
 Finally  for
$\lambda\in  I_2\cup I_4$, let
$$\mathcal C_{\lambda}^{\,\rm{orb}}=\left\{C_{\alpha}\subset H_{\lambda}
\set -Q-\sqrt{f}\leq\alpha\leq -Q+\sqrt{f},\hspace{1.5mm}
 C_{\theta} {\mbox{ is defined by (\ref{eqn-c})}}\right\}$$
be
the set of real touching conics of  orbit type on
$H_{\lambda}$.
Note that the restriction on $\alpha$ implies that the two irreducible
components of the inverse image are respectively real (Proposition
\ref{prop-inv-c} (v)). 
$\mathcal C_{\lambda}^{\,\rm{gen}}$ and $\mathcal
C_{\lambda}^{\,\rm{sp}}$ are parametrized by a circle on which $U(1)$
naturally acts transitively, whereas $\mathcal C_{\lambda}^{\,\rm{orb}}$
is parametrized by a closed interval on which $U(1)$ acts trivially.

\subsection{The case of generic type}\label{ss-generic}
This is the easiest case since touching conics of generic type do not pass through
the singular points of $B_0$ and hence resolutions of $Z_0$ do not have effect on 
the normal bundles.
Suppose $\lambda\in I_2\cup I_4^-\cup I_4^+$ and take
$C_{\theta}\in\mathcal C_{\lambda}^{\,\rm{gen}}$.
First we  calculate the intersection of $C_{\theta}$ and
$l_{\infty}$, where $l_{\infty}$ is the real line defined by $y_0=y_1=0$
as before.
 Let
$x_2=y_2/y_3$ be a non-homogeneous coordinate on $l_{\infty}$ (around
$P_{\infty}=(0:0:0:1)$).

\begin{lemma}\label{lemma-int-a} The set $\{C_{\theta}\cap l_{\infty}\set
C_{\theta}\in \mathcal C_{\lambda}^{\rm{gen}} \}$  consists of  disjoint two
circles about
$P_{\infty}$ in $l_{\infty}$, whose radiuses 
(with respect to the coordinate $x_2$ above)
 are
 given by
$$
h_0(\lambda):=\frac{Q+\sqrt{Q^2-f}}{\sqrt{f}}\hspace{3mm}{\rm{and}}\hspace{3mm}
h_0(\lambda)^{-1}=\frac{Q-\sqrt{Q^2-f}}{\sqrt{f}}$$  respectively, where
we put
$Q=Q(\lambda,1)$ and $f=f(\lambda)$ as before. 
\end{lemma}

Note that  we have $Q^2-f>0$ and $Q>\sqrt{f}$ by Proposition \ref{prop-necessa},
and therefore $h_0>1>h_0^{-1}>0$ holds.
Moreover, $h_0$ and $h_0^{-1}$ are differentiable on
 $I_2\cup I_4^-\cup I_4^+$.

\vspace{3mm}
\noindent Proof.
 On $H_{\lambda}=\{(y_0:y_1:y_2)\}$, $l_{\infty}$ is defined by 
$y_1=0$.
 Therefore by (\ref{eqn-a-0}) we readily have 
\begin{equation}\label{eqn-intersec100}C_{\theta}\cap l_{\infty}=
\left\{x_2=\frac{-Q\pm\sqrt{Q^2-f}}{\sqrt{f}}\cdot
e^{-i\theta}\right\}.\end{equation} This directly implies the claim of
the lemma. 
\proofend

\vspace{3mm} By Proposition \ref{prop-inv-a}, $\Phi^{-1}(C_{\theta})$
consists of two irreducible components, both of which are real rational
curves.
 We denote these components by $L^+_{\theta}$ and $L^-_{\theta}$,
although there is no canonical way of making a distinction of these two.
Again by Proposition \ref{prop-inv-a}, 
$L^+_{\theta}$ and $L^-_{\theta}$ respectively form 
disjoint families 
$$\mathcal L_{\lambda}^+=\{L^+_{\theta}\set\theta\in\mathbf
R\}\hspace{3mm} {\rm{and}}\hspace{3mm}\mathcal
L_{\lambda}^-=\{L_{\theta}^-\set\theta\in\mathbf R\}$$ of (real and
smooth)  rational curves on
$S_{\lambda}=\Phi^{-1}(H_{\lambda})$.
 These are real members 
of real pencils on $S_{\lambda}$ and each member has no real point by 
Proposition \ref{prop-a} (iii). 
 Because $U(1)$  acts also on the parameter spaces ($=S^1$) of
$\mathcal L^+_{\lambda}$ and
$\mathcal L^-_{\lambda}$, the normal bundles of $L^+_{\theta}$ and
$L^-_{\theta}$ inside
$Z_0$ are  independent of the choice of $\theta$. 
The following proposition plays a
key role in determining the normal bundle:

\begin{prop}\label{prop-nbc} For any $L\in \mathcal L^+_{\lambda}
\cup \mathcal L^-_{\lambda}$, the normal bundle of $L$ in
$Z_0$ is isomorphic to either $ O(1)^{\oplus 2}$ or $ O\oplus O(2)$.
 Further, the
latter holds iff
$\lambda$ is a critical point of
$h_0(\lambda)$ defined in Lemma \ref{lemma-int-a}.
\end{prop}

In particular, members of $\mathcal L^+_{\lambda}$ and 
$ \mathcal L^-_{\lambda}$ have the same normal bundle in $Z_0$.

\vspace{3mm}
\noindent Proof.
$L$ is contained in the
smooth surface
$S_{\lambda}=\Phi^{-1}(H_{\lambda})$ and therefore we have an exact
sequence
$0\ra N_{L/S_{\lambda}}\ra N_{L/Z_0}\ra N_{S_{\lambda}/Z_0}|_L\ra 0$. 
By  Proposition \ref{prop-inv-a} (ii), we have $N_{L/S_{\lambda}}\simeq 
O_L$.
On the other hand, $S_{\lambda}$ is a smooth member of $|(-1/2)K_Z|$.
Therefore by adjunction formula we have
$K_{S_{\lambda}}\simeq K_Z|_{S_{\lambda}}\otimes N_{S_{\lambda}/Z}\simeq
K_Z |_{S_{\lambda}}\otimes (-1/2)K_Z|_{S_{\lambda}}$ and hence 
$N_{S_{\lambda}/Z}\simeq (-1/2)K_Z|_{S_{\lambda}}\simeq -K_{S_{\lambda}}$.
Hence we get $N_{S_{\lambda}/Z}|_L\simeq-K_{S_{\lambda}}|_L\simeq
-K_L\otimes N_{L/S}\simeq O_L(2)$.
Therefore by the short exact sequence above, $N_{\lambda}:=N_{L/Z}$ is
isomorphic to either
$ O\oplus O(2)$ or $O(1)^{\oplus 2}$.
 Thus we get the first claim of
the proposition.
 
In order to show the second claim, we recall three facts; the first one
is about the 
 the natural real structure on $\Gamma(N_{\lambda})$,
the space of sections of $N_{\lambda}$.
Since $L$ is real,  $\sigma$ naturally acts on $\Gamma(N_{\lambda})$ as the
complex conjugation.
 For $s\in\Gamma(N_{\lambda})$ we denote by ${\rm Re} s$ 
and ${\rm Im} s$ the real part and the imaginary part of $s$ respectively.
Namely, ${\rm Re} s=(s+\ol{\sigma(s)})/2$ and 
${\rm Im} s=(s-\ol{\sigma(s)})/2$.
Secondly,  recall that $L\cap \Phi^{-1}(l_{\infty})$ consists of a conjugate pair
of points, one of which corresponds to the point of $l_{\theta}$ satisfying $x_2=h_2(\lambda)e^{i\theta}$, and the other one corresponds to
to the point of $l_{\theta}$ satisfying $x_2=h_2^{-1}(\lambda)e^{i\theta}$.
(See (\ref{eqn-intersec100})).
 Let $z_{\lambda}$ and $\ol{z}_{\lambda}$
be the former and the latter point respectively.
Thirdly recall that any one-parameter family of holomorphic  deformation of
$L$ in $Z$ naturally gives rise to a holomorphic section of
$N_{\lambda}$:
roughly this section is obtained by taking the tangent vector of 
the `orbit' for each point of $L$.
More precisely, take a neighborhood $U$ of each point of $L$ and
a holomorphic coordinate $(z_1,z_2,z_3)$ on $U$, such that 
each $L_x$, $|x|<\epsilon$ with $L_0=L$ of the given one-parameter
family  is defined by
$z_2=f(z_1,x)$ and $z_3=g(z_1,x)$ satisfying $f(z_1,0)\equiv g(z_1,0)\equiv 0$.
Then a representative of the section of $N_{\lambda}$ is given by
$(\partial f/\partial x,\partial g/\partial y)$.

Now we have the following two one-parameter families of deformations of
$L$ in
$Z$: the first one is obtained by moving $L$ by $\mathbf C^*$-action,
where the $\mathbf C^*$-action is
the complexification of the $U(1)$-action.
 The second one  is
obtained by moving the parameter $\lambda$ in $\mathbf C$, while fixing 
$\theta$.
 Let $s\in\Gamma(N_{\lambda})$ and
 $t\in\Gamma(N_{\lambda})$ be
the holomorphic sections associated to the former and the latter  family
respectively.
 These are linearly independent sections,
since each representative of $s$ is tangent to $S_{\lambda}$ (by the
$\mathbf C^*$-invariance of $S_{\lambda}$ ), whereas that of $t$ is not.
 Because the
$\mathbf C^*$-action preserves
$S_{\lambda}$, it follows from Proposition \ref{prop-inv-a} (ii) and (iii) that
each of the curves of the former family are disjoint.
 This implies that $s$ is
nowhere vanishing.
  Next we consider the latter family.
   First, noting
$H_{\lambda}\cap H_{\lambda'}=l_{\infty}$ for $\lambda\neq\lambda'$,
 $t$ can be zero only on $\Phi^{-1}(l_{\infty})$
 (cf. Figure \ref{fig-trope} in Section \ref{s-defeq}).
  Suppose $\lambda$ is a critical point of $h_0$; namely
$h_0'(\lambda)=0$, where the derivative is with respect to real $\lambda$,
of course.
Then, since $h$ is a holomorphic function of $\lambda$,
it can be easily derived from the Cauchy-Riemann equation that
$(\partial h_0/\partial \lambda)(\lambda)=0$.
By the way how we take a representative of the section of $N$
explained in the previous paragraph,
this directly implies that $t$ vanishes at $z_{\lambda}$.
On the other hand, $h_0'(\lambda)=0$ implies
$(h_0^{-1})'(\lambda)=0$.
  Then in the same manner as above, we have
$(\partial h_0^{-1}/\partial \lambda)(\lambda)=0$.
 These imply that $t$ also vanishes at $\ol{z}_{\lambda}$ and
we have obtained $t(z_{\lambda})=t(\ol{z}_{\lambda})=0$.
  Therefore by Lemma \ref{lemma-nb}, we get
$N_{\lambda}\simeq O\oplus O(2)$. 

Next suppose $h_0'(\lambda)\neq 0$, so
that  
$(h_0^{-1})'(\lambda)\neq 0$.
 We claim that Re$(bs+t)$ cannot vanish
at $z_{\lambda}$ and $\ol{z}_{\lambda}$ simultaneously,
for any $b\in\mathbf C$. 
Because $C_{\theta}$ intersects
$l_{\infty}$ transversally, $t$ also becomes a nowhere vanishing section
under our assumption.
 Hence the claim is true for $b=0$.
Putting $b=b_1+ib_2$, $b_1,b_2\in\mathbf R$,
we easily get
\begin{equation}\label{eqn-reim}
{\rm{Re}} (bs+t)=b_1{\rm{Re}}s+({\rm{Re}} t-b_2{\rm{Im}} s).
\end{equation}
Since $s$ comes from the
$\mathbf C^*$-action, and since its real part corresponds to the 
$U(1)$-action,
$({\rm{Re}}s)(z_{\lambda})$ is  represented by the
tangent vector of the $U(1)$-orbit going through $z_{\lambda}$.
On the other hand, (\ref{eqn-intersec100}) implies that  
$({\rm{Re}}t)(z_{\lambda})$ is represented by a tangent vector
which is parallel to $({\rm{Im}}s)(z_{\lambda})$.
Hence by (\ref{eqn-reim}), we can deduce that 
${\rm{Re}} (bs+t)(z_{\lambda})=0$ implies $b_1=0$ and 
\begin{equation}\label{eqn-reim2}
({\rm{Re}}t)(z_{\lambda})=b_2({\rm{Im}} s)(z_{\lambda}).
\end{equation}
Similarly,  ${\rm{Re}} (bs+t)(\ol{z}_{\lambda})=0$ implies  $b_1=0$ and 
\begin{equation}\label{eqn-reim3}
({\rm{Re}}t)(\ol{z}_{\lambda})=b_2({\rm{Im}} s)(\ol{z}_{\lambda}).
\end{equation}
Suppose $b_2>0$.
 Then since $\{{\rm{Re}}s(z_{\lambda}),
{\rm{Im}}s(z_{\lambda})\}$ is an oriented basis of 
$T_{z_{\lambda}}(\Phi^{-1}(l_{\infty}))$
from the beginning,
(\ref{eqn-reim2}) implies that $\{{\rm{Re}}s(z_{\lambda}),
{\rm{Re}}t(z_{\lambda})\}$ is an oriented basis of
$T_{z_{\lambda}}(\Phi^{-1}(l_{\infty}))$.
 Further, we have
$$({\rm{Re}}s)(\ol{z}_{\lambda})=\sigma_*(({\rm{Re}}s)(z_{\lambda}))
\hspace{3mm}{\mbox{and}}\hspace{3mm}
({\rm{Re}}t)(\ol{z}_{\lambda})=\sigma_*(({\rm{Re}}t)(z_{\lambda})).$$
Hence we get by (\ref{eqn-reim3})
 $$({\rm{Im}} s)(\ol{z}_{\lambda})
=\frac{1}{b_2}({\rm{Re}}t)(\ol{z}_{\lambda})=
\frac{1}{b_2}\sigma_*\left(({\rm{Re}}t)(z_{\lambda})\right).$$
So we have $$\{({\rm{Re}}s)(\ol{z}_{\lambda}),({\rm{Im}} s)(\ol{z}_{\lambda})\}=
\{\sigma_*(({\rm{Re}}s)(z_{\lambda})),
 \sigma_*(({\rm{Re}}t)(z_{\lambda}))/b_2\}.$$
But since $\sigma$ is anti-holomorphic, $\sigma$ is orientation reversing.
Further, as is already seen, $\{{\rm{Re}}s(z_{\lambda}),
{\rm{Re}}t(z_{\lambda})\}$ is an oriented basis of  $T_{z_{\lambda}}
(\Phi^{-1}(l_{\infty}))$
 (if $b_2>0$).
  This implies that $\{\sigma_*(({\rm{Re}}s)(z_{\lambda})),
 \sigma_*(({\rm{Re}}t)(z_{\lambda}))/b_2\}$ is an anti-oriented basis
of  $T_{\ol{z}_{\lambda}}(\Phi^{-1}(l_{\infty}))$.
  This contradicts to the fact that
$\{({\rm{Re}}s)(\ol{z}_{\lambda}),$ $({\rm{Im}} s)(\ol{z}_{\lambda})\}$ is an
oriented basis of $T_{\ol{z}_{\lambda}}(\Phi^{-1}(l_{\infty}))$.
Therefore, Re$(bs+t)$ cannot vanish
at $z_{\lambda}$ and $\ol{z}_{\lambda}$ simultaneously, provided $b_2>0$.
Parallel arguments show the same claim holds for the case $b_2<0$.
Thus we have shown that  Re$(bs+t)$ cannot vanish
at $z_{\lambda}$ and $\ol{z}_{\lambda}$ at the same time, as claimed.
On the other hand, it is obvious that 
$bs+t$ does not vanish except $\{z_{\lambda},\ol{z}_{\lambda}\}$.
Therefore, the zero locus of $bs+t$ consists of at most one point
for any $b\in\mathbf C$.
Since $s$ is a nowhere vanishing section, Lemma \ref{lemma-nb}
implies $N_{\lambda}\simeq O(1)^{\oplus 2}$.
\proofend

\begin{figure}[htbp]
\includegraphics{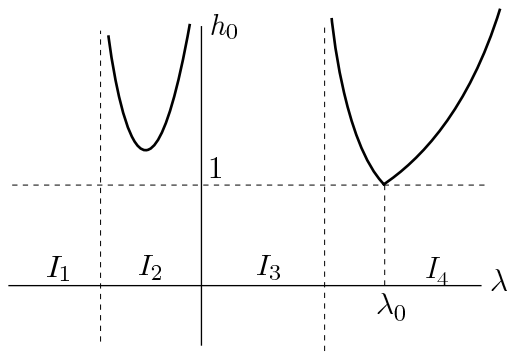}
\caption{behavior of $h_0$}
\label{fig-h0}
\end{figure}

\begin{lemma}\label{lemma-bhv-1} Let $h_0=h_0(\lambda)$ be the positive
valued function on $I_2\cup I_4$ defined in Lemma \ref{lemma-int-a},
which is differentiable on
$I_2\cup I_4\backslash\{\lambda_0\}$.
 Then 
$h_0$  has a unique critical point on $I_2$, and has no critical point on
$I_4\backslash\{\lambda_0\}$.
 (See Figure \ref{fig-h0}.)
\end{lemma}

\noindent Proof.
 Although  elementary, we include the proof since it is not
so easy (at least for the author) and it needs the previous result (Proposition 
\ref{prop-realpoint}).
We have $Q(-1)>0$, $Q(0)>0$ and $Q(a)>0$ by
Proposition \ref{prop-realpoint} (i), and  
\begin{equation}\label{eqn-g-1}
h_0=\frac{Q+\sqrt{Q^2-f}}{\sqrt{f}}=
\frac{Q}{\sqrt{f}}+\sqrt{\frac{Q^2}{f}-1}.
\end{equation} From these,  it follows that 
$\lim_{\lambda\downarrow -1}h_0(\lambda)=\lim_{\lambda\uparrow
0}h_0(\lambda)=+\infty$.
Therefore $h_0$ has at least one critical
point on $I_2$, since $h_0$ is differentiable on $I_2$.
 So to prove the lemma it
suffices to show that this is a unique critical point on $I_2\cup
I_4\backslash\{\lambda_0\}$.

We consider the real valued function
$\gamma:=Q^2/f$ defined on
$I:=I_1\cup I_2\cup I_3\cup I_4$, which is clearly differentiable
on $I$.
Then $h_0=\sqrt{\gamma}+\sqrt{\gamma-1}$ on $I_2\cup I_4$, and we have 
$$h_0'=\gamma'\cdot\left(\frac{1}{2\sqrt{\gamma}}+
\frac{1}{2\sqrt{\gamma-1}}\right),$$
provided
$\lambda\neq\lambda_0$.
 Therefore on $I_2\cup
I_4\backslash\{\lambda_0\}$,  $h_0'(\lambda)=0$ iff
$\gamma'(\lambda)=0$. 
It is readily seen  that
$\lim_{\lambda\da-\infty}\gamma(\lambda)=
\lim_{\lambda\ua-1}\gamma(\lambda)=-\infty$,
$\lim_{\lambda\da-1}\gamma(\lambda)=\lim_{\lambda\ua0}\gamma(\lambda)=\infty$,
$\lim_{\lambda\da 0}\gamma(\lambda)=\lim_{\lambda\ua
(a)}\gamma(\lambda)=-\infty$, and
$\lim_{\lambda\da
(a)}\gamma(\lambda)=\lim_{\lambda\ua\infty}\gamma(\lambda)=\infty$.  
Therefore
$\gamma$ has at least one critical point on each $I_j$, $1\leq j\leq 4$.
We also have  
\begin{equation}\label{eqn-diff}
 \gamma'=Q(2Q'f-Qf')/f^2.
\end{equation}
Suppose that  critical points of  $\gamma$ on $I_2$ are not unique.
 Then
$\gamma$ has at least three critical points on $I_2$.
This implies that $\gamma$ has at least four critical points on $I_2\cup I_4$. 
Because $Q>0$  on $I_2\cup I_4$ (Proposition
\ref{prop-realpoint} (i)), 
these critical points must correspond to zeros of
$2Q'f-Qf'$ whose degree is just four.
By (\ref{eqn-diff})
this implies that the other  critical points of $\gamma$
on
$I_1$ and $I_3$ must correspond to zeros of $Q$.
But this cannot happen since $Q>0$  on
$I_2\cup I_4$ and since $Q$ is degree two.
Therefore, our assumption fails and it follows that critical points on $\gamma$ 
on $I_2$ is unique.
 Hence  critical points of $h_0$ on $I_2$ is also
unique.
 Exactly the same argument shows that
$\gamma$  has a unique critical point on $I_4$.
This critical point must be $\lambda_0$, since 
$\gamma$ attains the minimal value $(=1)$ there.
 This implies that
$g$ has no critical point on $I_4\backslash\{\lambda_0\}$.
 Thus we
obtain the claims of the lemma.
\proofend

\vspace{3mm} We use the lemma to prove the   main result of this subsection: 

\begin{prop}\label{prop-nottl} (i) If $\lambda\in I_4$ and if
$\lambda\neq
\lambda_0$, we  have $N_{L/Z}\simeq O(1)^{\oplus 2}$ for any
$L\in\mathcal L_{\lambda}^+\cup\mathcal L_{\lambda}^-$.
 (ii) There is a
unique
$\lambda\in I_2$ such that
$N_{L/Z}\simeq O\oplus O(2)$ for any $L\in\mathcal
L_{\lambda}^+\cup\mathcal L_{\lambda}^-$.
 For any other $\lambda\in
I_2$,  we  have $N_{L/Z}\simeq O(1)^{\oplus 2}$ for arbitrary $L\in\mathcal
L_{\lambda}^+\cup\mathcal L_{\lambda}^-$.
  (iii) If $\lambda\in I_2$, 
any  member of $\mathcal L_{\lambda}^+\cup\mathcal L_{\lambda}^-$ is not
a twistor line in $Z$ (even if
$Z$ is actually a twistor space).
\end{prop}

\noindent Proof.
 (i) and (ii) are direct consequences of Proposition
\ref{prop-nbc} and Lemma
\ref{lemma-bhv-1}.
  To show (iii), let $\lambda'\in I_2$ be the unique
critical point of $g$.
 Then by (ii), any $L'\in\mathcal
L_{\lambda'}^+\cup\mathcal L_{\lambda'}^- $ is not a twistor line. 
We can see that for any
$\lambda\in I_2$ and for any $L\in\mathcal L_{\lambda}^+\cup\mathcal
L_{\lambda}^-$, $L$ can be deformed into some $L'\in\mathcal
L_{\lambda'}^+\cup\mathcal L_{\lambda'}^- $ preserving the real
structure.
 In fact, we have $\Phi(L)=C_{\theta}$ for some
$C_{\theta}\in\mathcal C_{\lambda}^{\,\rm{gen}}$.
  Then since $I_2$ is an interval in
$\mathbf R$,
$C_{\theta}$ can be canonically deformed into some
$C_{\theta}'\in \mathcal C_{\lambda'}^{\,\rm{gen}}$.
 (The point is that we take
a constant
$\theta$ for any $\lambda\in I_2$.)
  Correspondingly, we obtain
deformation of
$L$ into  $L'\in\mathcal L_{\lambda'}^+\cup\mathcal L_{\lambda'}^- $
such that
$\Phi(L')=C'_{\theta}$.
 Thus we get an explicit real one-dimensional
 family of  rational curves in
$Z$ containing $L$ and $L'$ as its members,
 as claimed.
  Any member of this family is
 real by  Proposition \ref{prop-inv-a} (iv).
   Since any
deformation of twistor line preserving the real structure is still a
twistor line, it follows that
$L$ is not a twistor line.
\proofend
 
\vspace{3mm} We note the proof of (iii) does not work for $I_4$, since as
$\lambda$ goes to
$\lambda_0$, the curve
$C_{\theta}$ (defined by (\ref{eqn-a-0})) degenerates into a double line.
This is an important point for our global construction of arbitrary twistor lines in $Z$.

\begin{cor}\label{cor-image1} If $\lambda\in I_2$, only the members of
$\mathcal C_{\lambda}^{\,\rm{orb}}$ can be the image of twistor lines.
Namely  over $I_2$, members of $\mathcal C_{\lambda}^{\,\rm{gen}}$ cannot be the
images of twistor lines.
\end{cor}

\noindent Proof.
By Proposition \ref{prop-type1} (i), the image of a twistor line contained
in
$H_{\lambda}$ is either a touching conic of generic type or that of orbit
type for $\lambda\in I_2\cup I_4$.
 But by Proposition \ref{prop-nottl}
(iii) the former cannot be the image of a twistor line if $\lambda\in
I_2$.
\proofend

\begin{cor}\label{cor-chatl} A twistor line of a self-dual 4-manifold (i.e.\! a
fiber of the twistor fibration) is not in general characterized by  the property
that it is a real smooth rational curve without real point whose normal bundle is
isomorphic to $O(1)^{\oplus 2}$. 
More concretely, the twistor space of any non-LeBrun self-dual metric on
$3\mathbf{CP}^2$ of positive scalar curvature with a non-trivial 
Killing field always possesses such a real rational curve.
\end{cor}

\noindent Proof. 
Let $Z$ be a twistor space as in the corollary.
Then $Z$ has a
structure as in  Proposition \ref{prop-def-B}, where $Q$ and $a$
satisfy the conditions in Proposition \ref{prop-necessa}.
 By (ii) and (iii)
of Proposition
\ref{prop-nottl}, $Z$ always has a real smooth rational
curve $L$ satisfying
$N_{L/Z}\simeq O(1)^{\oplus 2}$, but which is not a twistor line.
This $L$ has no real point by Proposition \ref{prop-a} (iii) and
the reality of $\Phi$.
 \proofend

\vspace{3mm} 
Next we give another geometric proof for the fact that $L$ cannot be a
twistor line for
$\lambda\in I_2$ (although we will not need this result in the sequel).

\begin{prop} \label{prop-break} If
$\lambda\in I_2$  is not a critical point of $h_0$, there exists a unique
$\mu\in I_2$ with
$\lambda\neq\mu$ satisfying the following: for any $L\in\mathcal
L_{\lambda}^+$ (resp.\,$L\in\mathcal L^-_{\lambda}$) there exists
$L'\in\mathcal L_{\mu}^+$ (resp.\,$L'\in\mathcal L^-_{\mu}$) such that $L\cap
L'\neq\phi$. 
\end{prop}

\noindent Proof.
 Let $\lambda'\in I_2$ be the unique critical point of
$h_0$ as before.
 By our proof of Lemma \ref{lemma-bhv-1} we have
$\lim_{\lambda\downarrow -1}h_0(\lambda)=\lim_{\lambda\uparrow
0}h_0(\lambda)=+\infty$ and
$g$ is strictly decreasing on $(-1,\lambda')$ and strictly
increasing on
$(\lambda',0)$.
 Suppose $\lambda<\lambda'$.
  Let $\Xi$ and $\ol{\Xi}$ be the 
conjugate pair of rational curves which are mapped biholomorphically onto
$l_{\infty}$.
 (See Proposition \ref{prop-image}.)
   Then by
Lemma
\ref{lemma-int-a}, $L\cap \Xi$ is a point  which is  either
$x_2=h_0(\lambda) e^{i\theta}$ or
$x_2=h_0(\lambda)^{-1}e^{i\theta}$ for some 
$\theta\in\mathbf R$,  where we identify $\Xi$ and
$l_{\infty}$ via $\Phi$ and use  $x_2=y_2/y_3$  as an affine coordinate on
$l_{\infty}$ as before.
 If $x_2=h_0(\lambda)^{-1}e^{i\theta}$, 
$\ol{\Xi}\cap L$ is a point having $x_2=h_0(\lambda)e^{i\theta}$.
 Thus (by
a possible exchange of $\Xi$ and
$\ol{\Xi}$)  we may suppose that $\Xi\cap L$ is a point satisfying
$x_2=h_0(\lambda)e^{i\theta}$.
 Then by the behavior of $h_0$ mentioned above,
 there exists a unique
$\mu>\lambda'$, $\mu\in I_2$  such that $h_0(\lambda)=h_0(\mu)$.
 On the
other hand, by our choice of $L$ we have $\Phi(L)=C_{\theta}$ for some
$C_{\theta}\in\mathcal C_{\lambda}^{\,\rm{orb}}$.
 Then take $L'\in\mathcal L_{\mu}^+$
such that
$\Phi(L')=C_{\theta}\in \mathcal C_{\mu}^{\,\rm{orb}}$.
 (Although we use the same symbol
$C_{\theta}$, they represent different conics since $\lambda\neq \mu$.
The point is that we take the same $\theta$ for different $\lambda$'s.)
Then
$L\cap L'\cap \Xi$ is  a point satisfying
$x_2=h_0(\lambda)e^{i\theta}$. 
(We also have $L\cap L'\cap \ol{\Xi}$ is a point satisfying
 $x_2=h_0(\lambda)^{-1}e^{i\theta}$.)
Thus we have proved the claim for
$\lambda<\lambda'$. Of course, the case $\lambda>\lambda'$
and the case $L\in \mathcal L^{-}_{\lambda}$ are similar.
\proofend
 
\vspace{3mm} The proposition shows that when $\lambda\in I_2$ passes
through the critical point ($=\lambda'$) of
$h_0$,
 the local  twistor fibration arising from $L\in\mathcal L_{\lambda}^+\cup
\mathcal L^-_{\lambda}$, $\lambda\neq \lambda'$ breaks down.
 Note also
that Proposition
\ref{prop-break} also holds for $I_4$ without any change of the
proof, and it implies that if members of $\mathcal L_{\lambda}^+$ (resp.
\!$\mathcal L_{\lambda}^-$) are twistor lines for $\lambda\in I_4^-$, 
 members of $\mathcal L_{\lambda}^-$ (resp. $\mathcal L_{\lambda}^+$)
must be  twistor lines for $\lambda \in I_4^+$.

\subsection{The case of special type}\label{ss-special} In this subsection we
calculate the normal bundle of $L^+$ and $L^-$ in $Z$, where
$L^+$ and $L^-$ are  curves which are mapped biholomorphically onto a real
touching conic of special type.
 Compared to the case of generic type, the
problem becomes  harder and the result becomes  more complicated, since
touching conics of special type go through the singular point
$P_{\infty}$ and
$\ol{P}_{\infty}$ of
$B$, so that  the situation, and hence  the result also,  depend on how
we resolve the corresponding singularities of $Z_0$.

\begin{figure}[htbp]
\includegraphics{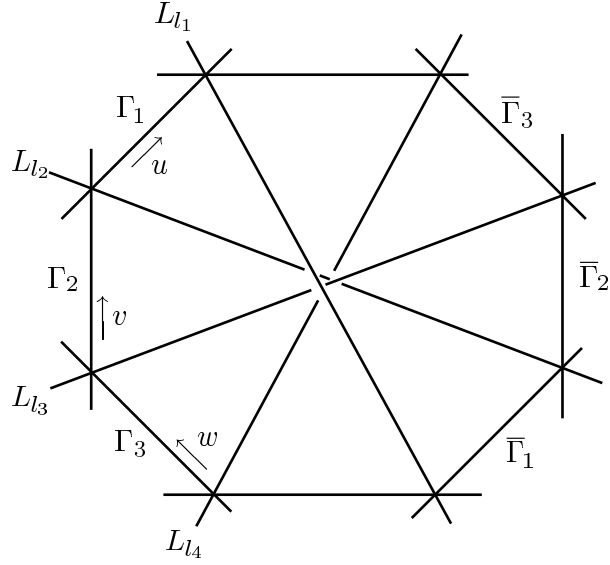}
\caption{the inverse image of $l_{\infty}$ and four tropes}
\label{fig-cycle3}
\end{figure}

First we recall the situation and fix notations.
 Let $\Phi_0:Z_0\ra
\mathbf{CP}^3$ be the double covering branched along $B$.
 Put
$p_{\infty}:=\Phi_0^{-1}(P_{\infty})$.
 In a neighborhood of $P_{\infty}=(0:0:0:1)$, we use 
$(x_0,x_1,x_2)$ as an affine coordinate by setting $x_i=y_i/y_3$.
 Then around
$P_{\infty}=(0,0,0)$,
$B$ is given by the equation
$(x_2+Q(x_0,x_1))^2-x_0x_1(x_0+x_1)(x_0-ax_1)=0$.
  Let $z$ be a fiber
coordinate of $O(2)\ra\mathbf{CP}^3$.
 Then $Z_0$ is given by the
equation 
\begin{equation}\label{eqn-Z_0-2}
z^2+\left(x_2+Q(x_0,x_1)\right)^2-x_0x_1(x_0+x_1)(x_0-ax_1)=0.
\end{equation} This can be also written as
$\{z+i(x_2+Q(x_0,x_1))\}\{z-i(x_2+Q(x_0,x_1))\}=x_0x_1(x_0+x_1)(x_0-ax_1)$.
Setting $\xi=z+i(x_2+Q(x_0,x_1))$ and $\eta=z-i(x_2+Q(x_0,x_1))$, we get
\begin{equation} Z_0:\hspace{3mm}\xi\eta=x_0x_1(x_0+x_1)(x_0-ax_1).
\end{equation} Thus $p_{\infty}=\{(x_0,x_1,\xi,\eta)=(0,0,0,0)\}$ is a compound
$A_3$-singularity.
 Small resolutions of $p_{\infty}$ are explicitly
given by steps as follows: first we choose ordered three linear forms
$\{\ell_1,\ell_2,\ell_3\}\subset\{x_0,x_1,x_0+x_1,x_0-ax_1\}$.
 Next
blow-up
$Z_0$ along a surface
$\{\xi=\ell_1=0\}$.
 Since this surface is contained in the threefold $Z_0$,
and since it goes through $p_{\infty}$, the exceptional locus is 
a smooth rational curve.
 Concretely,
 by setting $\xi=u\ell_1$ we get a threefold
\begin{equation}\label{eqn-0001}
 u\eta=x_0x_1(x_0+x_1)(x_0-ax_1)/\ell_1
\end{equation} which has a compound $A_2$-singularity at the origin.
The exceptional curve  is given by
$\Gamma_1:=\{(u,\eta,x_0,x_1)\set \eta=x_0=x_1=0\}$.
 We can use $u$
as an affine coordinate on $\Gamma_1$.
 Next, as the center of the  second blow-up, we choose $\{u=\ell_2=0\}$  
 which is also contained in the threefold (\ref{eqn-0001}).
 Setting
$u=v\ell_2$, we get
\begin{equation} \label{eqn-0002}
v\eta=x_0x_1(x_0+x_1)(x_0-ax_1)/\ell_1\ell_2
\end{equation} 
which still has an ordinary double point at the origin.
The exceptional curve of this second blow-up is given by
$\Gamma_2:=\{(v,\eta,x_0,x_1)\set \eta=x_0=x_1=0\}$, on which we can
use $v$ as an affine coordinate.
 Finally by blowing up along $\{v=\ell_3=0\}$ and setting $v=w\ell_3$,
we get 
\begin{equation} w\eta=x_0x_1(x_0+x_1)(x_0-ax_1)/\ell_1\ell_2\ell_3,
\end{equation} which is clearly smooth in a neighborhood of the origin.
 The
exceptional curve of the last small resolution is given by
$\Gamma_3:=\{(w,\eta,x_0,x_1)\set
 \eta=x_0=x_1=0\}$, on which we can use $w$ as an affine coordinate. 
Thus the composition of these three blow-ups provides a small resolution 
of $p_{\infty}$, whose exceptional curve is $\Gamma:=\Gamma_1+\Gamma_2
+\Gamma_3$ that is a chain of smooth rational curves.
(See Figure \ref{fig-cycle3}.)

Note that this small resolution of the singularity $p_{\infty}$ is 
easily seen to be 
$U(1)$-equivariant, since our $U(1)$-action is written as
$$(\xi, \eta,x_0,x_1)\mapsto (e^{2i\theta}\xi,e^{2i\theta}\eta,
e^{i\theta} x_0, e^{i\theta}x_1).$$

Once a resolution of $p_{\infty}$ is given, it naturally determines that of 
$\ol{p}_{\infty}$
by reality.
 Let $\ol{\Gamma}=
\ol{\Gamma}_1+\ol{\Gamma}_2+\ol{\Gamma}_3$ be the exceptional curve
over $\ol{p}_{\infty}$.
Let $\mu_{\infty}:Z'_0\ra Z_0$ be the small resolution of $p_{\infty}$ and $\ol{p}_{\infty}$
 produced in this way
(for some choice of $\ell_1$, $\ell_2$ and $\ell_3$).
Note that $Z'_0$ yet has a real  ordinary double point over the unique real double 
point of $B$ (which was denoted by $P_0$).
Each choice of $\ell_1,\ell_2,\ell_3$ (and
$\ell_4$) determines a small resolution of $p_{\infty}$
and $\ol{p}_{\infty}$ and there are
$4!=24$ ways of resolutions in all.

\vspace{3mm}
Next we prove the following lemma which is promised in the proof of
Proposition \ref{prop-class}.
Let $H_{\ell_i}\subset\mathbf{CP}^3$ ($1\leq i\leq 4$) be a real plane defined by $\ell_i=0$.
Let $T_{\ell_i}\subset H_{\ell_i}$ be the tropes of $B$. 
 $\Phi^{-1}(T_{\ell_i})$ is defined by $\xi=\eta=\ell_i=0$ around $p_{\infty}$.
 $\Phi_0^{-1}(H_{\ell_i})$ consists of two irreducible components $H_{\ell_i}'$ and $\ol{H}_{\ell_i}'$, both of which are $\mathbf C^*$-equivariantly biholomorphic to $H_{\ell_i}$ ($\simeq \mathbf{CP}^2$).
 (These are defined by $\xi=\ell_i=0$ and $\eta=\ell_i=0$ around $p_{\infty}$.
 We have $\Phi^{-1}_0(T_{\ell_i})=H_{\ell_i}'\cap \ol{H}_{\ell_i}'$.)
Therefore $\mu_{\infty}^{-1}( \Phi_0^{-1}(H_{\ell_i}))$ also consists of two irreducible components, which we write $D_{\ell_i}$ and $\ol{D}_{\ell_i}$.

%

\begin{lemma}\label{lemma-promise-1}
Consider a\, $\mathbf C^*$-action on $\mathbf{CP}^2$ defined by
$(y_1,y_2,y_3)\mapsto (y_1,ty_2,t^{-1}y_3)$, $t\in\mathbf C^*$, which has smooth conics as the closure of  general orbits.
group of $\mathbf{CP}^2$.) Let $C$ be  one of such conics. 
Choose any one of the  two $\mathbf C^*$-fixed points on $C$
and let $D\ra \mathbf{CP}^2$ be the composition of three blow-ups
whose centers are always  the  $\mathbf C^*$-fixed point on 
(the strict transform of) $C$ just chosen.
Then the natural $\mathbf C^*$-equivariant morphism $D_{\ell_4}\ra H_{\ell_4}$ is 
equivariantly biholomorphic to the morphism $D\ra\mathbf{CP}^2$.
\end{lemma}

In particular, $Z$ always contains a $\mathbf C^*$-invariant divisor $D_{\ell_4}$ which
is $U(1)$-equivariantly biholomorphic to $D$, and such that $D_{\ell_4}+\ol{D}_{\ell_4}$ is a fundamental divisor.

\vspace{3mm}
\noindent Proof of the lemma.
First, notice that for any $\mathbf C^*$-invariant plane $H\supset l_{\infty}$, our $\mathbf C^*$-action coincides the $\mathbf C^*$-action given in the lemma.
We prove the claim for $D_{\ell_4}$ by direct calculation using the local coordinate
$(\xi,\eta,x_0,x_1)$ introduced above.
$D_{\ell_4}$ can be assumed to be the inverse image (by $\mu$) of the surface
$H_{\ell_4}'=\{\eta=\ell_4=0\}\subset Z_0$.
As is described above, we first blow-up $Z_0$ along $\{\xi=\ell_1=0\}$.
Then since we can use $(\xi,\ell_1)$ as a local coordinate on $H_{\ell_4}'$ around $p_{\infty}$(= the origin), 
$H_{\ell_4}'$ is actually blown-up at  $p_{\infty}$.
Then (the inverse image of) $H_{\ell_4}'$ is still defined by $\{(u,\eta,x_0,x_1)\set \eta=\ell_4=0\}$,
where $u=\xi/\ell_1$ is an affine coordinate on the exceptional curve as before.
 Hence for the second blow-up whose center is
$\{u=\ell_2=0\}$, the origin of the surface $\{\eta=\ell_4=0\}$ is again blown-up.
The situation is the same for the third blow-up.
 Thus the original surface 
$H'_{\ell_4}\subset Z_0$ is blown-up three times at the $\mathbf C^*$-fixed point.
Moreover, the blown-up points are always over the trope on $H_{\ell_4}$, 
as is easily verified from the fact that
on the surface $H'_{\ell_4}$, $\Phi^{-1}_0(T_{\ell_4})$ is given by $\xi=0$.
 (Note that our coordinate change is not linear and
a conic on $H_{\ell_i}$ is transformed into a `line' on $H'_{\ell_4}$ with
respect to the coordinate $(\xi, \ell_1)$.)

So far we have only seen what happens around $p_{\infty}$.
Around $\ol{p}_{\infty}$, we claim that on the surface $H'_{\ell_4}$ nothing change through our blow-ups.
 To see this, set $x'_i=y_i/y_2$ for $i=0,1,3$.
 Then in a neighborhood of  $\ol{p}_{\infty}$, $Z_0$ is defined by
$$
z^2+(x'_3+Q(x_0',x'_1))^2-x'_0x'_1(x'_0+x'_1)(x'_0-ax'_1)=0.
$$
Thus setting $\xi'=z+i(x'_3+Q(x_0',x'_1))$ and $\eta'=z-i(x'_3+Q(x_0',x'_1))$,
we can write
$$
Z_0: \hspace{3mm}\xi'\eta'=x'_0x'_1(x'_0+x'_1)(x'_0-ax'_1).
$$
Then since we have $\sigma^*z=\ol{z}$, $\sigma^*x_0=x_0'$, $\sigma x_1=\ol{x}_1'$  and $\sigma^*x_2=\ol{x}_3'$, we have
$$
\sigma^*\xi =\sigma^*\left(z+i(x_2+Q(x_0,x_1))\right)
=\ol{z}+i(\ol{x}_3+Q(\ol{x}_0,\ol{x}_1))$$
$$
=\ol{z}+i(\ol{x}_3+Q(\ol{x}_0,\ol{x}_1))=
\ol{z-i(x_3+Q(x_0,x_1)) }=\ol{\eta}'.
$$
Similarly we have $\sigma^*\eta=\ol{\xi}'$.
On the other hand, we easily get $\sigma^*\ell_i=\ol{\ell}_i$ for $1\leq i\leq 4$.
Hence $\ol{H}_{\ell_4}'=\sigma(H_{\ell_4}')$  is defined by $\xi'=\ell_4=0$ in a neighborhood of $\ol{p}_{\infty}$.
This implies that $H_{\ell_4}'$ is defined by $\eta'=\ell_4=0$ around $\ol{p}_{\infty}$.
  On the other hand, by reality,
the center of the first blow-up must be $\{\eta'=\ell_1=0\}$, which intersect $H_4'$
transversally along a smooth curve.
 Therefore, $H'_{\ell_4}$ has no effect under the first
blow-up of $Z_0$.
 Similarly, nothing happens for the remaining two blow-ups.
Combined with what we have seen around $p_{\infty}$, we have shown that 
$D_{\ell_4}\ra H_{\ell_4}$ is $\mathbf C^*$-equivariantly biholomorphic to $D$ in the lemma.
\proofend

\vspace{3mm}
By similar direct calculations, we can show the following 

\begin{lemma}\label{lemma-degree01}
Let $\mu^{-1}_{\infty}(\Phi^{-1}_0(H_{\ell_i}))=D_{\ell_i}+\ol{D}_{\ell_i}$  and $T_{\ell_i}\subset H_{\ell_i}$ $(1\leq i\leq 4)$ be as above. 
Then $L_{\ell_i}=D_{\ell_i}\cap \ol{D}_{\ell_i}$ is a real smooth rational curve in $Z'_0$ satisfying the following properties: (i) $L_{\ell_i}$ are mapped biholomorphically onto  $T_{\ell_i}$, (ii) the  normal bundle of $L_{\ell_i}$ in $Z'_0$ is isomorphic to $O(1)^{\oplus 2}$,
(iii) the intersection of $L_{\ell_i}$ with $\Phi^{-1}(l_{\infty})$ is as in Figure \ref{fig-cycle3}.
\end{lemma}

\noindent Sketch of the proof.
We use the local coordinates in the previous proof.
$\Phi_0^{-1}(H_{\ell_i})$ consists of two irreducible components $H_{\ell_i}'$ and $\ol{H}_{\ell_i}'$. 
These are $\mathbf C^*$-invariant surfaces.
Since the trope $T_{\ell_i}$ is a conic on $H_{\ell_i}$, the normal bundle of  the intersection $H_{\ell_i}'\cap \ol{H}_{\ell_i}'$ in $H_{\ell_i}'$ and $\ol{H}_{\ell_i}'$ are isomorphic to $O(4)$.
We may suppose that $H_{\ell_i}'$ is defined by $\eta=\ell_i=0$ as in the previous proof.
Then if we blow-up along $\xi=\ell_1=0$, $H_{\ell_i}'$ is actually blown-up at $p_{\infty}$ for $2\leq i\leq 4$, whereas nothing happens for $H_{\ell_1}'$.
On the other hand, around $\ol{p}_{\infty}$, nothing happens on $H_{\ell_i}'$ for $2\leq i\leq 4$, whereas $H_{\ell_1}'$ is blown-up at $\ol{p}_{\infty}$.
(Here note that $p_{\infty}$ and $\ol{p}_{\infty}$ are just the fixed locus of $\mathbf C^*$-action on $H_{\ell_i}'\cap \ol{H}_{\ell_i}'$.)
Consequently, the resulting space has two divisors $H_{\ell_i}''$ and $\ol{H}_{\ell_i}''$ that are the inverse images of $H_{\ell_i}'$ and $\ol{H}_{\ell_1}'$ respectively.
The intersection $H_{\ell_i}''\cap\ol{H}_{\ell_i}''$ is a smooth $\mathbf C^*$-invariant rational curves whose normal bundles in $H_{\ell_i}''$ and $\ol{H}_{\ell_i}''$ are $O(3)$. 
  For the second  blow-ups, $H_{\ell_i}''$, $i=1,3,4$ (resp.\,$\ol{H}_{\ell_i}''$, $i=1,3,4$), are blown-up at the $\mathbf C^*$-fixed point on $H_{\ell_i}''\cap\ol{H}_{\ell_i}''$ lying over $p_{\infty}$ (resp.\,$\ol{p}_{\infty}$), and $H_{\ell_2}''$ (resp.\,$\ol{H}_{\ell_2}''$)  is blown-up at the $\mathbf C^*$-fixed point on $H_{\ell_i}''\cap\ol{H}_{\ell_i}''$ lying over $\ol{p}_{\infty}$ (resp.\, $p_{\infty}$).
 Denoting  $H_{\ell_i}^{(3)}$ and $\ol{H}_{\ell_i}^{(3)}$ the inverse image of $H_{\ell_i}''$ and $\ol{H}_{\ell_i}''$ respectively,  the normal bundles of $H_{\ell_i}^{(3)}\cap\ol{H}_{\ell_i}^{(3)}$ in $H_{\ell_i}^{(3)}$ and $\ol{H}_{\ell_i}^{(3)}$ become $O(2)$.
The third (=the last) blow-up yields $D_{\ell_i}$ and $\ol{D}_{\ell_i}$ as the inverse images of $H_{\ell_i}^{(3)}$ and $\ol{H}_{\ell_i}^{(3)}$ respectively.
 $D_{\ell_i}\ra H_{\ell_i}^{(3)}$, $i=1,2,4$ (resp.\,$\ol{D}_{\ell_i}\ra \ol{H}_{\ell_i}^{(3)}$, $i=1,2,4$) is the blow-up at the $\mathbf C^*$-fixed point of $H_{\ell_i}^{(3)}\cap\ol{H}_{\ell_i}^{(3)}$  lying over $p_{\infty}$ (resp.\,$\ol{p}_{\infty})$, whereas $D_{\ell_3}\ra H_{\ell_3}^{(3)}$(resp.\,$\ol{D}_{\ell_3}\ra \ol{H}_{\ell_i}^{(3)}$) is the blow-up at the $U(1)$-fixed point lying over $\ol{p}_{\infty}$ (resp.\,$p_{\infty}$).
Hence the normal bundles of $L_{\ell_i}=D_{\ell_i}\cap \ol{D}_{\ell_i}$ in $D_{\ell_i}$ and $\ol{D}_{\ell_i}$ become isomorphic to $O(1)$. 
Moreover, the intersection of $D_{\ell_i}$ and $\ol{D}_{\ell_i}$ in $Z'_0$ is transversal.
Therefore we obtain $N_{L_{\ell_i}/Z'_0}\simeq O(1)^{\oplus 2}$.
Thus we get (i) and (ii) of the lemma.
(iii) is rather easily verified by direct calculations using the local coordinates above, and we omit the detail. (The case $L_{\ell_4}$ is actually proved in the last lemma.)\proofend

\vspace{3mm}
Next we prove that,
as promised in the proof of Proposition \ref{prop-str_of_S},
 our small resolution $\mu_{\infty}$ gives a resolution  of 
the surface $\Phi_0^{-1}(H_{\lambda})$ which has 
$p_{\infty}$ and $\ol{p}_{\infty}$ as its $A_3$-singularities:
\begin{lemma}\label{lemma-promise-2}
If $\lambda\neq -1,0,a,\infty$,
 $S_{\lambda}=\Phi^{-1}(H_{\lambda})$ is non-singular. 
\end{lemma}

\noindent
Proof. This can be also seen by direct calculation using coordinate 
$(\xi,\eta,x_0,x_1)$ around $p_{\infty}$ above.
 As $H_{\lambda}$ is
defined by $x_0=\lambda x_1$, $\Phi_0^{-1}(H_{\lambda})$ is locally defined
by $\xi\eta=f(\lambda)x_1^4$ around $p_{\infty}$. 
  Substituting $\xi=u\ell_1$ for
the first blow-up, the inverse image becomes 
$u\eta=(f(\lambda)/\ell_1(\lambda,1))x_1^3$.
 Next substituting
$u=v\ell_2$ for the second blow-up, we get 
$v\eta=(f(\lambda)/\ell_1(\lambda,1)\ell_2(\lambda,1))x_1^2$.
Finally substituting $v=w\ell_3$ for the third blow-up, we get
$w\eta=(f(\lambda)/\ell_1(\lambda,1)\ell_2(\lambda,1)
\ell_3(\lambda,1))x_1$.
 This is smooth. 
By reality, the conjugate singular point $\ol{p}_{\infty}$ of $S_{\lambda}$ is 
also resolved by our blow-ups.
 Thus both $p_{\infty}$ and 
$\ol{p}_{\infty}$  are resolved. 
These are of course minimal resolution of the $A_3$-singularity
and the self-intersection number of each irreducible component in $S_{\lambda}$ is $-2$.

If $\lambda\neq \lambda_0$
(and $\lambda\neq  -1,0,a,\infty$), $B\cap H_{\lambda}$ has no singular point
other than $p_{\infty}$ and $\ol{p}_{\infty}$, and it follows that 
$S_{\lambda}$ is smooth provided $\lambda\neq\lambda_0$.
Because $B$ has a real ordinary double point at $P_0=(\lambda_0:1:0:0)$,
$\Phi_0^{-1}(H_{\lambda_0})$ has a real ordinary double point over $P_0$.
But this singular point is also resolved through small resolution of $Z_0$.
This is of course rather simpler than the above case and we omit the calculation.\proofend
\vspace{3mm}

Next in order to calculate the intersection $L^+$ and $L^-$ (cf. the beginning of this subsection) with $\Gamma$,
we  need a one-parameter presentation of $C_{\theta}$, in a neighborhood of
$p_{\infty}$:

\begin{lemma}\label{lemma-one-rep} Let $C_{\theta}\subset H_{\lambda}$
be a real touching conic of special type whose equation is given by
(\ref{eqn-b-0}), and $(x_0,x_1,x_2)$ the affine coordinate around
$P_{\infty}$ as above.
 Then in  a neighborhood of
$P_{\infty}$,
$C_{\theta}$ has a one-parameter presentation of the following form:
\begin{equation}\label{eqn-solve}
\left\{
\begin{array}{l}x_0=\lambda x_1\\ \displaystyle
x_2=-Be^{-i\theta}x_1-\frac{\sqrt{Q^2-f}+Q}{2}x_1^2+
\frac{\sqrt{Q^2-f}+Q}{2}Be^{i\theta}x_1^3+O(x_1^4),
\end{array}\right.
\end{equation}
where we put
$$B:=B(\lambda)=\left(\frac{\sqrt{Q^2-f}-Q}{2}\right)^{\frac{1}{2}}.$$
\end{lemma}

Note again that $f<0$ guarantees $\sqrt{Q^2-f}-Q>0$.

\vspace{3mm}
\noindent Proof. By  solving (\ref{eqn-b-0}) with respect to $x_2$, we
get
\begin{equation}\label{eqn-solv} x_2=-g(x_1)\cdot x_1,\hspace{3mm}
g(x_1):=
\frac{Be^{-i\theta}+\sqrt{Q^2-f}\,\,x_1}{1+Be^{i\theta}x_1}.
\end{equation} Calculating the Maclaurin expansion of $g(x_1)$, we get
(\ref{eqn-solve}).
 This is a routine work and we omit the
detail.\proofend

\begin{lemma}\label{lemma-tlr} In a neighborhood of $p_{\infty}$, each of the two
irreducible components of
$\Phi^{-1}_0(C_{\theta})$ has a one-parameter presentation with
respect to
$x_1$ in the following forms respectively:
\begin{equation}\label{eqn-L1}
\xi=-2iBe^{-i\theta}x_1+O(x_1^2),\hspace{3mm}
\eta=\frac{ie^{i\theta}f}{2B}x_1^3+O(x_1^4),\hspace{3mm} x_0=\lambda x_1,
\end{equation} and
\begin{equation}\label{eqn-L2}
\xi=-\frac{ie^{i\theta}f}{2B}x_1^3+O(x_1^4),\hspace{3mm} 
\eta=2iBx_1e^{-i\theta}+O(x_1^2),\hspace{3mm} x_0=\lambda x_1.
\end{equation}

\end{lemma}

\noindent Proof.
 First by substituting $x_0=\lambda x_1$ into
(\ref{eqn-Z_0-2}), we get
$$ z^2=(f-Q^2)x_1^4-2Qx_1^2x_2-x_2^2.
$$ Substituting (\ref{eqn-solv}) into this, we get 
$$z^2=\left\{(f-Q^2)x_1^2+2Qg(x_1)x_1-g(x_1)^2\right\}\,x_1^2.$$
 Hence we have
$$z=\pm k(x_1)\,x_1,\hspace{3mm}
k(x_1)=\left\{(f-Q^2)x_1^2+2Qg(x_1)x_1-g(x_1)^2
\right\}^{\frac{1}{2}}.$$
 From this
we deduce 
$$
\xi=z+i(x_2+Qx_1^2)=\left(\pm k(x_1)-ig(x_1)\right)\,x_1+iQ x_1^2
$$  and 
$$\eta=z-i(x_2+Qx_1^2)=\left(\pm k(x_1)+ig(x_1)\right)\,x_1-iQ x_1^2.$$ Then we get
the  desired equations by calculating the Maclaurin expansions of 
$\pm k(x_1)-ig(x_1)$ and $\pm k(x_1)+ig(x_1)$.
 These are also  routine works and we omit the detail.
\proofend

\vspace{3mm}

\begin{lemma}\label{lemma-interse} Let $L^+_{\theta}$ and
$L^-_{\theta}$  be the curves in $Z'_0$ which are the proper transforms of the
curves (\ref{eqn-L1}) and (\ref{eqn-L2}) respectively.
  Then we have: (i)
$L^+_{\theta}$ and $\Gamma_1$ intersect transversally at  a unique point satisfying 
\begin{equation}\label{eqn-inters}
u=-2iBe^{-i\theta}\cdot\frac{x_1}{\ell_1}
\end{equation}
and $L^+_{\theta}\cap
\Gamma_2$ and $L^+_{\theta}\cap\Gamma_3$ are empty, (ii)
$L^-_{\theta}\cap\Gamma_1$ and $L^-_{\theta}\cap\Gamma_2$ are empty and
$L^-_{\theta}$ and $\Gamma_3$ intersect transversally at  a unique point satisfying
\begin{equation}\label{eqn-inter2}
w=-\frac{ie^{i\theta}f}{2B}\cdot
\frac{x_1^3}{\ell_1\ell_2\ell_3}
\end{equation} 
where we use $u$ and
$w$ as local coordinates on $\Gamma_1$ and $\Gamma_3$ respectively as
explained before, and $B=B(\lambda)$ is as in Lemma
\ref{lemma-one-rep}.  
\end{lemma}

Here note that
$x_1/\ell_1$ and $x_1^3/\ell_1\ell_2\ell_3$ do not depend on $x_1$, and
depend  on $\lambda$ only.
 Further, we have $B>0$ since $f<0$.

\vspace{3mm}
\noindent Proof.
 By substituting $\xi=u\ell_1$ into (\ref{eqn-L1}), we
get the inverse image of (\ref{eqn-L1}) to be
$$ u\ell_1=-2iBe^{-i\theta}x_1+O(x_1^2),\hspace{3mm}
\eta=\frac{ie^{i\theta}f}{2B}x_1^3+O(x_1^4),\hspace{3mm} x_0=\lambda x_1.
$$ Excluding the equation of  $\Gamma_1=\{
\eta=x_0=x_1=0\}$, we get the equation of the
proper transform in the threefold (\ref{eqn-0001}) to be
$$ u=-2iBe^{-i\theta}\frac{x_1}{\ell_1}+O(x_1),\hspace{3mm}
\eta=\frac{ie^{i\theta}f}{2B}x_1^3+O(x_1^4),\hspace{3mm} x_0=\lambda x_1.
$$ By setting $x_1=0$, we get $u=-2iBe^{-i\theta}\cdot x_1/\ell_1$,
and the intersection (of $\Gamma_1$ and $L^+_{\theta}$) is transversal.
Finally as remarked above, $B$ is non-zero.
 Therefore the remaining two blow-ups 
 do not have effect on the intersection.
  Hence we get
(\ref{eqn-inters}).
 Similar calculations show (\ref{eqn-inter2}).
  (Note
that we have
$\xi=w\,\ell_1\ell_2\ell_3$.)
 \proofend

\vspace{3mm} As in the previous subsection we put $\mathcal
L_{\lambda}^+=\{L^+_{\theta}\set\theta\in\mathbf R\}$ and $\mathcal
L_{\lambda}^-=\{L^-_{\theta}\set\theta\in\mathbf R\}$.
 (Note that this
time we explicitly specified $L^+_{\theta}$ and $L^-_{\theta}$
respectively in Lemma \ref{lemma-tlr}.)
$U(1)$ again acts transitively  on the parameter spaces of these
families.
  By Lemma
\ref{lemma-interse},
$L^+_{\theta}\cap \Gamma_1$ is a point, and 
$\{L^+_{\theta}\cap \Gamma_1\set\theta\in\mathbf R\}$ is a circle in
$\Gamma_1$ whose radius is
$$h_1(\lambda):=2B\cdot |x_1/\ell_1|.$$
 Similarly,
$\{L^-_{\theta}\cap
\Gamma_3\set\theta\in\mathbf R\}$ is a circle in $\Gamma_3$ whose radius
is
$$h_3(\lambda):= (-f/2B)\cdot |x_1^3/\ell_1\ell_2\ell_3|.$$ 
($h_2$ will appear in the next subsection.)
 Then  we
have the following proposition, which implies that the normal bundles of
$L^+_{\theta}$ and
$L^-_{\theta}$  in $Z'_0$ are  determined by the behavior of
$h_1$ and $h_3$ respectively.

\begin{prop}\label{prop-nbc2} Assume $\lambda\in I_1\cup I_3$ and take
any  $L=L^+_{\theta}\in\mathcal L_{\lambda}^+$.
 Then either
$N_{L/Z'_0}\simeq O(1)^{\oplus 2}$ or
$N_{L/Z'_0}\simeq O\oplus O(2)$ holds.
 Further, the latter holds iff
$\lambda$ is a critical point of
$h_1$ above.
 The same claims  also hold  for $\mathcal L^-_{\lambda}$ if
we replace $h_1$ by $h_3$. 
\end{prop}

\noindent Proof. The first claim can be proved in the same way as in 
Proposition \ref{prop-nbc}.
The other claims can also be  proved in the same manner as
in Proposition \ref{prop-nbc}:
take any $L\in \mathcal L_{\lambda}^+$.
 Then the two
one-parameter families of $L$ in $Z'_0$  in the previous proof make senses
also in this case, so that  we again have two linearly independent
sections $s$ and $t$ of $N_{\lambda}$, $N_{\lambda}=N_{L/Z'_0}$.
Then the previous proof works if we replace $h_0$ by $h_1$, 
$\Phi^{-1}(l_{\infty})$ by $\Gamma_1\cup \ol{\Gamma}_1$, and
(\ref{eqn-intersec100}) by (\ref{eqn-inters}).
For $L\in \mathcal L_{\lambda}^-$,
replace $h_0$ by $h_3$, 
$\Phi^{-1}(l_{\infty})$ by $\Gamma_3\cup \ol{\Gamma}_3$, and
(\ref{eqn-intersec100}) by (\ref{eqn-inter2}).
\proofend

\vspace{3mm} By definition, $h_1$ and $h_3$ depend on the choice of
$\ell_1,\ell_2$ and $\ell_3$.
 Therefore, Proposition \ref{prop-nbc2}
implies that 
 the normal bundles of $L^+_{\theta}$ and $L^-_{\theta}$  in $Z'_0$ depend on how we resolve
$p_{\infty}$.
 More precisely, the normal bundle of $L^+_{\theta}$ depends on the
choice of $\ell_1$ only, whereas  the normal bundle  of $L^-_{\theta}$ depends on
that of
$\{\ell_1,\ell_2,\ell_3\}$.

Thus we need to know the critical points of $h_1$ and $h_3$ for every
choices of
$\ell_1,\ell_2$ and $\ell_3$.
  At first glance there may seem to be too
many functions to be investigated, but it is easily seen that $h_3$ is a
reciprocal of $h_1$ (up to a constant) for some other choice of
$\ell_1,\ell_2$ and $\ell_3$.
 Consequently, what we need to know is the
behavior of $h_1$ for  the  four choices of $\ell_1$.
(Behavior of these functions near the endpoints of $I_1$ and $I_3$ below will be needed  in \S \ref{ss-con}.)

\begin{figure}[htbp]
\includegraphics{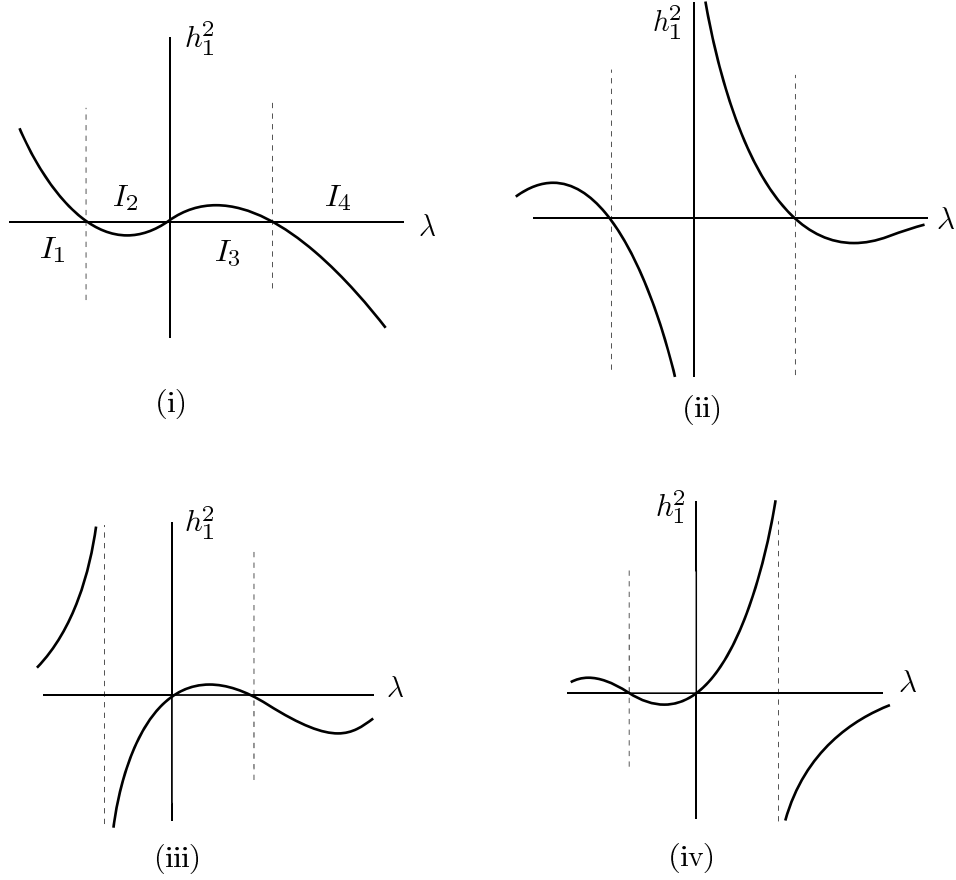}
\caption{behaviors of $h_1^2$}
\label{fig-h1}
\end{figure}

\begin{lemma}\label{lemma-crtcl-sp} (i) If $\ell_1=x_1$,
$h_1$ has no critical point on
$I_1$, and has a unique critical point on $I_3$.
  Further, we have
$\lim_{\lambda
\da-\infty}h_1(\lambda)=+\infty$ and $h_1(-1)=0$.
 (ii)  If $\ell_1=x_0$,
$h_1$ has a unique critical point on
$I_1$, and no critical point on $I_3$.
 Further, we have $\lim_{\lambda
\da 0}h_1(\lambda)=+\infty$ and $h_1(a)=0$.
 (iii) If
$\ell_1=x_0+x_1$,
$h_1$ has no critical point on
$I_1$, and has a unique critical point on $I_3$.
  Further, we have
$\lim_{\lambda
\da-\infty}h_1(\lambda)=0$ and $\lim_{\lambda\uparrow
-1}h_1(\lambda)=+\infty$.
 (iv)  If $\ell_1=x_0-ax_1$, $h_1$ has a
unique critical point on
$I_1$, and no critical point on $I_3$.
 Further, we have $h_1(0)=0$ and
$\lim_{\lambda\uparrow a}h_1(\lambda)=+\infty$.
(See Figure \ref{fig-h1}.)
\end{lemma}

\noindent Proof.  (i) If
$\ell_1=x_1$, we have $h_1^2=2(\sqrt{Q^2-f}-Q)$.
Since $h_1>0$ on $I_1\cup I_3$, the critical points of 
$h_1^2$ and $h_1$ coincide on $I_1\cup I_3$.
We think of $h^2_1$  as a real valued function defined on the whole
of $\mathbf R$, but which is not differentiable at $\lambda=\lambda_0$
in general.
 It is immediate to see that
$h_1^2(-1)=h_1^2(0)=h_1^2(a)=0$,
$\lim_{\lambda\da-\infty}h_1^2(\lambda)=+\infty$ and 
$\lim_{\lambda\ua\infty}h_1^2(\lambda)=-\infty$.
Then because $h_1^2$ is differentiable on $\lambda\neq\lambda_0$,
$h_1^2$ has a critical point on $I_2$ and $I_3$ respectively.
On the other hand, we have
$$
\left(\sqrt{Q^2-f}-Q\right)'=\frac{2QQ'-f'-2Q'\sqrt{Q^2-f}}{2\sqrt{Q^2-f}},
$$ and it follows that $(\sqrt{Q^2-f}-Q)'=0$ implies
\begin{equation}\label{eqn-diff00}
(2QQ'-f')^2=4Q'\
\!^2(Q^2-f).
\end{equation} It is readily seen that the degree of both hand sides
of (\ref{eqn-diff00}) are six, and that both have  $(\lambda-\lambda_0)^2$ as a factor.
Since we have already got two critical points of $h_1^2$
other than $\lambda=\lambda_0$, 
there are at most two solutions of (\ref{eqn-diff00})  remaining.
 
We set $g:=2(-\sqrt{Q^2-f}-Q)$ which is also defined on $\mathbf R$ and
possibly not differentiable at $\lambda=\lambda_0$.
Note that if we replace $h_1^2$ by $g$ on $\lambda\geq\lambda_0$,
then the resulting function is differentiable at $\lambda=\lambda_0$.
It is  easily verified that 
$g'=0$ also implies (\ref{eqn-diff00}) and it gives a solution
not coming from $(h_1^2)'=0$.
Further, we readily have 
$\lim_{\lambda\da-\infty}g(\lambda)=
\lim_{\lambda\ua\infty}g(\lambda)=-\infty$.

Suppose that $g$ has a critical point.
Then together with the above two critical points of $h_1^2$ on $I_2\cup
I_3$, we have three solutions of (\ref{eqn-diff00}) other than
$\lambda=\lambda_0$.
 If $h_1^2$ has  critical points on $I_1$, its
number is at least two.
 This implies that (\ref{eqn-diff00}) has five
solutions other than
$\lambda=\lambda_0$ and this is a contradiction.
Therefore $h_1^2$  has no critical points on $I_1$,
if $g$ has a critical point.
Similarly, if the number of the critical points on $I_3$ is not one,
then it must be at least three.
This is again
a contradiction.
 Thus if $g$ has a critical point,
$h_1^2$, and hence $h_1$ has no critical point on $I_1$ and a
unique critical point on $I_3$.

So suppose that $g$ has no critical point.
This happens exactly when $g$ attains the maximal
value at $\lambda=\lambda_0$.
 Then we have $\lim_{\lambda\ua\lambda_0}
g'(\lambda)>0$,
since otherwise $g$ has a critical point on $\lambda<\lambda_0$.
Because we have $\lim_{\lambda\ua\lambda_0}
g'(\lambda)=\lim_{\lambda\da\lambda_0}(h_1^2)'(\lambda)$,
we get $\lim_{\lambda\da\lambda_0}(h_1^2)'(\lambda)>0$.
Since $\lim_{\lambda\ua\infty}h_1^2(\lambda)=-\infty$,
it follows that $h_1^2$ has a critical point on $I_4$.
 Thus we get
three solutions of (\ref{eqn-diff00}) other than $\lambda_0$.
Then the same argument in the case that $g$ has a critical point as above,
we can deduce that $h_1$ has no critical point on $I_1$ and
a unique critical point on $I_3$.
Thus we get the claim of (i) concerning  critical points of $h_1$.
The remaining claims of (i) immediately follows from the definition of
$h_1$.

Claims of 
(ii), (iii) and (iv) about critical points can be obtained by applying a
projective transformation $\lambda\mapsto 1/\lambda$ for the case (ii), 
$\lambda\mapsto 1/(\lambda+1)$ for the case (iii), and 
$\lambda\mapsto 1/(\lambda-a)$ for the case (iv) respectively.
The other claims are immediate to see.
\proofend

\vspace{3mm} As is already mentioned, the behavior of $h_3$ can be
easily seen from that of $h_1$ for some other choice of $\ell_1,\ell_2$
and $\ell_3$.
 The result is the following:

\begin{lemma}\label{lemma-crtcl-sp2} (i) If
$\{\ell_1,\ell_2,\ell_3\}=\{x_0,x_0+x_1,x_0-ax_1\}$, $h_3$ has no
critical point on
$I_1$, and has a unique critical point on $I_3$.
  Further, we have
$\lim_{\lambda
\da-\infty}h_3(\lambda)=0$ and $\lim_{\lambda\ua
-1}h_3(\lambda)=+\infty$.
 (ii)  If
$\{\ell_1,\ell_2,\ell_3\}=\{x_1,x_0+x_1,x_0-ax_1\}$,
 $h_3$ has a unique critical point on
$I_1$, and no critical point on $I_3$.
 Further, we have $h_3(0)=0$ and
$\lim_{\lambda\ua a}h_3(\lambda)=+\infty$.
   (iii) If
$\{\ell_1,\ell_2,\ell_3\}=\{x_0,x_1,x_0-ax_1\}$, $h_3$ has no critical
point on
$I_1$, and has a unique critical point on $I_3$. 
 Further, we have $\lim_{\lambda
\da-\infty}h_3(\lambda)=+\infty$ and $h_3(-1)=0$.
 (iv)  If
$\{\ell_1,\ell_2,\ell_3\}=\{x_0,x_1,x_0+x_1\}$, $h_3$ has a unique
critical point on
$I_1$, and no critical point on $I_3$.
 Further, we have $\lim_{\lambda\da
0}h_3(\lambda)=+\infty$ and
$h_3(a)=0$.
\end{lemma}

Note that in the lemma we do not specify the ordering of $\ell_1,\ell_2$ and $\ell_3$; for instance  $\{\ell_1,\ell_2,\ell_3\}=\{x_0,x_0+x_1,x_0-ax_1\}$ in (i) does not imply $\ell_1=x_0,\ell_2=x_0+x_1$ and $\ell_3=x_0-ax_1$.

\begin{cor}\label{cor-oneof} For any choice of $\ell_1,\ell_2$ and
$\ell_3$,  the following (i) and (ii) hold: (i) 
members of $\mathcal L^+_{\lambda}$ ($\lambda\in I_1$) and  $\mathcal
L^+_{\lambda}$ ($\lambda\in I_3$) cannot be twistor lines at the
same time,  (ii)
the same claim holds also for $\mathcal L^-_{\lambda}$.
\end{cor}

\noindent Proof.
 Suppose that $\lambda\in I_1\cup I_3$ is a critical
point of
$h_1$.
 Then by Proposition \ref{prop-nbc2}, any member of $\mathcal
L^+_{\lambda}$ is not a twistor line because its normal bundle in $Z'_0$ is
$O\oplus O(2)$.
 Then just as in the proof of Proposition
\ref{prop-nottl}, any member of $\mathcal L^+_{\mu}$ cannot be a twistor
line provided that $\mu\in I_1\cup I_3$ and $\lambda$ belong to the same 
interval ($I_1$ or $I_3$).
 By Lemma \ref{lemma-crtcl-sp},
 $h_1$  necessarily  has  a critical point on just one of $I_1$ and
$I_3$.
 Hence (i) holds.
  The proof is the same for $\mathcal
L^-_{\lambda}$ if we use  Lemma \ref{lemma-crtcl-sp2} instead.
\proofend

\vspace{3mm} Thus together with Proposition \ref{prop-nbc2},
 we have obtained  new families of real smooth rational curves which have
$O(1)^{\oplus 2}$ as their normal bundles, but which are not  twistor lines.

 Proposition
\ref{prop-nbc2} and Lemmas
\ref{lemma-crtcl-sp} and
\ref{lemma-crtcl-sp2} enable us to determine the normal bundles of
$L_{\theta}^+$ and
$L_{\theta}^-$ in $Z'_0$ for every choices of small resolutions of
$p_{\infty}$.
 In particular, the normal bundles of $L_{\theta}^+$ and
$L_{\theta}^-$ in $Z'_0$,
and also which component has to be chosen as
 candidates of twistor lines, depend on the choice made.

\subsection{The case of orbit type}\label{ss-orbit} Suppose $\lambda\in I_2\cup I_4$.
In this subsection we calculate the normal bundles of
$L^+_{\alpha}$ and
$L^-_{\alpha}$ in $Z'_0$, where
$L^+_{\alpha}$ and $L^-_{\alpha}$ are  curves which are mapped
biholomorphically onto a real touching conic $C_{\alpha}\in 
\mathcal C_{\lambda}^{\,{\rm{orb}}}$
defined by (\ref{eqn-c}).
   Note again that $C_{\alpha}$
and $L^{\pm}_{\alpha}$ depend not only on $\alpha$, but also on
$\lambda\in I_2\cup I_4$.
 Compared to  generic type and special type,
calculations are much easier since the equations of touching conics of
orbit type are much simpler.

First we make a distinction of
$L^+_{\alpha}$ and 
$L^-_{\alpha}$.
  We use  local   coordinates $(x_0,x_1,x_2,z)$
and
$(x_0,x_1,\xi,\eta)$ as in the previous subsection.
 Recall that 
$\Phi_0^{-1}(H_{\lambda})$ is defined by
$\xi\eta=fx_1^4$, and that  the equation of irreducible components of
$\Phi_0^{-1}(C_{\alpha})$ is  given by $z=\pm
(f-(\alpha+Q)^2)^{\frac{1}{2}}\, x_1^2$ ((\ref{eqn-inv-gen})).
 Then we
denote by $L^+_{\alpha}$ (resp.
$L^-_{\alpha}$)  the components corresponding to 
$z=(f-(\alpha+Q)^2)^{\frac{1}{2}}\, x_1^2$ (resp.
$z=-(f-(\alpha+Q)^2)^{\frac{1}{2}}\, x_1^2$).
 $L^+_{\alpha}$ and
$L^-_{\alpha}$ are curves in $Z'_0$.
Of course we have $L^+_{\alpha}$ and
$L^-_{\alpha}$ coincide if $\alpha=-Q\pm\sqrt{f}$.

Recall that in the previous subsection we have introduced an affine
coordinate
$v=\xi/\ell_1\ell_2$  on the exceptional curve $\Gamma_2$.
 Points on
$\Gamma_2$ are indicated by using this $v$.

\begin{lemma}\label{lemma-intersection3} Let $L^+_{\alpha}$ and
$L^-_{\alpha}$ be as above.
 Then $ L^+_{\alpha}\cap \Gamma_1,
L^+_{\alpha}\cap\Gamma_3,L^-_{\alpha}\cap\Gamma_1$ and
$L^-_{\alpha}\cap\Gamma_3$ are empty, and $L^+_{\alpha}\cap\Gamma_2$ and
$L^-_{\alpha}\cap\Gamma_2$ are points satisfying respectively
$$ L^+_{\alpha}\cap \Gamma_2=\left\{v=\left(\sqrt{f-(\alpha+Q)^2}+
i(\alpha+Q)\right)\,\frac{x_1^2}{\ell_1\ell_2}\right\}
$$ and
$$ L^-_{\alpha}\cap \Gamma_2=\left\{v=\left(-\sqrt{f-(\alpha+Q)^2}+
i(\alpha+Q)\right)\,\frac{x_1^2}{\ell_1\ell_2}\right\},
$$ 
where each intersection is transversal.
Moreover, $L^+_{\alpha}$ and $L^-_{\alpha}$ do not intersect
provided $\alpha\neq-Q\pm\sqrt{f}$ (namely 
$L^+_{\alpha}\neq L^-_{\alpha}$).
\end{lemma}

Here, note again that $\,x_1^2/\ell_1\ell_2$ does not depend on
$x_1$.

\vspace{3mm}
\noindent Proof.
 Substituting $x_2=\alpha x_1^2$, we have
$$\xi=z+i(x_2+Qx_1^2)=
\left\{\pm \sqrt{f-(\alpha+Q)^2}+i(\alpha+Q)\right\}\,x_1^2$$ and $$\eta=
z-i(x_2+Qx_1^2) =\left\{\pm
\sqrt{f-(\alpha+Q)^2}-i(\alpha+Q)\right\}\,x_1^2$$ over $C_{\alpha}$.
($\pm$ corresponds to $L_{\alpha}^{\pm}$.)
From these and from the explicit resolutions of the
previous subsection, we can easily see  that for any choice of
$\ell_1,\ell_2$ and
$\ell_3$,
$L_{\alpha}^{\pm}\cap \Gamma_1$ and
$L_{\alpha}^{\pm}\cap\Gamma_3$ are empty and that $L_{\alpha}^+\cap
\Gamma_2$ and
$L_{\alpha}^-\cap
\Gamma_2$ are points satisfying 
$$v=\xi/\ell_1\ell_2=\left\{\pm\sqrt{f-(\alpha+Q)^2}+
i(\alpha+Q)\right\}\,\frac{x_1^2}{\ell_1\ell_2},$$ where $\pm$
corresponds to
$L_{\alpha}^+$ and $L_{\alpha}^-$ respectively. 
The transversality is evident from this representation.
Finally, suppose that $\alpha\neq-Q\pm\sqrt{f}$.
 Then
since $C_{\alpha}$ intersect $B$ only at $P_{\infty}$ and 
$\ol{P}_{\infty}$, $L_{\alpha}^+$ and $L_{\alpha}^-$ intersect at most on
$\Gamma\cup\ol{\Gamma}$.
 But this does not happen because
$f-(\alpha+Q)^2$ is non-zero and hence the values of $v$ we have
already got are different.
 Thus we have obtained
all of the claims of the lemma. \proofend

\vspace{3mm} 
Since $L_{\alpha}^{\pm}$ and $\Gamma_2$ are $U(1)$-invariant,
$L_{\alpha}\cap \Gamma_2$ must be  $U(1)$-fixed point.
In particular, any points on $\Gamma_2$ is $U(1)$-fixed.
From these lemmas, we immediately get the following

\begin{lemma}\label{lemma-intersection4} Fix $\lambda\in I_2\cup I_4$.
 Then the set $\{ (L_{\alpha}^+\cup L_{\alpha}^-)\cap\Gamma_2\set
-Q-\sqrt{f}\leq\alpha\leq -Q+\sqrt{f} \}$ is a circle in
$\Gamma_2$ whose center is 
$\Gamma_2\cap
\Gamma_3$ ($=\{v=0\}$) and whose radius is 
$\sqrt{f}\,|x_1^2/\ell_1\ell_2|$. 
\end{lemma}

\vspace{3mm} The following proposition, which corresponds to Propositions
\ref{prop-nottl} (generic type) and
\ref{prop-nbc2} (special type), can be proved by using the same idea as
in  Proposition
\ref{prop-nbc2}.
 So we omit the proof.

\begin{prop}\label{prop-criet} Set
$h_{2}(\lambda)=\sqrt{f}\,(x_1^2/\ell_1\ell_2)$, which is clearly
differentiable on
 $I_2\cup I_4$.
  Let $N$ denote the normal bundle of $L^+_{\alpha}$ in
$Z'_0$.
 Then we have either 
$N\simeq O(1)^{\oplus 2}$ or $N\simeq O\oplus O(2)$, and the latter
holds iff
$\lambda$ is  a critical point of $h_2$.
  The same claim holds also for
$L^-_{\alpha}$.
\end{prop}

\begin{figure}[htbp]
\includegraphics{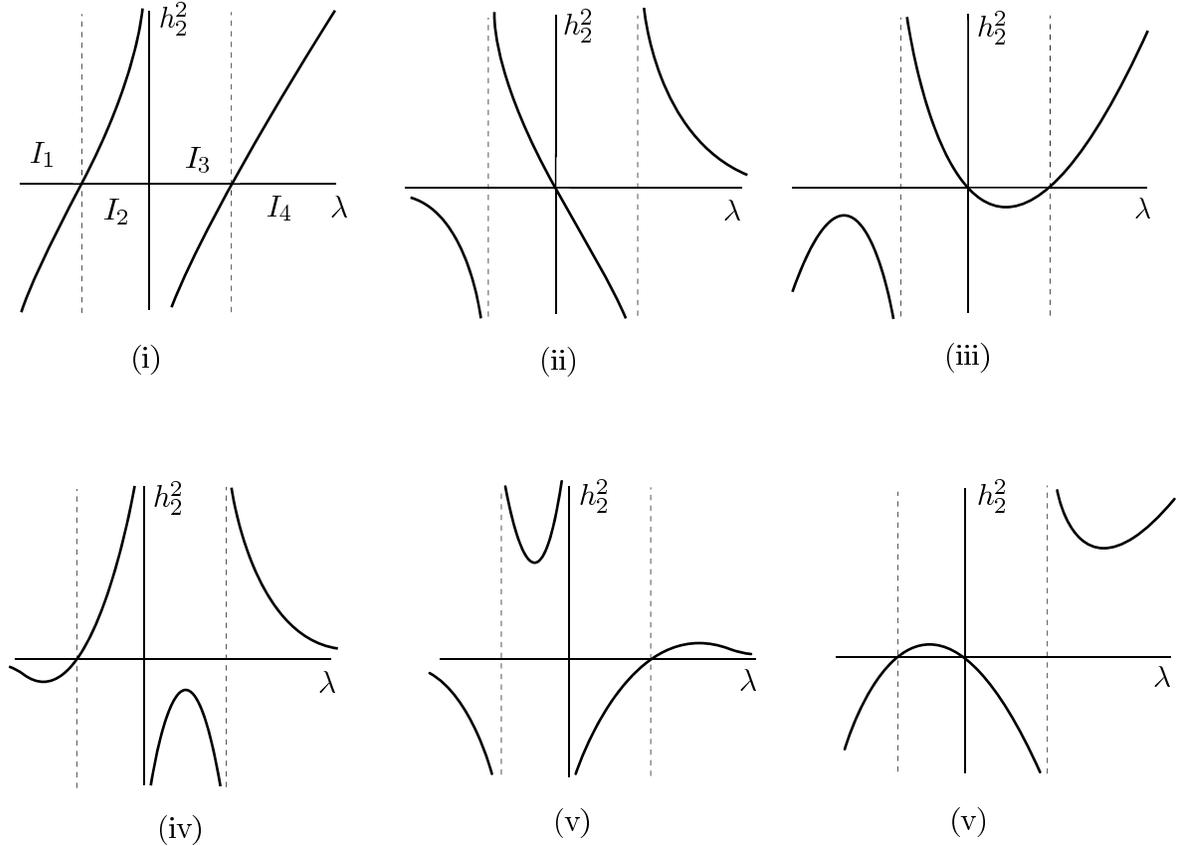}
\caption{behaviors of $h_2^2$}
\label{fig-h2}
\end{figure}

Needless to say, $h_2$ depends on the choice of $\ell_1$ and $\ell_2$.
Thus as in the case of special type, the normal bundles of $L_{\alpha}^+$
and
$L_{\alpha}^-$ depend on the choice of small resolution of
$p_{\infty}$.
  In view of Proposition \ref{prop-criet}, we need to know
the critical point of
$h_2$ for each choice of
$\{\ell_1,\ell_2\}$.
 There are $4!/(2!2!)=6$  choices of
$\{\ell_1,\ell_2\}$.
   If we take $\{\ell_1,\ell_2\}=\{x_0,x_1\}$ for
instance, we have
$h_2(\lambda)^2=(\lambda+1)(\lambda-a)/\lambda$, and it is elementary 
to determine the critical points of this function.
 For any other
choices, we always get
$h_2$ in explicit form and it is easy to determine their critical
points.
  So here we only present the result:

\begin{lemma}\label{lemma-beha} (i) If $\{\ell_1,\ell_2\}=\{x_0,x_1\}$,
$h_2$ has no critical point on $I_2\cup I_4$.
 Further,
$h_2(-1)=0,\lim_{\lambda\ua 0}h_2(\lambda)=+\infty,h_2(a)=0$ and
$\lim_{\lambda\ua\infty}h_2(\lambda)=+\infty$.
  (ii) If
$\{\ell_1,\ell_2\}=\{x_0+x_1,x_0-ax_1\}$, $h_2$ has no critical point on
$I_2\cup I_4$.
 Further, $\lim_{\lambda\da
-1}h_2(\lambda)=+\infty,h_2(0)=0,\lim_{\lambda\da
a}h_2(\lambda)=+\infty$ and
$\lim_{\lambda\ua\infty}h_2(\lambda)=0$.
 (iii) If
$\{\ell_1,\ell_2\}=\{x_1,x_0+x_1\}$, $h_2$ has no critical point on
$I_2\cup I_4$.
 Further, $\lim_{\lambda\da
-1}h_2(\lambda)=+\infty,h_2(0)=0, h_2(a)=0$ and
$\lim_{\lambda\ua\infty}h_2(\lambda)=\infty$.
 (iv) If
$\{\ell_1,\ell_2\}=\{x_0,x_0-ax_1\}$, $h_2$ has no critical point on
$I_2\cup I_4$.
 Further,  $h_2(-1)=0,\lim_{\lambda\ua
0}h_2(\lambda)=+\infty,
\lim_{\lambda\da a}h_2(\lambda)=+\infty$.
 (v) If
$\{\ell_1,\ell_2\}=\{x_0,x_0+x_1\}$, or if 
$\{\ell_1,\ell_2\}=\{x_1,x_0-ax_1\}$, $h_2$ has a unique critical point
on
$I_2$ and
$I_4$ respectively.
(See Figure \ref{fig-h2}.)
\end{lemma}
 
By Corollary \ref{cor-image1}, if $\lambda\in I_2$, images of
twistor lines in
$H_{\lambda}$ must be of orbit type.
 Therefore by Proposition
\ref{prop-criet},  if a small resolution of $Z_0$ yields a twistor space, 
$h_2$ does not have critical points on $I_2$.
  Hence by Lemma
\ref{lemma-beha}, we can conclude that
$\{\ell_1,\ell_2\}\neq\{x_0,x_0+x_1\}$ and 
$\{\ell_1,\ell_2\}\neq\{x_1,x_0-ax_1\}$.
 Namely, our investigation
decreases the possibilities of small resolutions.
 We postpone further consequences until the next subsection.

\subsection{Consequences of the results in \S
\ref{ss-generic}--\ref{ss-orbit}}
\label{ss-con} Before stating the results, we again recall our setup.
  Let
$B$ be a quartic surface defined by (\ref{eqn-B}) and assume that $Q$ and $f$
 satisfy  the necessary conditions as in Proposition
\ref{prop-necessa}.
  Let
$\Phi_0:Z_0\ra\mathbf{CP}^3$ be the double covering branched along $B$.
On $\mathbf{CP}^3$ there is a pencil of $U(1)$-invariant planes
$\{H_{\lambda}\}$, where  $H_{\lambda}$ is defined by $x_0=\lambda x_1$
which is real iff
$\lambda\in\mathbf R\cup\{\infty\}$. 
As explained in \S \ref{ss-special}, an ordering $\ell_1,\ell_2,\ell_3,
\ell_4$ of the four forms $\{x_0,x_1,x_0+x_1, x_0-ax_1\}$ determines
a small resolution of $p_{\infty}$ and $\ol{p}_{\infty}$ which are
compound $A_3$-singularities of $Z_0$.
 $\mu_{\infty}:
Z'_0\ra Z_0$ denotes this small resolution.
 Then 
$\mu_{\infty}$ is only a partial resolution of $Z_0$; namely
$Z'_0$ still has a (unique and real) ordinary double point $p_0$ corresponding
to the ordinary double point $P_0$ of $B$.
For any small resolution $\nu:Z\ra
Z'_0$ of $p_0$,  we put $\Phi:=\Phi_0\mu_{\infty}\nu$, and
let
 $\{S_{\lambda}=\Phi^{-1}(H_{\lambda})\}$ be (the real part of) a pencil
of
$U(1)$-invariant divisors on
$Z$, where we put  $S_{\lambda}=\Phi^{-1}(H_{\lambda})$ as before.

We start with the following proposition, which uniquely determines the type
of  real touching conics which can be the images of  twistor lines
contained in $U(1)$-invariant fundamental divisors.

\begin{prop}\label{prop-type2}
 Suppose that there is a small resolution  $\nu:Z\ra Z'_0$
 of the real ordinary double point $p_0$ such that $Z$ is a
twistor space.
  Let $L$ be a twistor line of $Z$ contained in
$S_{\lambda}$ for some
$\lambda\in\mathbf R$.
 Then $\Phi(L)$ is a real touching conic of:  (i) 
special type if
$\lambda\in I_1\cup I_3$, (ii)  orbit type if $\lambda\in I_2$,  (iii)  
generic type if
$\lambda\in I_4$ and if
$\lambda\neq\lambda_0$.
\end{prop}

Note that by Proposition \ref{prop-image}, $\Phi(L)\subset H_{\lambda}$
is a line if $\lambda=\lambda_0$.
 
\vspace{3mm}
\noindent Proof.
 (i) immediately follows from (ii) of Proposition
\ref{prop-type1} and  (v) of Proposition \ref{prop-inv-c}.
 (ii) is just
Corollary \ref{cor-image1}.
Finally we show   (iii).
  By (i) of Proposition \ref{prop-type1} it
suffices to show that if
$\lambda\in I_4$, the image cannot be of orbit type.
 In view of Lemma
\ref{lemma-beha}, we have $h_2(I_2)=(0,\infty)$ and
$h_2(I_4)=(0,\infty)$ for any of the cases (i)--(iv) of the lemma.
 (We
have already seen that the case (v) can be eliminated.)
 This implies that  the
circles
 appeared in Lemma
\ref{lemma-intersection3} sweep out
$\Gamma_2\backslash\{\Gamma_2\cap\Gamma_1,\Gamma_2\cap\Gamma_3\}$.
Therefore, $L_{\alpha}^{\pm}\subset S_{\lambda}$ with $\lambda\in
I_2$, and 
$L_{\alpha}^{\pm}\subset S_{\lambda}$ with $\lambda\in I_4$
cannot be the images of twistor lines at the same time.
 Therefore, if
$\lambda\in I_4$ and if $\lambda\neq \lambda_0$,
the images of twistor lines must be of generic type, as required.
\proofend

\vspace{3mm} The following is the main result of this section.
Recall that 
$p_{\infty}$ is a compound $A_3$-singularity of $Z_0$, and there are
$4!=24$ choices of  small resolutions of $p_{\infty}$,
each one corresponding to a choice of $\ell_1,\ell_2,\ell_3$ (and $\ell_4$)
(see  Section \ref{ss-special}).
 Recall also that
once a resolution of
$p_{\infty}$ is given, it naturally induces that of $\ol{p}_{\infty}$ by
reality.

\begin{prop}\label{prop-sr} Among 24 ways of possible small resolutions of
$p_{\infty}$, 22 resolutions do not yield a twistor space.
 The
remaining two  resolutions are given by the following two choices of
linear forms:
$$(I) \hspace{3mm}\ell_1=x_1,\hspace{2mm}\ell_2=x_0+x_1,\hspace{2mm}\ell_3=x_0,$$ and
$$(II)\hspace{3mm}\ell_1=x_0-ax_1,\hspace{2mm}\ell_2=x_0,\hspace{2mm}\ell_3=x_0+x_1.$$
\end{prop}

Here we do not yet claim that the threefolds obtained by these two
resolutions are actually twistor spaces.

\vspace{3mm}
\noindent Proof.
 By Proposition \ref{prop-type2} (i), if $\lambda\in I_1$, the
images of twistor lines in
$S_{\lambda}$ are real touching conics  of special type.
 As in Section
\ref{ss-special}, there are two families $\mathcal L_{\lambda}^+$ and
$\mathcal L_{\lambda}^-$ of real rational curves which are candidates of
 twistor lines in $S_{\lambda}$.
  As we have already remarked in
Section \ref{s-inv},
$L^+_{\theta}\in\mathcal L_{\lambda}^+$ and 
$L^-_{\theta}\in\mathcal L_{\lambda}^-$ 
cannot be twistor lines simultaneously.
 Suppose first that (any of the)
members of $\mathcal L_{\lambda}^+$  are twistor lines.
 Then by
Proposition \ref{prop-nbc2}, the function $h_1$ does not have critical
points on $I_1$.
 By Lemma
\ref{lemma-crtcl-sp}, this implies that  we have either
\begin{equation}\label{eqn-candi1}
\ell_1=x_1 \hspace{3mm}{\rm{or}}\hspace{3mm}
\ell_1=x_0+x_1.
\end{equation} On the other hand, by Corollary
\ref{cor-oneof} (i), under our assumption, members of $\mathcal
L_{\lambda}^-$ are twistor lines for $\lambda\in I_3$.
 Therefore again
by Proposition
\ref{prop-nbc2},
$h_3$ does not have critical points on
$I_3$.
 Then by Lemma \ref{lemma-crtcl-sp2}, the cases (i) and (iii)
of the lemma are
eliminated and  we have either 
\begin{equation}\label{eqn-candi2}
\{\ell_1,\ell_2,\ell_3\}=\{x_1,x_0+x_1,x_0-ax_1\}
\hspace{3mm}{\rm{or}}\hspace{3mm}
\{\ell_1,\ell_2,\ell_3\}=\{x_0,x_1,x_0+x_1\}.
\end{equation}(Note again that we do not specify the order.) 

Next we consider twistor lines in $S_{\lambda}$ for $\lambda\in I_2$.
 By 
Proposition \ref{prop-type2} (ii) the images are real touching conics of orbit
type.
 Then  Proposition \ref{prop-criet} implies that $h_2$ has no
critical point on
$I_2$.
 Hence by Lemma \ref{lemma-beha},  we have either 
\begin{equation}\label{eqn-candi3}
\{\ell_1,\ell_2\}=\{x_0,x_1\}
\hspace{2mm}{\rm{or}}\hspace{2mm}
\{x_0+x_1,x_0-ax_1\}
\hspace{2mm}{\rm{or}}\hspace{2mm}
\{x_1,x_0+x_1\}
\hspace{2mm}{\rm{or}}\hspace{2mm}
\{x_0,x_0-ax_1\}.
\end{equation} 

 Now we note other
restrictions: namely,  when
$\lambda$ increases to pass from $I_1$ to
$I_2$,  twistor lines in $S_{\lambda}$ must vary continuously, so that  we
have 
\begin{equation}\label{eqn-limit-1}
\lim_{\lambda\ua-1}h_1(\lambda)=\left(\lim_{\lambda\da-1}h_2(\lambda)\right)^{-1}.
\end{equation} (Here the inverse of the right hand side is a consequence
of the fact that 
$\Gamma_1\cap \Gamma_2=\{u=0\}=\{v=\infty\}$.)
 Similarly, moving
$\lambda$ from
$I_2$ to $I_3$, we have 
\begin{equation}\label{eqn-limit-2}
\lim_{\lambda\ua 0}h_2(\lambda)=\left(\lim_{\lambda\da
0}h_3(\lambda)\right)^{-1}.
\end{equation} Take $\ell_1=x_1$ for the first example.
 Then by Lemma
\ref{lemma-crtcl-sp} (i) we have
$h_1(-1)=0$.
 Hence it follows from (\ref{eqn-limit-1}) that  
$\lim_{\lambda\da-1}h_2(\lambda)=\infty$.
 Then the cases (i) and (iv)
of Lemma \ref{lemma-beha} fail and we have 
$\{\ell_1,\ell_2\}=\{x_0+x_1,x_0-ax_1\}$ ((ii)) or
$\{\ell_1,\ell_2\}=\{x_1,x_0+x_1\}$ ((iii)).
 The former clearly fails and
we get
$\ell_2=x_0+x_1$.
 This appears  in (\ref{eqn-candi3}).
  Then we have from
Lemma
\ref{lemma-beha} (iii) that $h_2(0)=0$.
 Hence by (\ref{eqn-limit-2}), we
have
$\lim_{\lambda\da 0}h_3(\lambda)=\infty$.
 It then follows from Lemma
\ref{lemma-crtcl-sp2} that $\ell_3=x_0$.
 Thus we get
$\ell_1=x_1,\ell_2=x_0+x_1,\ell_3=x_0$.

Next take $\ell_1=x_0+x_1$.
 Then we have $\lim_{\lambda\uparrow
-1}h_1(\lambda)=+\infty$ (Lemma \ref{lemma-crtcl-sp} (iii)), so that
$\lim_{\lambda\da-1}h_2(\lambda)=0$.
 Then looking (i)--(iv) of Lemma
\ref{lemma-beha}, this possibility fails.
 Namely, 
we have $l_1\neq x_0+x_1$.
 Thus we can conclude that if
$L^+_{\theta}\in\mathcal L^+_{\lambda}$ is supposed to be a twistor line over $I_1$,  we have to choose $\ell_1=x_1,\ell_2=x_0+x_1$ and $\ell_3=x_0$.
  This is the item (I) of the theorem. 

Next suppose  that  $L^-_{\theta}\in\mathcal L^-_{\lambda}$ is a twistor
line over $I_1$ and repeat similar argument above.
 By Proposition
\ref{prop-nbc2}, $h_3$ has no critical point on $I_1$.
 It then follows
from Lemma
\ref{lemma-crtcl-sp2} that  either
$\{\ell_1,\ell_2,\ell_3\}=\{x_0,x_0+x_1,x_0-ax_1\}$ ((i)) or
$\{\ell_1,\ell_2,\ell_3\}=\{x_0,x_1,x_0-ax_1\}$ ((iii)) holds.
On the other hand,  (\ref{eqn-candi3}) is valid also in this case.
Further we have as before
\begin{equation}\label{eqn-fsfa}
\lim_{\lambda\ua-1}h_3(\lambda)=\left(\lim_{\lambda\da-1}h_2(\lambda)\right)^{-1}
\hspace{2mm}{\rm{and}}\hspace{3mm}
\lim_{\lambda\ua 0}h_2(\lambda)=\left(\lim_{\lambda\da
0}h_1(\lambda)\right)^{-1}.
\end{equation}

If $\{\ell_1,\ell_2,\ell_3\}=\{x_0,x_0+x_1,x_0-ax_1\}$, then
$\lim_{\lambda\ua -1}h_3(\lambda)=+\infty$ (Lemma \ref{lemma-crtcl-sp2}
(i)), so that we have 
$\lim_{\lambda\da -1}h_2(\lambda)=0$ by (\ref{eqn-fsfa}).
 Hence by Lemma
\ref{lemma-beha} we have 
$\{\ell_1,\ell_2\}=\{x_0,x_0-ax_1\}$, which implies 
$\lim_{\lambda\ua 0}h_2(\lambda)=\infty$ ((iv) of Lemma
\ref{lemma-beha}) and
$\ell_3=x_0+x_1$.
 Hence $\lim_{\lambda\da 0}h_1(\lambda)=0$
by (\ref{eqn-fsfa}).
  It follows
from Lemma \ref{lemma-crtcl-sp} that $\ell_1=x_0-ax_1$, which means
$\ell_2=x_0$.
   (iv) of Lemma \ref{lemma-crtcl-sp} says that $h_1$ has
no critical point on
$I_3$, which is consistent with the fact that 
$L^+_{\theta}$ is a twistor line over $I_3$.

If $\{\ell_1,\ell_2,\ell_3\}=\{x_0,x_1,x_0-ax_1\}$, 
$h_3(-1)=0$ (Lemma \ref{lemma-crtcl-sp2} (iii)), so that we have 
$\lim_{\lambda\da -1}h_2(\lambda)=\infty$ by (\ref{eqn-fsfa}).
 Therefore
we get the two possibilities (ii) and (iii) of  Lemma
\ref{lemma-beha}, but both contain $x_0+x_1$ which is not compatible
with our choice of $\{\ell_1,\ell_2,\ell_3\}$.
 Thus we have
$\{\ell_1,\ell_2,\ell_3\}\neq \{x_0,x_1,x_0-ax_1\}$.
 This implies that  if $L^-_{\theta}\in\mathcal L^-_{\lambda}$ is a
twistor line for $\lambda\in I_1$,  then 
$\ell_1=x_0-ax_1,\ell_2=x_0$, and $\ell_3=x_0+x_1$.
  This is the case (II) of the theorem, and we have completed the proof.
\proofend

\vspace{3mm}
Next we summarize which irreducible component of the inverse image of touching conics have to be chosen as candidates of twistor lines. 
Note that for touching conic of special type, we have made in Section \ref{ss-special} distinction of the two components $L^+_{\theta}$ and $L^-_{\theta} $ by 
the property that $L^+_{\theta}$ intersects $\Gamma_1$ and $\ol{\Gamma}_1$, and $L^-_{\theta} $ intersects $\Gamma_3$ and $\ol{\Gamma}_3$ (Lemma \ref{lemma-interse}).
For touching conics of generic type, we have not made distinction of $L^+_{\theta}$ and $L^-_{\theta} $ so far.
To make a distinction, we write $\Phi^{-1}(l_{\infty})-\Gamma-\ol{\Gamma}=\Xi+\ol{\Xi}$, and $\Phi$ gives isomorphism of $\Xi\simeq l_{\infty}$ and $\ol{\Xi}\simeq l_{\infty}$ (cf.\,Proposition \ref{prop-octagon}).
We can suppose that $\Xi$ intersects $\Gamma_1$ (and hence $\ol{\Gamma}_3$)
as in Figure \ref{fig-cycle3}. 
Then  $\ol{\Xi}$ intersects $\Gamma_3$ (and $\ol{\Gamma}_1$).
Recall that $C_{\theta}$ intersects $l_{\infty}$ at two points indicated by $e^{i\theta}h_0(\lambda)$ and $e^{-i\theta}h_0^{-1}(\lambda)$ (Lemma \ref{lemma-int-a} and its proof).
Thus $L^+_{\theta}$ intersects $\Xi$ at a point indicated by $e^{i\theta}h_0$ or $e^{-i\theta}h_0^{-1}$.
Now we make a distinction of $L_{\theta}^+$ and $ L^-_{\theta}$ by declaring that $L^+_{\theta}\cap \Xi$ is indicated by $e^{i\theta}h_0$ for $\lambda<\lambda_0$, and $e^{-i\theta}h_0^{-1}$ for $\lambda>\lambda_0$. 
It then follows that  $L^-_{\theta}\cap \Xi$ is indicated by $e^{-i\theta}h_0^{-1}$ for $\lambda<\lambda_0$ and $e^{i\theta}h_0$ for $\lambda>\lambda_0$.
The following result is essentially showed in the proof of Proposition \ref{prop-sr}.

\begin{prop}\label{prop-choice01}
(i) If we take the small resolution determined by (I) of Proposition \ref{prop-sr}, members of $\mathcal L_{\lambda}^+$ must be chosen as twistor lines on $S_{\lambda}$ for $\lambda\in I_1$, both members $L_{\alpha}^+$ and $ L_{\alpha}^-$ must be chosen for $\alpha\in I_2$, members of $\mathcal L_{\lambda}^-$ must be chosen for $\lambda\in I_3$,  and  $\mathcal L_{\lambda}^+$ must be chosen for $\lambda\in I_4^-\cup I_4^+$.
(ii) If we take the small resolution determined by (II) of Proposition \ref{prop-sr}, members of $\mathcal L_{\lambda}^-$ must be chosen as twistor lines on $S_{\lambda}$ for $\lambda\in I_1$, both members $L_{\alpha}^+$ and $ L_{\alpha}^-$ must be chosen for $\lambda\in I_2$, members of $\mathcal L_{\lambda}^+$ must be chosen for $\lambda\in I_3$,  and $\mathcal L_{\lambda}^-$ must be chosen for $\lambda\in I_4^-\cup I_4^+$.
\end{prop}

\noindent Proof.
We only verify (i) because (ii) can be seen by parallel argument. 
The claims for $\lambda\in I_1\cup I_3$ are already seen in the proof of Proposition \ref{prop-sr}.
(ii) immediately follows form Proposition \ref{prop-inv-c} (and its subsequent comment) and Corollary \ref{cor-image1}.
It remains to see the claim for $\lambda\in I_4^{\pm}$.
Over $I_3$, members of
$\mathcal L^-_{\lambda}$ must be chosen as above.
 Further, we have
$h_3(a)=0$ by Lemma \ref{lemma-crtcl-sp2} (iv).
 This implies that as $\lambda$ goes to $a$, the  intersection circle 
$\cup\{\Gamma_3\cap L^-_{\theta}\set
 L^-_{\theta}\in\mathcal L^-_{\lambda}\}$ shrinks
to be the point $\Gamma_3\cap\ol{\Xi}$.
Therefore the intersection circle of twistor lines in $S_{\lambda}$ and $\ol{\Xi}$ also has to shrink to be $\Gamma_3\cap\ol{\Xi}$.
 On $\ol{\Xi}$, one can use
$x_2=y_2/y_3$ as an affine coordinate whose center is the intersection
point with $\Gamma_3$.
By the rule explained just before the proposition, the radius (with respect to $x_2$) of  intersection circle of $L^+_{\theta}\cap \ol{\Xi}$  and $L^-_{\theta}\cap \ol{\Xi}$ are respectively indicated by $h_0^{-1}$ and $h_0$.
We have $\lim_{\lambda\da a}h_0(\lambda)=+\infty$ and $\lim_{\lambda\da a}h_0^{-1}(\lambda)=0$ (cf.\,Figure \ref{fig-h0}).
Hence we conclude that $L^+_{\theta}\cap \ol{\Xi}$  must be chosen for $\lambda\in I_4^-$.
Similar argument shows that for $\lambda\in I_4^-$, $L^+_{\theta}\cap \ol{\Xi}$  must still be chosen.\proofend

\begin{figure}[htbp]
\includegraphics{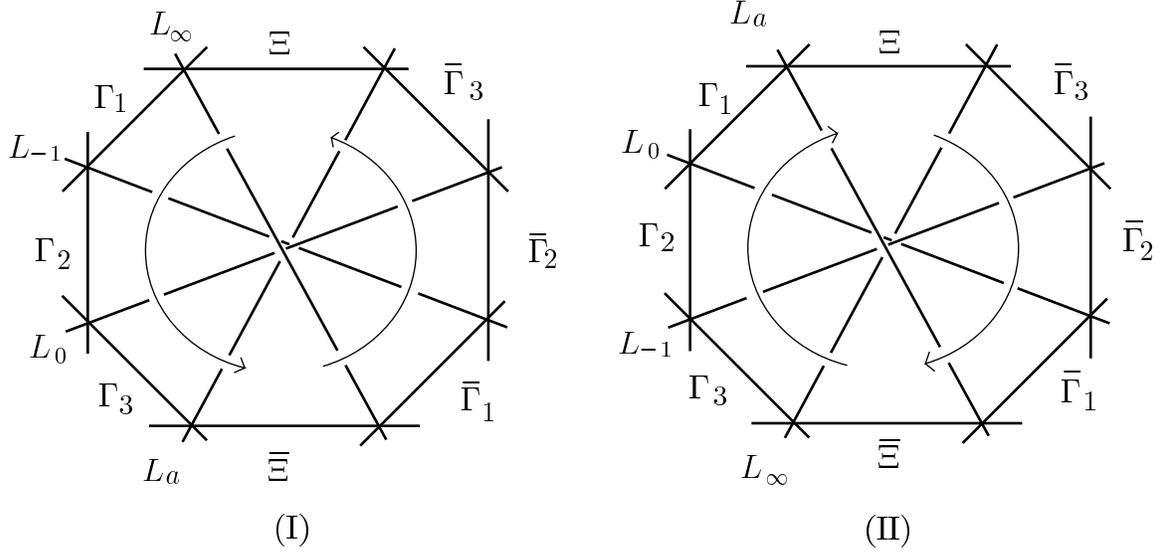}
\caption{The intersection circles rotate as $\lambda$ increases.}
\label{fig-cycle4}
\end{figure}

\vspace{3mm} Differences of the two resolutions (I) and (II) of Proposition \ref{prop-sr} can be displayed as in Figure \ref{fig-cycle4}. 
Briefly speaking, exchanging the two resolutions reverses the direction of the moving of the intersection circles as $\lambda$ increases.

\subsection{Connectedness of the families of real lines obtained in \S \ref{ss-con}}\label{ss-conn}
In the last subsection we have determined the families of twistor lines contained in $S_{\lambda}=\Phi^{-1}(H_{\lambda})$, where $\lambda$ is in the intervals $I_{j}$, $1\leq j\leq 4$. 
Namely, we have detected two $S^1$-families of real lines on $S_{\lambda}$ for $\lambda\in \cup_{j=1}^4 I_{j}$, $\lambda\neq \lambda_0$.
It is obvious from our explicit description that, for each $1\leq j\leq 4$, the two $S^1$-families respectively forms a connected 2-dimensional family of real lines whose parameter space is a cylinder $I_{j}\times S^1$.
As stated in Proposition \ref{prop-type2}, these families have different description depending on $1\leq j\leq 4$.
In this subsection, we show that, any member of these families convergent to
 the inverse image of the four tropes of $B$, when the plane moves to the four planes on which the tropes lie.
 This implies that the parameter spaces $I_j\times S^1$ are joined (connected) by adding four points corresponding to the tropes.
 This connected result will be needed in our proof of main theorem.

Recall that  the branch quartic surface $B$ is defined by 
$$\left(y_2y_3+Q(y_0,y_1)\right)^2-y_0y_1(y_0+y_1)(y_0-ay_1)=0.$$
Let $T_{-1},T_{0},T_{a}$ and $T_{\infty}$ be the four tropes which are intersections of $B$ with the planes $y_0=-y_1,y_0=0,y_0=ay_1$ and $y_1=0$ respectively.
The inverse image of the tropes in $Z_0$ is of course a smooth rational curves which are mapped biholomorphically onto the tropes. 
These four curves on $Z_0$ define, as the strict transforms, smooth rational curves on $Z'_0$ for any small resolution $Z'_0\ra Z_0$ of $p_{\infty}$ and $\ol{p}_{\infty}$.
We denote these curves in $Z'_0$ by $L_{-1}, L_0,L_{a}$ and $L_{\infty}$ respectively.
These are $\mathbf C^*$-invariant real rational curves naturally biholomorphic to the tropes.
By Lemma \ref{lemma-degree01}, the normal bundles in $Z'_0$ are isomorphic to $O(1)^{\oplus 2}$ by Lemma \ref{lemma-degree01}.
(In the lemma these were written $L_{\ell_i}$.)

In the following we investigate limits of real lines on $S_{\lambda}$ when $\lambda$ approaches the endpoints of $I_{j}$. 
We begin with the case of $\lambda\in I_1\cup I_3$.
In this case we know by Proposition \ref{prop-type2} that the images of twistor lines on $S_{\lambda}$ is a touching conic of special type.
Moreover, Proposition \ref{prop-choice01} shows that if we choose the small resolution determined by (I)  
of Proposition \ref{prop-sr}, then $L_{\theta}^+$ 
 must be selected as twistor lines on $S_{\lambda}$ for $\lambda\in I_1$, whereas 
$L_{\theta}^-$  must be selected for $\lambda\in I_3$.
Similarly, if we choose the small resolution determined by (II) 
of Proposition \ref{prop-sr}, then $L_{\theta}^-$ 
 must be selected as twistor lines on $S_{\lambda}$ with $\lambda\in I_1$, whereas 
$L_{\theta}^+$  must be selected for $\lambda\in I_3$.
The next two propositions in particular imply that when the real plane $H_{\lambda}$ moves to the endpoints of $I_1$ and $I_3$, these real lines convergent to the real lines $L_{\lambda_i}$, $\lambda_i=-1,0,a$ or $\infty$.

\begin{prop}\label{prop-limit-special1}
Assume $\lambda\in I_1\cup I_3$ and let $L_{\theta}^{\pm}$ be the real lines in $S_{\lambda}$ which are the irreducible components of the inverse image of real touching conic $C_{\theta}\subset H_{\lambda}$ of special type as in Section \ref{ss-special}. 
If we choose the small resolution (I)  of Proposition \ref{prop-sr}, we have the following:
$$\left\{
\begin{array}{ll}\displaystyle\lim_{\lambda\da -\infty}L_{\theta}^+=L_{\infty}, &
\displaystyle\lim_{\lambda\da -\infty}L_{\theta}^-=L_{\infty}+\Gamma_1+\Gamma_2+\ol{\Gamma}_1+\ol{\Gamma}_2,\\
\displaystyle\lim_{\lambda\ua -1}L_{\theta}^+=L_{-1},&
\displaystyle\lim_{\lambda\ua-1}L_{\theta}^-=L_{-1}+\Gamma_2+\ol{\Gamma}_2,
\end{array}
\right.
$$
and
$$ \left\{
\begin{array}{ll}\displaystyle\lim_{\lambda\da 0}L_{\theta}^+=L_{0}+\Gamma_2+\ol{\Gamma}_2, &
\displaystyle\lim_{\lambda\da 0}L_{\theta}^-=L_{0},\\
\displaystyle\lim_{\lambda\ua\,a}L_{\theta}^+=L_{a}+\Gamma_2+\Gamma_3+\ol{\Gamma}_2+\ol{\Gamma}_3,&
\displaystyle\lim_{\lambda\ua\,a}L_{\theta}^-=L_{a}.
\end{array}
\right.$$
\end{prop}

\noindent Proof.
It is immediate to verify from the explicit equation of touching conic of $C_{\theta}$ obtained in Proposition \ref{prop-b} that $\lim_{\lambda\da -\infty}C_{\theta}=T_{\infty}.$
Since $\Phi^{-1}(C_{\theta})=L_{\theta}^++L_{\theta}^-$ and $\Phi^{-1}(T_{\infty})=L_{\infty}+\Gamma+\ol{\Gamma}$, we obtain
$$
\lim_{\lambda\da -\infty}(L_{\theta}^++L_{\theta}^-)=L_{\infty}+\Gamma+\ol{\Gamma}.
$$
Moreover, both $\lim_{\lambda\ra-\infty}L^+_{\theta}$ and $\lim_{\lambda\ra-\infty}L^-_{\theta}$ contain $L_{\infty}$ as their irreducible components, since $L^+_{\theta}$ and $L_{\theta}^-$ are mapped biholomorphically onto $C_{\theta}$.
First we consider $\lim L_{\theta}^+$ as $\lambda\da-\infty$.
By Lemma \ref{lemma-interse}, $L_{\theta}^+$ intersects $\Gamma_1$ at a unique point $$u=-2iBe^{-i\theta}\cdot x_1/\ell_1=ie^{-i\theta}h_1,$$ where $u$ is the affine coordinate on $\Gamma_1$ introduced in the beginning of Section \ref{ss-special}, and $h_1=h_1(\lambda)$ is a function introduced before Proposition \ref{prop-nbc2}.
Now since  we are choosing the small resolution determined by (I) of Proposition \ref{prop-sr},  we have $\ell_1=x_1$ and the graph of $h_1^2$ looks like (i) of Figure \ref{fig-h1}.
(See Lemma \ref{lemma-crtcl-sp} (i) for precisely.)
In particular, we have $\lim_{\lambda\da-\infty}h_1(\lambda)=\infty$. 
Because $L_{\infty}$ is as in Figure \ref{fig-cycle4} in this case,  this implies that the intersection point $L^+_{\theta}\cap \Gamma_1$ goes to the point $L_{\infty}\cap \Gamma_1$.
Hence $\lim_{\lambda\da-\infty}L^+_{\theta}$ does not contain any components of $\Gamma$ and $\ol{\Gamma}$ and we get $\lim_{\lambda\da-\infty}L^+_{\theta}=L_{\infty}$.
$\lim_{\lambda\ua-1}L^+_{\theta}=L_{-1}$ can be proved in a similar way.
Next we look into the limit of $L^-_{\theta}$ as $\lambda$ approaches $-\infty$ and $-1$.
We know by Lemma \ref{lemma-interse} that $L^-_{\theta}$ does not intersect $\Gamma_1\cup\Gamma_2$ and that $L^-_{\theta}\cap \Gamma_3$ is a unique point satisfying 
$$w=-\frac{ie^{i\theta}f}{2B}\cdot\frac{x_1^3}{\ell_1\ell_2\ell_3}=-ie^{i\theta}h_3,$$ where $w$ is the affine coordinate on $\Gamma_3$ introduced in  Section \ref{ss-special}, and $h_3$ is a function introduced before Proposition \ref{prop-nbc2}.
Moreover, by Lemma \ref{lemma-crtcl-sp2} (iv) (plus some easy calculations), we have $\lim_{\lambda\da-\infty}h_3(\lambda)=\infty$ and $\lim_{\lambda\ua-1}h_3(\lambda)=\infty$.
This implies that the point $L^-_{\theta}\cap \Gamma_3$ approaches $L_0\cap\Gamma_3$ when $\lambda\da-\infty$ and $\lambda\ua-1$.
On the other hand, as is already mentioned, $\lim_{\lambda\da-\infty}L_{\theta}^-$ contains $L_{\infty}$ as its irreducible components. 
Then since $\lim_{\lambda\da-\infty}L_{\theta}^-$ must be connected,
it must contain $\Gamma_1$ and $\Gamma_2$ (see Figure \ref{fig-cycle4}). 
By reality, we conclude $\lim_{\lambda\da-\infty}L_{\theta}^-=L_{\infty}+\Gamma_1+\Gamma_2+\ol{\Gamma}_1+\ol{\Gamma}_2$.
Similarly, $\lim_{\lambda\ua-1}L_{\theta}^-$ must contain $\Gamma_2$ and $\ol{\Gamma}_2$ and we get $\lim_{\lambda\ua-1}L_{\theta}^-=L_{-1}+\Gamma_2+\ol{\Gamma}_2$, as claimed.

For other limits $\lambda\da 0$ and $\lambda\ua a$, we can prove the claim of the proposition in a similar way  by using Lemmas \ref{lemma-crtcl-sp} and \ref{lemma-crtcl-sp2}.  We omit the detail.
\proofend

\vspace{3mm}
For another small resolutions, we have the following proposition.
We omit the proof since it is completely parallel to the above case.

\begin{prop}\label{prop-limit-special2}
In Proposition \ref{prop-limit-special1}, if we choose the small resolution (II)  of Proposition \ref{prop-sr} instead, we have the following:
$$\left\{
\begin{array}{ll}\displaystyle\lim_{\lambda\da -\infty}L_{\theta}^+=L_{\infty}+\Gamma_2+\Gamma_3+\ol{\Gamma}_2+\ol{\Gamma}_3, &
\displaystyle\lim_{\lambda\da -\infty}L_{\theta}^-=L_{\infty},\\
\displaystyle\lim_{\lambda\ua -1}L_{\theta}^+=L_{-1}+\Gamma_2+\ol{\Gamma}_2,&
\displaystyle\lim_{\lambda\ua-1}L_{\theta}^-=L_{-1},
\end{array}
\right.
$$
and
$$ \left\{
\begin{array}{ll}\displaystyle\lim_{\lambda\da 0}L_{\theta}^+=L_{0}&
\displaystyle\lim_{\lambda\da 0}L_{\theta}^-=L_{0}+\Gamma_2+\ol{\Gamma}_2, ,\\
\displaystyle\lim_{\lambda\ua\,a}L_{\theta}^+=L_{a},&
\displaystyle\lim_{\lambda\ua\,a}L_{\theta}^-=L_{a}+\Gamma_1+\Gamma_2+\ol{\Gamma}_1+\ol{\Gamma}_2.
\end{array}
\right.$$
\end{prop}

Next we investigate the limit of twistor lines on $S_{\lambda}$ when $\lambda$ approaches the endpoints of  $I_2$. 
In this case we know by  Proposition \ref{prop-type2}  that the image of twistor line is a touching conic of orbit type  denoted by $C_{\alpha}\subset H_{\lambda}$, and that by  Proposition \ref{prop-inv-c} there are two irreducible components $L_{\alpha}^+$ and $L_{\alpha}^-$ of $\Phi^{-1}(C_{\alpha})$ which are biholomorphic to $C_{\alpha}$ by $\Phi$.
As remarked at the end of Section \ref{s-inv}, both $L_{\alpha}^+$ and $L_{\alpha}^-$ must be twistor lines (because they are members of the same 
pencil on $S_{\lambda}$).

\begin{prop}\label{prop-lim2}
Assume $\lambda\in I_2$ and let $L_{\alpha}^{\pm}$ be the real lines in $S_{\lambda}$ which are  irreducible components of the inverse image of real touching conic $C_{\alpha}\subset H_{\lambda}$ of orbit type as in Section \ref{ss-orbit}.
Then regardless of which small resolution (I) or (II) of Proposition \ref{prop-sr} we take, we have the following
$$\displaystyle\lim_{\lambda\da -1}L_{\theta}^+=\lim_{\lambda\da -1}L_{\theta}^-=L_{-1},\hspace{3mm}\displaystyle\lim_{\lambda\ua 0}L_{\theta}^+=
\displaystyle\lim_{\lambda\ua 0}L_{\theta}^-=L_{0}.
$$
\end{prop}

\noindent Proof.
This can be proved by the same argument as in Proposition \ref{prop-limit-special1}, if we replace Proposition \ref{prop-b} by Proposition \ref{prop-c}, Lemma \ref{lemma-interse} by Lemma \ref{lemma-intersection3}, and Lemmas \ref{lemma-crtcl-sp} and \ref{lemma-crtcl-sp2} by Lemma \ref{lemma-beha}.
\proofend

\vspace{3mm}
Finally we study the case $\lambda\in I_4$. 
Again by Proposition \ref{prop-type2} the image of twistor lines are real touching conics of generic type. 
We have $\Phi^{-1}(C_{\theta})=L_{\theta}^++L^-_{\theta}$.
Recall that in Section \ref{ss-con} (the explanation before Proposition \ref{prop-choice01})
we have made distinction of
 $L_{\theta}^+$ and $ L^-_{\theta}$ by considering the intersection with $\Xi$ and $\ol{\Xi}$.

\begin{prop}\label{prop-conn-4}
Assume $\lambda\in I_4^-\cup I_4^+$ and let $L_{\theta}^{\pm}$ be as above.
Then  the following (I') and (II') respectively hold depending on which small resolution (I) and (II) of Proposition \ref{prop-sr} we take:
$$(I')\cdots \left\{
\begin{array}{ll}\displaystyle\lim_{\lambda\da \,a}L_{\theta}^+=L_{a}, &
\displaystyle\lim_{\lambda\da \,a}L_{\theta}^-=L_{a}+\Gamma+\ol{\Gamma},\\
\displaystyle\lim_{\lambda\ua\infty}L_{\theta}^+=L_{\infty},&
\displaystyle\lim_{\lambda\ua\infty}L_{\theta}^-=L_{\infty}+\Gamma+\ol{\Gamma}.
\end{array}
\right.$$
$$(II')\cdots \left\{
\begin{array}{ll}\displaystyle\lim_{\lambda\da \,a}L_{\theta}^+=L_{a}+\Gamma+\ol{\Gamma}, &
\displaystyle\lim_{\lambda\da \,a}L_{\theta}^-=L_{a},\\
\displaystyle\lim_{\lambda\ua\infty}L_{\theta}^+=L_{\infty}+\Gamma+\ol{\Gamma},&
\displaystyle\lim_{\lambda\ua\infty}L_{\theta}^-=L_{\infty}.
\end{array}
\right.$$
\end{prop}

\noindent Proof.
It is immediate to verify from the explicit equation of touching conic of generic type obtained in Proposition \ref{prop-a} that $\lim_{\lambda\da a}C_{\theta}=T_{a}$
and $\lim_{\lambda\ua \infty}C_{\theta}=T_{\infty}$.
As $\Phi^{-1}(C_{\theta})=L_{\theta}^++L_{\theta}^-$ and $\Phi^{-1}(T_{a})=L_{a}+\Gamma+\ol{\Gamma}$ and $\Phi^{-1}(T_{\infty})=L_{\infty}+\Gamma+\ol{\Gamma}$, we obtain
$$
\lim_{\lambda\da a}(L_{\theta}^++L_{\theta}^-)=L_{a}+\Gamma+\ol{\Gamma},\hspace{3mm}\lim_{\lambda\ua \infty}(L_{\theta}^++L_{\theta}^-)=L_{\infty}+\Gamma+\ol{\Gamma}.
$$
By Lemma \ref{lemma-bhv-1} (plus easy calculation; see Figure \ref{fig-h0}), we have $\lim_{\lambda\da a} h_0(\lambda)=\infty$ and $\lim_{\lambda\ua\infty} h_0^{-1}(\lambda)=0$.
By the rule for making  distinction of $L^+_{\theta}$ and $L^-_{\theta}$, these imply that $L^+_{\theta}\cap \Xi$ goes to the point $\Xi\cap\ol{\Gamma}_3$ as $\lambda\da a$, and goes to the point $ \Xi\cap \Gamma_1$ as $\lambda\ua\infty$.
On the other hand, we have $L_{a}\cap \Xi=\ol{\Gamma}_3\cap \Xi$ and $L_{\infty}\cap \Xi=\Gamma_1\cap \Xi$ for the small resolution determined by (I) of Proposition \ref{prop-sr} (Figure \ref{fig-cycle4}).
Namely as $\lambda\da a$ and $\lambda\ua\infty$ respectively, $\Xi\cap L_{\theta}^-$ approaches  $\Xi\cap L_{a}$ and $\Xi\cap L_{\infty}$. 
Therefore we have $\lim_{\lambda\da \,a}L_{\theta}^+=L_{a}$ and $\lim_{\lambda\ua\infty}L_{\theta}^+=L_{\infty}$ if we take the small resolution coming from (I) of Proposition \ref{prop-sr}.

Next we look into the limit of $L_{\theta}^-$ as $\lambda\da a$ and $\lambda\ua \infty$ for type (I) resolution.
It is immediate from Lemma \ref{lemma-bhv-1} and Figure \ref{fig-h0} again to get that the point $\Xi\cap L^-_{\theta}$ goes to the points $\Xi\cap L_{\infty}$ as $\lambda\da a$ and goes to the points $\Xi\cap L_{a}$ as $\lambda\ua\infty$.
Then since $\lim_{\lambda\da a} L_{\theta}^-$ is connected and contains $L_{a}$ as its irreducible component (and since moreover the limit does not contain $\Xi$), it follows that $\lim_{\lambda\da a} L_{\theta}^-$ contains $\Gamma$.
Then by reality $\lim_{\lambda\da a} L_{\theta}^-$ also contains $\ol{\Gamma}$.
Thus we get $\lim_{\lambda\da a} L_{\theta}^-=L_{a}+\Gamma+\ol{\Gamma}$.
Similarly we get $\lim_{\lambda\ua \infty} L_{\theta}^-=L_{\infty}+\Gamma+\ol{\Gamma}$.

Finally, we can immediately obtain the limit for the type (II) resolution, if we  note that exchanging resolution of type (I) and type (II) interchanges $L_{a}$ and $L_{\infty}$ (Figure \ref{fig-cycle4}). \proofend

\vspace{3mm}
Propositions \ref{prop-limit-special1}--\ref{prop-conn-4} can be summarized as follows:

\begin{prop}\label{prop-limit02}
Let $Z'_0\ra Z_0$ be the small resolution of the conjugate pair of singularities determined by (I) or (II) of Proposition \ref{prop-sr}, and consider the families of real lines on $S_{\lambda}$ ($\lambda\in I_j$) for each $1\leq j\leq 4$ determined in Proposition \ref{prop-choice01}, (i) or (ii) respectively.
Then if we move $H_{\lambda}$ into the planes on which the tropes of $B$ lie, the real lines on $S_{\lambda}$ converge to the real line which is (an irreducible component of) the inverse image of the trope.
\end{prop}

\section{Twistor lines whose images are lines and a resolution of the ODP}\label{s-line}

In Sections \ref{s-how}--\ref{s-nb} we have intensively studied real lines in $Z$ whose images are touching conics.
But as Proposition \ref{prop-image} says, we cannot obtain arbitrary twistor lines by considering such real lines only.
In the first half of this section we study real lines whose images become lines.
These real lines are needed for compactifying the space of real lines obtained in Sections \ref{s-how}--\ref{s-nb}.

On the other hand we showed in the last section that among many (twenty-four) ways of possible small resolutions of the conjugate pair of singularities of $Z_0$, there are just two resolutions which can yield a twistor space (Proposition \ref{prop-sr}). 
In the second half of this section we show that once a small resolution of the unique ordinary double point of $Z_0$ is given, 
it uniquely determines which one of  the above two resolutions have to be taken for the conjugate pair of singularities.
We keep the notations and assumptions in Sections \ref{s-how}--\ref{s-nb}.

As showed in Proposition \ref{prop-image}, if
$L\subset Z$ is a real line  intersecting $\Gamma_0$, then $\Phi(L)$ is a line
 going through $P_0$.
 The next proposition is its
 converse.

\begin{prop}\label{prop-imline} Let $l\subset\mathbf{CP}^3$ be any real line
going through
$P_0$. Then $\Phi^{-1}(l)$ has just two irreducible components, both of which are
real, smooth and rational. One of the components is the exceptional
curve $\Gamma_0$ and  another component is  mapped (2 to 1) onto $l$.
Further, the normal bundle of the latter component in $Z$ is isomorphic
to $O(1)^{\oplus 2}$. 
\end{prop}

\noindent Proof. First we note that if $l$ is a  real line,  $B\cap
l$ consists of just three points, one of which is
$P_0$. This follows from the facts that, $B$ is a quartic,  $B\cap l$ is
real,
$P_0$ is the unique real point of $B$ (Proposition
\ref{prop-realpoint}), and $P_0$ is a double point. Therefore
$\Phi_0^{-1}(l)\ra l$ is  two-to-one covering branched at
three points. Let $P$ and
$\ol{P}$ be the two branch points other than $P_0$. 
Because $l$ intersects $B$ transversally at these two points,
$\Phi_0^{-1}(P)$ and $\Phi_0^{-1}(\ol{P})$ are smooth points of
$\Phi_0^{-1}(l)$. Further, since $P_0$ is an ordinary double point of $B$,
$p_0=\Phi_0^{-1}(P_0)$ is a node of
$\Phi_0^{-1}(l)$. From these it  follows that  
$\Phi^{-1}(l)$ has just two irreducible components, one of which is
$\Gamma_0$. Let
$L$ be the irreducible component different from $\Gamma_0$. 
Then $L$ is  smooth and $L\ra l$ is two-to-one covering whose branch points are $P$ and $\ol{P}$.  
 Therefore (by Hurwitz) $L$ is a  rational curve.
$L$ is  real since $\Phi_0^{-1}(l)$ is real.

It remains to show that $N_{L/Z}\simeq  O(1)^{\oplus 2}$. The idea is
similar to   Propositions \ref{prop-nbc},
\ref{prop-nbc2}, and
\ref{prop-criet}. We first show that $N_{L/Z}\simeq O(1)^{\oplus 2}$
 or $N_{L/Z}\simeq O\oplus O(2)$. By Bertini, $H\cap B$ is smooth
outside $P_0$ for general plane $H$ containing $ l$. Further, since
$H\cap B$ is a quartic,
$S:=\Phi^{-1}(H)$ is a smooth rational surface with $c^2_1=2$. Moreover,
$\Phi^{-1}(l)$ is an anticanonical curve of
$S$ so that we have $(\Gamma_0+L)^2=2$ on $S$. Furthermore, it is readily
seen that
$\Gamma_0^2=-2$ and $\Gamma_0\cdot L=2$ on $S$. Therefore we have $L^2=0$ on $S$.
Then the argument in the proof of Proposition \ref{prop-nbc} implies 
$N_{L/Z}\simeq O(1)^{\oplus 2}$ or $N_{L/Z}\simeq O(2)\oplus O$.
 To show that the latter does
not hold, we first see that $\Gamma_0/\langle\sigma\rangle$ is canonically
identified with the projective space of real lines going through $P_0$. Concretely,
for each real line $l\ni P_0$, we associate the intersection $\Gamma_0\cap
(\Phi^{-1}(l)-\Gamma_0)$ which is a conjugate pair of points.
We show by explicit calculation that this correspondence,
which we will denote by $\psi$, is actually
an isomorphism. The problem being local, we use a local coordinate
$(w_1,w_2,w_3)$  (around $P_0$) defined in (\ref{coord}).  Then in a
neighborhood of
$p_0=\Phi_0^{-1}(P_0)$, $Z_0$ is given by the equation 
\begin{equation}\label{eqn-Z_0}z^2+w_1^2+w_2w_3=0.\end{equation}  Small
resolutions of the double point $p_0\in Z_0$ are
explicitly obtained by blowing-up along $\{z+iw_1=w_3=0\}$ or
$\{z+iw_1=w_2=0\}$. In the former case, we can use
$(z+iw_1:w_3)=(-w_2:z-iw_1)$ as a homogeneous coordinate on $\Gamma_0$,
whereas in the latter case we can use
$(z+iw_1:w_2)=(-w_3:z-iw_1)$ instead. We see only in the former case,
since the calculation is identical. Let
$(w_1:w_2:w_3)$ be a real line through $P_0$.  Namely, we assume
$w_1\in\mathbf R$,  $\ol{w}_2=w_3$, and $w_1^2+|w_2|^2\neq 0$.  Then by
(\ref{eqn-Z_0}), we have $z=\pm i(w_1^2+|w_2|^2)^{1/2}$. Hence we get
$(z+iw_1:w_3)=(i(w_1\pm (w_1^2+|w_2|^2)^{1/2}):\ol{w}_2)$. Namely,
$\psi$ is explicitly given by 
\begin{equation}\label{eqn-correspon} 
\psi:(w_1:w_2:\ol{w}_2)\longmapsto 
\left(i\left(w_1\pm \sqrt{w_1^2+|w_2|^2}\right):\ol{w}_2\right).
\end{equation} 
(Note that the image of (\ref{eqn-correspon}) is considered as a point of 
$\Gamma_0/\langle\sigma\rangle$.)
First suppose 
$w_1\neq 0$. It is readily seen that we can suppose $w_1=1$. Then  in
(\ref{eqn-correspon})  the image becomes 
$(i(1\pm(1+|w_2|^2)^{1/2}):\ol{w}_2)$. Taking the sign `$+$', 
(\ref{eqn-correspon}) can be rewritten as
\begin{equation}\label{eqn-correspon2}
\psi:\mathbf C\ni w_2\mapsto 
\frac{-i\ol{w}_2}{1+\sqrt{1+|w_2|^2}},
\end{equation} where we use (the second entry)$/$(the first entry) as
an affine coordinate on $\Gamma_0$. The image of (\ref{eqn-correspon2}) is
clearly contained in the unit disk
$\{u\in\mathbf C\set |u|<1\}$. We show that (\ref{eqn-correspon2})
give a diffeomorphism between
$\mathbf C=\mathbf R^2$ and the unit disk. Putting $w_2=re^{i\theta}$,
(\ref{eqn-correspon2})  is rewritten as
\begin{equation}\label{eqn-correspon3}
\psi: re^{i\theta}\longmapsto
\frac{-ire^{-i\theta}}{1+\sqrt{1+r^2}}.
\end{equation} It is elementary to show that
$k(r):=r/(1+\sqrt{1+r^2})$ is differentiable on
$\{r>0\}$ and  its derivative is always positive, and that
$\lim_{r\ua\infty}k(r)=1$ and $\lim_{r\da 0}k(r)=0$ hold.
Hence $k$ gives  a bijection between $\{r\geq 0\}$ and
$\{0\leq s<1\}$.
It follows that (\ref{eqn-correspon2}) gives a bijection between
$\mathbf C$ and the unit disk.
 Moreover, the positivity of $k'$ implies 
that (\ref{eqn-correspon3}) is a diffeomorphism on $\mathbf C^*$.
For $w_2=0$, it can be easily checked that 
$(\partial w_2/\partial \ol{w}_2)(0)\neq 0$. Therefore (\ref{eqn-correspon2})
is a diffeomorphism on $\mathbf C$.

Next consider the case $w_1=0$. Then  we have $w_2\neq 0$, and the
image becomes
$(\pm i|w_2|:\ol{w}_2)=(1:\pm i\ol{w}_2/|w_2|)$. From this, it easily
follows that (\ref{eqn-correspon}) gives a diffeomorphism between the two
subsets
$\{(0:1:w)\set w\in U(1)\}$ and 
$\{(1:u)\in \Gamma_0\set u\in U(1)\}/\langle\sigma\rangle$. 
Moreover on $\mathbf{RP}^2\backslash\mathbf R^2$ we can use
$(1/r,\theta)$ as a local coordinate on $\mathbf{RP}^2$.
Then we can readily show that $(d/ds)(k(1/s))|_{s=0}\neq 0$.
This implies that $\psi$ is diffeomorphic also on a neighborhood of
$\mathbf{RP}^2\backslash\mathbf R^2$.
Note that the bijectivity of $\psi$ 
 implies that if
$l\neq l'$, then the corresponding rational curves in $Z$ are
disjoint. 

Next take any  real plane $H$ containing $l$.
On $H$ there is a one-dimensional family of lines through
$P_0$. Taking the inverse image, we obtain a one-dimensional 
holomorphic family $\mathcal L_{H}$
of  rational curves in $Z$, containing $L$ as a real member.
Any real member of $\mathcal L_{H}$ defines a conjugate pair of points
as the intersection with $\Gamma_0$. 
Consequently, real members of $\mathcal L_H$ determine a real circle $\mathcal C_H$
 in $\Gamma_0$.
If $s$ denotes the section of $N=N_{L/Z}$ associated to $\mathcal L_H$ (viewed as a deformation of $L$ in $Z$),
then Re$s(z)$ is non-zero by the diffeomorphicity of $\psi$, and
is represented by a tangent vector of
$\mathcal C_{H}$ at $z$, where we put $\{z,\ol{z}\}=\Gamma_0\cap L$.

Let $\{v_1,v_2\}$ be any oriented orthogonal basis of $T_z\Gamma_0$,
where we take the complex orientation and orthogonality.
Then since we have the isomorphism $\psi$,
there is a unique real plane $H_i$ ($i=1,2$) containing $l$ such that
$v_i$ is tangent to $\mathcal C_{H_i}$.
Let $s$ (resp.\,$t$)   be the global section of $N=N_{L/Z}$
associated to $\mathcal L_{H_1}$ (resp.\,$\mathcal L_{H_2}$).
 We now claim that
$as+bt$ does not vanish at $z$ and $\ol{z}$ simultaneously, 
unless $(a,b)= (0,0)$.
Putting $a=a_1+ia_2$ and $b=b_1+ib_2$, we readily have
\begin{equation}\label{eqn-reim4}
{\rm{Re}}(as+bt)=(a_1{\rm{Re}}s-b_2{\rm{Im}}t)+(b_1{\rm{Re}}t-a_2{\rm{Im}}s).
\end{equation}
Since (Re$s)(z)$ and (Re$s)(\ol{z})$ 
(resp.\,(Re$t)(z)$ and (Re$t)(\ol{z})$) are represented by  tangent vectors of
$\mathcal C_{H_1}$ (resp. $\mathcal C_{H_2}$), our choice of $H_1$ and $H_2$
implies that $({\rm{Re}}s)(z)$ is parallel to $({\rm{Im}}t)(z)$ and
$({\rm{Re}}t)(z)$ is parallel to $({\rm{Im}}s)(z)$.
The same is true at $\ol{z}$.
Hence by (\ref{eqn-reim4})
 if ${\rm{Re}}(as+bt)(z)=0$, then $a_1{\rm{Re}}s(z)=b_2{\rm{Im}}t(z)$
and $b_1{\rm{Re}}t(z)=a_2{\rm{Im}}s(z)$,
and ${\rm{Re}}(as+bt)(\ol{z})=0$ implies similar equalities.
Since ${\rm{Re}}s, {\rm{Re}}t,{\rm{Im}}s$ and ${\rm{Im}}t$ do not
be zero at both of $z$ and $\ol{z}$  as is already mentioned,
$a_1=0$ iff $b_2=0$ and $a_2=0$ iff $b_1=0$.
Therefore either $a_1b_2\neq 0$ or $b_1a_2\neq 0$ holds.
Suppose $a_1b_2\neq 0$. Then we show that 
$a_1{\rm{Re}}s(z)=b_2{\rm{Im}}t(z)$ and 
$a_1{\rm{Re}}s(\ol{z})=b_2{\rm{Im}}t(\ol{z})$ cannot hold simultaneously:
suppose that $a_1b_2>0$. Then 
${\rm{Re}}s(z)$ and ${\rm{Im}}t(z)$ have the same direction and
it follows that $\{{\rm{Re}}t(z),{\rm{Re}}s(z)\,
(=(b_2/a_1){{\rm{Im}}t(z)})\} $ is an
oriented basis of
$T_z\Gamma_0$. On the other hand, we have  
 ${\rm{Re}}s(\ol{z})=\sigma_*({\rm{Re}}s(z))$ and
 ${\rm{Re}}t(\ol{z})=\sigma_*({\rm{Re}}t(z))$. 
Therefore $\{{\rm{Re}}t(\ol{z}), {\rm{Re}}s(\ol{z})\}$ is an anti-oriented
basis of $T_{\ol{z}}\Gamma_0$ because $\sigma$ is orientation reversing.
On the other hand, $a_1{\rm{Re}}s(\ol{z})=b_2{\rm{Im}}t(\ol{z})$ and
$a_1b_2\neq 0$ imply that 
$\{{\rm{Re}}t(\ol{z}),{\rm{Re}}s(\ol{z})\} $ is an oriented basis of 
 $T_{\ol{z}}\Gamma_0$. This is a contradiction.
The case $b_1a_2>0$ is similar.
Therefore ${\rm{Re}}(as+bt)$ cannot be zero at 
$z$ and $\ol{z}$ simultaneously  provided $a_1b_2\neq 0$.
If $b_1a_2\neq 0$, then $b_1{\rm{Re}}t(z)=a_2{\rm{Im}}s(z)$ and 
$b_1{\rm{Re}}t(\ol{z})=a_2{\rm{Im}}s(\ol{z})$ do not hold at the same time.
 Thus we have shown that ${\rm{Re}}(as+bt)$ cannot be zero at 
$z$ and $\ol{z}$ simultaneously for any $(a,b)\neq (0,0)$. 
Hence so does $as+bt$. 
Therefore we get
$N\simeq O(1)^{\oplus 2}$ by Lemma \ref{lemma-nb} and complete a proof of Proposition \ref{prop-imline}.
\proofend

\vspace{3mm} Thus in our complex manifold $Z$ there actually exists a
connected family  of real lines. 
 Obviously this family is  $U(1)$-invariant, although
general members are not  $U(1)$-invariant:

\begin{prop} Among this family of real lines in $Z$, just one member is
$U(1)$-invariant. Further, the member is fixed  by $U(1)$ pointwisely.
\end{prop}

\noindent Proof. 
Recall that in a neighborhood of $p_0$,
$Z_0$ is defined by the equation
$z^2+w_1^2+w_2w_3=0$ ((\ref{eqn-Z_0})). It is immediate to see that the
$U(1)$-action looks like $(w_1,w_2,w_3)\mapsto
(w_1,tw_2,t^{-1}w_3)$ for
$t\in U(1)$. Thus using  homogeneous coordinates used in the last
proof, the $U(1)$-action on $\Gamma_0$ is given by 
$(u:v)\mapsto (u:tv)$ or $(u:v)\mapsto (u:t^{-1}v)$, depending on the
choice of a small resolution of $p_0$. Therefore only the real line
corresponding to 
$[(1:0)](=[(0:1)])\in \Gamma_0/\langle\sigma\rangle$ is $U(1)$-fixed. 
In view of (\ref{eqn-correspon}) and (\ref{coord}), the equation of this
line is explicitly given by
$y_2=y_3=0$, which is   pointwisely  $U(1)$-fixed by Proposition \ref{prop-def-B}.
 Since $\Phi:Z\ra\mathbf{CP}^3$ is 
$U(1)$-equivariant, it follows that the corresponding real line in $Z$ is also pointwise fixed. \proofend

\begin{definition}{\rm{
We will call the real lines in $Z$ obtained by Proposition \ref{prop-imline} real lines {\em  at  infinity}. Namely, a real line is said to be at infinity if its image in $\mathbf{CP}^3$ (by $\Phi$) is a line (necessarily going through $P_0$).
}}
\end{definition}

Because there are $\mathbf{RP}^2$'s worth of real lines in $\mathbf{CP}^3$ going through $P_0$,  real lines at infinity are parametrized by $\mathbf{RP}^2$. 
In the rest of this section we will study  deformations in $Z$ of real lines at infinity.

For any real plane $H$ going through $P_0$, we define a line bundle $\mathcal L_H$ over the smooth surface $S_H=\Phi^{-1}(H)$ (cf.\,Proposition \ref{prop-str_of_S}) by
$$\mathcal L_H=\Phi^*O_H(1)-\Gamma_0.$$
Of course, we consider $\Gamma_0$ as a divisor on $S_H$. 
Clearly $|\mathcal L_H|$ is a real pencil on $S_H$ and its real part is precisely the family of real lines at infinity lying on $S_H$.
Since $S_H$ is  rational, $\mathcal L_H\ra S_H$
  uniquely extends to a line bundle $\mathcal L_{H'}\ra S_{H'}$ for any real plane $H'$ sufficiently near to $H$, and that $|\mathcal L_{H'}|$ is still a real pencil whose real members are smooth rational curves parametrized by $S^1$.
Let us consider the particular case that the above planes $H$ and $H'$  contain $l_{\infty}$. 
Namely  $H=H_{\lambda_0}$ and  $H'=H_{\lambda}\in\linfty^{\sigma}$ with $\lambda\in I_4^-\cup I_4^+$ in the notations of Section \ref{s-nb}.
As introduced in Section \ref{ss-generic}, there are two families $$\mathcal L_{\lambda}^+=\{L^+_{\theta}\subset S_{\lambda}\set e^{i\theta}\in U(1)\},\hspace{2mm}\mathcal L_{\lambda}^-=\{L^-_{\theta}\subset S_{\lambda}\set e^{i\theta}\in U(1)\}$$ of real lines lying on $S_{\lambda}=S_{H_{\lambda}}$;
$L^+_{\theta}$ and $L^-_{\theta}$ are irreducible components of the inverse images of touching conics of generic type on $H_{\lambda}$.
It is  expected that either $\mathcal L_{\lambda}^+$ or $\mathcal L_{\lambda}^-$ is precisely the family obtained by deforming  real lines at infinity.
This is true,  since   the image of real line not at infinity is always a real touching conic (Proposition \ref{prop-image}) and since touching conics of generic type are the only one which do not go through $P_{\infty}$ and $\ol{P}_{\infty}$ on $H_{\lambda}$ (Propositions \ref{prop-a}, \ref{prop-b} and \ref{prop-c}), and since being real line  is preserved under small deformation.
Now recall that the distinction of $\mathcal L_{\lambda}^+$ and $\mathcal L_{\lambda}^-$ ($\lambda\in I_4^-\cup I_4^+$) made in Section \ref{ss-conn}.
We have been assumed 
$\Xi\cap L^-_{\theta}=\{x_2=e^{i\theta}h_0(\lambda)\}$ for $\lambda<\lambda_0$ and $\Xi\cap L^-_{\theta}=\{x_2=e^{-i\theta}h_0^{-1}(\lambda)\}$ for $\lambda>\lambda_0$, where $x_2=y_2/y_3$ is used as a coordinate on $\Xi\simeq l_{\infty}$ and $h_0$ is a function on $I_4$ introduced in Section \ref{ss-generic} that behaves as in Figure \ref{fig-h0}.
Assume for instance that members of $\mathcal L^+_{\lambda}$ are obtained as a deformation of real lines at infinity for $\lambda<\lambda_0$.
We claim that this assumption implies that the same is true for $\lambda>\lambda_0$.
Suppose not. 
Then members of $\mathcal L^-_{\lambda}$ must be deformations of real lines at infinity.
However, the behavior of $h_0$ displayed in Figure \ref{fig-h0} shows that, for any $\lambda<\lambda_0$ and for any $L^+_{\theta}\in\mathcal L^+_{\theta}$ there exists $\lambda'>\lambda_0$ such that $L^-_{\theta}\in\mathcal L^-_{\lambda'}$ intersects $L_{\theta}^+$ on $\Xi$.
Namely $L_{\theta}^+\subset S_{\lambda}$ intersects $L_{\theta}^-\subset S_{\lambda'}$  (cf.\,Proposition \ref{prop-break}). 
On the other hand, we already know that real lines at infinity has the right normal bundle (Proposition \ref{prop-imline}) and therefore $Z$ has a structure of  twistor space, at least on a neighborhood of real lines at infinity. 
In particular, $L^+_{\theta}\subset S_{\lambda}$ and  $L_{\theta}^-\subset S_{\lambda'}$ must be disjoint, because they are assumed to be obtained by deforming real lines at infinity.
This is a contradiction and we obtain that members of $\mathcal L^+_{\lambda}$ is obtained as a deformation of real lines at infinity for $\lambda>\lambda_0$.
By the same reason, if members of $\mathcal L^-_{\lambda}$ are obtained as  deformation of real lines at infinity for $\lambda<\lambda_0$, then the same is true for  $\mathcal L^-_{\lambda}$  for $\lambda>\lambda_0$.
Thus we have proved that real lines at infinity `connects' $\mathcal L^+_{\lambda} (\lambda<\lambda_0)$ and  $\mathcal L^+_{\lambda} (\lambda>\lambda_0)$, or $\mathcal L^-_{\lambda} (\lambda<\lambda_0)$ and  $\mathcal L^-_{\lambda}(\lambda>\lambda_0)$. 
It is obvious that just one of these two situations happens.
The following proposition shows that which one happens depends on the choice of a small resolution of $p_0$.

\begin{prop}\label{prop-dfmofline}
Changing a small resolution of $p_0\in Z_0$ switches which one of $\mathcal L^+_{\lambda}$ and $\mathcal L^-_{\lambda}$ ($\lambda\in I_4$) is obtained as deformation of real lines at infinity . 
Namely, which irreducible component of $\Phi^{-1}(C_{\theta})$ ($C_{\theta}$ being a touching conic of generic type) is a deformation of real lines at infinity depends on the choice of small resolution of $p_0$.
\end{prop}

\noindent Proof.
Recall that by Proposition \ref{prop-a}  touching conic of generic type on $H_{\lambda}$ is defined by 
\begin{equation}\label{eqn-tcgen}
C_{\theta}:\hspace{3mm}2(Q^2-f)y_1^2+\sqrt{f}e^{i\theta}y_2^2+2Qy_2y_3+\sqrt{f}e^{-i\theta}y_3^2=0.
\end{equation}
(We always need to keep in mind that $C_{\theta}$ depends on $\lambda$.)
Also recall that $H_{\lambda}$ contains $P_0$ precisely when $\lambda=\lambda_0$ which satisfies $(Q^2-f)(\lambda_0)=0$.
Hence substituting $\lambda=\lambda_0$ into the above equation, we get
$
e^{i\theta}y_2^2+2y_2y_3+e^{-i\theta}y_3^2=0.
$
This can be written as $
(e^{\frac{\theta}{2}}y_2+e^{-\frac{\theta}{2}}y_3)^2=0.
$
Let $l_{\theta}$ be the line on $ H_{\lambda_0}$ defined by $e^{\frac{\theta}{2}}y_2+e^{-\frac{\theta}{2}}y_3=0$ (and $y_0=\lambda_0y_1)$.
We have seen that $\lim_{\lambda\ra\lambda_0}C_{\theta}=l_{\theta}.$
Taking the inverse image, we obtain  
$$\lim_{\lambda\ra\lambda_0}\Phi^{-1}(C_{\theta})=L_{\theta}+\Gamma_0,$$
 where $L_{\theta}$ denotes the irreducible component of $\Phi^{-1}(l_{\theta})$ that are mapped surjectively  onto $l_{\theta}$. 
Since $\Phi^{-1}(C_{\theta})=L_{\theta}^++L_{\theta}^-$, 
it follows that  either
\begin{equation}\label{eqn-limit01}
\lim_{\lambda\ra\lambda_0}L_{\theta}^+=L_{\theta},\hspace{5mm}
\lim_{\lambda\ra\lambda_0}L_{\theta}^-=L_{\theta}+\Gamma_0
\end{equation}
or 
\begin{equation}\label{eqn-limit02}
\lim_{\lambda\ra\lambda_0}L_{\theta}^+=L_{\theta}+\Gamma_0,\hspace{5mm}
\lim_{\lambda\ra\lambda_0}L_{\theta}^-=L_{\theta}
\end{equation}
holds.
To prove the proposition it suffices to show the following:

\vspace{3mm}
\noindent
($\ast$) Changing a choice of small resolution of $p_0$ interchanges which one of (\ref{eqn-limit01}) and (\ref{eqn-limit02}) happens.
\vspace{3mm}

By applying $U(1)$-action, it suffices to prove ($\ast$) for the case $\theta=0$.
We write $C_{\theta=0}$ and $L^{\pm}_{\theta=0}$ to mean the curves 
obtained by substituting $\theta=0$ into $C_{\theta}$ and $L_{\theta}$ respectively.
We note  that two small resolutions of $p_0$ given in the proof of Proposition \ref{prop-imline} are also obtained by first taking the blowing-up at $p_0$ (in the usual sense) and then blowing-down the exceptional divisor ($\simeq\mathbf{CP}^1\times\mathbf{CP}^1$) in the two directions.
Consider  $\cup_{\lambda} L^+_{\theta=0}$ and  $\cup_{\lambda} L^-_{\theta=0}$ which are surfaces in $Z_0$.
(Remind that $L^{\pm}_{\theta}$ are contained in $S_{\lambda}$ and depend on $\lambda$.) 
To prove $(\ast)$ we consider {\em the tangent cones}  at $p_0$ of these surfaces.
As in the proof of Proposition \ref{sing. of B} putting $v_i=y_i/y_1$ and $(Q^2-f)(v_0)=g(v_0)(v_0-\lambda_0)^2$, we can use three functions
$$w_1=\sqrt{g(v_0)}\cdot(v_0-\lambda_0),
w_2=\sqrt{2Q(v_0,1)+v_2v_3}\cdot v_2,
w_3= \sqrt{2Q(v_0,1)+v_2v_3}\cdot v_3
$$ 
as a local coordinate around $P_0$.
Then rewriting (\ref{eqn-tcgen}), the equation of $C_{\theta=0}\subset H_{\lambda}$ in a neighborhood of $P_0$ becomes
$$
w_1=\sqrt{g}(\lambda-\lambda_0),
$$
$$
\sqrt{f}\left(w_2^2+w_3^2\right)+2\left\{(Q^2-f)\sqrt{Q^2+w_2w_3}+Q\left(Q^2-f+w_2w_3\right)\right\}=0.
$$
Then since we have $Q^2-f=w_1^2$ and since $Q(\lambda_0,1)\neq 0$ and $f(\lambda_0)\neq 0$ (basically because we are assuming Condition (A) of Proposition \ref{prop-necessa}),  a surface $\cup_{\lambda}C_{\theta=0}$ is singular at $P_0$ and its tangent cone is
$$\sqrt{f}\left(w_2^2+w_3^2\right)+2Qw_1^2+2Q\left(w_1^2+w_2w_3\right)=0,
$$
where we consider $Q$ and $f$ as a function of $w_1$.
But since $Q(\lambda_0,1)=\sqrt{f(\lambda_0)}$ we get an equation of the tangent cone of $\cup_{\lambda}C_{\theta=0}$ to be 
$$
4w_1^2+2w_2w_3+w_2^2+w_3^2=0.
$$
Thus the tangent cone splits into a union of two planes
$2w_1=\pm i(w_2+w_3).$
On the other hand, in a neighborhood of $p_0$, $Z_0$ was defined by $z^2+w_1^2+w_2w_3=0$.
Hence a equation of the tangent cone of $\cup_{\lambda}\Phi^{-1}_0(C_{\theta=0})$ at $p_0$ becomes 
$$
\left\{2z-(w_2-w_3)\right\}\left\{2z+(w_2-w_3)\right\}=0,\,\,\,2w_1=\pm i(w_2+w_3).
$$
which is also a union of two planes.
On the exceptional divisor of blowing-up at $p_0$, these equations define a reduced reducible curve of bidegree $(2,2)$ constituting of four irreducible components, two of which are conjugate pair of curves of bidegree $(1,0)$ and the other two are conjugate pair of curves of bidegree $(0,1)$.
Hence we get that the tangent cone of  $\cup_{\lambda}\Phi_0^{-1}(C_{\theta=0})$ at $p_0$ defines, on the exceptional divisor $(\simeq \mathbf{CP}^1\times\mathbf{CP}^1)$, the reducible curve of bidegree $(2,2)$.
On the other hand,   we know that each irreducible component of $\Phi^{-1}_0(C_{\theta})$ is real (Proposition \ref{prop-inv-a}) (iv)).
Moreover, easy calculations using local coordinate shows that real structure on the exceptional divisor acts as a product of complex conjugations on each factor and does not exchange the factors.
Therefore, on the full-blowup,  either
$$\lim_{\lambda\ra\lambda_0}L_{\theta=0}^+=L_{\theta=0}+(2,0),\hspace{2mm}\lim_{\lambda\ra\lambda_0}L_{\theta=0}^-=L_{\theta=0}+(0,2)$$
or
$$\lim_{\lambda\ra\lambda_0}L_{\theta=0}^+=L_{\theta=0}+(0,2),\hspace{2mm}\lim_{\lambda\ra\lambda_0}L_{\theta=0}^-=L_{\theta=0}+(2,0)$$
holds, where $(2,0)$ (resp. $(0,2)$) means a reducible curve of bidegree $(2,0)$ (resp. $(0,2)$) that is a conjugate pair of curves of bidegree $(1,0)$ (resp. $(0,1)$).
If we blow down the exceptional divisor along the direction whose fibers are curves of bidegree $(2,0)$, we get, on the blown-down space $Z$, that either
$$\lim_{\lambda\ra\lambda_0}L_{\theta=0}^+=L_{\theta=0},\hspace{2mm}\lim_{\lambda\ra\lambda_0}L_{\theta=0}^-=L_{\theta=0}+\Gamma_0$$
or
$$\lim_{\lambda\ra\lambda_0}L_{\theta=0}^+=L_{\theta=0}+\Gamma_0,\hspace{2mm}\lim_{\lambda\ra\lambda_0}L_{\theta=0}^-=L_{\theta=0}$$
holds, depending on the above two cases respectively.
Also we have an alternative conclusion for another blow-down of the exceptional divisor.
This implies that just one of (\ref{eqn-limit01}) and (\ref{eqn-limit02}) happens, depending on the direction of blow-down of the exceptional divisor.
This implies the desired conclusion.\proofend

\vspace{3mm}
Combined Proposition \ref{prop-limit02} with Proposition \ref{prop-dfmofline}, we obtain the following proposition which will be needed in our proof of the main theorem.

\begin{prop}\label{prop-summary}
For each of the two small resolutions  $Z'_0\ra Z_0$ determined in Proposition \ref{prop-sr}, there exists a unique  small resolution $\nu:Z\ra Z'_0$ of the ordinary double point satisfying the following properties:

Let $\mu:Z\ra Z_0$ be the composition of the two resolutions, and put $\Phi=\Phi_0\mu$. Then 
(i) $\Phi$ preserves the real structure,
(ii) for each $H_{\lambda}\in \linfty^{\sigma}$, $\lambda\in I_1\cup I_2\cup I_3\cup I_4$, there exists a family of real lines on $S_{\lambda}=\Phi^{-1}(H_{\lambda})$ whose parameter space is $S^1$, thereby getting  families of real lines parametrized by $I_{j}\times S^1$, $1\leq j\leq 4$,
(iii) for $\lambda=-1,0,a,\infty$, the intersection of the two irreducible components of $\Phi^{-1}(H_{\lambda})$ is a real line in $Z$,
(iv) when $\lambda$ approaches the endpoints of $I_{j}$, the real lines in (ii) converge to the real lines in (iii),
(v) if $\lambda=\lambda_0$, the $S^1$-family of real lines on $S_{\lambda_0}$ in (ii) are real lines at infinity,
(vi) when $\lambda\in I_4^-\cup I_4^+$ approaches $\lambda_0$, the real lines on $S_{\lambda}$ obtained in (ii) converge to the real lines at infinity. 

Moreover, if we change the small resolution ((I) and (II) of Proposition \ref{prop-sr}), then the above resolution of the ordinary double point also changes.
\end{prop}

The point of this proposition is (iv) and (v) which imply that the $S^1$-family of real lines on $S_{\lambda}$, $\lambda\in I_1\cup I_2\cup I_3\cup I_4^-\cup I_4^+$ obtained in (ii) form a connected family, whose parameter space is a union of four spheres, where `North pole' of each sphere is identified with `South pole' of the next sphere.

\section{Global construction of twistor lines and their disjointness}\label{s-disj}

In Section \ref{ss-con} we have proved that there are only two (small) resolutions of the conjugate pair of singularities of $Z_0$ which can yield a twistor space (Proposition \ref{prop-sr}).
In the previous section we have shown that once one of the above two resolutions is given, a resolution   of the remaining singularity (the real ordinary double point) is automatically determined and that, then, real lines contained in $S_{\lambda}$, $\lambda\in\mathbf R\cup \{\infty\}$ naturally forms a connected family. 
Because the parameter space of this family is  real two-dimensional, and because we already know that the normal bundles  in $Z$ of the members are $O(1)^{\oplus 2}$, the family is not a complete family. 
In this section we first show that members of this family can be extended in $Z$ to give (real) four-dimensional family of real rational curves (Proposition \ref{prop-globalext}).
Next we show that these real rational curves are disjoint and further cover $Z$; namely they foliate $Z$ (Proposition \ref{prop-disj}).

For these purposes, we briefly recall notations and  results in Section \ref{s-how}.
Let $(\mathbf{RP}^3)^{\vee}$ be the dual projective planes of real planes in $\mathbf{CP}^3$,  and $\mathbf{RP}^2_{\infty}\subset(\mathbf{RP}^3)^{\vee} $ the set of real planes going through $P_0$. 
Then $\linfty^{\sigma}$, the set of real planes containing $l_{\infty}$,  becomes a line in $(\mathbf{RP}^3)^{\vee}$, intersecting $\mathbf{RP}^2_{\infty}$ transversally at a point (for which we have denoted $H_{\lambda_0}$).
Clearly we have $(\mathbf{RP}^3)^{\vee}\backslash \mathbf{RP}^2_{\infty}=\mathbf R^3$. 
By Lemma \ref{lemma-singsection}, $\mathbf R^3\backslash \linfty^{\sigma}$ is precisely the set of real planes intersecting $B$ smoothly; 
namely we have $U^{\sigma}=\mathbf R^3\backslash \linfty^{\sigma}(\simeq\mathbf C^*\times\mathbf R)$.
%
We have constructed  a fibration ${\mathcal C}_{U^{\sigma}}\ra U^{\sigma}$  whose fiber over $H\in U^{\sigma}$ is the set of touching conics of $H\cap B$. 
 ${\mathcal C}_{U^{\sigma}}\ra U^{\sigma}$ is a fiber bundle whose fiber is a union of 63 copies of $\mathbf{CP}^1$ (Proposition \ref{prop-ftc4}).
Recalling that the inverse image of touching conics splits into two curves and taking the inverse image by the double cover $Z_0\ra\mathbf{CP}^3$, we get a fiber bundle $\tilde{\mathcal C}_{U^{\sigma}}\ra U^{\sigma}$ whose fiber is a union of 128 copies of $\mathbf{CP}^1$. 
Each $\mathbf{CP}^1$ of the fibers of $\tilde{\mathcal C}_{U^{\sigma}}\ra U^{\sigma}$  is a parameter space of pencils on  $S_H$, whose image is a one-dimensional family of touching conics.
The following proposition gives all twistor lines in $Z$:

\begin{prop}\label{prop-globalext}
For any real plane $H$ different from the four planes on which the tropes of $B$ lie, we can find a $S^1$-family of real smooth  rational curves on $S_H=\Phi^{-1}(H)$ satisfying the following properties:
(i) if $H\in\linfty^{\sigma}$, the $S^1$-family on $S_H$ is just the family of real lines obtained in (ii) of Proposition \ref{prop-summary}, 
(ii) if $H\in\mathbf{RP}^2_{\infty}$, the $S^1$-family is just the family of real lines at infinity whose images are real lines in $H$, 
(iii)  there is a real connected component of the total space of $\tilde{\mathcal C}_{U^{\sigma}}\ra U^{\sigma}$ such that the real circles of each fiber ($\mathbf{CP}^1$) is just the $S^1$-family of real curves on $S_H$  for any $H\in U^{\sigma}$,
(iv) for any two real planes $H$ and $H'$, and for any members $L\subset S_H$ and $L'\subset S_{H'}$ of the $S^1$-families, $L$ and $L'$ can be connected by deformation in $Z$ preserving the real structure,
(v) when a real plane moves to any one of the four planes on which the tropes of $B$ lie, members of the $S^1$-family converge to the real lines over the tropes.
\end{prop}

\begin{figure}[htbp]
\includegraphics{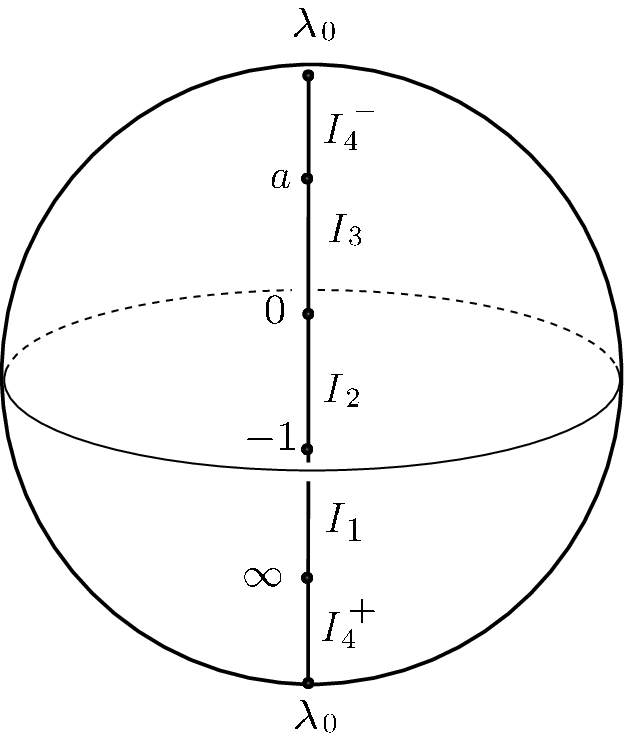}
\caption{$\linfty^{\sigma}$ can be considered as an `axis' in $\mathbf U$.}
\label{fig-sphere}
\end{figure}

\noindent Proof.
This is a compilation of the results we have obtained in this paper so far.
In this proof we put $\mathbf U=(\mathbf{RP}^3)^{\vee}\backslash$(4 points), where the four points represent the four planes on which the tropes of $B$ are lying.
Note that $U^{\sigma}$ is an open subset of $\mathbf U$.
For any $H\in\mathbf U$, $S_H=\Phi^{-1}(H)$ is a smooth rational surface with $c_1^2=2$ (Proposition \ref{prop-str_of_S}).
So we can consider the local system $\mathcal H^2$ over $\mathbf U$ whose fiber over $H$ is the cohomology group $H^2(S_H,\mathbf Z)$ $(\simeq\mathbf Z^8)$.
Note that $\mathcal H^2\ra\mathbf U$ is not necessarily trivial, since $\pi_1(\mathbf U)\simeq\pi_1(\mathbf{RP}^3)\simeq\mathbf Z_2$.
The four intervals $I_j$, $1\leq j\leq 4$, can naturally be considered as a subset of $\mathbf U$ (Figure \ref{fig-sphere}). 
For each $I_j$  consider the $S^1$-family of real lines on $S_{\lambda}=S_{H_{\lambda}}$ for $\lambda\in I_j$ obtained in Proposition \ref{prop-summary} (ii).
As in Section \ref{s-inv}, this $S^1$-family is the real part of a real pencil on $S_{\lambda}$.
Let $\mathcal L_{\lambda}\ra S_{\lambda}$ be the associated real line bundle.
(In Section \ref{s-nb}, $\mathcal L_{\lambda}$ was used to mean the $S^1$-family of real lines on $S_{\lambda}$. 
But there should be no confusion since the real lines are simply the real members of $|\mathcal L_{\lambda}|$ for the present $\mathcal L_{\lambda}\ra S_{\lambda}$.)
Because these $S^1$-families form a connected family (Proposition \ref{prop-summary}), the cohomology class determined by $\mathcal L_{\lambda}$ is independent of $I_j$ to which $\lambda$ belongs. 
Therefore, by associating the cohomology class of $\mathcal L_{\lambda}$, we obtain a section (of $\mathcal H^2$) over $\linfty^{\sigma}\backslash$(the 4 points). 
This section  can be naturally extended to a section over $\mathbf{U}\backslash\mathbf{RP}^2_{\infty}$,
since $\mathbf U\backslash\mathbf{RP}^2_{\infty}$ is contractible to $\linfty^{\sigma}$.
On the other hand, we have a section over $\mathbf{RP}^2_{\infty}$ of $\mathcal H^2$, which associates the cohomology class of the line bundle $\Phi^*O_H(1)-\Gamma_0$ for each $H\in \mathbf{RP}^2_{\infty}$.
These two sections fit together to give a global section of $\mathcal H^2\ra\mathbf U$, since when $\lambda$ approaches $\lambda_0$, $\mathcal L_{\lambda}$ converges to the line bundle $\Phi^*O_H(1)-\Gamma_0$ by Proposition \ref{prop-summary} (vi).
 This section must be real, since $\mathcal L_{\lambda}$ is real.
 Thus we have obtained a real global section of $\mathcal H^2\ra\mathbf U$.
 Recalling that $S_H$ is a smooth rational surface (implying $H^2(S_H,\mathbf Z)\simeq$ Pic $S_H$), the section uniquely determines a real line bundle $\mathcal L_H\ra S_H$ for any $H\in\mathbf U$.
Now we claim that $|\mathcal L_H|$ is a (real) pencil for any $H\in\mathbf U$.
First,  if $H\in \linfty^{\sigma}$ or $H\in\mathbf{RP}^2_{\infty}$, the claim is obvious from the above construction of the section.
(This  implies (i) and  (ii) of the proposition.)
So it suffices to prove the claim for $H\in U^{\sigma}$.
Since we have $H^0(S_{\lambda},\mathcal L_{\lambda})\simeq\mathbf C^2$ and $H^i(S_{\lambda},\mathcal L_{\lambda})=0$ for $i\geq 1$, it is easily verified that  $|\mathcal L_H|$ is also a (real) pencil for any real plane $H$ sufficiently near to the original planes $H_{\lambda}$.
This implies that there exists a neighborhood $W$ of $\linfty^{\sigma}$ in $\mathbf U$ such that $|\mathcal  L_H|$ is a pencil for any $H\in W$.
Real members of this pencil can be assumed to be real lines, since real members of $|\mathcal L_{\lambda}|$ are actually  real lines, and since  $N\simeq O(1)^{\oplus 2}$ is an open condition.
Then by Proposition \ref{prop-image}, the image of these real lines on $S_H$ must be  real touching conics on $H$ and these form a $S^1$-family of real touching conics on $H$.
If $H$ is moved in $U^{\sigma}$, then the $S^1$-family of real touching conics automatically survives since as in Proposition \ref{prop-ftc3} any smooth quartic has 63 families of touching conics.
Taking the inverse image, this implies that $|\mathcal L_H|$ is a pencil for any $H\in U^{\sigma}$, as claimed.

Thus we have specified  a real pencil $|\mathcal L_H|$ for any $H\in \mathbf U$.
We claim that any real member of this pencil must be irreducible.
For $H\in \linfty^{\sigma}$ this is obvious from our explicit construction of the family of touching conics.
For $H\in\mathbf{RP}^2_{\infty}$, this is also immediate from Proposition \ref{prop-imline}.
For $H\in U^{\sigma}$, $H\cap B$ is a smooth quartic curve.
It follows that reducible members of the  pencil $|\mathcal L_H|$ must be a sum of two smooth rational curves intersecting transversally at a unique point (cf. the proof of Proposition \ref{prop-ftc3}).
Therefore real reducible member of the pencil gives a real point, contradicting the fact that $Z$ has no real point (Proposition \ref{prop-realpoint}).
Thus $|\mathcal L_H|$ has no real reducible member  for any $H\in \mathbf U$.

Then for each $H\in \mathbf U$ we associate $|\mathcal L_H|^{\sigma}$ which forms $S^1$-family of  smooth real rational curves on $S_H$.
Then it is obvious from the above construction that (i)--(iv) of the proposition hold.
Further, (v) is also immediate from Proposition \ref{prop-summary} (iv).
\proofend

\begin{rmk}{\rm
In view of our proof of Proposition \ref{prop-globalext}, one may expect that the family of  line bundles $\mathcal L_H\ra S_H$, $H\in \mathbf U$ are further extended to $\Phi^{-1}(H_i)$, $1\leq i\leq 4$, where $H_i$ are the four planes on which the tropes of $B$ lie.
However this is not true, since the real line lying over the tropes are not Cartier divisors on $\Phi^{-1}(H_i)$. 
(Recall that $\Phi^{-1}(H_i)$ is a reducible surface and that the two irreducible components intersect transversally along the real lines over the tropes.)

}
\end{rmk}

The following proposition, saying that the real rational curves in Proposition \ref{prop-globalext} foliate $Z$, is the final main step in our proof of the main theorem  given in the next section:
\begin{prop}\label{prop-disj}
For any point of $Z$ there uniquely exists a unique real smooth rational curves obtained in Proposition \ref{prop-globalext} going through the point.
\end{prop}

\noindent Proof.
Because our proof is somewhat long, we give an outline.
First we take any real line $l$ in $\mathbf{CP}^3$.
Then  $E:=\Phi^{-1}(l)$ is a real curves in $Z$ which becomes either a smooth elliptic curve or a cycle of smooth rational curves consisting of two or eight irreducible components, depending on how $l$ intersects $B$.
Next we show the following claim:

\vspace{2mm}\noindent
 $(\sharp)$ For any point of $E$ there uniquely exists a real plane $H$ containing $l$, on which there exists a member $L$ of the $S^1$-family (obtained in Proposition \ref{prop-globalext}) going through the given point.
 
 \vspace{2mm}\noindent
This directly implies the claim of the proposition.
In the following, let $\mathbf U$ and   $\mathcal L_H\ra S_H$ ($H\in\mathbf U$) have the same meaning in the previous proof.
In particular,  the $S^1$-family of real smooth rational curves on $S_H$ of Proposition \ref{prop-globalext} is just the real part $|\mathcal L_H|^{\sigma}$ of $|\mathcal L_H|$.
Moreover, for a real line $l\subset\mathbf{CP}^3$, $\langle l\rangle$ denotes the pencil of planes containing $l$ and $\langle l\rangle^{\sigma}$ denotes its real members parametrized by $S^1$.

First we consider the case that the real line $l$ intersects $B$ transversally (at four points).
In this case $E=\Phi^{-1}(l)$ is a real smooth elliptic curve.
For each real plane $H$ containing $l$ we define a subset of $E$ by
$$\mathcal T_H=\{L\cap E\set L\in |\mathcal L_H|^{\sigma}\}.$$
Since $(E,L)_{S_H}=2$ (by adjunction), and since $E$ has no real point, $L\cap E$ is a pair of conjugate points for any $L\in|\mathcal L_H|^{\sigma}.$
It is readily verified that the restriction map $H^0(S_H,\mathcal L_H)\ra H^0(E,\mathcal L_H|_E)$ is isomorphic.
Therefore $\mathcal T_H$ is also the union of the zero locus of all real sections of the line bundle $\mathcal L_H|_E$.
This implies that if $\mathcal T_H\cap \mathcal T_{H'}$ is non-empty, then $\mathcal T_H=\mathcal T_{H'}$ (as a subset of $E$).

Next we investigate structure of $\mathcal T_H$.
We show in this paragraph that $\mathcal T_H$ is a disjoint union of two circles $\mathcal T_H^+$ and $\mathcal T_H^-$ such that (a) $\sigma(\mathcal T_H^+)=\mathcal T_H^-$, (b) $\mathcal T_H^+$ and $\mathcal T_H^-$ are homologous but not zero-homologous.
Since the parameter space of $|\mathcal L_H|^{\sigma}$ is a circle, it is obvious from the beginning that $\mathcal T_H$ is a (connected) circle or a union of two disjoint circles. 
To show that the former case cannot happen, consider the double covering $E\ra l$.
Since $E$ has no real point (because $Z$ has no real point by Proposition \ref{prop-realpoint}), the four branch points of $E\ra l$ are  not real. 
Since $l$ is  real, we can decompose $l$ as $l^+\cup l^{\sigma}\cup l^-$, where $l^{\sigma}$ is the real locus (a circle) of $l$ and $l^{\pm}$ are hemi-spheres satisfying $\sigma(l^+)=l^-$.
First we consider the case that  $H\in\langle l\rangle^{\sigma}$ does not go through $P_0$.
In this case $\Phi$ gives an isomorphism between $L\in |\mathcal L_H|^{\sigma}$ and $\Phi(L)$ which is a real touching conic. 
Hence (because $L$ does not have real point) $\Phi(L)\cap l$ consists of two points, one of which belongs to $l^+$ and the other belongs to $l^-$. 
Therefore the set $\{\Phi(L)\cap l\set L\in|\mathcal L_H|^{\sigma}\}\subset l$ must be two circles. 
It follows that $\mathcal T_H$ is a union of two circles, which we denote by $\mathcal T_H^+$ and $\mathcal T_H^-$. 
The decomposition $l=l^+\cup l^{\sigma}\cup l^-$ yields the decomposition $E=E^+\cup \Phi^{-1}(l^{\sigma})\cup E^-,$ where we set $E^+=\Phi^{-1}(l^+)$ and $E^-=\Phi^{-1}(l^-)$.
We can suppose $\mathcal T^+_H\subset E^+$ and $\mathcal T_H^-\subset E^-.$
Then obviously $\sigma(\mathcal T_H^+)=\mathcal T_H^-$ and we get the claim (a) for $H$ not going through $P_0$. 
We postpone (b) for planes of this kind and consider a unique real plane $H_0\in\langle l\rangle^{\sigma}$ going through $P_0$.
It is readily seen that $\mathcal T_{H_0}=\Phi^{-1}(l^{\sigma})$.
Because the branch points of $E\ra l$ are non-real, $\Phi^{-1}(l^{\sigma})$ is a circle or a disjoint union of two circles.
Suppose the former.
Then since $E\backslash \Phi^{-1}(l^{\sigma})$ is evidently disconnected, $\Phi^{-1}(l^{\sigma})$ becomes homologous to zero and hence bounds a disk.
$\sigma$ clearly has to act on this disk, so has a fixed point by the Brouwer fixed point theorem.
This is a contradiction and we conclude that $\Phi^{-1}(l^{\sigma})$ is a disjoint union of two circles.
We also denote these circles by $\mathcal T^+_{H_0}$ and $\mathcal T_{H_0}^-.$
Then since $\Phi(\mathcal T^+_{H_0})=\Phi(\mathcal T^-_{H_0})=l^{\sigma}$, we have $\sigma(\mathcal T_{H_0}^+)=\mathcal T_{H_0}^-$ and we get (a) for $H=H_0$.
Further, since $\sigma$ gives an automorphism of $E$ and since $\mathcal T^+_{H_0}$ and $\mathcal T_{H_0}^-$ are disjoint, $\mathcal T^+_{H_0}$ and $\mathcal T_{H_0}^-$ must be mutually homologous. 
Furthermore, they are not homologous to zero since otherwise $E\backslash \Phi^{-1}(l^{\sigma})$ would have three connected components, which is not the case.
Thus we obtained (b) for $H=H_0$.
This immediately implies  (b) for $H\neq H_0$, since when $H$ continuously moves in $\langle l\rangle^{\sigma}$, $\mathcal T_H$ also moves in $E$ continuously. 

Next we show that $\{\mathcal T_H\set H\in\langle l\rangle^{\sigma}\}$ gives a foliation of $E$. 
Namely we show 
\begin{equation}\label{eqn-alpha1}
E=\displaystyle\bigcup_{H\in \langle l\rangle^{\sigma}}\mathcal T_H
=\displaystyle\bigcup_{H\in \langle l\rangle^{\sigma}}(\mathcal T^+_H\cup\mathcal T^-_H) \hspace{5mm}{\mbox{(disjoint union)}}
\end{equation}

\begin{figure}[htbp]
\includegraphics{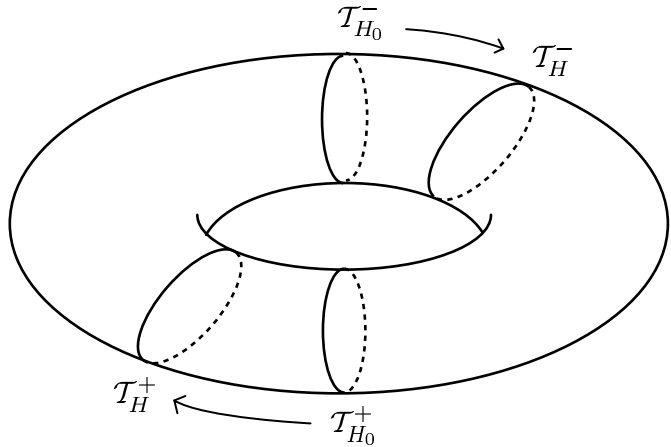}
\caption{foliation on $E$}
\label{fig-torus1}
\end{figure}

\noindent
(Figure \ref{fig-torus1}) To prove (\ref{eqn-alpha1}) let Pic$^2 E$ be the set of line bundles on $E$ of degree two.
Of course, Pic$^2 E$ can be identified with Pic$^0 E$ (though there is no canonical isomorphism).
Then we define a map $$\alpha_l: \langle l\rangle^{\sigma}\ra\mbox{Pic}^2 E;
\hspace{3mm} H\mapsto \mathcal L_H|_E.$$
Then $\alpha_l(H)=\alpha_l(H')$ if and only if $\mathcal T_H=\mathcal T_{H'}$.
Further, as is already seen, $\mathcal T_H=\mathcal T_{H'}$ if and only if they are not disjoint.
Hence to prove (\ref{eqn-alpha1}) it suffices to show that $\alpha$ is diffeomorphic.
Since the line bundle $\mathcal L_H\ra S_H$  can be holomorphically extended   to  that on $S_{H'}$ for any $H'\in\langle l\rangle$ sufficiently near to $H$, $\alpha_l$ can be extended to a holomorphic map $\tilde{\alpha}_l:A\ra$ Pic$^2E$, where $A$ is an annulus which is an open neighborhood of  $\langle l\rangle^{\sigma}$ in  $\langle l\rangle$. 
In particular, the ramification index of $\alpha$ makes sense for any point of $\langle l\rangle^{\sigma}$.
Let $H_0\in\langle l\rangle^{\sigma}$ be the unique plane satisfying $P_0\in H_0$ as in the last paragraph.
We then have $\alpha_l(H_0)=\Phi^*O_l(2)$.
Moreover, since $L\in|\mathcal L_H|^{\sigma}$ is mapped biholomorphically onto a touching conic for any $H\neq H_0$, $\alpha_l(H)\neq\Phi^*O_l(2)$ for any $H\neq H_0$.
In particular, $\alpha_l$ is not a constant map.
Now recall that we have already seen in Proposition \ref{prop-imline} that $N_{L/Z}\simeq O(1)^{\oplus 2}$ for $L\in |\mathcal L_{H_0}|^{\sigma}$ and therefore $Z$ has a structure of twistor space in a neighborhood of $L$.
This implies that the ramification index of $\alpha_l$ at $H_0$ is one, since 
otherwise $L$ can be deformed in $Z$ in such a way that its Kodaira-Spencer class in $H^0(N_{L/Z})$ vanishes at two points $(=L\cap E)$, which means that $N_{L/Z}$ contains $O(2)$ as a subbundle.
On the other hand, the image of $\alpha_l$ is contained in (Pic$^2 E)^{\sigma}$. 
Because Pic$^2 E$ is an elliptic curve, its real locus is, if not empty, either a circle or a union of two circles. 
In any way, the image of $\alpha_l$ is contained in the connected component of (Pic$^2 E)^{\sigma}$ containing $\Phi^*O_l(1)$.
In particular the degree of $\alpha_l$ as a differential map from $S^1$ to $S^1$ makes sense.
Now since we have already shown that the ramification index of $\alpha_l$  at $H_0$ is one and that $\alpha_l^{-1}(\Phi^*O_l(1))=\{H_0\}$, it follows that the degree of $\alpha_l$ is one. 
Hence $\alpha_l$ must be surjective.
It remains to show the injectivity of $\alpha_l$, or more precisely, that the differential $d\alpha_l$ does not vanishes everywhere on $\langle l\rangle^{\sigma}$.
Suppose that $H\in\langle l\rangle^{\sigma}$ is a critical point of $\alpha_l$.
This implies, as in the proof of Proposition \ref{prop-nbc}, that $N_{L/Z}$ contains $O(k)$ as a holomorphic subbundle, where $k$ denotes the ramification index of $\tilde{\alpha}_l$ at $H\in\langle l\rangle$. 
It is readily seen that $N_{L/Z}$ is either $O\oplus O(2)$ or $O(1)^{\oplus 2}$.
Therefore we have $k=2$ for any critical points (if any).
Under the assumption of the existence, 
 the number of critical points of $\alpha_l$ is obviously greater than one.
On the other hand, our choice of real lines over $H\in\linfty^{\sigma}$ (done in Section \ref{ss-generic}--\ref{ss-con}) implies that for any $L\in\mathcal L_{H_{\lambda}}$, $H_{\lambda}\in\linfty^{\sigma}$, we have $N_{L/Z}\simeq O(1)^{\oplus 2}$
(because we have chosen the touching conics on $H_{\lambda}$ in such a way that the function $h_i$ ($0\leq i\leq 3$) have no critical points).
Therefore, (because $N\simeq O(1)^{\oplus 2}$ is an open condition) we have  
$N_{L/Z}\simeq O(1)^{\oplus 2}$ for  any $L\in|\mathcal L_{H}|^{\sigma}$, 
where $H\in\langle l'\rangle^{\sigma}$, $l'$  being a real line sufficiently near to $l_{\infty}$.
Hence for this real line $l'$, $\alpha_{l'}$ must be diffeomorphic and does not have critical points.
Now consider a smooth real one-dimensional family of real lines $\{l_t\set0\leq t\leq 1\}$ such that $l_0=l'$ and $l_1=l$, and such that $l_t\cap B$ is transversal for any $t$.
Then the map $\alpha_t=\alpha_{l_t}$ clearly varies differentiably. 
Because $\alpha_0$ is diffeomorphic, and because $\alpha_1$ was assumed to have at least two critical points, there must be some $t$, $0<t<1$, such that $\alpha_t$ has a critical point whose ramification index is greater than 2.
This implies that for $L_t\in |\mathcal L_{H_t}|^{\sigma}$, $H_t\in\langle l_t\rangle^{\sigma}$, $N_{L_t/Z}$ contains $O(k)$, $k\geq 3$ as a subbundle.
This is a contradiction and we get that $\alpha_l$ does not have critical points, as claimed.
%
%
Thus we obtain (\ref{eqn-alpha1}), which of course imply that $(\sharp)$ is true for  our $l$ intersecting $B$ transversally.

Next we consider the case that the real line $l$ does not intersect $B$ transversally.
Then we have $l\not\subset B$ since $B$ does not have real point except $P_0$.
Hence $l\cap B$ consists of  one, two or three points.
If $l\cap B$ is one or three point, then it must contain $P_0$ by reality. 
But we have already shown in the proof of Proposition \ref{prop-imline} that if a real line $l$ goes through $P_0$, then $l\cap B$ consists of three points.
Hence we consider the case that $l\cap B$ consists of three points, one of which is $P_0$.
In this case $\Phi^{-1}(l)$ is a union of the exceptional curve $\Gamma_0$ and the real line $L$ over $l$, intersecting transversally at two points (cf. the proof of Proposition \ref{prop-imline}). 
As proved in Proposition \ref{prop-imline}, real lines $L$ of this form (which were called `at infinity') are disjoint. 
Furthermore, the isomorphism $\psi$ between $\Gamma_0/\langle\sigma\rangle$ and the projective space of real lines through $P_0$ in the proof means that, for any point of $\Gamma_0$, there uniquely exists a real line at infinity through the point. 
This directly implies $(\sharp)$  for real lines $l$ going through $P_0$.

\begin{figure}[htbp]
\includegraphics{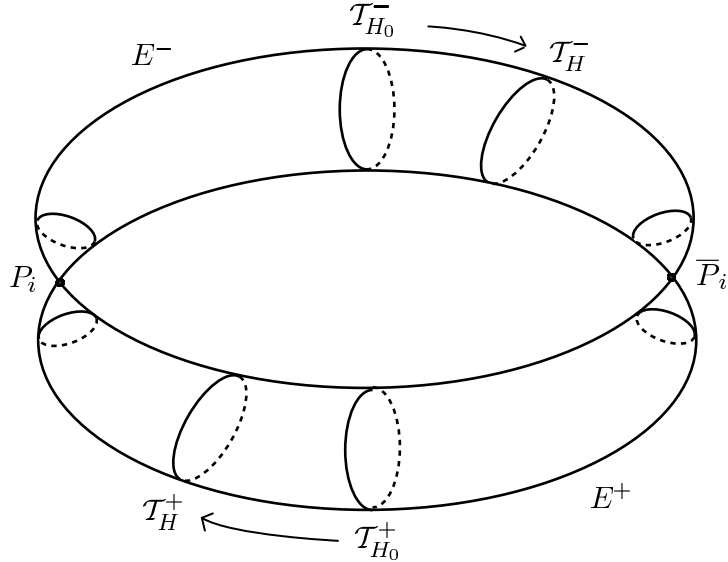}
\caption{foliation on $E=E^++E^-$}
\label{fig-torus2}
\end{figure}

So it remains to see the case that $l\cap B$ consists of a conjugate pair of points; namely $l$ being a real bitangent.
By Proposition \ref{prop-norealbitan}, such a bitangent $l$ must be always trivial, which means that $l$ lies on one of the four real planes on which the tropes lie.
Let $H_i\supset l$ denote this real plane.
First we assume $l\neq l_{\infty}$.
In this case $\Phi^{-1}(l)$ becomes a cycle of two smooth rational curves intersecting transversally over the tangent points.
We write $\Phi^{-1}(l)=E=E^++E^-$ with $E^-=\ol{E}^+$, and $P_i$ and $\ol{P}_i$ denote the intersection points of $E^+$ and $E^-$.
Further, write $\Phi^{-1}(H_i)=D_i+\ol{D}_i$ and $L_i=D_i\cap \ol{D}_i$, which is the only candidate of twistor line lying on  $\Phi^{-1}(H_i)=D_i+\ol{D}_i$.
Then we have $L_i\cap E=\{P_i,\ol{P}_i\}$.
Next, for any real plane $H\supset l$ different from $H_i$, we set $\mathcal T_H=\{L\cap  E\set L\in|\mathcal L_H|^{\sigma}\}$.
Then as in the case of $E$ being a smooth elliptic curve, each $L\cap E$ consists of a conjugate pair of points and $\mathcal T_H$ is a disjoint union of two circles $\mathcal T_H^+$ and $\mathcal T_H^-$, one of which is contained in $E^+$ and the other is in $E^-$.
We claim that 
\begin{equation}\label{eqn-fol2}
E\backslash \{P_i,\ol{P}_i\}=\bigcup_{H\in\langle l\rangle^{\sigma}\backslash\{ H_i\}}
\mathcal T_H
\end{equation}
holds; namely $E\backslash \{P_i,\ol{P}_i\}$
 is foliated by $\mathcal T_H$ (Figure \ref{fig-torus2}).
 
To show this, let Pic$^{1,1}E$ denote the subset of the Picard group of $E$ whose elements represent holomorphic line bundles on $E$ whose degree of their restriction onto $E^+$ and $E^-$ are both one. 
Take an affine coordinate $z$ (resp. $w$) on $E^+$ (resp. $E^-$), whose origin  represents $P_i$ and infinity represents $\ol{P}_i$.

\vspace{2mm}\noindent
{\textbf{Claim.}} Suppose that  $F\in$ Pic$^{1,1}E$ is real  (i.e. $\sigma^*F\simeq\ol{F}$). 
Then the union of the zero locus of real sections of $F$ consists of two circles $\mathcal T^+\subset E^+$ and $\mathcal T^-\subset E^-$ having the radiuses $|z|=r>0$ and $|w|=1/r$ respectively.
Moreover, $F\simeq F'$ iff the two circles coincides respectively.

\vspace{2mm}
We postpone a proof of this claim and continue the proof of Proposition \ref{prop-disj}.
We define a map
$$\alpha_l:\langle l\rangle^{\sigma}\backslash\{H_i\}\ra \hspace{2mm}{\mbox{Pic}}^{1,1}E;
\hspace{3mm} H\mapsto \mathcal L_H|_E.$$
Then by definition $\mathcal T_H$ is  again the union of the zero locus of real sections of $\alpha_l(H)=\mathcal L_H|_E$. 
Since $L\cap E$ belongs to the linear system $|\mathcal L_H|_E|$ for any $L\in |\mathcal L_H|^{\sigma}$,  $\mathcal T_H=\mathcal T_{H'}$ holds iff $\mathcal T_H$ and $\mathcal T_{H'}$ intersect.
Thus by the above claim, to prove (\ref{eqn-fol2}), it suffices to show that $\alpha_l$ is a bijection  onto (Pic$^{1,1}E)^{\sigma}$.
First we show the surjectivity.
We have seen in Proposition \ref{prop-globalext} (v) that as the real plane $H$ goes to $H_i$, members of the $S^1$-family $|\mathcal L_H|^{\sigma}$ converge to $L_i$.
Hence the two circles $\mathcal T_H^+$ and $\mathcal T_H^-$ will shrink to the two points $P_i$ and $\ol{P}_i$ as $H$ goes to $H_i$.
On the other hand, we have $N_{L_i/Z}\simeq O(1)^{\oplus 2}$ (Lemma \ref{lemma-degree01}).
Therefore $Z$ has a structure of twistor space in a neighborhood of $L_i$.
In particular, all small deformations of $L_i$ in $Z$ keeping the reality must be disjoint.
This implies that, when a real plane $H$ moves in $\langle l\rangle^{\sigma}$ to pass $H_i$ from one side to another side, then the intersection circle $\mathcal T_H^+$ and $\mathcal T_H^-$ moves to pass (via shrinking) from $E^+$ to $E^-$ (or $E^-$ to $E^+$) without turning the direction.
This implies that $\mathcal T_H$ sweeps the whole of $E\backslash\{P_i,\ol{P}_i\}$ when $H$ runs through $\langle l\rangle^{\sigma}\backslash\{H_i\}$.
Hence $\alpha_l$ must be surjective onto (Pic$^{1,1}E)^{\sigma}$.
The injectivity can be proved in a similar way as in the case that $l\cap B$ is transversal.
Therefore $\alpha_l$ must be injective and  (\ref{eqn-fol2}) holds.
Hence $(\sharp)$ is proved for a real bitangent $l\neq l_{\infty}$.

Finally, we prove that $(\sharp)$ is true for $l=l_{\infty}$, or equivalently, 
\begin{equation}\label{eqn-fol3}
E\backslash \{P_i,\ol{P}_i\}_{i=1}^4=\bigcup_{H\in\langle l_{\infty}\rangle^{\sigma}\backslash\{ H_i\}_{i=1}^4}
\mathcal T_H,
\end{equation}
where $E=\Phi^{-1}(l_{\infty})$, $\mathcal T_H$ have the same meaning as before, and $\{P_i,\ol{P}_i\}_{i=1}^4$ denote the singular points of $E$.
This is possible thanks to the explicit construction of real lines on $S_{\lambda}$ given in Sections \ref{s-detc}--\ref{s-nb}.
By Proposition \ref{prop-octagon}, (see Figure \ref{fig-cycle3}), $\Phi^{-1}(l_{\infty})$ is a cycle of eight rational curves.
First consider the case when we take the small resolution determined by (I) of Proposition \ref{prop-sr}.
If $H$ is of the form $H_{\lambda}=\{y_0=\lambda y_1\}$ where $\lambda\in I_1$ (resp.\,$\lambda\in I_3$), then $L^+_{\theta}$ (resp.\,$L_{\theta}^-$) must be chosen  by Proposition \ref{prop-choice01} (i), where $L^{\pm}_{\theta}$ is of course the irreducible components of the inverse images of touching conics of special type. 
Hence by Lemma \ref{lemma-interse},  $\mathcal T_H$ is two circles, one of which is contained in $\Gamma_1$ (resp.\,$\Gamma_3$) and the other is contained in $\ol{\Gamma}_1$ (resp.\,$\ol{\Gamma}_3$). 
Further, their radiuses are given by the function $h_1$ (resp.\,$h_3$) defined after Lemma \ref{lemma-interse}.
For our choice of $\ell_1(=x_1)$, Lemma \ref{lemma-crtcl-sp} (i) implies that the function $h_1$ gives a bijection between $I_1$ and $\mathbf R^+$.
This implies that $\mathcal T_{H_{\lambda}}$ foliates $\Gamma_1\cup\ol{\Gamma}_1$ (with the $U(1)$-fixed points removed) when $\lambda$ runs through $I_1$.
In the same manner, $\mathcal T_{H_{\lambda}}$ foliates $\Gamma_3\cup\ol{\Gamma}_3$ (with the $U(1)$-fixed points removed) when $\lambda$ runs through $I_3$.
Next, if $H$ is of the form $H_{\lambda}=\{y_0=\lambda y_1\}$ where $\lambda$ is in $I_2$ then by Lemma \ref{lemma-intersection3}, $\mathcal T_H$ is two circles, one of which is contained in $\Gamma_2$ and the other is contained in $\ol{\Gamma}_2$. 
Moreover, their radiuses are given by the function $h_2$  defined in Proposition \ref{prop-criet}.
For our choices of $\ell_1(=x_1)$ and $\ell_2(=x_0+x_1)$, Lemma \ref{lemma-beha} (iii) implies that $h_2$ gives a bijection between $I_2$ and $\mathbf R^+$.
This implies that $\mathcal T_{H_{\lambda}}$ foliates $\Gamma_2\cup\ol{\Gamma}_2$ (with the two endpoints removed) when $\lambda$ runs through $I_2$.
Finally when $H$ is of the form $H_{\lambda}=\{y_0=\lambda y_1\}$ where $\lambda$ is in $I_4^-\cup I_4^+$, then $\mathcal T_H$ is two circles, one of which is contained in $\Xi$ and the other is contained in $\ol{\Xi}$.
Further, their radiuses are given by $h_0$ and $h^{-1}_0$. 
Thus if we recall that $L^+_{\theta}$ has to be chosen (Proposition \ref{prop-choice01}), and if we recall the rule of distinction of $L^+_{\theta}$ and $L^-_{\theta}$ made in the explanation just before Proposition \ref{prop-choice01}, $\mathcal T_{H_{\lambda}}$ covers $\Xi\cup \ol{\Xi}$ with the unit circles in $\Xi$ and $\ol{\Xi}$, and the four $U(1)$-fixed points removed.
The unit circles (in $\Xi$ and $\ol{\Xi}$) are covered by real lines in $H_{\lambda_0}$ whose images are real lines going through $P_0$.
Thus $\mathcal T_{H_{\lambda}}$ foliates $l\cup\ol{l}$ (with the four $U(1)$-fixed points removed) when $\lambda$ runs through $ I_4$.
Finally, the remaining points, that is, the eight singular points of $E$ are passed by $L_i$, $1\leq i\leq 4$.
Thus we have verified that $\mathcal T_H$ foliate $\Phi^{-1}(l_{\infty})$ as $H$ runs through $\linfty^{\sigma}$, and we have completed the claim of Proposition \ref{prop-disj}.
\proofend

\vspace{2mm}
\noindent{\textbf{Proof of the claim.}
Let $F\in$ Pic$^{1,1}E$ be a line bundle which is not necessarily real.
Because $F|_{E^+} \simeq O(1)$ and $F|_{E^-} \simeq O(1)$, $F$ must be the line bundle
associated to the effective divisor $a+b$ on $E$, where $a\in\mathbf C^*$ (resp. $b\in\mathbf C^*$) represents a point on $E^+\backslash\{P_i,\ol{P}_i\}$ (resp. $E^-\backslash\{P_i,\ol{P}_i\}$) satisfying $z=a$ (resp. $w=b$).
Namely, we have $F=[a+b]$ in the usual notation.
Then it is elementary to see that for $(a',b')\in\mathbf C^*\times\mathbf C^*$ the divisor $a'+b'$ is linearly equivalent to $a+b$ iff $(a',b')=(ca,cb)$ for some $c\in\mathbf C^*$.
Indeed, the condition of linear equivalence is the same thing as the existence of a meromorphic function $f$ on $E$ having $z=a$ and $w=b$ as simple zeros, and $z=a'$ and $w=b'$ as simple poles.
This $f$ must be of the form $c(z-a)/(z-a')$ on $E^+$ for some $c\in\mathbf C^*$ and $c'(w-b)/(w-b')$ on $E^-$ for some $c'\in \mathbf C^*$.
Further, the values of these two functions must coincide at $P_i$ (i.e. $z=w=0$) and $\ol{P}_i$( i.e. $z=w=\infty$).
This implies that $a/b=a'/b'$. 
Now consider a map from Div$_+^{1,1}E$ to $\mathbf C^*$ sending $a+b$ to $a/b$, where Div$_+^{1,1}E$ denotes the set of effective Cartier divisors whose degrees on $E^{\pm}$ are one.
Then the above argument implies that this map descends to give an isomorphism Pic$^{1,1}E\simeq\mathbf C^*$.

Next assume that $F\in $Pic$^{1,1}E$ is real and $s$ is a real section of $F$.
Recall that there is a (branched) double covering $E\ra l$ commuting with the real structures, where the real structure on $E$ interchanges $E^+$ and $E^-$, and the real structure on $l$ (a real line) is a complex conjugation.
Therefore we can write $F=[a+1/\ol{a}]$ for some $a\in\mathbf C^*$
Thus the image of $F$ in $\mathbf C^*$ by the above isomorphism Pic$^{1,1}E\ra\mathbf C^*$ is $(1/\ol{a})/a=1/|a|^2$. 
This implies the image of real line bundles is $\mathbf R_{>0}$.
(Although not needed in the sequel, this calculation also implies that if we change the real structure on $l$ from the complex conjugation to the antipodal map (and yet assuming $\sigma(E^+)=E^-$ and that $E\ra l$ commutes with the real structures), then $F$ becomes real iff the image is in $\mathbf R_{<0}$.)
It is immediate to see that for a real divisor $a+1/\ol{a}\in$Div$_+^{1,1}E $,  $\alpha a+(\alpha/\ol{a})$ remains to be real iff $\alpha\in U(1)$.
This implies that the union of the zero locus of a real line bundle $F=[a+1/\ol{a}]$  is two circles $|z|=|a|$ (in $E^+$) and $|w|=1/|a|$ (in $E^-$), which implies the claim.
\proofend

\section{The main theorems
}\label{s-thms}
In this section, based on the results in the previous section, we prove our main result  that our quartic surfaces always determine twistor spaces of $3\mathbf{CP}^2$ whose self-dual metrics are non-LeBrun while admitting a non-trivial Killing field (Theorem \ref{thm-inv}). 
Then by using this, we determine the moduli space of these self-dual metrics on $3\mathbf{CP}^2$ (Theorem \ref{thm-mod}).

In order to state our result precisely, we now recall the setting.
Let $(y_0,y_1,y_2,y_3)$ be a homogeneous coordinate on $\mathbf{CP}^3$ and consider a real structure defined by 
\begin{equation}\label{eqn-finalreal}
(y_0,y_1,y_2,y_3)\mapsto(\ol{y}_0,\ol{y}_1,\ol{y}_3,\ol{y}_2).
\end{equation}
Let $B$ be a quartic surface defined by
\begin{equation}\label{eqn-finalbranch}
(y_2y_3+Q(y_0,y_1))^2-y_0y_1(y_0+y_1)(y_0-ay_1)=0,\end{equation}
where  $a$ is a positive real numbers, and $Q$ is a quadratic form of $y_0$ and $y_1$ with real coefficients. 
Then $B$ is preserved by the real structure.
Suppose that $Q$ and $a$  satisfy the following properties which are necessary conditions for getting a twistor space (Proposition \ref{prop-necessa}):

\vspace{3mm}\noindent {\textbf{Condition (A)}}:{\em{
 If $\lambda$ satisfies $\lambda(\lambda+1)(\lambda-a)\geq 0$ (i.e. if $\lambda\in I_2\cup I_4$), then 
 \begin{equation}\label{eqn-finalcond}
 Q(\lambda,1)\geq\sqrt{\lambda(\lambda+1)(\lambda-a)}
 \end{equation}
holds. Moreover, there exists a unique $\lambda_0$, $\lambda_0>a$ (i.e. $\lambda_0\in I_4$) such that  the equality of (\ref{eqn-finalcond}) holds for $\lambda=\lambda_0$.}}
 
 \vspace{3mm}\noindent
 (We can assume $\lambda_0>a$ by possible application of a projective transformation  $(y_0,y_1)\mapsto(ay_1,-y_0)$ as explained in Section \ref{ss-prel}.
 This assumption will be important afterwards.)
 Then $B$ has three singular points $P_0$, $P_{\infty}$ and $\ol{P}_{\infty}$, where $P_0$ is a real ordinary double point and $\{P_{\infty},\ol{P}_{\infty}\}$ is a pair of conjugate points which are simple elliptic singularities of type $\tilde{E}_7$
(Proposition \ref{sing. of B}).
Let $Z_0\ra\mathbf{CP}^3$ be the double covering branched along $B$.
Since $B$ is real, the real structure on $\mathbf{CP}^3$ naturally induces that on $Z_0$.
$Z_0$ has three singular points $p_0$, $p_{\infty}$ and $\ol{p}_{\infty}$ over the singularities of $B$ respectively.
Then the following result, saying in particular that a quartic surface (\ref{eqn-finalbranch}) canonically determines (via twistor space) a self-dual metric on $3\mathbf{CP}^2$, is the main result of this paper:

\begin{thm}\label{thm-inv}
(i) There are precisely two small resolutions of $Z_0$ such that the resulting compact complex threefolds are  twistor spaces of $3\mathbf{CP}^2$.
(ii) The two self-dual metrics on $3\mathbf{CP}^2$ defined by the two twistor spaces are mutually conformally isometric.
(iii) The self-dual metrics on $3\mathbf{CP}^2$ determined by the above twistor spaces have an isometric $U(1)$-action but are not conformally isometric to LeBrun metric.
(iv) For  any non-LeBrun self-dual metric on $3\mathbf{CP}^2$ of positive scalar curvature admitting an isometric $U(1)$-action, there exists a quartic surface (\ref{eqn-finalbranch}) satisfying Condition (A) such that the resulting twistor space is the twistor space of the given self-dual metric.
\end{thm}

(Concerning (iv) of the theorem, there are two choices of quartic surface representing given self-dual metric in general,  as will be shown in Lemma \ref{lemma-projequiv}.)

\vspace{3mm}
\noindent Proof of Theorem \ref{thm-inv}.
By Propositions \ref{prop-sr} and  \ref{prop-summary}  there are precisely two small resolutions of $Z_0$  such that the resulting non-singular threefolds $Z$ can have a structure of twistor space.
We  show that these two threefolds  are actually  twistor spaces of $3\mathbf{CP}^2$.
By Proposition \ref{prop-globalext}, $Z$ has a connected family of real smooth rational curves. 
Moreover, by   Proposition \ref{prop-disj}, different members of this family are disjoint and  $Z$ is foliated by the curves of this family. 
Therefore $Z$ must be a twistor space of some four-dimensional manifold $M$, having this family  as the set of twistor lines.
By Lemma \ref{lemma-promise-1} or \ref{lemma-degree01}, $Z$ always possesses a divisor $D+\ol{D}$ such that $D$ is biholomorphic to a three points blown-up of $\mathbf{CP}^2$. 
In particular, $D$ is diffeomorphic to $\mathbf{CP}^2\#(\ol{\mathbf{CP}}^2\#\ol{\mathbf{CP}}^2\#\ol{\mathbf{CP}}^2)$, where $\ol{\mathbf{CP}}^2$ denotes the complex projective plane with the complex orientation reversed.
Then $L=D\cap \ol{D}$ is a real line in $Z$ (Lemma \ref{lemma-degree01}) and $L$ is a member of our family of real lines in $Z$ (Proposition \ref{prop-globalext} (v)).
Since $D$ and $\ol{D}$ intersect transversally along $L$,
we have $(N_{D/Z})|_L\simeq N_{L/\ol{D}}$.
Hence $(D,L)_Z=\deg([D]|_L)=\deg(N_{D/Z}|_L)=\deg(N_{L/\ol{D}})$, and the last term is one.
Thus we get $(D,L)_Z=1$.
It follows that $(D,L')_Z=1$ for any member $L'$ of our family.
Moreover, our construction of the family implies that $L$ is the unique member which is contained in $D+\ol{D}$.
Therefore for any $L'\neq L$, $L'$ intersects $D$ (and $\ol{D}$) transversally at a unique point.
This implies that the parameter space of our family can be identified with a one-point compactification of $D\backslash L$.
It is easily seen that the latter space is diffeomorphic to $3\mathbf{CP}^2$ with the usual complex orientation reversed.
Thus we have shown that $Z$ is a twistor space of $3\mathbf{CP}^2$ and we get (i) of the theorem.

Next we show (ii).
We remark that our family of real lines is a unique family which gives $Z$ a structure of twistor space, since by Proposition \ref{prop-image} real lines in $Z$ intersecting $\Gamma_0$ must be of the form $\Phi^{-1}(l)-\Gamma_0$ (which is actually contained in our family of real lines), where $l$ is a real line through $P_0$.
Hence to prove (ii)  it suffices to construct a biholomorphic map $\phi$ between the two twistor spaces, such that $\phi$ commutes with the real structures on the twistor spaces,  because we have the above uniqueness of the family of twistor lines.
Let $\iota:\mathbf{CP}^3\ra\mathbf{CP}^3$ be a holomorphic involution defined by $\iota(y_0,y_1,y_2,y_3)= (y_0,y_1,y_3,y_2)$. 
It is immediate to verify that $\iota$ commutes with our real structure on $\mathbf{CP}^3$,  $\iota$ preserves $B$, and that $\iota(P_0)=P_0,\iota(P_{\infty})=\ol{P}_{\infty}$ hold.
Moreover, $\iota$ naturally lifts on the line bundle $O(2)\ra\mathbf{CP}^3$ in such  a way that $\iota^*x_0=x'_0,\iota^* x_1=x'_1,\iota^*x_2=x'_3$ and $\iota^*z=z$ hold,  where $(x_0,x_1,x_2)=(y_0/y_3,y_1/y_3,y_2/y_3)$ and $(x'_0,x'_1,x'_3)=(y_0/y_2,y_1/y_2,y_3/y_2)$ are non-homogeneous coordinates on $y_3\neq 0$ and $y_2\neq 0$ respectively as in the proof of Lemma \ref{lemma-promise-1}, and $z$ is a fiber coordinate on the bundle $O(2)\ra \mathbf{CP}^3$.
If we still denote $\iota$ for this lift on $O(2)$, $\iota$ preserves $Z_0\subset O(2)$ invariant, and $\iota$ commutes with the real structure on $O(2)$.
Moreover, by (\ref{eqn-Z_0}), $\iota$ interchanges two factors of $\mathbf{CP}^1\times\mathbf{CP}^1$ that is the projectified tangent cone of $Z_0$ at $P_0$.
Hence $\iota$ is lifted to give a  a biholomorphic map between two small resolutions of $Z_0$, and we get the desired biholomorphic map $\phi$.
Thus we obtain (ii).

To prove (iii) let $\rho$ be  a $U(1)$-action  on $\mathbf{CP}^3$ defined by 
$(y_0,y_1,y_2,y_3)\mapsto$ 
$(y_0,y_1,e^{i\theta}y_2$, $e^{-i\theta}y_3)$, $e^{i\theta}\in U(1).$
It is immediate to see that $\rho$ commutes with the real structure on $\mathbf{CP}^3$ and preserves $B$.
Hence it is naturally lifted on $Z_0$  to be a (holomorphic) $U(1)$-action  commuting with the real structure.
Moreover, since  $\rho$ fixes the singular points of $B$, this $U(1)$-action on $Z_0$ is automatically lifted on a small resolution $Z$.
This action on $Z$ is clearly a holomorphic $U(1)$-action commuting with the real structure. 
Therefore it induces an isometric $U(1)$-action of the corresponding self-dual metric on  $3\mathbf{CP}^2$.
Moreover, the self-dual metrics are not conformally isometric to LeBrun metrics, since $|(-1/2)K_Z|$ of our twistor spaces do not have base point (cf. Proposition \ref{prop-pic}) but LeBrun twistor spaces have. 
Thus we get (iii).

Finally, (iv) is obvious from Proposition \ref{prop-def-B}.
\proofend

\vspace{3mm}
We use Theorem \ref{thm-inv} to determine the moduli space of self-dual metrics on $3\mathbf{CP}^3$ of this kind. 
We begin with the following

\begin{lemma}\label{lemma-moduli03}
Let $\tilde{\mathcal M}$ be the set of the quartic surface (\ref{eqn-finalbranch}) satisfying Condition (A).
Then $\tilde{\mathcal M}$ is diffeomorphic to $\mathbf R^3$.
\end{lemma}

\noindent Proof.
Since the quartic (\ref{eqn-finalbranch}) is determined by a quadratic form $Q(y_0,y_1)$ and $a>0$, we represent a point of $\tilde{\mathcal M}$ by a pair $(a,Q)$.
Let $\mathcal N$ be an open subset of $\mathbf R^2$ defined by 
$$\mathcal N=\{(x,y)\in\mathbf R^2\set x>0,\,y>x\}.$$
We consider the map $\varpi:\tilde{\mathcal M}\ra\mathcal N$ defined by
$\varpi(a,Q)=(a,\lambda_0)$,
where $\lambda_0$ is the unique real number attaining the equality of (\ref{eqn-finalcond}). 
Then it is obvious that $\varpi$ is surjective.
Fix $(a,\lambda_0)\in\mathcal N$ and
set $f(\lambda)=\lambda(\lambda+1)(\lambda-a)$ as before.

\begin{figure}[htbp]
\includegraphics{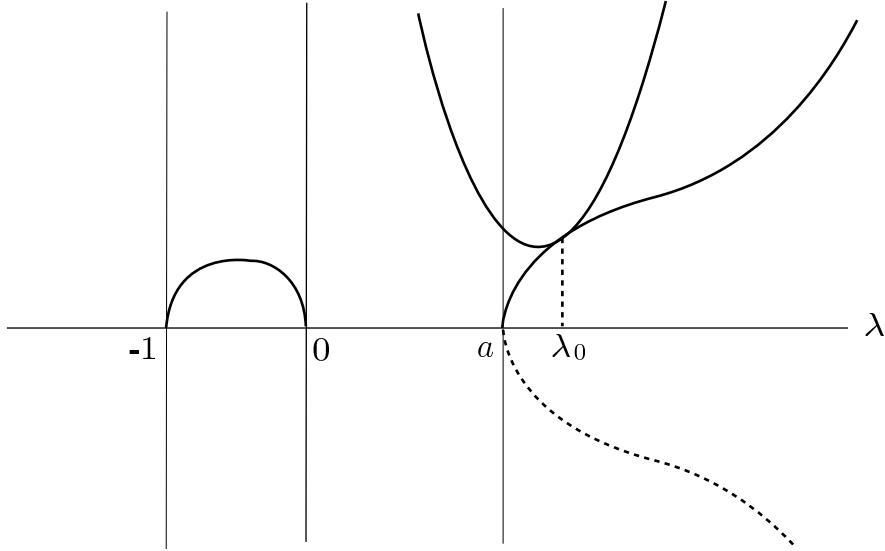}
\caption{graphs of $\sqrt{f}$ and $Q$}
\label{fig-sqrt}
\end{figure}

The graph of $\sqrt{f(\lambda)}$ is as in Figure \ref{fig-sqrt} and $\varpi^{-1}(a,\lambda_0)$ is the set of $Q(y_0,y_1)$ such that the graph of $Q(\lambda,1)$ is over that of $\sqrt{f}$ ( on $I_2\cup I_4$) and is tangent (from above) at the point $(\lambda_0,\sqrt{f(\lambda_0)})$.
Writing $Q(y_0,y_1)=by_0^2+cy_0y_1+dy_1^2$, the last condition is equivalent to the two equalities
\begin{equation}\label{eqn-line001}
b\lambda_0^2+c\lambda_0+d=\sqrt{f(\lambda_0)},\,\,\,\,2b\lambda_0+c=\left.\left(\sqrt{f(\lambda)}\right)'\right|_{\lambda=\lambda_0}.
\end{equation}
These define a line in $\mathbf R^3$.
Moreover, an elementary argument shows that there is a positive constant $b_0=b_0(\lambda_0)>0$ such that Condition (A) (for fixed $\lambda_0$) is equivalent to (\ref{eqn-line001}) together with the inequality $b>b_0(\lambda)$.
Thus we get that $\varpi^{-1}(a,\lambda_0)$ is a half-line in $\mathbf R^3$.
Namely $\varpi:\tilde{\mathcal M}\ra\mathcal N$ is a half-line bundle on $\mathcal N$. On the other hand it is readily seen that $\mathcal N$ is diffeomorphic to $\mathbf R^2$. 
Hence $\tilde{\mathcal M}$ is the total space of $\mathbf R$-bundle over $\mathbf R^2$.
Therefore $\tilde{\mathcal M}$ must be diffeomorphic to $\mathbf R^3$, as desired.
\proofend

\vspace{3mm}
At first sight $\tilde{\mathcal M}$ may be thought of a moduli space of self-dual metrics. 
However, there is a possibility that different quartics of the form (\ref{eqn-finalbranch})  determines the same (i.e. conformally isometric) self-dual metric.
Actually, we have the following:

\begin{lemma}\label{lemma-projequiv}
Let $(a,Q)$ and $(a',Q')$ be points of $\tilde{\mathcal M}$, and $B$ and $B'$ the associated quartic surfaces defined by (\ref{eqn-finalbranch}) respectively.
Let $g$ and $g'$ be the self-dual metrics on $3\mathbf{CP}^2$ canonically associated to $B$ and $B'$ respectively by Theorem \ref{thm-inv}.
Then $g$ and $g'$ are conformally isometric if and only if the following two conditions are satisfied.
(i) $a=a'$, (ii) $Q'(ay_0+ay_1, y_0-ay_1)=a(a+1)Q(y_0,y_1)$.
\end{lemma}

\noindent Proof.
If there is a conformal isometry between $g$ and $g'$, 
it naturally defines a biholomorphic map between the twistor spaces. 
This automorphism also defines an automorphism $G$ of $\mathbf{CP}^3$ commuting with the real structure (\ref{eqn-finalreal}), because $\mathbf{CP}^3$ is canonically identified with the dual projective space of $H^0(Z,(-1/2)K_Z)\simeq\mathbf C^4$, and the automorphism of the twistor space naturally lifts on the line bundle $(-1/2)K_Z$.
We have clearly $G(B)=B'$.
So suppose that $G$ is an automorphism of $\mathbf{CP}^3$ commuting with the real structure (\ref{eqn-finalreal}), and satisfies $G(B)=B'$.
Of course, $G$ is a projective transformation.

Assume first that $G(P_{\infty})=P_{\infty}$. 
Then it follows $G(\ol{P}_{\infty})=\ol{P}_{\infty}$ and it can be easily verified that $G$ must be of the form
\begin{equation}\label{eqn-mat}
G=\left(\begin{array}{cc}
A_1&O\\
O&A_2\\
\end{array}\right),\hspace{3mm}A_1\in GL(2,\mathbf{R}),\hspace{3mm}
A_2=\left(\begin{array}{cc}
a_2&0\\
0&\ol{a}_2\\
\end{array}\right),\hspace{2mm}a_2\neq0.
\end{equation}
Moreover, we have $G(B)=B'$ if and only if the following two conditions are satisfied:
\begin{equation}\label{eqn-101}
Q'(y_0',y_1')=|a_2|^2Q(y_0,y_1)
\end{equation}
and 
\begin{equation}\label{eqn-102}
y_0'y_1'(y'_0+y'_1)(y'_0-a'y'_1)=|a_2|^4y_0y_1(y_0+y_1)(y_0-ay_1),
\end{equation}
where we put $^t(y_0',y_1')=A_1\cdot\,^t(y_0,y_1)$.
Because the cross-ratio is invariant under projective transformations, 
(\ref{eqn-102}) implies $a=a'$.
It is elementary to show that if $a\neq 1$, there are just three linear transformations of $\mathbf{CP}^1$ preserving the four points $\{-1,0,a,\infty\}$, where we put $-1=(-1,1), 0=(0,1),a=(a,1)$ and $\infty=(1,0)$ and that they are concretely given by
$$
A_1^{(1)}=\left(\begin{array}{rr}
0&a\\
-1&0\\
\end{array}\right),\hspace{3mm}
A_1^{(2)}=\left(\begin{array}{rr}
a&a\\
1&-a\\
\end{array}\right),\hspace{3mm}
A_1^{(3)}=\left(\begin{array}{rr}
1&-a\\
-1&-1\\
\end{array}\right).
$$
Note that these define involutions on $\mathbf{CP}^1$, $A_1^{(1)}$ preserves orientation of $\mathbf{RP}^1\simeq S^1$,   $A_1^{(2)}$ and $A_1^{(3)}$ reverse it, $ A_1^{(1)}$ and $ A_1^{(3)}$ interchange $I_2=(-1,0)$ and $I_4=(a,\infty)$, $ A_1^{(2)}$ preserves $I_2$ and $I_4$, and that $A_1^{(1)}= A_1^{(2)} A_1^{(3)}$ in $PGL(2,\mathbf{R})$. (Thus $ A_1^{(1)},  A_1^{(2)}$ and $ A_1^{(3)}$ generate a subgroup of $PGL(2,\mathbf{R})$ which is isomorphic to $\mathbf Z_2\times \mathbf Z_2$.)
Thus $ A_1^{(1)}$ and $ A_1^{(3)}$ do not preserves Condition (A) (since they interchange $I_2$ and $I_4$ and hence does not preserve the property $\lambda_0\in I_4$).
Therefore, at least if $a\neq 1$, $ A_1^{(2)}$ is the unique one which can determine $G$ satisfying $G(B)=B'$.
It is immediate to check that under $A_1^{(2)}$, the right hand side of $(\ref{eqn-102})$ is just multiplied by $a^2(a+1)^2$.
Hence we obtain $|a_2|^2=a(a+1)$.
Substituting this into (\ref{eqn-101}), we get $Q'(y_0',y_1')=a(a+1)Q(y_0,y_1)$.
This is the equation in (ii) of the lemma.
If $a=1$, there arise other $A_1$'s preserving the four points.
But all of them satisfies $y_0'y_1'(y'_0+y'_1)(y'_0-y'_1)=cy_0y_1(y_0+y_1)(y_0-y_1)$ with $c<0$ (or more precisely $c=-1$ or $-4$), and  is not compatible with (\ref{eqn-101}).
Thus even if $a=1$, we have $A_1=A_1^{(2)}$.

Next, assume that $G(P_{\infty})=\ol{P}_{\infty}$.
Then it can be also easily seen that
\begin{equation}\label{eqn-mat2}
G=\left(\begin{array}{cc}
A_1&O\\
O&A_2\\
\end{array}\right),\hspace{3mm}A_1\in GL(2,\mathbf{R}),\hspace{3mm}
A_2=\left(\begin{array}{cc}
0&a_2\\
\ol{a}_2&0\\
\end{array}\right),\hspace{2mm}a_2\neq0.
\end{equation}
Once this is obtained, it is completely parallel to the last paragraph to conclude that $a=a'$, $A_1=A_1^{(2)}$ and $|a_2|^2=a(a+1)$.
Hence we again get  $Q'(ay_0+ay_1, y_0-ay_1)=a(a+1)Q(y_0,y_1)$ and we have proved the necessity.

Conversely, suppose that $(a,Q)$ and $(a',Q')\in \tilde{\mathcal M}$ satisfy $a=a'$ and $Q'(ay_0+ay_1, y_0-ay_1)=a(a+1)Q(y_0,y_1)$.
Then for the real projective transformation
$$
G:(y_0,y_1,y_2,y_3)\mapsto \left(ay_0+ay_1, y_0-ay_1,\sqrt{a(a+1)}y_2,\sqrt{a(a+1)}y_3\right),
$$
we readily have $G(B)=B'$. 
This $G$ (automatically) maps singular points of $B$ onto $B'$.
Therefore $G$ naturally lifts on a small resolution of the double covers.
Thus $G$ gives rise to an isomorphism of the two twistor spaces.
Hence the converse is also proved.
\proofend

\vspace{3mm}
We use these two lemmas to obtain the following result which describes a global structure of the moduli space:

\begin{thm}\label{thm-mod}
Let $\mathcal M$ be the set of conformal classes of self-dual metrics $g$ on $3\mathbf{CP}^2$ satisfying the following properties: (i) the scalar curvature of $g$ is positive, (ii) $g$ admits a non-trivial Killing field, (iii) $g$ is not conformally isometric to LeBrun's explicit self-dual metrics \cite{LB91}.
Then $\mathcal M$ is naturally identified with an orbifold\, $\mathbf{R}^3/G$,  where $G$ is an involution of\, $\mathbf R^3$ having two-dimensional fixed locus.
\end{thm}

\noindent Proof.
By Lemma \ref{lemma-projequiv}, $\mathcal M$ is naturally identified with $\tilde{\mathcal M}/G$, where $G$ is an involution on $\tilde{\mathcal M}$ defined by 
\begin{equation}\label{eqn-involution2}
\left(a,Q(y_0,y_1)\right)\mapsto \left(a,\,\,\frac{1}{a(a+1)}{Q(ay_0+ay_1,y_0-ay_1)}\right).
\end{equation}
By Lemma \ref{lemma-moduli03}, $\tilde{\mathcal M}$ is diffeomorphic to $\mathbf R^3$. 
Therefore, in order to prove the theorem, it suffices to show that $G$ acts on $\tilde{\mathcal M}\simeq\mathbf R^3$ having two-dimensional fixed locus.
To see this, write $Q(y_0,y_1)=by_0^2+cy_0y_1+dy_1^2$ with $b,c,d\in\mathbf R$, as in the proof of Lemma \ref{lemma-moduli03}.
Then as in the proof, we can use $(a,\lambda_0,b)$ as a coordinate on $\tilde{\mathcal M}$, because $c$ and $d$ are uniquely determined from $b$ by (\ref{eqn-line001}).
Write $G(a,\lambda_0,b)=(a',\lambda_0',b')$. 
Then (\ref{eqn-involution2}) means that 
\begin{equation}\label{eqn-invol3}
a'=a,\hspace{3mm}\lambda_0'=\frac{a\lambda_0+a}{\lambda_0-a}.
\end{equation}
Moreover, (\ref{eqn-involution2}) also implies that $b'$ is the coefficient of $y_0^2$ of the polynomial 
\begin{equation}\label{eqn-invol4}
\frac{1}{a(a+1)}{Q(ay_0+ay_1,y_0-ay_1)}.
\end{equation}
Now we assert that  the fixed locus of $G$ is precisely the set $\{(a,\lambda_0,b)\in\tilde{\mathcal M}\set \lambda_0=a+\sqrt{a^2+a}\}$.
This directly implies the  theorem.
To show the assertion, first it is obvious that $\lambda_0$ must be the fixed point of the fractional transformation in (\ref{eqn-invol3}), which is easily known to be $a+\sqrt{a^2+a}$.
On the other hand,  by (\ref{eqn-invol4}), we have
$b'=(ba^2+ca+d)/a(a+1)$, where $c$ and $d$ are determined by (\ref{eqn-line001}).
Substituting $\lambda_0=a+\sqrt{a^2+a}$,  various cancellation occurs and we finally  get the last one becomes just $b$.
Therefore, if $\lambda_0=a+\sqrt{a^2+a}$, then $b$ gets no effect from $G$.
This proves the assertion and we obtain the theorem.

\proofend

\vspace{3mm}
By the above proof,  a self-dual metric on $3\mathbf{CP}^2$ of positive scalar curvature with a non-trivial Killing field admits an isometry not contained in the $U(1)$-isometry generated by the Killing field, iff the corresponding point of the moduli space $\mathcal M$ is a $G$-fixed point. 
Moreover,   in such a case, the isometry must be an involution.

We also note that by construction of LeBrun \cite{LB91} and his classification theorem \cite{LB93} the moduli space of self-dual metrics on $3\mathbf{CP}^2$ whose isometry group is just $U(1)$ (i.e. cannot extends to $U(1)^2$-isometry) acting semi-freely on $3\mathbf{CP}^2$ is identified with the set of different non-collinear three points on the hyperbolic three-space  modulo the isometry group.
In particular, it has dimension  $3\times 3-6=3$.
It will be really interesting to try to find an explicit description of the self-dual metrics associated to our quartics.

\small
\vspace{13mm}
\hspace{5.5cm}
$\begin{array}{l}
\mbox{Department of Mathematics}\\
\mbox{Graduate School of Science and Engineering}\\
\mbox{Tokyo Institute of Technology}\\
\mbox{2-12-1, O-okayama, Meguro, 152-8551, JAPAN}\\
\mbox{{\tt {honda@math.titech.ac.jp}}}
\end{array}$

\end{document}